%% file: Main.tex
\documentclass[english,a4paper,10pt]{book}
\usepackage{babel,amsmath,amssymb,amsthm,mathrsfs,graphicx}
\usepackage[top=3cm, bottom=3cm, outer=3cm, inner=4.5cm, heightrounded, marginparwidth=3cm, marginparsep=0.5cm]{geometry}
\usepackage[latin1]{inputenc}
\usepackage[T1]{fontenc}
\usepackage{indentfirst}
\usepackage{url,amsfonts,epsfig,eso-pic}
\usepackage{tikz-cd}
\usepackage{leftidx}
\usepackage{marginnote}
\theoremstyle{definition}
\newtheorem*{theorem*}{Theorem}
\newtheorem{theorem}{Theorem}[chapter]
\newtheorem{prop}[theorem]{Proposition}
\newtheorem{lemma}[theorem]{Lemma}
\newtheorem{defn}[theorem]{Definition}
\newtheorem{obs}[theorem]{Observation}
\newtheorem{cor}[theorem]{Corollary}

\newtheorem{rmk}[theorem]{Remark}

\allowdisplaybreaks

\newcommand{\V}{\varnothing}
\newcommand{\E}{\exists}
\renewcommand{\P}{\mathbb{P}}
\newcommand{\C}{\mathbb{C}}
\renewcommand{\E}{\mathbb{E}}
\renewcommand{\V}{\mathbb{V}}

\newcommand{\Q}{\mathbb{Q}}

\newcommand{\Z}{\mathbb{Z}}
\newcommand{\A}{\mathbb{A}}

\renewcommand{\L}{\mathbb{L}}

\renewcommand{\H}{\mathbb{H}}
\newcommand{\id}{\mathbb{I}}

\newcommand{\stless}{\operatornamewithlimits{<}_{(=)}}
\newcommand{\stgreat}{\operatornamewithlimits{>}_{(=)}}

\newcommand{\git}{\mathbin{
  \mathchoice{/\mkern-6mu/}% \displaystyle
    {/\mkern-6mu/}% \textstyle
    {/\mkern-5mu/}% \scriptstyle
    {/\mkern-5mu/}}}% \scriptscriptstyle

\DeclareMathOperator{\im}{Im}
\DeclareMathOperator{\rk}{rk}

\DeclareMathOperator{\supp}{supp}
\DeclareMathOperator{\s}{s}

\DeclareMathOperator{\GL}{GL}

\DeclareMathOperator{\Hom}{Hom}
\DeclareMathOperator{\Ext}{Ext}
\DeclareMathOperator{\End}{End}
\DeclareMathOperator{\Aut}{Aut}

\DeclareMathOperator{\Sym}{Sym}

\DeclareMathOperator{\Spec}{Spec}
\DeclareMathOperator{\Gr}{Gr}
\DeclareMathOperator{\hilb}{Hilb}
\DeclareMathOperator{\R}{R}

\def\quotient#1#2{%
    \raise1ex\hbox{$#1$}\Big/\lower1ex\hbox{$#2$}%
}
\begin{document}
\frontmatter
\date{\today}
\title{\huge{\textbf{Motivic invariants of moduli spaces\\ of rank 2 Bradlow - Higgs triples}}}
\author{Riccardo Grandi}
\maketitle
\tableofcontents
\input{abstract}
\input{Acknowledgements}
\input{Intro}
\mainmatter
\input{Background}
\input{BHT}
\input{Motives}
\input{Poles}

\input{MSMY}
\input{Higher}
\backmatter
\nocite{*}
\bibliographystyle{plain}
\bibliography{biblioRePlan}
\end{document}

%% file: abstract.tex
\chapter{Abstract}
In the present thesis we study the geometry of the moduli spaces of Bradlow-Higgs triples on a smooth projective curve $C$. $(E,\phi, s)$ is a Bradlow-Higgs triple if $(E,\phi)$ is a Higgs bundle and $s$ is a non-zero global section of $E$. There is a family of stability conditions for triples that depends on a positive real parameter $\sigma$. The moduli spaces $\mathcal{M}_\sigma^{r,d}$ of $\sigma$-semistable triples of rank $r$ and degree $d$ vary with $\sigma$. The phenomenon arising from this is known as wall-crossing.

In the first half of the thesis we will examine how the moduli spaces $\mathcal{M}_\sigma^{r,d}$ and their universal additive invariants change as $\sigma$ varies, for the case $r=2$. In particular we will study the case of $\sigma$ very close to 0, for which $\mathcal{M}_\sigma^{r,d}$ relates to the moduli space of stable Higgs bundles, and $\sigma$ very large, for which $\mathcal{M}_\sigma^{r,d}$ is a relative Hilbert scheme of points for the family of spectral curves. Some of these results will be generalized to Bradlow-Higgs triples with poles.

In the second half we will prove a formula relating the cohomology of $\mathcal{M}_\sigma^{2,d}$ for small $\sigma$ and $d$ odd and the perverse filtration on the cohomology of the moduli space of stable Higgs bundles. The formula is not far from the generalized Macdonald formulas found in \cite{migliorini2011support}, \cite{maulik2014macdonald} and \cite{migliorini2015support}. We will also partially generalize this result to the case of rank greater than 2.\\[10pt]

\textbf{Keywords:} Moduli spaces, wall-crossing, Bradlow-Higgs triples, Macdonald formula, Hilbert scheme.
\chapter{R\'esum\'e}
Dans la pr\'esente th\`{e}se, nous \'etudions la g\'eom\'etrie des espaces de modules de triplets de Bradlow-Higgs sur une courbe lisse projective $C$. $(E, \phi, s)$ est un triplet de Bradlow-Higgs si $ (E, \phi) $ est un fibr\'e de Higgs et $s$ est une section globale de $ E $ non-nulle. Il y a une famille de conditions de stabilit\'e pour les triplets qui d\'epend d'un param\`{e}tre positif r\'eel $ \sigma $. Les espaces de modules $ \mathcal {M}_\sigma^{r, d} $ de triplets $\sigma$-semistable de rang $r$ et degr\'e $d$ varient avec $ \sigma $. Le ph\'enom\`{e}ne r\'esultant est connu comme wall-crossing.

Dans la premi\`{e}re moiti\'e de la th\`{e}se, nous examinerons comment les espaces de modules $\mathcal{M}_\sigma^{r, d} $ et leurs invariants additifs universels changent en fonction de $ \sigma $, pour le cas $ r = 2 $. En particulier, nous allons \'etudier le cas de $ \sigma $ tr\`{e}s proche de 0, pour lequel $ \mathcal {M}_\sigma^{r, d} $ est reli\'e \`{a} l'espace des modules de fibr\'es de Higgs stables, et $ \sigma $ tr\`{e}s grand, pour lequel $ \mathcal{M}_\sigma^{r, d} $ est un sch\'ema de Hilbert de points relatif pour la famille de courbes spectrales. Certains de ces r\'esultats seront g\'en\'eralis\'es aux triplets de Bradlow-Higgs avec des p\^oles.

Dans la seconde moiti\'e, nous allons prouver une formule concernant la cohomologie de $ \mathcal{M}_\sigma^{2, d} $ pour $ \sigma $ petit et $d$ impair et la filtration perverse sur la cohomologie de l'espace des modules de fibr\'es de Higgs stables de degr\'e impair. La formule est proche des formules de Macdonald g\'en\'eralis\'ees trouv\'ees dans \cite{migliorini2011support}, \cite{maulik2014macdonald} et \cite{migliorini2015support}. Nous allons aussi partiellement g\'en\'eraliser ce r\'esultat au cas du rang sup\'erieur \`{a} 2.\\[10pt]

\textbf{Mots-cl\'e:} Espaces de modules, wall-crossing, triplets de Bradlow-Higgs, formule de Macdonald, sch\'ema de Hilbert.

%% file: Acknowledgements.tex
\chapter{Acknowledgements}
First of all, I would like to thank my thesis advisor, Tam\'as Hausel, who made all this possible by pointing out the project to me, supervising my activity through these years and forming a lively geometry group in which I grew a lot as a person and as a mathematician.

I would also like to thank our beloved Pierrette, that helped and took care of me even before my official start at EPFL. She always says that it is her job, but one can notice the passion, kindness and care for others that she provides every day.

I have a huge debt of gratitude towards Daniele Boccalini, Alexandre Peyrot, Dimitri Wyss and Martina Rovelli, my closest "academic relatives" and good friends, for providing countless hours of fun, company and mathematical discussions as well.

A big thank you goes to all present and former members of the Geometry Group: Michael Groechenig, Michael Lennox Wong, Ben Davison, Martin Mereb, Zongbin Chen, Alexander Noll. I was very lucky to share this time with all of you and thank you for your useful suggestions.

Many visitors came to Lausanne to exchange thoughts on research for brief periods of time. Among those, I would like to thank Luca Migliorini and Jochen Heinloth for taking an interest in my work and answering many of my questions. Another big thank you for accepting to review my final work. I would like to thank Kathryn Hess and Donna Testerman as well for being so kind to spend their time being part of the jury for my PhD defense.

Last but not least I would like to thank my friends and family whose moral support was fundamental through these years.

During these years, my research was supported by the Project Funding number 144119: \emph{Arithmetic harmonic analysis on Higgs, character and quiver varieties} issued by the Swiss National Science Foundation.

%% file: Intro.tex
\chapter{Introduction}
\section*{Description of the project}
The main goal of the thesis is to understand the geometry of the moduli spaces of Bradlow-Higgs triples on a smooth projective curve $C$, with particular focus on the case of rank 2. These objects consist of a vector bundle $E$ on $C$, a morphism $\phi: E \rightarrow E \otimes K$ and a non-zero global section of $E$. Since there is a Higgs field $\phi$ and a section $s$, Bradlow-Higgs triples relate to both Higgs bundles and Bradlow pairs. See \cite{nitsure1991moduli} and \cite{thaddeus1994stable} for an introduction on the two moduli problems. Features from both original moduli problems are then inherited by the moduli problem of Bradlow-Higgs triples.

As for Bradlow pairs, there exists a family of stability conditions that depends on a positive real parameter $\sigma$. This will allow to vary the stability and hence produce several moduli spaces $\mathcal{M}_\sigma^{r,d}$ of $\sigma$-semistable triples of rank $r$ and degree $d$. For fixed $r$ and $d$, the set of positive real numbers is then partitioned into a finite number of intervals, the last one of which is unbounded. The defining property of these intervals is that varying $\sigma$ in the interior of each one of these intervals will yield the same moduli space $\mathcal{M}_\sigma^{r,d}$, but different intervals have non-isomorphic associated moduli spaces. The endpoints of these intervals are called \emph{critical values} and $\mathcal{M}_\sigma^{r,d}$ will change whenever $\sigma$ crosses one. This phenomenon is well known as \emph{wall-crossing}.

The moduli spaces $\mathcal{M}_\sigma^{r,d}$ will inherit properties from the moduli problem of Higgs bundles as well. On the other hand they are not all smooth, and this makes their geometry richer but more complicated. As we will see, each of them admits a proper \emph{Hitchin map}
$$\chi_\sigma^{r,d}: \mathcal{M}_\sigma^{r,d} \rightarrow \mathcal{A}^r$$
whose target, the \emph{Hitchin base}, is an affine space. Furthermore, $\C^*$ acts on each of the $\mathcal{M}_\sigma^{r,d}$ by scaling the Higgs field and there is an action of $\C^*$ on $\mathcal{A}^r$ for which the Hitchin map is equivariant. In particular we will prove that for every $(E,\phi,s) \in \mathcal{M}_\sigma^{r,d}$ there exists the limit
$$\lim_{\lambda \rightarrow 0} \lambda \cdot (E,\phi,s)$$
and this has important consequences on the geometry of our moduli spaces.

When $\sigma$ crosses one of the critical values, the geometry of the moduli space will change, in the sense that some triples will become unstable and they will be erased from $\mathcal{M}_\sigma^{r,d}$, while some others will become stable and hence will be added to $\mathcal{M}_\sigma^{r,d}$. The loci that are added and removed are often known as \emph{flip loci}. Furthermore there are two distinguished stability conditions. Namely, for $\sigma=\varepsilon$ smaller than the least critical value, the stability of the triple will imply the semistability of the underlying Higgs bundle and so we have a forgetful map, also known as \emph{Abel-Jacobi map},
$$\mathcal{M}_\varepsilon^{r,d} \rightarrow \mathcal{M}^{r,d}$$
whose target is the moduli space of semistable Higgs bundles of rank $r$ and degree $d$. If $r$, $d$ are coprime and $d$ is very large compared to $r$ then the Abel-Jacobi map is a projective bundle, providing a direct relation between $\mathcal{M}_\varepsilon^{r,d}$ and $\mathcal{M}^{r,d}$.

Instead, when $\sigma$ is large enough to be bigger than the last critical value then the stability condition is equivalent to requiring that the section $s$ is a cyclic vector for $\phi$. In this case, the map $\chi_\sigma^{r,d}$ can be proved to be a relative Hilbert scheme of points. It is also related to moduli spaces of stable pairs studied by Pandharipande and Thomas in \cite{pandharipande2009curve}. The connection between Bradlow-Higgs triples and Hilbert schemes of points is also an important part of the project and a consequence of the study of the geometry of Bradlow-Higgs triples.

\section*{Content of the thesis and main results}
In chapter 1 we collect some known results from the literature about the main topics that are needed to understand the thesis. First we introduce the invariants that we will compute throughout the work i.e. the Grothendieck motives. We recall some results about Higgs bundles and Bradlow pairs that will set the starting point for the definition of Bradlow-Higgs triples. We will also recollect some classical results about non-abelian Hodge theory and the relation between moduli spaces of Higgs bundles and character varieties. Some attention will be given to the wall-crossing for the moduli problem of Bradlow pairs, which is a good starting point for the understanding of the wall-crossing for Bradlow-Higgs triples. Finally we give some insight about the notion of Hilbert scheme, especially in the case of curves, since moduli spaces of Bradlow-Higgs triples are connected to relative Hilbert schemes.

Chapter 2 will contain the basic definition of Bradlow-Higgs triples and their $\sigma$-stabilty. We also outline the construction of the moduli spaces $\mathcal{M}_\sigma^{r,d}$ as open subsets of moduli spaces of coherent systems on the surface $\P(\mathcal{O}_C \oplus K)$ that were originally studied in \cite{le1995faisceaux}. This approach will also lead to the proof of the properness of the Hitchin maps $\chi_\sigma^{r,d}$. There will be a study of the deformation theory of Bradlow-Higgs triples, in order to understand if the moduli spaces are singular. We will prove that after the first critical value is crossed the moduli spaces are in fact singular. The $\C^*$-action on the moduli spaces of Bradlow-Higgs triples is understood and for rank 2 we provide a decomposition of the moduli spaces in terms of the connected components of the fixed point loci. For some combinations of $\sigma$, $r$ and $d$ we can prove that $\mathcal{M}_\sigma^{r,d}$ is smooth and this will allow to compute its cohomology from the motives. Here is the main outcome of the chapter:
\begin{theorem*}
Bradlow-Higgs triples $(E,\phi,s)$ on $C$ correspond to pairs $(\mathcal{F},s)$ where $\mathcal{F}$ is a rank one pure one dimensional sheaf on the surface $\P(\mathcal{O}_C \oplus K)$ whose support does not intersect the divisor at infinity and $s$ is a non-zero global section of $\mathcal{F}$. Under this identification the Hitchin base $\mathcal{A}^r$ is a family curves in $\P(\mathcal{O}_C \oplus K)$ whose support does not intersect the divisor at infinity. Furthermore the Hitchin maps $\chi_\sigma^{r,d}: \mathcal{M}_\sigma^{r,d} \rightarrow \mathcal{A}^r$ correspond to taking the scheme theoretic support of the sheaves $\mathcal{F}$.

For $\sigma \rightarrow 0$ there is an Abel-Jacobi map
$$AJ: \mathcal{M}_\sigma^{r,d} \rightarrow \mathcal{M}^{r,d}$$
forgetting the section and whose target is therefore the moduli space of semistable Higgs bundles.

For $\sigma \rightarrow \infty$, $\mathcal{M}_\sigma^{r,d} \rightarrow \mathcal{A}^r$ is a relative Hilbert scheme of points for the family of curves parametrized by $\mathcal{A}^r$.

For $r=2$, $\sigma \rightarrow 0$ and either $d<0$ or $d>4g-4$ odd, $\mathcal{M}_\sigma^{r,d}$ is smooth. In all other cases it is singular.

$\mathcal{M}_\sigma^{r,d}$ can be decomposed into attracting sets for a $\C^*$-action, explicitly for $r=2$:
$$\mathcal{M}_\sigma^{2,d}=F_{(2),\sigma}^{(d),1+} \sqcup \bigsqcup_{I_1} F_{(1,1),\sigma}^{(d_1,d_2),1+}\sqcup \bigsqcup_{I_2} F_{(1,1),\sigma}^{(d_1,d_2),2+}.$$

Each of the $F^+$ contains a connected component $F$ of the fixed point locus of $\mathcal{M}_\sigma^{2,d}$ and is characterized by the property that the limit as $\lambda \rightarrow 0$ of $\C^*$ acting on points in $F^+$ will belong to $F$.
\end{theorem*}

In chapter 3 we compute the motivic invariants of the $\mathcal{M}_\sigma^{r,d}$ for $r=2$. There are three main parts. First, we examine the flip loci and we understand what triples are added and erased from $\mathcal{M}_\sigma^{r,d}$ as a critical value $\bar\sigma$ is crossed. We can prove that the two flip loci $\mathcal{W}^{d,+}_{\bar\sigma}$ and $\mathcal{W}^{d,-}_{\bar\sigma}$ admit maps:
$$
\pi_{\bar\sigma}^{d,+}: \mathcal{W}_{\bar\sigma}^{d,+} \rightarrow X^d_{\bar\sigma}
$$
and
$$
\pi_{\bar\sigma}^{d,-}: \mathcal{W}_{\bar\sigma}^{d,-} \rightarrow X^d_{\bar\sigma}
$$
where $X^d_{\bar\sigma}$ is the cartesian product of a symmetric power of the curve $C$, of its Jacobian and of two copies of an affine space. The flip loci themselves contain triples that can be described as extensions of Higgs bundles for which the canonical subobject and quotient have prescribed degrees. The fibers of the maps $\pi_{\bar\sigma}^{d,+}$ and $\pi_{\bar\sigma}^{d,-}$ are projective spaces, although not of constant dimension, and this suffices to compute the motive of the flip loci. We will also point out how the flip loci interact with a decomposition of the $\mathcal{M}_\sigma^{2,d}$ that is obtained by exploiting the $\C^*$-action. Second we compute $[\mathcal{M}_\varepsilon^{2,d}]$ for small $\varepsilon$. The strategy relies once again on the $\C^*$-action. Last, we will point out a strategy to compute $[\mathcal{M}_\sigma^{2,d}]$ for $\sigma$ bigger than the last critical value. The main results of the chapter can be formulated as follows:
\begin{theorem*}
There are maps:
$$
\pi_{\bar\sigma}^{d,+}: \mathcal{W}_{\bar\sigma}^{d,+} \rightarrow X^d_{\bar\sigma}
$$
and
$$
\pi_{\bar\sigma}^{d,-}: \mathcal{W}_{\bar\sigma}^{d,-} \rightarrow X^d_{\bar\sigma}
$$
whose fibers are projective spaces. This allows to compute motives:
\begin{equation*}
[\mathcal{W}_{\bar\sigma}^{d,+}]=\L^{2g} \cdot [\C\P^{2g-3}]\cdot [S^{(d-{\bar\sigma})/2}( C)] \cdot [J( C)]+\L^{3g-2} \cdot [S^{(d-{\bar\sigma})/2}( C)] \cdot [S^{\bar\sigma}( C)],
\end{equation*}
\begin{equation*}
[\mathcal{W}_{\bar\sigma}^{d,-}]=\L^{2g}\cdot[S^{(d-{\bar\sigma})/2}( C)]\cdot[J^{(d+{\bar\sigma})/2}( C)]\cdot [\C\P^{(d+\bar\sigma)/2+g-2}].
\end{equation*}

The motive of $\mathcal{M}_\varepsilon^{2,d}$ for $0<\varepsilon<1$ and $d >0$ odd can be computed from the decomposition into attracting sets:
\begin{align*}
[\mathcal{M}_\varepsilon^{2,d}]=&\L^{1+4(g-1)}[M_\varepsilon^{2,d}]+\sum_{(d_1,d_2) \in I_1^o(d)} \L^{1+3(g-1)+d_2}[S^{d_1}( C)][S^{d_1-d_2+2g-2}]+\\
&+\sum_{(d_1,d_2) \in I_2^o(d)}\L^{1+4(g-1)} [S^{d_2}( C)][S^{d_1-d_2+2g-2}( C)]+\\
&+\sum_{(d_1,d_2) \in I_1^o(d)} (\L^{4g-2}-\L^{4g-3}) [S^{d_1}( C)] \left ([S^{d_2}( C)]-[J( C)] \frac{\L^{d_2+1-g}-1}{\L-1}\right).
\end{align*}

A similar formula also holds for $d \geq 0$ even.
\begin{align*}
[\mathcal{M}_\varepsilon^{2,d}]&=\L^{1+4(g-1)}  [M^{2,d}_\varepsilon]+\sum_{(d_1,d_2) \in I_1^e(d)} \L^{1+3(g-1)+d_2} [S^{d_1}( C)][S^{d_1-d_2+2g-2}( C)] +\\
&+\sum_{(d_1,d_2) \in I_2^e(d)}\L^{1+4(g-1)} [S^{d_2}( C)][S^{d_1-d_2+2g-2}( C)]+\\
&+(\L-1) \L^{1+4(g-1)}[\Sym^2( S^{d/2} ( C) )] +\L^{1+4(g-1)} [S^{d/2}( C)][J^{d/2}( C)][\C\P^{g-2}]+\\
&+[S^{d/2}( C)] \L^{3g-2} \left ( \L^{d/2+g-1}+\L^{2g-2}-1\right )+\\
&+\sum_{(d_1,d_2) \in I_1^e(d)} (\L^{4g-2}-\L^{4g-3}) [S^{d_1}( C)] \left ([S^{d_2}( C)]-[J( C)] \frac{\L^{d_2+1-g}-1}{\L-1} \right).
\end{align*}

For $\sigma > d$, $[\mathcal{M}_\infty^{2,d}]$ can be computed either by combining the above formulas and the motive of the flip loci or directly.
\end{theorem*}

Chapter 4 contains some comments about the generalization of the results of the previous chapters to the case of Bradlow-Higgs triples with poles. Fixed an integer $\gamma \geq 1$ and a point $P \in C$, a Bradlow-Higgs $\gamma$-triple is the datum $(E,\phi,s)$ where $E$ is a vector bundle on $C$, $s $ is a non-zero global section of $E$ and $\phi: E \rightarrow E \otimes K(\gamma P)$ is a Higgs field that is allowed to have poles at $P$. Most of the properties of the $\mathcal{M}_\sigma^{2,d}$ are still valid for their analogues with poles $\mathcal{M}_\sigma^{2,d}(\gamma)$. We also prove that when $\gamma$ is chosen large enough, $\mathcal{M}_\sigma^{2,d}(\gamma)$ are smooth, the flip loci in the wall crossing are actually projective bundles and their motive can be computed in an easier way. The two main results are:
\begin{theorem*}
$\mathcal{M}_\varepsilon^{2,d}(\gamma)$ is smooth for all $\gamma \geq 1$. In particular $\mathcal{M}_\varepsilon^{2,d}(\gamma)$ is always semiprojective.

If $\gamma > d$ then $\mathcal{M}_\sigma^{2,d}(\gamma)$ is smooth regardless of $\sigma$ as long as it is different from a critical value. In this case then $\mathcal{M}_\sigma^{2,d}(\gamma)$ is semiprojective.

We have:
\begin{align*}
&\lim_{\gamma \rightarrow \infty} P(\mathcal{M}_{\varepsilon}^{2,d}(\gamma),t) = \lim_{\gamma \rightarrow \infty} P(\mathcal{M}_{\infty}^{2,d}(\gamma),t)=\\
&=\frac{(1+t^3)^{2g}(1+t)^{2g}}{(1-t^2)^2(1-t^4)}=P(\C \P^{\infty},t) P(B \overline{\mathcal{G}},t)
\end{align*}
where $B \overline{\mathcal{G}}$ is the classifying space of the group $\overline{\mathcal{G}}$ mentioned in \cite[section 7.2]{hausel2001geometry}.
\end{theorem*}

In chapter 5 we explore the relation between the moduli spaces of Bradlow-Higgs triples and relative Hilbert schemes of points on curves. We already mentioned that for very large $\sigma$, $\chi_{\sigma}^{2,d}$ is a Hilbert scheme relative to the family of locally planar curves parametrized by the Hitchin base $\mathcal{A}^2$. A lot of recent work by several authors, see \cite{migliorini2011support}, \cite{maulik2014macdonald} and \cite{migliorini2015support}, relates the cohomology of the Hilbert scheme of points on a curve to the cohomology of the compactified Jacobian of the same curve. The key assumption is that the curve has to be locally planar and at least reduced. Note that these formulas are valid for certain families of curves and their relative Hilbert scheme. In the chapter we exploit the connection between Bradlow-Higgs triples and relative Hilbert schemes of points on curves to prove a partial formula relating the cohomology of $\mathcal{M}_\varepsilon^{2,d}$ for $d$ odd to the perverse filtration on the cohomology of $\mathcal{M}^{2,1}$. The main formula is as follows:
\begin{theorem*}
Setting
\begin{equation*}
F^{sh}(q)= \sum_{n \geq 1-g} \R (\chi_\varepsilon^{2,2n+1})_* (IC_{\mathcal{M}_\varepsilon^{2,2n+1}}) q^{2n+2g-1},
\end{equation*}
\begin{equation*}
F^{vir}(q,t)= \sum_{n \geq 1-g} P^{vir}(\mathcal{M}_\varepsilon^{2,2n+1},t) q^{2n+2g-1}
\end{equation*}
and
\begin{equation*}
G(q,t)=\text{odd}_q \left (\frac{PH(\mathcal{M}^{2,1},q,t)}{(1-q)(1-qt^2)} \right ),
\end{equation*}
\begin{equation*}
G^{sh}(q)=\text{odd}_q \left (\frac{\bigoplus_{i=0}^{8g-6} IC\left (\bigwedge^i\R^1 \right) }{(1-q\Q)(1-q\Q[-2](-1))} \right )
\end{equation*}
we can prove that $F^{sh}$ and $G^{sh}$ coincide for $\deg q \leq 2g-3$. $P(F^{sh}(q),t)$ and $G$ coincide for $\deg q \leq 2g-3$ and for $\deg q \geq 6g-5$. Furthermore $P(F^{sh}(q),t)-G(q,t)$ is a polynomial with non-negative coefficients.
\end{theorem*}

In the last chapter we will discuss the issues that arise when trying to generalize the content of the thesis to Bradlow-Higgs triples of rank bigger than 2. In the first part, we focus on the fact that while there exists an approach to understand the wall-crossing of Bradlow pairs for arbitrary rank and degree, see \cite{mozgovoy2013moduli}, this will fail for Bradlow-Higgs triples and we point out why. In the second part we will partly generalize the content of chapter 5 to higher rank.

%% file: Background.tex
\chapter{General facts and background results}
%
%\begin{quote}
%{\small
%Summary of chapter contents\\
%\begin{itemize}
%\item Cohomology, weights, Grothendieck ring of varieties
%\item Generalities on Higgs bundles (def, moduli and stability, deformations, Hitchin map, properness, BNR, $\C^*$-action and fixed points, cohomology, ???twisted Higgs bundles)
%\item Bradlow pairs (def, moduli and stability, deformations, Abel-Jacobi map, cohomology and motives)
%\item Hilbert schemes (definiton, moduli, deformations, families of Hilbert schemes)
%\end{itemize}
%}
%\end{quote}
%
\section{Grothendieck ring of varieties over $\C$ and motivic invariants}
In the rest of the thesis we will use a ring containing universal additive invariants of varieties over $\C$. Good references for this are \cite[Section 2]{behrend2007motivic} and \cite[Section 1]{garcia2011motives}.
\begin{defn}
Define $K_0(\text{Var}_\C)$ to be the ring with the following presentation:
\begin{itemize}
\item the generators are the isomorphism classes $[X]$ of varieties over $\C$
\item the relations are generated by $[X]-[Y]-[X \setminus Y]$ whenever $Y$ is a closed subvariety of $X$
\item the product is defined by $[X]\cdot [Y]=[X \times Y]$ where $\times$ is the product of $\C$-schemes over $\Spec(\C)$ (or in this case the Cartesian product of varieties).
\end{itemize}
$K_0(\text{Var}_\C)$ is called \emph{Grothendieck ring of varieties over $\C$}.
\end{defn}
A few remarks are necessary. The relations extend to $[X]=[Y]+[X \setminus Y]$ when $Y$ is a locally closed subvariety of $X$. The unit of this ring is the class of a point $[pt]$, while the $0$ is the class of the empty variety. The class of the affine line $\A^1$ is denoted by $\L$.

Here we mention some relevant examples of relations in $K_0(\text{Var}_\C)$ some of which will be used in the rest of the work.
\begin{itemize}
\item If $X \rightarrow Y$ is a Zariski locally trivial fibration with fiber $F$, then $[X]=[Y][F]$.
\item If $X \rightarrow Y$ is bijective on closed points then $[X]=[Y]$ and this is also true under the more general condition that $X$ and $Y$ can be written as disjoint unions $X= \sqcup X_i$, $Y=\sqcup Y_j$ with a bijection between the index sets such that $[X_i]=[Y_i]$. An example of this occurs if there are isomorphisms $X_i \rightarrow Y_i$ that do not extend to maps with larger domain than the $X_i$.
\item We could have defined the Grothendieck ring using isomorphism classes of schemes of finite type. However this is not important for our purposes and also, using this last definition it is easy to see that $[X]=[X_{red}]$.
\item Pick a smooth projective curve $C$ of genus $g \geq 2$. We can consider the symmetric powers of $C$, $S^n( C)$ with the Abel-Jacobi map:
\begin{align*}
AJ_n: S^n( C) & \rightarrow J^n( C)\\
D & \mapsto \mathcal{O}(D).
\end{align*}
It is a well known fact that the fibers of $AJ_n$ are projective spaces, not necessarily of constant dimension. More precisely $AJ_n^{-1}(\{L\})=\P H^0(L)$. It is also known that $J^n ( C)$ can be stratified in such a way that $AJ_n$ is a projective bundle on each of the strata. We can be more precise. Define the \emph{Brill-Noether locus}:
\begin{equation*}
V^n_i=\{L \in J^n( C) : \dim H^0(L)=i\}
\end{equation*}
which are locally closed subvarieties of $J^n( C)$. Then we have the motivic relations:
\begin{equation*}
[S^n ( C)]=\sum_{i=1}^{n+1} [V_i^n] [\C \P^{i-1}]
\end{equation*}
and
\begin{equation*}
[J^n ( C)]=\sum_{i=0}^{n+1} [V_i^n].
\end{equation*}
Note that, in general, the strata $V_i^n$ are very complicated, maybe singular, and they depend on the complex structure of $C$, i.e. on the specific curve and not only on the genus. The sum of their motives, as well as the sum of the motives weighted with the appropriate projective space, can be expressed nicely as we have seen above.\\
For $n>2g-2$, $AJ_n$ is a $\C\P^{n-g}$ bundle since $H^1$ will vanish for line bundles of degree bigger than $2g-2$ and hence the dimension on $H^0$ is constant.
\end{itemize}

$K_0(\text{Var}_\C)$ is often completed to the ring $\widehat{K_0(\text{Var}_\C)}$ in the following way. First invert $\L$ in $K_0(\text{Var}_\C)$. Then inside the ring $K_0(\text{Var}_\C)_\L=K_0(\text{Var}_\C)[\L^{-1}]$ consider the filtration $F^\bullet$ defined by imposing that $F^m$ is generated by the elements:
\begin{equation*}
\frac{[X]}{\L^n}
\end{equation*}
with $\dim X -n \leq -m$. Then $\widehat{K_0(\text{Var}_\C)}$ is obtained by completing $K_0(\text{Var}_\C)[\L^{-1}]$ with respect to this filtration. We omit most of the details, but we just remark that in $\widehat{K_0(\text{Var}_\C)}$, $\L^n-1$ is invertible for all $n$ since:
\begin{equation*}
\frac{1}{\L^n-1}=\L^{-n}(1+\L^{-n}+\dots).
\end{equation*}

Also $[\GL_n]$ is invertible for all $n$ since:
\begin{equation*}
[\GL_n]=(\L^n-1)(\L^n-\L)\dots(\L^n-\L^{n-1}).
\end{equation*}

The idea behind this completion is that we would like to be able to define the motive of a certain class of stacks. Namely Artin stacks, locally of finite type, whose geometric stabilizers are linear algebraic groups. A result by Kresch \cite{kresch1999cycle} states that all of these stacks admit a stratification whose strata are quotients of varieties by $\GL_n$ for different $n$. Note that even though the stratification is not necessarily finite, the properties of these stacks imply that the motive is well defined in the completion.

As we said in the beginning, $K_0(\text{Var}_\C)$ is the ring of universal additive invariants associated to varieties. Recall that an \emph{additive invariant} $\theta$ with values in a ring $R$ is a function $K_0(\text{Var}_\C)$ that associates to the isomorphism class of a variety $X$ the invariant $\theta(X) \in R$. $\theta$ must also satisfy the scissor relations, i.e. $\theta(X)=\theta(X \setminus Y)+\theta(Y)$ whenever $Y$ is a closed subvariety of $X$, and the product rule, i.e. $\theta(X \times Y)=\theta(X)\cdot \theta(Y)$. This is equivalent to asking that $\theta$ is a ring homomorphism $K_0(\text{Var}_\C) \rightarrow R$.

Some examples are:
\begin{itemize}
\item the Euler characteristic $\chi: K_0(\text{Var}_\C) \rightarrow \Z$
\item the $E$-polynomial defined as:
\begin{equation*}
E(X, x,y)=\sum_{p,q,j} (-1)^j h_c^{p,q,j}(X) x^p y^q
\end{equation*}
where $h_c^{p,q,j}(X)=\dim \left ( Gr^F_p Gr^W_{p+q} H_c^j(X)\right )$ are the compactly supported Hodge numbers of $X$. The additive nature of the $E$-polynomial follows from the exact sequence in compactly supported cohomology arising from the decomposition $X= U \cup X \setminus U$ of a variety into an open subset and its complement. For more details see \cite[Section 2]{hausel2008mixed}. Here $E$ takes values in $R=\Z[x,y]$.
\end{itemize}

As a last remark, we recall that if the cohomology of a variety $X$ is pure, meaning that the weight filtration on $H^i(X)$ is limited to the $i$-th weight, then the $E$ polynomial (which is additive) will determine both the Hodge polynomial:
\begin{equation*}
H(X, x,y,t)=\sum_{p,q,j} h^{p,q,j}(X) x^p y^q t^j
\end{equation*}
by the substitution:
\begin{equation*}
H(X, x,y,t)=(xyt)^{2 \dim (X)} E\left ( X, -\frac{1}{xt}, -\frac{1}{yt}\right )
\end{equation*}
and the Poincar\'e polynomial by the substitution:
\begin{equation*}
P(X,t)=t^{2 \dim (X)} E\left ( X, -\frac{1}{t}, -\frac{1}{t}\right ).
\end{equation*}
\section{Moduli spaces of Higgs bundles}
\label{higgs}
Let us fix a smooth projective curve $C$ over $\C$ of genus $g \geq 2$. Denote by $K$ the canonical bundle of $C$.
\begin{defn}[Higgs bundles and stability]
A \emph{Higgs bundle} on $C$ is a pair $(E,\phi)$ consisting of a vector bundle $E$ and a twisted endomorphism $\phi: E \rightarrow E \otimes K$. $\phi$ is also referred to as the \emph{Higgs field}.

The slope of a Higgs bundle $(E,\phi)$ is defined to be the slope of $E$, i.e. $\mu(E)=\deg E / \rk E$.\\
We say that a Higgs bundle is (semi)stable if for every $F \subset E$ proper $\phi$-invariant subbundle of $E$ we have:
\begin{equation*}
\mu(F) \stless \mu(E).
\end{equation*}
\end{defn}
After fixing the rank $r$ and the degree $d$, we can define the moduli space $\mathcal{M}^{r,d}$ of S-equivalence classes of (semi)stable Higgs bundles. We summarize here the main properties of $\mathcal{M}^{r,d}$. Good references for this section are \cite{hitchin1987self}, \cite{nitsure1991moduli} among many others. For a nice summary of the case $r=2$ see \cite{hausel2001geometry}.

$\mathcal{M}^{r,d}$ is a complex variety of dimension $2+2 r^2(g-1)$ and it is smooth when $(r,d)=1$. The Zariski tangent space to a point $(E,\phi) \in \mathcal{M}^{r,d}$ is given by the hypercohomology of the following two step complex:
\begin{align*}
\End E &\rightarrow \End E \otimes K\\
f & \mapsto [f, \phi].
\end{align*}

Note that the complex is invariant under Serre's duality and this is responsible, in the smooth case, for $\mathcal{M}^{r,d}$ being a symplectic variety.

There are several useful maps between these moduli spaces, as $d$ varies. Here we outline just a few that will be used in the following chapters.\\
Tensoring with a line bundle $L$:
\begin{align*}
\mathcal{M}^{r,d} &\rightarrow \mathcal{M}^{r,d+r \deg L}\\
(E,\phi) & \mapsto (E \otimes L, \phi)
\end{align*}
via the identification $\End E \cong \End (E \otimes L)$. This is an isomorphism with inverse given by tensoring with $L^*$. In particular $\mathcal{M}^{r,d}$ is isomorphic to $\mathcal{M}^{r,d'}$ when $d-d'$ is divisible by $r$. This for example implies that in rank 2, up to isomorphism, there are only two moduli spaces of (semi)stable Higgs bundles, namely $\mathcal{M}^{2,0}$ which is singular and $\mathcal{M}^{2,1}$ which is smooth.

There is also a form of Serre's duality for Higgs bundles:
\begin{align*}
\mathcal{M}^{r,d} &\rightarrow \mathcal{M}^{r,2r(g-1)-d }\\
(E,\phi) & \mapsto (E^* \otimes K, \phi)
\end{align*}
using the identification $\End E \cong \End (E^* \otimes K)$. This is also an isomorphism since it is an involution. As an example, for rank $3$, this implies that $\mathcal{M}^{3,1}$ and $\mathcal{M}^{3,2}$ are smooth and isomorphic, while $\mathcal{M}^{3,0}$ is singular.

$\mathcal{M}^{r,d}$ always contains a copy of the moduli space of S-equivalence classes of rank $r$ degree $d$ (semi)stable vector bundles over $C$, denoted by $N^{r,d}$, embedded as Higgs bundles with zero Higgs field. More precisely, let $(r,d)=1$ so that both $N^{r,d}$ and $\mathcal{M}^{r,d}$ are smooth. Then the fiber of the cotangent bundle $T^*N^{r,d}$ over $E \in N^{r,d}$ is canonically identified with $H^1(\End E)^*$ which, in turn, is canonically identified with $H^0(\End E \otimes K)$.

Observe that if $E \in N^{r,d}$ is a stable vector bundle, then certainly $(E,\phi)$ is a stable Higgs bundle regardless of $\phi$. This means that the embedding:
\begin{align*}
N^{r,d} &\rightarrow \mathcal{M}^{r,d}\\
E &\mapsto (E,0)
\end{align*}
extends to an embedding:
\begin{align*}
T^*N^{r,d} \rightarrow \mathcal{M}^{r,d}
\end{align*}
whose image is an open dense symplectic subset of $\mathcal{M}^{r,d}$.

$\mathcal{M}^{r,d}$ is not proper but admits a proper map to an affine space. The map is known as \emph{Hitchin morphism} and is defined as follows. Let $\mathcal{A}^r = H^0(K) \oplus H^0(K^2) \oplus \dots \oplus H^0(K^r)$ be the so called \emph{Hitchin base}. Define:
\begin{align*}
h^{r,d}: \mathcal{M}^{r,d} &\rightarrow \mathcal{A}^r\\
(E,\phi) &\mapsto \text{char poly}(\phi).
\end{align*}
The dimension of $\mathcal{A}^r$ is $1+r^2(g-1)$ which is exactly half of $\dim \mathcal{M}^{r,d}$. The map $h^r$ is proper  complete integrable system whose generic fibers, as we will see later, are (torsors over) Abelian varieties. The proof of the properness can be found in \cite{nitsure1991moduli}, while the proof of second assertion can be found in \cite{hitchin1987stable}.
Beauville, Narasimhan, Ramanan (and many others extended their initial results) proved in \cite{beauville1989spectral} that $\mathcal{M}^{r,d}$ can be identified with a moduli space of stable sheaves on the cotangent bundle $T^*C=\text{Tot}(K)$ of $C$. Indeed we have the following equivalence:
\begin{equation*}
\left \{
\begin{array}{c}
\mathcal{F} \text{ pure one dimensional}\\
\text{sheaves of rank one on } T^*C\\
\text{whose support does not}\\
\text{intersect the divisor at infinity}
\end{array} \right \}
\longleftrightarrow
\left \{
\begin{array}{c}
(E,\phi) \text{ Higgs bundles on } C\\
\end{array} \right \}
\end{equation*}
given by the pushforward $\pi_*$ with respect to the projection map $\pi: T^*C \rightarrow C$. More precisely, given a sheaf $\mathcal{F}$ as above we can see that, since $\pi$ is affine, $\supp \mathcal{F}$ does not intersect the divisor at infinity and $\mathcal{F}$ is pure dimensional, $E=\pi_* \mathcal{F}$ is a torsion free sheaf on $C$ and therefore a vector bundle. With the hypotheses above, $\supp\mathcal{F}$ must be an $r$ to $1$ cover of the zero section in $T^*C$ and the rank of $E$ is equal to $r$. Since $E=\pi_* \mathcal{F}$ it also carries an action of the sheaf of algebras $\pi_* \mathcal{O}_{T^*C}=\Sym^* K^*$. This is the same as a morphism $K^* \otimes E \rightarrow E$ which can be identified with a Higgs field $\phi: E \rightarrow E \otimes K$.

Observe that $\mathcal{A}^r$ is the parameter space of a flat family of curves in $T^*C$ given by their equations. Namely the point $(\psi_1, \dots, \psi_r) \in \mathcal{A}_r$ corresponds to the curve of equation $y^r+\psi_r(x) y^{r-1}+\dots+\psi_1(x)=0$ where $x$ is the coordinate on the curve $C$ and $y$ is the coordinate on the fibers of the cotangent bundle of $C$. These curves are all $r$ to $1$ covers of the zero section of $T^*C$ and are generically smooth projective curves of genus $1+r^2(g-1)$. The curves parametrized by $\mathcal{A}^r$ are called \emph{spectral curves}. The family however contains curves that are not irreducible and even non-reduced. Under the identification of Higgs bundles on $C$ with one dimensional sheaves on $T^*C$, the Hitchin morphism will send a sheaf $\mathcal{F}$ to $\supp(\mathcal{F})$ which is a curve in $T^*C$. Therefore the statement of the BNR correspondence can be slightly improved by adding that if $(E,\phi)=\pi_* \mathcal{F}$ then $\supp (\mathcal{F})$ is the curve whose equation is given by $h^{r,d}(E,\phi)$. It is also important to note that the generic curve in the Hitchin base is integral and that if $(E,\phi)$ has a $\phi$-invariant proper subbundle $F$, then the characteristic polynomial of $\phi$ will factor and will be divisible by the characteristic polynomial of $\phi_{|F}$. These two facts together imply that the generic Higgs bundle in $\mathcal{M}^{r,d}$ will have no $\phi$-invariant subbundles and therefore will be trivially stable.

Let $X$ be a spectral curve. The fiber of $h^{r,d}$ over $X$ is easily seen to be (a torsor over) the Jacobian of $X$ if $X$ is smooth. If $X$ is singular but integral, then the Jacobian is replaced by the compactified Jacobian of $X$ which parametrizes torsion free sheaves on $X$. When $X$ is only reduced, meaning it could be singular and have more than one irreducible component, then the fiber of $h^{r,d}$ is a fine compactified Jacobian of $X$ for a specific choice of the stability condition.%\marginpar{\footnotesize Spoke with Luca about this, can't find the reference yet}. This is the reason why $\mathcal{M}^{r,d} \rightarrow \mathcal{A}^r$ could be regarded as a relative fine compactified Jacobian even though it is not clear what the fiber of $h^{r,d}$ is above a non-reduced spectral curve.

There is an action of $\C^*$ on $\mathcal{M}^{r,d}$ defined by:
\begin{align*}
\C^* \times \mathcal{M}^{r,d} &\rightarrow \mathcal{M}^{r,d}\\
(\lambda, (E, \phi)) & \mapsto (E, \lambda \phi).
\end{align*}

After we endow the Hitchin base with the action:
\begin{align*}
\C^* \times \mathcal{A}^r &\rightarrow \mathcal{A}^r\\
(\lambda, (\psi_1, \dots ,\psi_r)) & \mapsto (\lambda\psi_1, \lambda^2 \psi_2, \dots ,\lambda^r\psi_r),
\end{align*}
we can see that the Hitchin map has the property that $h^{r,d}(\lambda \cdot (E,\phi)=\lambda \cdot h^{r,d}(E,\phi)$.

This, together with the fact that $h^{r,d}$ is proper, implies that for all $(E,\phi) \in \mathcal{M}^{r,d}$:
\begin{equation*}
\lim_{\lambda \rightarrow 0} \lambda \cdot (E,\phi)
\end{equation*}
exists and lies in $\mathcal{M}^{r,d}$. Clearly, the limit points for the $\C^*$-action are the fixed points for the action. In general the fixed points have the following form:
\begin{prop}
\label{split}
Suppose that $(E,\phi)$ is a semistable Higgs bundle of rank $r$ such that $\lambda \cdot (E,\phi) \cong (E,\phi)$ for some $\lambda \in \C^*$ that is not a root of unity, then $E=E_1 \oplus \dots \oplus E_m$, for some $m \leq r$, $\sum \rk E_i=r$ and $\sum \deg E_i =\deg E$. Furthermore, $\phi$ has the property that $\phi(E_i) \subseteq E_{i-1} \otimes K$, $\phi(E_1)=0$.
\begin{proof}
See for example \cite[Lemma 4.1]{simpson1992higgs}.
\end{proof}
\end{prop}
For example, if $m=1$ in the proposition above, then $E$ is (semi)stable as a vector bundle and $\phi=0$ so that we find the moduli space of semistable vector bundles as one of the fixed point components for the $\C^*$-action. For a discussion about the cohomology of $\mathcal{M}^{r,d}$ see section \ref{semipr}.\\
\section{Semiprojective varieties and their cohomology}
\label{semipr}
The moduli space of Higgs bundles, as well as some of the moduli spaces of Bradlow-Higgs triples satisfy some properties that have strong implications on their cohomology. In this section we give a short overview of semiprojective varieties. All the material can be found in greater detail in \cite{hausel2013cohomology}.
\begin{defn}
Let $X$ be a quasi-projective complex variety with an action of $\C^*$. We call $X$ \emph{semiprojective} if the following two properties are satisfied:
\begin{itemize}
\item[(i)] the fixed point set $X^{\C^*}$ is proper
\item[(ii)] for all $x \in X$, the limit $\lim_{\lambda \rightarrow 0} \lambda \cdot x$ for $\lambda \in \C^*$ exists in $X$.
\end{itemize}

Let $X^{\C^*}=\bigsqcup_i F_i$ be the decomposition of the fixed point locus of $X$ into connected components, where $i$ ranges over some index set $I$. Define $U_i$ to be the set of points $x \in X$ for which $\lim_{\lambda \rightarrow 0} \lambda \cdot x $ lies in $F_i$. Define also $D_i$ to be the set of points $x \in X$ for which $\lim_{\lambda \rightarrow \infty} \lambda \cdot x $ lies in $F_i$.
\end{defn}
The $U_i$ are often referred to as \emph{attracting sets (or cells, or loci)} and the $D_i$ are referred to as \emph{downward flows}. Both the $U_i$ and the $D_i$ are locally closed subsets of $X$.\\
Note that if $X$ is semiprojective then condition (ii) implies the decomposition:
\begin{equation*}
X=\bigsqcup_i U_i.
\end{equation*}

The decomposition of $X$ into attracting sets is referred to as \emph{Bia\l{}ynicki-Birula decomposition}, after \cite[Theorem 4.1]{bialynicki1973some}.\\
We also give an important definition.
\begin{defn}
Let $X$ be a semiprojective variety. The \emph{core} of $X$ is:
\begin{equation*}
\mathcal{C}:= \cup_i D_i.
\end{equation*}
\end{defn}
According to \cite[Corollary 1.2.2]{hausel2013cohomology}, if $X$ is semiprojective then the core $\mathcal{C}$ is a proper subvariety of $X$. Moreover, according to \cite[Theorem 1.3.1]{hausel2013cohomology} the embedding $\mathcal{C} \rightarrow X$ induces an isomorphism $H^*(\mathcal{C},\Z)\cong H^*(X,\Z)$ when $X$ is smooth.

If we assume that $X$ is smooth, in addition to being semiprojective, then \cite[Corollary 1.3.2]{hausel2013cohomology} states that $X$ has pure cohomology. In fact, on one hand, $H^*(X)=H^*(\mathcal{C})$ is the cohomology of a proper variety so the weight filtration $W_\bullet$ on $H^*(X)$ has the property that $W_r H^i(X)=W_i H^i(X)=H^i(X)$ for $r \geq i$. On the other hand, since $X$ is smooth we have $W_r H^i(X)=0$ for $r<i$ and therefore $H^i(X)$ is pure of weight $i$. If $X$ is smooth, according to \cite[Corollary 1.3.6]{hausel2013cohomology} the core $\mathcal{C}$ is a deformation retract of $X$. The smoothness of $X$ also implies, see \cite[Theorem 4.1]{bialynicki1973some}, that both the $U_i$ and the $D_i$ are Zariski locally trivial affine fibrations over the $F_i$.

Thanks to the purity of $H^*(X)$ we see that the $E$-polynomial of $X$ will determine the Poincar\'e polynomial and hence we can deduce the Betti numbers from the motive, if we are able to compute it.

Let us now apply the previous discussion to $\mathcal{M}^{r,d}$. As we already observed, the limits as $\lambda \rightarrow 0$ of the $\C^*$ action applied to any point of $\mathcal{M}^{r,d}$ exist. This is a consequence of the equivariance and properness of the Hitchin map. Also, since the fixed point components of the $\C^*$-action are closed subvarieties of the fiber $(h^{r,d})^{-1}(\{0\})$ which is proper, the fixed point locus is proper. Finally, if $r$ and $d$ are coprime, then $\mathcal{M}^{r,d}$ is smooth.\\
Given an ordered partition $\underline{r}$ of $r$ and an ordered partition $\underline{d}$ of $d$, we can denote by $F_{\underline{r}}^{\underline{d}}$ the component of the fixed point locus in $\mathcal{M}^{r,d}$ whose points split as in proposition \ref{split}. Furthermore we can denote by $F_{\underline{r}}^{\underline{d},+}$ and by $F_{\underline{r}}^{\underline{d},-}$ the loci of $\mathcal{M}^{r,d}$ of those $(E,\phi)$ for which $\lim_{\lambda \rightarrow 0}\lambda \cdot (E,\phi) \in F_{\underline{r}}^{\underline{d}}$ and $\lim_{\lambda \rightarrow \infty} \lambda \cdot (E,\phi) \in F_{\underline{r}}^{\underline{d}}$ respectively. The notations given here are consistent with \cite{heinloth2014intersection}.

The previous discussion on smooth semiprojective varieties then implies the following.
\begin{prop}
Let $(r,d)=1$. Then:
\begin{itemize}
\item The fixed point components of the $\C^*$-action on $\mathcal{M}^{r,d}$ are the $F_{\underline{r}}^{\underline{d}}$ (when they are non-empty), which are all contained in the fiber $(h^{r,d})^{-1}(\{0\})$ of the Hitchin map.
\item The limit maps $F_{\underline{r}}^{\underline{d},+} \rightarrow F_{\underline{r}}^{\underline{d}}$ and $F_{\underline{r}}^{\underline{d},+} \rightarrow F_{\underline{r}}^{\underline{d}}$ are Zariski locally trivial fibrations whose fibers are isomorphic to affine spaces.
\item $\mathcal{M}^{r,d}=\bigcup F_{\underline{r}}^{\underline{d},+}$ and $(h^{r,d})^{-1}(\{0\})=\bigcup F_{\underline{r}}^{\underline{d},-}$. In particular the closure of the $F_{\underline{r}}^{\underline{d},-}$ are the irreducible components of $(h^{r,d})^{-1}(\{0\})$ and $(h^{r,d})^{-1}(\{0\})$ is the core of $\mathcal{M}^{r,d}$.
\item The dimension of the fibers of $F_{\underline{r}}^{\underline{d},+} \rightarrow F_{\underline{r}}^{\underline{d}}$ is constant and equal to $1+r^2(g-1)$.
\item There is a motivic equality $\left [\mathcal{M}^{r,d} \right ]=\L^{1+r^2(g-1)} \cdot \sum \left [ F_{\underline{r}}^{\underline{d}}\right ]$.
\end{itemize}
\begin{proof}
Follows from the previous discussion. A proof of this is also found in \cite[Proposition 2.1]{garcia2011motives}.
\end{proof}
\end{prop}
Let us briefly recall how to derive the structure of the attracting sets $F_{\underline{r}}^{\underline{d},+}$ in $\mathcal{M}^{r,d}$. Fix $r,d$ coprime. Let us consider a fixed point $(E,\phi)\in F_{\underline{r}}^{\underline{d}}$. If $r=(r)$ then $\phi=0$ and $E$ is a stable vector bundle. Since $\dim H^0(K \End E)=1/2 \dim \mathcal{M}^{r,d}$ then we know that the tangent space at $(E,\phi)$ decomposes into two parts of the same dimension, one of which is acted on with weight 0 by $\C^*$ and the other one with positive weight.

In all other cases, when $E=E_1 \oplus \dots \oplus E_m$ is split and has the form above. By the definition of hypercohomology we see that $T_{(E,\phi)} \mathcal{M}^{r,d}$ is generated by the kernel of the map:
\begin{align*}
C^1(\End E) \oplus C^0(K\End E) &\rightarrow C^1(K \End E)\\
(\tau, \nu) & \mapsto [\tau,\phi]+d \nu.
\end{align*}

The action of $\C^*$ induced on $T_{(E,\phi)} \mathcal{M}^{r,d}$ is then given by the pullback with respect to $f$, where $f: E \rightarrow E$ is the (one of the) automorphism of $E$ giving the isomorphism $(E,\lambda \phi) \cong (E,\phi)$. More explicitly if $E=E_1 \oplus \dots \oplus E_m$ and $\phi(E_i) \subseteq E_{i-1} \otimes K$ then $f$ is diagonal and defined by $f_{E_i}=\lambda^i \id_{E_i}$. Therefore we see that:
\begin{align*}
\lambda(\tau,0)=\lambda^{j-i}(\tau,0) &\text{ if } \tau \in C^1(E_i^* E_j)\\
\lambda(0,\nu)=\lambda^{j-i+1}(0,\nu) &\text{ if } \nu \in C^0(E_i^* E_j K).
\end{align*}

An easy way to remember the exponents above is to imagine that $\C^*$ acts with weight $i$ on $E_i$ and with weight $1$ on $K$. We can then split the deformation complex:
\begin{equation*}
\End E \rightarrow K \End E
\end{equation*}
into parts of positive, zero and negative weights according to the exponents above.

The splitting is as follows:
\begin{equation*}
\bigoplus_{i=1}^{m} \Hom(E_i,E_i) \rightarrow \bigoplus_{i=2}^{m} \Hom(E_i,E_{i-1} \otimes K)
\end{equation*}
gives the 0 weight part,
\begin{equation*}
\bigoplus_{i<j} \Hom(E_i,E_j) \rightarrow \bigoplus_{i\leq j} \Hom(E_i,E_{j} \otimes K)
\end{equation*}
gives the positive weight part and
\begin{equation*}
\bigoplus_{i>j}^{m} \Hom(E_i,E_j) \rightarrow \bigoplus_{i>j+1}^{m} \Hom(E_i,E_{j} \otimes K)
\end{equation*}
gives the negative weight part.
We can, for example, compute the positive weight part of the tangent space at $(E,\phi)$ by writing the long exact sequence associated to the positive weight complex:
\begin{align*}
0 &\rightarrow \bigoplus_{i<j} \Hom(E_i,E_j) \rightarrow \bigoplus_{i\leq j} \Hom(E_i,E_{j} \otimes K) \rightarrow T^+_{(E,\phi)} \mathcal{M}^{r,d} \rightarrow\\
&\rightarrow \bigoplus_{i<j} \Ext^1(E_i,E_j) \rightarrow \bigoplus_{i\leq j} \Ext^1(E_i,E_{j} \otimes K) \rightarrow \C \rightarrow 0.
\end{align*}
From this we can compute the dimension of the positive weight part using the formula:
\begin{equation*}
\dim T^+_{(E,\phi)} \mathcal{M}^{r,d} =1+ \sum_{i \leq j} \chi( \Hom(E_i,E_{j} \otimes K) ) - \sum_{i<j}  \chi( \Hom(E_i,E_{j}))=1+(g-1)r^2.
\end{equation*}

Here we see that, as it should be, the dimension of the positive weight part is constant, i.e. it does not depend on $\underline{r}$ and $\underline{d}$, and is always half of the dimension of the moduli space.

Since $\mathcal{M}^{r,d}$ is smooth, \cite[Theorem 4.1]{bialynicki1973some} implies that $F_{\underline{r}}^{\underline{d},+} \rightarrow F_{\underline{r}}^{\underline{d}}$ is a Zariski locally trivial affine fibration whose fibers have dimension equal to the dimension of the positive weight part of the $\C^*$ induced on the tangent space at one of the fixed points, which is what we computed with the above complexes.
\section{Character varieties and P = W}
\label{charvar}
For all details in the forthcoming discussion we refer to \cite{de2010topology} and \cite{hausel2008mixed}. Define the following variety.

\begin{defn}
Let $C$ be a smooth projective curve of genus $g \geq 2$. Let $\mathcal{M}^{2}_B$ be the variety defined by the following GIT quotient:
\begin{equation*}
\mathcal{M}^{2}_B = \left \{ A_1,\dots, A_g, B_1,\dots, B_g \in \GL_2(\C) : [A_1, B_1] \dots [A_g, B_g]=-\id\right \} \git \GL_2(\C)
\end{equation*}
where $\GL_2(\C)$ acts by conjugation. Such a variety is known as (rank 2) \emph{twisted character variety}.
\end{defn}
Note that $\mathcal{M}_B^2$ parametrizes isomorphism classes of representations of the fundamental group of $C$ into $\GL_2(\C)$ and is a purely topological invariant of $C$, i.e. it does not depend on the complex structure on the curve but only the genus. Also, $\mathcal{M}^2_B$ is a smooth affine variety since it is obtained from a GIT quotient of an affine space.

From the non-Abelian Hodge theorem (see for example \cite[proposition 10]{simpson1990nonabelian}) follows that $\mathcal{M}^{2,1}$ is canonically diffeomorphic to $\mathcal{M}^2_B$.

The cohomology of $\mathcal{M}^2_B$ is not pure and carries a nontrivial weight filtration. Let $Y$ be a complex variety. Recall that, if we denote by $F^\bullet$ the decreasing Hodge filtration on $H^*(Y)$ and by $W_\bullet$ the increasing weight filtration on $H^*(Y)$, we can define the mixed Hodge numbers:
\begin{equation*}
h^{p,q,j}(Y)=\dim (Gr_p^F Gr_{p+q}^W H^j(Y))
\end{equation*}
and the Hodge polynomial:
\begin{equation*}
H(Y,x,y,t)=\sum_{p,q,j} h^{p,q,j}(Y)x^p y^q t^j.
\end{equation*}

From \cite[corollary 4.1.11]{hausel2008mixed} we see that $\mathcal{M}^2_B$ has the property that $h^{p,q,j}(\mathcal{M}^2_B)=0$ if $p \neq q$. In particular $Gr_{2i+1}^W H^*(\mathcal{M}^2_B)=0$. We can define a variant of the mixed Hodge polynomial:
\begin{equation*}
H(\mathcal{M}^2_B,q,t)=H(\mathcal{M}^2_B,\sqrt{q},\sqrt{q},t).
\end{equation*}

On the other hand $\mathcal{M}^{2,1}$ has pure cohomology, so carries trivial weight filtration, but the Hitchin map $h^{2,1}: \mathcal{M}^{2,1} \rightarrow \mathcal{A}^2$ allows to define a \emph{perverse filtration} (see \cite[section 1.4]{de2010topology}). 

Consider the object $\R h^{2,1}_* \Q \in D_c^b(\mathcal{A}^2)$. Similarly to the case of the perverse filtration on the cohomology of the compactified Jacobian, $\R h^{2,1}_* \Q$ is filtered by its perverse truncations $\leftidx{^p}{\tau}{^{\leq i}}\R h^{2,1}_* \Q$ and the perverse filtration consists of the images:
\begin{equation*}
P_i H^*(\mathcal{M}^{2,1})= \im \left ( H^*\left (\leftidx{^p}{\tau}{^{\leq i}}\R h^{2,1}_* \Q \right) \rightarrow H^*\left (\R h^{2,1}_* \Q \right) \right ).
\end{equation*}

The following theorem is central:
\begin{theorem}[{\cite[theorem 1.1.1]{de2010topology}}]
The diffeomorphism given by non-abelian Hodge theory $\mathcal{M}^{2,1} \rightarrow\mathcal{M}^2_B$ induces an isomorphism in cohomology $H^*(\mathcal{M}^{2,1}) \rightarrow H^*(\mathcal{M}^2_B)$ with the property that:
\begin{equation*}
W_{2k}H^*(\mathcal{M}^2_B)=W_{2k+1}H^*(\mathcal{M}^2_B)=P_k H^*(\mathcal{M}^{2,1}).
\end{equation*}
\end{theorem}

The previous result is also known as "P = W", from the name of the filtrations involved. As of now, it is proved only for rank 2, while for higher rank it remains a conjecture. If we define the perverse Hodge polynomial by:
\begin{equation*}
PH(\mathcal{M}^{2,1},q,t)=\sum_{k,j} \dim P_k H^j(\mathcal{M}^{2,1}) q^k t^j
\end{equation*}
then the P = W result implies:
\begin{cor}
\begin{equation*}
H(\mathcal{M}^2_B,q,t)  = PH(\mathcal{M}^{2,1},q,t).
\end{equation*}
\end{cor}

In \cite[Theorem 1.1.3]{hausel2008mixed} we find the computation of the mixed Hodge polynomial of $\mathcal{M}_B^2$:
\begin{align*}
H(\mathcal{M}_B^2,q,t)&=
( 1+qt )^{2g} \left ({\frac { ( 1+q^2t^3)^{2g}}{ ( 1-q^2t^2 )( 1-q^2t^4) }}+{\frac {q^{2g-2}t^{4g-4}( 1+q^2t)^{2g}}{ ( 1-q^2)( 1-q^2t^2) }} \right )+\\
&-\frac{1}{2} q^{2g-2} t^{4g-4}( 1+qt)^{2g}  \left ( {\frac {( 1+qt )^{2g}}{( 1-qt^2 )( 1-q ) }}+ \frac{( 1-qt )^{2g}}{( 1+qt^2)( 1+q ) }  \right ).
\end{align*}
which, in turn, is equal to $PH(\mathcal{M}^{2,1},q,t)$.
\section{Higgs bundles with poles}
\label{hbpoles}
There is a variant of Higgs bundles that is also of interest for us. Choose once and for all a point $P \in C$.
\begin{defn}
\label{gammahiggs}
A pair $(E,\phi)$ is a \emph{Higgs $\gamma$-bundle} for some integer $\gamma \geq 0$ if $E$ is a vector bundle on $C$ and $\phi: E \rightarrow E \otimes K(\gamma P)$ is a twisted endomorphism of $E$.
\end{defn}

The definition is very similar to that of a Higgs bundle, and actually is the same when $\gamma=0$. Stability is the defined in the same way as for Higgs bundles and we denote by $\mathcal{M}^{r,d}(\gamma)$ the moduli space of S-equivalence classes of semistable Higgs $\gamma$-bundles of rank $r$ and degree $d$. Higgs $\gamma$-bundles can be thought as Higgs bundles whose Higgs field is allowed to have a pole of order at most $\gamma$ at $P$. The choice of putting $\gamma$ in parentheses instead of at the subscript will be apparent when we introduce Bradlow-Higgs triples. The construction of Higgs $\gamma$-bundles (and more general twists) can be found in \cite{nitsure1991moduli}.

Most of the properties of the moduli space of Higgs bundles carry over to Higgs $\gamma$-bundles. In particular we have semiprojectivity, smoothness for $r$ and $d$ coprime, the fixed points of the $\C^*$-action have the same form. There is a slight difference in the BNR correspondence, since Higgs $\gamma$-bundles will now correspond to sheaves on the total space of $K(\gamma P)$. This difference also carries to the Hitchin base $\mathcal{A}^{r,d}(\gamma)= H^0(K(\gamma P)) \oplus H^0(K^2(2\gamma P)) \oplus \dots \oplus H^0(K^r(r\gamma P))$ and to the Hitchin map $h^{r,d}_\gamma$.

It is worth observing that since we have a canonical map $\mathcal{O} \rightarrow \mathcal{O}(P)$ we also have embeddings $\mathcal{M}^{r,d}(\gamma) \rightarrow \mathcal{M}^{r,d}(\gamma+1)$. This approach, for $r=2$, is carried over in \cite[Chapter 7]{hausel2001geometry}. There we find the definition of $ \mathcal{M}^{2,d}(\infty)$ as the direct limit of the increasing embeddings $\mathcal{M}^{2,d}(\gamma) \rightarrow \mathcal{M}^{2,d}(\gamma+1)$.

The Poincar\'{e} polynomial of $\mathcal{M}^{2,d}(\infty)$ is then computed in \cite[Section 7.2]{hausel2001geometry}:
\begin{equation*}
P(\mathcal{M}^{2,d}(\infty),t) =\lim_{\gamma \rightarrow \infty} P(\mathcal{M}^{2,d}(\gamma) ,t)=\frac{(1+t)^{2g}(1+t^3)^{2g}}{(1-t^2)(1-t^4)}=P(B\overline{\mathcal{G}},t)
\end{equation*}
where $B\overline{\mathcal{G}}$ is the classifying space of a gauge group $\mathcal{G}$ modulo scalars. In the same section it is also proved that $\mathcal{M}^{2,d}(\infty)$ and $B\overline{\mathcal{G}}$ are homotopically equivalent, but we will not attempt the same with Bradlow-Higgs triples.
\section{Moduli spaces of Bradlow pairs and wall crossing}
\label{Bradlow}
Let us first define Bradlow pairs on a smooth projective curve $C$ over $\C$. The main reference for rank 2 Bradlow pairs will be \cite{thaddeus1994stable}. In \cite{mozgovoy2013moduli} it is possible to find the computation of the motives of the moduli spaces of Bradlow pairs for arbitrary rank and degree. Further comments on this approach will be made in the last chapter.
\begin{defn}[Bradlow pairs and stability]
A \emph{Bradlow pair} $(E,s)$ is a pair consisting of a vector bundle $E$ on $C$ and a nonzero section $s \in H^0(C,E)$.

Let $\sigma >0$ be a real number. We say $(E,s)$ is \emph{$\sigma$-(semi)stable} if for every $F \subset E$ proper subbundle we have:
\begin{align*}
&\frac{\deg(F)+\sigma}{\rk(F)}\stless \frac{\deg(E)+\sigma}{\rk(E)} &\text{if } s \in H^0(F)\\
&\frac{\deg(F)}{\rk(F)} \stless\frac{\deg(E)+\sigma}{\rk(E)} &\text{if } s \notin H^0(F).\\
\end{align*}
\end{defn}
There are several other equivalent notions of Bradlow pairs and stability. Here we outline some:
\begin{itemize}
\item we say $(E,s)$ is $\tau$-(semi)stable if for all $F \subseteq E$ subbundles we have:
\begin{align*}
&\mu(F) \stless \tau & \\
&\mu(E/F) \stgreat \tau &\text{if } s \in H^0(F).\\
\end{align*}
This is an equivalent notion of stability ($\tau = (d+\sigma)/r$) mentioned in \cite{mozgovoy2013moduli}.
\item A Bradlow pair can be equivalently defined as the datum of a map $s: \mathcal{O} \rightarrow E$. This interpretation falls into a more general framework of triples $\mathcal{E}=(f:E_0 \rightarrow E_1)$. There is a notion of $\gamma$-slope for this kind of triples:
\begin{equation*}
\mu_\gamma(\mathcal{E})=\frac{\deg E_0 +\deg E_1+\gamma \rk E_0}{\rk E_0+\rk E_1}
\end{equation*}
and $\mathcal{E}$ is called $\gamma$-(semi)stable if for all subtriples $\mathcal{F} \subset \mathcal{E}$ we have 
\begin{equation*}
\mu_\gamma(\mathcal{F}) \stless \mu_\gamma(\mathcal{E}).
\end{equation*}
This notion of stability is equivalent to the original one via $\gamma= \sigma+ (d+\sigma)/r$. This stability is mentioned in \cite{mozgovoy2013moduli}.
\end{itemize}
If we fix a $\sigma>0$, a rank $r$ and a degree $d$ we can define the moduli space of S-equivalence classes of $\sigma$-(semi)stable Bradlow pairs $M_\sigma^{r,d}$.

Note that for $d < 0$ or for $\sigma > d/(r-1)$ the moduli space is empty since the line bundle generated by the section $s$ is destabilizing. In the interval $[0,d/(r-1)]$ there is a finite set of values of $\sigma$ for which the stability inequalities can possibly become equalities. We call these values \emph{critical values}. For $\sigma$ different from the critical values, semistability and stability coincide and $M_\sigma^{r,d}$ is a smooth projective variety of dimension $d+(r^2-r)(g-1)$. See e.g. \cite{thaddeus1994stable} for the deformation theory of the rank 2 case.

It is also clear that, for $\sigma$ lying in the interval between two consecutive critical values, the inequalities for $\sigma$-stability will all be equivalent and therefore $M_\sigma^{r,d}$ only depends on the interval and not on the specific value of $\sigma$.

As we already noticed for $\sigma>d/(r-1)$ the moduli space of Bradlow pairs is empty. Near the other extremal value, meaning for $\sigma$ very close to $0$, $\sigma$-stability of a Bradlow pair $(E,s)$ implies semistability of the underlying vector bundle $E$. Therefore we have an \emph{Abel-Jacobi map}:
\begin{align*}
AJ: M_\sigma^{r,d} & \rightarrow N^{r,d}\\
(E,s) &\mapsto E
\end{align*}

In this notation we avoid as much as possible putting too many decorations on $AJ$ and it should be clear from the context what ranks and degrees we are referring to. Note that if $r$ and $d$ are coprime and $\sigma$ is small, then $\sigma$-stability implies stability of $E$ and $AJ^{-1}(E)=\P H^0(E)$. In particular, in this case, $AJ$ is a projective bundle for $d > r(2g-2)$.

For small $\sigma$ and $r$, $d$ not coprime the fibers of $AJ$ are more complicated. This is due to the fact that if $(E,s)$ is $\sigma$-stable, then it is possible that $E$ is strictly semistable and therefore $AJ$ will not only forget about the section, but will lose extra extension information about $E$ keeping only its S-equivalence class.

Let us discuss the role of critical values. First we introduce a notation that will be used throughout the whole thesis.
\begin{defn}
Fix a critical value $\bar \sigma$. Denote by $\bar\sigma_+$ any value of $\sigma$ bigger than $\bar\sigma$ but smaller than the consecutive critical value. Similarly $\bar\sigma_-$ will denote a value of $\sigma$ that is smaller than $\bar\sigma$ but bigger than the previous critical value.
\end{defn}
If $\sigma$ crosses $\bar\sigma$ and is increasing, then the first inequality for the stability becomes stronger, while the second one will become weaker. This means that as $\sigma$ crosses $\bar\sigma$, some pairs that were $\bar\sigma_-$-stable will become $\bar\sigma_+$-unstable and viceversa.

Let us see the rank 2 case, which is also the most relevant for the present thesis. The main reference here is of course \cite{thaddeus1994stable} even though some of the notations might vary. For rank 2 it is easy to see that the critical values are all the integers in $[0,d]$ with the same parity as $d$. Let us pick one of the critical values $\bar\sigma$.

There are two families $\P W^{d,+}_{\bar\sigma}$ and $\P W^{d,-}_{\bar\sigma}$ that parametrize the so called \emph{flip loci}, i.e. where $M^{2,d}_{\bar\sigma_-}$ and $M^{2,d}_{\bar\sigma_+}$ differ. More precisely, $\P W^{d,+}_{\bar\sigma}$ parametrizes pairs $(E,s)$ that are $\bar\sigma_-$-stable but not $\bar\sigma_+$-stable. This is equivalent to say that $E$ is a nonsplit extension
\begin{equation*}
0 \rightarrow L \rightarrow E \rightarrow M \rightarrow 0
\end{equation*}
with $L$ and $M$ two line bundles satisfying $\deg L =(d-\bar\sigma)/2$ and $\deg M =(d+\bar\sigma)/2$. Furthermore, $s \in H^0(L)$. Here $L$ is the canonical $\bar\sigma_+$-destabilizing subbundle of $(E,s)$. The condition of the extension being non-split is equivalent to the pair being $\bar\sigma_-$-stable.\\
There is a map:
\begin{align*}
\P W^{d,+}_{\bar\sigma} & \rightarrow S^{(d-\bar\sigma)/2} ( C) \times J^{(d+\bar\sigma)/2} ( C) \\
(E,s) & \mapsto (L,s, M)
\end{align*}
sending $(E,s)$ to the divisor of $s \in H^0(L)$ (which also retrieves $L$ itself) and to $M$, the quotient of $E$ by $L$. The fibers of the previous map are the projectivized extension spaces $\Ext^1(M,L)$. It is easy to check that the previous map is actually a projective bundle of rank $\bar\sigma+g-2$. In particular if we want to compute the motive of $\P W^{d,+}_{\bar\sigma}$ we get:
\begin{equation*}
[\P W^{d,+}_{\bar\sigma}]=[S^{(d-\bar\sigma)/2} ( C)] [J^{(d+\bar\sigma)/2} ( C)][\C\P^{\bar\sigma+g-2}].
\end{equation*}
$\P W^{d,-}_{\bar\sigma}$ instead parametrizes pairs $(E,s)$ that are $\bar\sigma_+$-stable but not $\bar\sigma_-$-stable. This is equivalent to say that $E$ is an extension
\begin{equation*}
0 \rightarrow M \rightarrow E \rightarrow L \rightarrow 0
\end{equation*}
where $L$ and $M$ are two line bundles with $\deg L =(d+\bar\sigma)/2$, $\deg M =(d-\bar\sigma)/2$ and $s \notin H^0(M)$. Note that, in this case, $s \notin H^0(M)$ is equivalent to $p(s)=\bar{s} \in H^0(L)$ being nonzero, where $p: E \rightarrow L$ is the projection. Once again, there is a map:
\begin{align*}
\P W^{d,-}_{\bar\sigma} & \rightarrow S^{(d-\bar\sigma)/2} ( C) \times J^{(d+\bar\sigma)/2} ( C) \\
(E,s) & \mapsto (L,\bar{s}, M)
\end{align*}
but this time $M$ is the canonical $\bar\sigma_-$-destabilizing subobject of $(E,s)$ and the divisor in $S^{(d-\bar\sigma)/2} ( C)$ is the divisor $D$ of the projection of $s$ on the canonical quotient $E/M$. The fibers of this map are the projectivization of the vector spaces $H^0(M \mathcal{O}_D)$ that parametrize extensions as above plus a lift of the section $\bar{s} \in H^0(L)$ to $s \in H^0(E)$. Again, this map is a projective bundle of rank $(d-\bar\sigma)/2 -1$. From this follows the relation:
\begin{equation*}
[\P W^{d,-}_{\bar\sigma}]=[S^{(d-\bar\sigma)/2} ( C)] [J^{(d+\bar\sigma)/2} ( C)][\C\P^{(d-\bar\sigma)/2-1}].
\end{equation*}

Observe that we have the motivic relation
\begin{equation*}
[M^{2,d}_{\bar\sigma_-}] - [\P W^{d,+}_{\bar\sigma}]=[M^{2,d}_{\bar\sigma_+}] - [\P W^{d,-}_{\bar\sigma}]
\end{equation*}
since outside the two flip loci, the two moduli spaces are isomorphic.

The previous phenomenon, involving a varying stability condition together with sets of critical values whose crossing makes the moduli problem change geometry, is known as \emph{wall crossing}.
\section{Hilbert schemes of points}
Here we introduce a geometric object that will be of interest in the last chapter of the thesis, namely the Hilbert scheme of points. A general introduction can be found in \cite{gottsche1994hilbert} and \cite{grothendieck1960techniques}, for more specific results about Hilbert schemes of points on a surface we refer to \cite{nakajima1999lectures}, while for a recollection on the properties of Hilbert schemes of singular curves and their relative version we refer to \cite{migliorini2011support}, \cite{maulik2014macdonald} and \cite{migliorini2015support}.

The most general definition of the Hilbert scheme, due to Grothendieck, passes through the functor of points. Since we don't need this much generality we just quickly sketch it.
Let $X$ be a scheme and $P$ a polynomial with integer coefficients. $\hilb^P(X)$ is a scheme (note that the existence of such a scheme is already a theorem by Grothendieck) such that maps $U \rightarrow \hilb^P(X)$ from another scheme $U$ are canonically identified with closed subschemes $Z \subset X \times U$ that are flat over $U$ and such that for every $u \in U$ the fiber $Z_u$ has Hilbert polynomial equal to $P$. In other words $\hilb^P(X)$ is the moduli space of subschemes of $X$ whose Hilbert polynomial is $P$.

For the purpose of some of the work in this thesis we just need to outline some of the properties of $\hilb^P(X)$ when $P$ is a constant polynomial and $X$ is either a smooth surface or a projective curve over $\C$.

When $P=n$ is a constant, $\hilb^P(X)$ is also commonly denoted by $X^{[n]}$. In this case the Hilbert scheme parametrizes subschemes of $X$ which are $0$ dimensional and whose space of global sections has dimension $n$. We can think of these subschemes as $n$-tuples of points and the Hilbert scheme in this case is called \emph{Hilbert scheme of points}.

For the case of $X$ a smooth surface we refer to \cite{nakajima1999lectures} for the proof of the main results and to the many works of G\"{o}ttsche \cite{gottsche1990betti,gottsche2000motive,gottsche1994hilbert} and G\"{o}ttsche, Soergel \cite{gottsche1993perverse} for the results on the cohomology.

When $X$ is a smooth surface then $X^{[n]}$ is a smooth variety of dimension $2n$ and is projective if $X$ itself is projective. Let $X^{(n)}$ or $\Sym^n X$ denote the $n$-fold product $X^n$ divided by the action of the symmetric group of the set of cardinality $n$ permuting the factors. There exists a map, called \emph{Hilbert-Chow morphism}:
\begin{equation*}
\pi_X^{n}: X^{[n]} \rightarrow X^{(n)}
\end{equation*}
sending a subscheme of length $n$ in $X$ to its support. This map is an isomorphism over the locus of $X^{(n)}$ where the points are all distinct, but has more complicated fibers when some points coincide. The fiber $(\pi_X^{n})^{-1}\{(P, \dots , P)\}$ for $P \in X$ is called the \emph{punctual Hilbert scheme of $n$-points} and denoted by $X^{[n]}_P$ or $\hilb_p^n(X)$, among the most common. It parametrizes subschemes of length $n$ whose support lies entirely in $P$. In other words it parametrizes the possible structures of subscheme of length $n$ that can be given to the point $P \in X$. The fibers of the Hilbert-Chow map are in general products of punctual Hilbert schemes. The fundamental property of the Hilbert-Chow morphism is that it is actually a resolution of the singular space $X^{(n)}$.

For the case of a projective curve we refer to \cite{migliorini2011support,maulik2014macdonald} for the integral case and to \cite{migliorini2015support} for the reduced case.

If $X$ is a projective locally planar curve then $X^{[n]}$ is a projective scheme over $\C$. The locally planar property for $X$ insures that it can be locally embedded in a surface and this implies some regularity for $X^{[n]}$. In general Hilbert schemes of points of varieties of dimension at least $3$ can be extremely complicated.

For the case of a curve $X$ we still have the Hilbert-Chow morphism $\pi_X^{n}: X^{[n]}_{red} \rightarrow X^{(n)}$, that is only defined on the reduction of $X^{[n]}$ which might be non-reduced if $X$ is non-reduced. If $X$ is smooth, then the Hilbert scheme of $n$-points is isomorphic to the symmetric power $X^{(n)}$ and the Hilbert-Chow morphism is an isomorphism.

If $X$ is integral then $X^{[n]}$ is integral of dimension $d$ and locally complete intersection. If we have a flat family of projective integral curves $f: \mathcal{X} \rightarrow B$ then there exists a variety $\mathcal{X}^{[n]}$ together with a proper map $f^{[n]}: \mathcal{X}^{[n]} \rightarrow B$ such that $(f^{[n]})^{-1}(b)=(f^{-1}(b))^{[n]}$ for all $b \in B$. $\mathcal{X}^{[n]}$ is called the \emph{Hilbert scheme of $n$-points relative to $B$}. If $X$ is singular, in general $X^{[n]}$ can be singular as well and the Hilbert-Chow morphism could have very complicated fibers. In \cite{ran2004note}, for example, it is proved that for the simplest possible planar singularity, i.e. the simple node $xy=0$, the punctual Hilbert scheme is a chain of $\C\P^1$'s touching at a point.

For reduced curves $X$ the Hilbert scheme behaves similarly to the case of integral curves but certainly has more than one component if the curve $X$ does. Not much is known about the Hilbert scheme of points of a non-reduced curve.
\section{The CKS complex}
\label{cks}
A more precise reference for this section is \cite[section 3.4]{migliorini2015support}. Let us assume that $\pi: \mathcal{C} \rightarrow B$ is a locally versal family of curves whose central fiber is the curve $\mathcal{C}_b$ which is nodal and reduced. We are not interested in the definition of locally versal family, as will become clear soon, just assume that it is a nice enough family.

The goal of this section is to introduce a tool to compute the stalk of
\begin{equation*}
IC(\bigwedge^i \R^1 \pi_{sm*} \Q)
\end{equation*}
where $\pi_{sm}$ is the restriction of $\pi$ to the locus of $B$ where the curves are smooth.

To the nodal curve $\mathcal{C}_b$ we can associate some combinatorial data, according to \cite[section 3]{migliorini2015support}. Denote by $\Gamma$ the graph of the curve, i.e. the set of vertices $V$ is determined by the components of $\mathcal{C}_b$ and the set of edges $E$ is formed by adding an edge between two vertices (not necessarily distinct) if the corresponding components intersect in a node. We will also denote by $\V$ and $\E$ the vector spaces spanned by the vertices and edges respectively.

Let us denote by $\mathcal{C}_\eta$ the versal deformation of our nodal curve. Then on $H^1(\mathcal{C}_\eta,\Q)$ is defined a monodromy weight filtration $W$ with the following identifications:
\begin{align*}
&Gr_0^W H^1(\mathcal{C}_\eta,\Q)=H^1(\Gamma)\\
&Gr_1^W H^1(\mathcal{C}_\eta,\Q)=H^1(\mathcal{C}_b^{\nu})\\
&Gr_2^W H^1(\mathcal{C}_\eta,\Q)=H_1(\Gamma) \otimes \L
\end{align*}
where $H^1(\Gamma)$ and $H_1(\Gamma)$ are the cohomology and homology of the graph $\Gamma$, $\mathcal{C}_b^{\nu}$ is the normalization of our nodal curve and $\L=\Q[-2](-1)$ as a Hodge structure. Note that the name $\L$ is no accident because it is in fact the image of $\L \in K_0(Var_\C)$ under the additive invariant taking a variety to its class in the Grothendieck ring of Hodge structures.

For each edge $e \in E$ we define a linear map:
\begin{equation*}
N_e : H^1(\mathcal{C}_\eta,\Q) \rightarrow H^1(\mathcal{C}_\eta,\Q) \otimes \L.
\end{equation*}

It is defined as the composition:
\begin{equation*}
H^1(\mathcal{C}_\eta,\Q) \rightarrow Gr_2^W H^1(\mathcal{C}_\eta,\Q) \rightarrow Gr_0^W H^1(\mathcal{C}_\eta,\Q)\otimes \L \rightarrow H^1(\mathcal{C}_\eta,\Q) \otimes \L
\end{equation*}
where the first map is the canonical projection, the last map is the canonical inclusion. The middle map instead is defined by knowing that $Gr_2^W H^1(\mathcal{C}_\eta,\Q)=H_1(\Gamma) \otimes \L$ and it will embed into $\E$, while $Gr_0^W H^1(\mathcal{C}_\eta,\Q)=H^1(\Gamma)$ and it is a quotient of the dual vector space $\E^*$. Therefore the middle map is defined by the composition:
\begin{equation*}
H_1(\Gamma) \rightarrow \E \rightarrow \E^* \rightarrow H^1(\Gamma)
\end{equation*}
where $\E \rightarrow \E^*$ is the natural duality map $t \mapsto \langle t, e^*\rangle e^*$. Here $e^*$ is the dual element of the edge $e$.

The operators $N_e$ extend to operators $N_e^{(i)}$ in a natural way:
\begin{align*}
N_e^{(i)}: \bigwedge^i H^1(\mathcal{C}_\eta,\Q) &\rightarrow \bigwedge^i H^1(\mathcal{C}_\eta,\Q) \otimes \L\\
c_1 \wedge \dots \wedge c_i & \mapsto \sum_{k=1}^i c_1 \wedge \dots \wedge N_e(c_k) \wedge \dots \wedge c_i.
\end{align*}

For a subset of edges $I \subseteq E$ we write $N_I^{(i)}$ to denote the iterated composition of the $N_e^{(i)}$ for $e \in I$. Note that this composition does depend on the order (in terms of a sign) but the image $\im N_I^{(i)} \subset H^1(\mathcal{C}_\eta,\Q) \otimes \L^{|I|}$ does not.

We are now finally ready to define the Cattani-Kaplan-Schmid complex associated to $\mathcal{C}_b$.
\begin{defn}
For a fixed $i$ define the following complex $CKS^i$:
\begin{equation*}
0 \rightarrow \bigwedge^i H^1(\mathcal{C}_\eta,\Q) \rightarrow \bigoplus_{\substack{I \subseteq E\\ |I|=1}} \im N_I^{(i)} \rightarrow \bigoplus_{\substack{I \subseteq E\\ |I|=2}} \im N_I^{(i)} \rightarrow \dots
\end{equation*}
the maps are defined by considering an ordering $e_1, e_2, \dots$ of the edges and then for $J \subseteq E$ with $|J|=j$:
\begin{align*}
\im N_J^{(i)} &\rightarrow \bigoplus_{\substack{I \subseteq E\\ |I|=j+1}} \im N_I^{(i)}\\
c &\mapsto (c_I)_{|I|=j+1}
\end{align*}
where $c_I=0$ if $J \nsubseteq I$ and if $I = J \cup \{e\}$ then $c_I= (-1)^a N_e^{(i)}(c)$. Here the sign is determined by ordering $J$ according to the order of $E$, then appending $e$ at the end of $J$ and counting the number of swaps necessary to order $J \cup \{e\}$.
\end{defn}

We have the following proposition.
\begin{prop}
The $i$-th CKS complex for the curve $\mathcal{C}_b$ is quasi isomorphic to the stalk of $IC(\bigwedge^i \R^1 \pi_{sm*} \Q)$ at $\mathcal{C}_b$.
\begin{proof}
See \cite[section 3.4]{migliorini2015support}.
\end{proof}
\end{prop}

As is also explained in \cite[section 3.4]{migliorini2015support} one might want to compute the weight polynomial of the CKS complexes, instead of trying to compute their cohomology. In this case it is enough to know the dimension of the $\im N_I^{(i)}$. Here we prove a lemma that we will use in one of the next chapters.
\begin{lemma}
\label{ratcurve}
Denote by $\Sigma$ a curve with two rational components meeting in $2g-2$ simple nodes and by $\overline\Sigma$ an integral curve with $2g-2$ simple nodes whose normalization is isomorphic to $\C\P^1$. Define:
\begin{equation*}
U(\Sigma)=\sum_{n = 0}^{2g-3} q^n [CKS^n [-n]]
\end{equation*}
and
\begin{equation*}
U(\overline\Sigma)=\sum_{n= 0}^{2g-3} q^n [\overline{CKS}^n [-n]]
\end{equation*}
where $[CKS^n]$ is the weight polynomial of the $n$-th CKS complex associated to $\Sigma$ and $[\overline{CKS}^n]$ is the analogous object for $\overline\Sigma$.
Then:
\begin{equation*}
U(\overline\Sigma)=(1-q\Q)(1-q\L)U(\Sigma) \qquad \text{mod } q^{2g-2}.
\end{equation*}
\begin{proof}
The dual graph $\Gamma$ of $\Sigma$ consists of two vertices and $2g-2$ edges connecting them, while the dual graph $\overline \Gamma$ of $\overline\Sigma$ consists of one vertex and $2g-2$ loops. For the rest of the proof we will crucially use the fact that $n$ is at most $2g-3$ and therefore there is no way we can disconnect the graphs of $\Sigma$ and $\overline\Sigma$ by removing $n$ edges. Choose $0 \leq n \leq 2g-3$, then:
\begin{equation*}
[CKS^n [-n]]=(-1)^n \sum_{i=0}^n (-1)^i \sum_{\substack{I \subseteq e(\Gamma)\\|I|=i}} \left [ \bigwedge^{n-i} \left( H^1(\Gamma \setminus I) \oplus H_1(\Gamma \setminus I) \otimes \L\right )\right ]
\end{equation*}
and
\begin{equation*}
[\overline{CKS}^n [-n]]=(-1)^n \sum_{i=0}^n (-1)^i \sum_{\substack{I \subseteq e(\overline\Gamma)\\|I|=i}} \left [ \bigwedge^{n-i} \left( H^1(\overline\Gamma \setminus I) \oplus H_1(\overline\Gamma \setminus I) \otimes \L\right )\right ].
\end{equation*}

Note that $H^1(\overline\Gamma \setminus I) = H^1(\Gamma \setminus I) \oplus \Q$ and $H_1(\overline\Gamma \setminus I) \otimes \L=H_1(\Gamma \setminus I) \otimes \L \oplus \L$. Therefore:
\begin{align*}
[\overline{CKS}^n [-n]]&=(-1)^n \sum_{i=0}^n (-1)^i \sum_{\substack{I \subseteq e(\overline\Gamma)\\|I|=i}} \left [ \bigwedge^{n-i} \left( H^1(\overline\Gamma \setminus I) \oplus H_1(\overline\Gamma \setminus I) \otimes \L\right )\right ]=\\
&=(-1)^n \sum_{i=0}^n (-1)^i \sum_{\substack{I \subseteq e(\Gamma)\\|I|=i}} \left [ \bigwedge^{n-i} \left( H^1(\Gamma \setminus I) \oplus H_1(\Gamma \setminus I) \otimes \L \oplus \Q \oplus \L \right )\right ]=\\
&=(-1)^n \sum_{i=0}^n (-1)^i \sum_{\substack{I \subseteq e(\Gamma)\\|I|=i}} \left [ \bigwedge^{n-i} \left( H^1(\Gamma \setminus I) \oplus H_1(\Gamma \setminus I) \otimes \L \right )\right ]+\\
&+(-1)^n \sum_{i=0}^n (-1)^i \sum_{\substack{I \subseteq e(\Gamma)\\|I|=i}} \left [ \bigwedge^{n-i-1} \left( H^1(\Gamma \setminus I) \oplus H_1(\Gamma \setminus I) \otimes \L \right )\otimes (\Q \oplus \L) \right ]+\\
&+(-1)^n \sum_{i=0}^n (-1)^i \sum_{\substack{I \subseteq e(\Gamma)\\|I|=i}} \left [ \bigwedge^{n-i-2} \left( H^1(\Gamma \setminus I) \oplus H_1(\Gamma \setminus I) \otimes \L \right )\otimes \L \right ]=\\
&=[CKS^n [-n]]-[CKS^{n-1} [-n+1]]\cdot (\Q+\L)+[CKS^{n-2} [-n+2]] \cdot \L.
\end{align*}
From this follows:
\begin{equation*}
U(\overline\Sigma)=(1-q\Q)(1-q\L)U(\Sigma) \qquad \text{mod } q^{2g-2}.
\end{equation*}
\end{proof}
\end{lemma}

%% file: BHT.tex
\chapter{Bradlow-Higgs triples}
%
%\begin{quote}
%{\small
%Summary of chapter contents\\
%\begin{itemize}
%\item  Intro triples (definition, stability, moduli spaces, BNR, extremal stability conditions)
%\item Properties of moduli of triples (Deformation theory, discussion of the singular points, $\C^*$-action, fixed points, attracting sets and properness of the Hitchin map)
%\end{itemize}
%}
%\end{quote}
%
\section{Definitions, properties and basic results}
Consider a smooth projective curve $C$ over $\C$ of genus $g\geq 2$ and the total space $X=T^*C$ of $K$, the canonical bundle of $C$. Let us also denote by $\P X$ the projectivization of $X$ obtained by adding the divisor at infinity i.e. $\P X= \P(\mathcal{O} \oplus K)$.

We want to study moduli spaces of objects that are a variant of both Higgs bundles and Bradlow pairs.
\begin{defn}[Bradlow-Higgs triples and $\sigma$-stability]
\label{sigmast}
A triple $(E,\phi,s)$ is said to be a \emph{Bradlow-Higgs triple} if $(E,\phi)$ is a Higgs bundle and $s \in H^0(C,E)$ is a nonzero section of the underlying vector bundle $E$.

Let $\sigma$ be a positive real number. We say $(E,\phi,s)$ is \emph{$\sigma$-(semi)stable} if, for all proper $\phi$-invariant subbundles $F \subset E$, we have:
\begin{align*}
\frac{\deg(F)}{\rk(F)} \stless \frac{\deg(E)+\sigma}{\rk(E)} \qquad & \text{if } s \notin H^0(C, F)\\
\frac{\deg(F)+\sigma}{\rk(F)} \stless \frac{\deg(E)+\sigma}{\rk(E)} \qquad & \text{if } s \in H^0(C, F)
\end{align*}

A morphism of Bradlow-Higgs triples is a morphism of the underlying Higgs bundles that sends the section of the source to a scalar multiple of the section of the target.
\end{defn}
The idea is to introduce a section of the underlying vector bundle so that the stability condition is allowed to vary together with the positive real parameter $\sigma$. This will produce more complicated moduli spaces but will also allow to study the relation among them as $\sigma$ varies. As we will see later (and as already seen about Bradlow pairs) there are two interesting extremal stability conditions, one of which relates the moduli space of Bradlow-Higgs triples with the moduli space of semistable Higgs bundles, and the other one relates the moduli space of Bradlow-Higgs triples with a relative Hilbert scheme of points.

Here is a basic but fundamental result on Bradlow-Higgs triples.
\begin{prop}
Let $\sigma>0$ and $(E,\phi,s)$ be a $\sigma$-stable triple. Let $f \in H^0(\End E)$ such that $[f, \phi]=0$ and $f\cdot s=0$. Then $f=0$.\\
In particular, the only endomorphisms of a $\sigma$-stable triple $(E,\phi,s)$ are scalar multiples of the identity.
\begin{proof}
Let $L$ denote the kernel of $f$ and suppose that $L \neq E$. Note that $L \neq 0$ because $s \in H^0(L)$. Then $L$ has rank $r' < r=\rk E$ and is $\phi$-invariant containing the section $s$, therefore we have $\mu(L) \leq d/r+\sigma(1/r-1/r')$. In particular $f$ induces an injective map $E/L \rightarrow E$ so that $E/L$ can be thought as a sub sheaf of $E$. Since it is $\phi$ invariant, its saturation will be $\phi$ invariant as well and $\mu( E/L) \geq d/r + \sigma/r$. This means that the saturation of $E/L$ will be $\sigma$-destabilizing unless $f=0$.

The second part follows by noting that if $f$ is an endomorphism of $(E,\phi,s)$ and $f(s)=\alpha s$ then $f-\alpha \id_E$ satisfies the hypotheses of the first part of the proposition and therefore is $0$.
\end{proof}
\end{prop}
\subsection{BNR correspondence and moduli space of coherent systems on a surface}
\label{secbnr}
Recall from section \ref{higgs} that due to the BNR correspondence, pure one dimensional sheaves on $X$ that do not intersect the divisor at infinity correspond to Higgs bundles on $C$ through the pushforward along the projection map $\pi: T^*C \rightarrow C$. The pushforward also canonically identifies the sections of such sheaves with the sections of their pushforward, i.e. the underlying vector bundle of the Higgs bundle obtained by applying $\pi_*$. We can therefore alternatively work with pairs $(\mathcal{F},s)$ on $X$, where $\mathcal{F}$ is a rank one pure one dimensional sheaf of degree $d'$ whose support $\supp (\mathcal{F})$ does not intersect the divisor at infinity and $s$ is a non-zero section in $H^0(X,\mathcal{F})$.

These objects are known as a particular case of (one dimensional) \emph{coherent systems} on the surface $\P X$. The moduli spaces of these objects were defined and studied in \cite{le1995faisceaux} (and in a slightly different context in \cite{pandharipande2009curve}). It is worth recollecting some of the results we need.

First we recall the definition of the stability condition \cite[d\'{e}finition 4.4]{le1995faisceaux} and \cite[section 1.1]{pandharipande2009curve}.
\begin{defn}
Fix a polarization $L$ on $\P X$ and denote by $\mathcal{F}(k)=\mathcal{F} \otimes L^k$ for $\mathcal{F}$ a sheaf on $\P X$. Let $\sigma$ be a positive real number, we say $(\mathcal{F},s)$ is \emph{$\sigma$-(semi)stable} if:
\begin{align*}
\frac{\chi(\mathcal{G}(k))}{a(\mathcal{G})} \stless \frac{\chi(\mathcal{F}(k))+\sigma}{a(\mathcal{F})} \qquad & \text{if } s \notin H^0(\P X, \mathcal{G})\\
\frac{\chi(\mathcal{G}(k))+\sigma}{a(\mathcal{G})} \stless \frac{\chi(\mathcal{F}(k))+\sigma}{a(\mathcal{F})} \qquad & \text{if } s \in H^0(\P X, \mathcal{G})
\end{align*}
for all proper subsheaves $\mathcal{G} \subset \mathcal{F}$. Both conditions on the Hilbert polynomials are meant to hold for big values of $k$.
Note that, since $\mathcal{F}$ is pure one dimensional, both Hilbert polynomials of $\mathcal{F}$ and $\mathcal{G}$ have degree 1. Here we denote by $a(\cdot)$ the corresponding leading terms and by $b(\cdot)$ the constant terms.
\end{defn}
Observe that, as explained in \cite[section 1.1]{pandharipande2009curve}, $\chi(\mathcal{F}(k))=a(\mathcal{F})k+b(\mathcal{F})$ with:
\begin{align*}
&a(\mathcal{F})=\int_\beta c_1(L)\\
&b(\mathcal{F})=\chi(\mathcal{F})
\end{align*}
with $\beta \in H_2(\P X)$ the class of the support of $\mathcal{F}$.\\
After this we can define moduli spaces of $\sigma$-(semi)stable coherent systems:
\begin{defn}
Fix $m$ an integer and the homology class $\beta=n[C] \in H_2(\P X)$ of the support for our sheaves (here $C$ embeds in $\P X$ as the zero section of $T^*C$).\\
We denote by $\text{Syst}_{\P X,\sigma}(\beta,m,1)$ the moduli space of $\sigma$-(semi)stable pure one dimensional rank one degree $d'$ pairs $(\mathcal{F},s)$ with $\chi(\mathcal{F})=m$ and support in class $\beta$ (compare with \cite[th\'{e}or\`eme 4.11]{le1995faisceaux}).
\end{defn}
The $1$ in the definition of $\text{Syst}_{\P X,\sigma}(\beta,m,1)$ denotes that we only have one section $s$ of the underlying sheaf.\\
From \cite[th\'{e}or\`eme 4.11]{le1995faisceaux} we deduce that the $\text{Syst}_{\P X,\sigma}(\beta,m,1)$ are projective varieties, since $\P X$ is projective.

In order to obtain Bradlow-Higgs triples from pairs on $\P X$ we need the extra condition that the support of the underlying sheaf does not intersect the divisor at infinity. Once we identify the stability condition for coherent systems with the stability condition for Bradlow-Higgs triples we can identify the moduli spaces of $\sigma$-(semi)stable Bradlow-Higgs triples with the open subset of $\text{Syst}_{\P X,\sigma}(\beta,m,1)$ of coherent systems whose support does not intersect the divisor at infinity of $\P X$.
\begin{prop}
Let $(\mathcal{F},s)$ be a coherent system (with one section) such that $\supp(\mathcal{F})$ does not intersect the divisor at infinity. Assume $\chi(\mathcal{F})=m$ and the class of the support of $\mathcal{F}$ is $\beta=n[C]$, where $C$ is the zero section of $\P X$. Denote by $(E,\phi,s)$ the Bradlow-Higgs triple obtained by pushing forward $(\mathcal{F},s) $ along $\P X \rightarrow C$. Then $\deg E=m+n(g-1)$ and $\rk E=n$. Furthermore, the stability condition for $(\mathcal{F},s)$ is the same as the stability condition for $(E,\phi,s)$, in particular this means that the stability condition does not depend on the polarization on $\P X$ as long as the hypotheses on the support are satisfied.
\begin{proof}
Since $\supp(\mathcal{F})$ does not intersect the divisor at infinity it is a branched $n$ to $1$ cover of $C$ embedded as the zero section of $\P X$. More precisely there exists an affine open subset $U=\Spec A$ of $C$ over which $T^*C$ trivializes and $T^*C_{|U}=\Spec A[t]$. Then $\supp(\mathcal{F}) \cap T^*C_{|U}$ is defined by a degree $n$ polynomial in $A[t]$. Therefore if $\mathcal{F}$ restricts to the rank one $A[t]$-module $M$ on $\Spec A[t]$, its pushforward is $M$ seen as an $A$-module and therefore it will be isomorphic to $M \oplus M t \oplus \dots \oplus M t^{n-1}$ so that it will have rank $n$.\\
Under our hypotheses:
\begin{align*}
&a(\mathcal{F})=\int_\beta c_1(L)=n \int_\beta c_1(L) >0\\
&b(\mathcal{F})=\chi(\mathcal{F})=m
\end{align*}
therefore $\deg E+n(1-g)=\chi(E)=\chi(\mathcal{F})=m$. Let $A=\int_\beta c_1(L)>0$.\\
For the stability condition let $\mathcal{G}$ be a subsheaf of $\mathcal{F}$ and call $G$ its pushforward, which is a $\phi$-invariant subbundle of rank $0<\rk G < n$. Note that the class of the support of $\mathcal{G}$ is then $\rk G \cdot A$. Therefore we have:
\begin{align*}
&\frac{\chi(\mathcal{G}(k))}{a(\mathcal{G})}=k+\frac{\deg G}{\rk G \cdot A}+(1-g),\\
&\frac{\chi(\mathcal{G}(k))+\sigma}{a(\mathcal{G})}=k+\frac{\deg G+ \sigma}{\rk G \cdot A}+(1-g),\\
&\frac{\chi(\mathcal{F}(k))}{a(\mathcal{F})}=k+\frac{\deg E}{\rk E \cdot A}+(1-g) \text{ and}\\
&\frac{\chi(\mathcal{F}(k))+\sigma}{a(\mathcal{F})}=k+\frac{\deg E+ \sigma}{\rk E \cdot A}+(1-g).
\end{align*}

Thus we see that the inequalities for stability are equivalent.
\end{proof}
\end{prop}
Therefore we can give the following definition, using the identification above we know that such a moduli space actually exists and is also a quasi-projective variety.
\begin{defn}
Fix integers $r \geq 1$ and $d$. We denote by $\mathcal{M}_\sigma^{r,d}$ the moduli space of S-equivalence classes of $\sigma$-(semi)stable Bradlow Higgs triples of degree $d$ and rank $r$.
\end{defn}
As we saw $\mathcal{M}_\sigma^{r,d}$ corresponds to an open subset of $\text{Syst}_{\P X,\sigma}(r[C],d+r(1-g),1)$. We conclude with an
\begin{obs}
\label{degs}
The sheaves in $\text{Syst}_{\P X,\sigma}(r[C],d+r(1-g),1)$ can also be regarded as sheaves that are rank one and degree $d'$ on their support. So if their support is a smooth projective curve of genus $g'$ we have $d+r(1-g)=\chi(E)=\chi(\mathcal{F}_{\supp \mathcal{F}})=d'+1-g'$. If we further observe that for the curves above the genus is $g'=1+r^2(g-1)$ (section \ref{higgs} of the introductory chapter) we have that $d'=d+(r^2-r)(g-1)$.
\end{obs}
%As for the moduli spaces of Bradlow pairs, there are a few values of $\sigma$ for which it is possible to have equality in the stability condition. These are once again called \emph{critical values}. For Bradlow-Higgs triples it is still true that we have a wall-crossing phenomenon and that the moduli spaces $\mathcal{M}_\sigma^{r,d}$ are constant for $\sigma$ lying in the intervals between two critical values.
\subsection{Hitchin maps}
\label{hitchinmaps}
An analogy with the moduli space of stable Higgs bundles is the existence, for every $\sigma$, of a Hitchin map, defined in the analogous way as for Higgs bundles.
\begin{defn}
Let $\mathcal{A}^r$ be the Hitchin base. Then for every $\sigma>0$ we can define the \emph{Hitchin map}:
\begin{align*}
\chi^{r,d}_\sigma: \mathcal{M}_\sigma^{r,d} & \rightarrow \mathcal{A}^r\\
(E,\phi,s) & \mapsto \text{char poly}(\phi).
\end{align*}
\end{defn}
As we already observed for Higgs bundles, the Hitchin map is the analogue, through the BNR correspondence, of the map taking $(\mathcal{F},s)$ to the equation of $\supp(\mathcal{F})$. Let us first make the following observation.
\begin{prop}
\label{opendense}
Let $\mathcal{A}^r_{int} \subset \mathcal{A}^r$ be the locus of integral spectral curves. Then the open subset $(\chi^{r,d}_\sigma)^{-1} (\mathcal{A}^r_{int})$ does not depend on $\sigma$.
\begin{proof}
This follows from what we already observed for Higgs bundles. If $(E,\phi,s) \in \mathcal{M}_\sigma^{r,d}$ and $E$ has a $\phi$-invariant subbundle, then the characteristic polynomial of $\phi$ will factor, and hence $\chi^{r,d}_\sigma(E,\phi,s) \in  \mathcal{A}^r \setminus \mathcal{A}^r_{int}$. It follows that $(\chi^{r,d}_\sigma)^{-1} (\mathcal{A}^r_{int})$ will consist of triples for which the Higgs bundle does not have any subobject and therefore are automatically $\sigma$-stable for all $\sigma$.
\end{proof}
\end{prop}
\subsection{Extremal stability conditions}
\label{extrstab}
We are interested in the stability for extremal values of $\sigma$ and for which combinations of $r ,d , \sigma$ the moduli spaces are non-empty.\\
Let us first define a filtration that will be useful in the following chapters.
\begin{theorem}
Let $(E,\phi,s)$ be a Bradlow-Higgs triple. Then there exists a filtration $0 \subset U_1 \subset \dots \subset U_l \subseteq E$ with $l \leq r$ and $U_r=E$ in case equality holds. The filtration satisfies the following conditions:
\begin{itemize}
\item[(i)] $\rk U_i=i$,
\item[(ii)] $\phi(U_i) \subseteq U_{i+1} \otimes K$ for $1 \leq i <l$,
\item[(iii)] for $1\leq i < l$, $\phi$ induces a nonzero map $U_i/U_{i-1} \rightarrow U_{i+1}/U_i \otimes K$,
\item[(iv)] $U_l$ is the smallest subbundle of $E$ that is $\phi$-invariant and contains the section $s$.
\end{itemize}
Furthermore, we have a bound $\deg U_i \geq i(i-1)(1-g)$.
\begin{proof}
We define the filtration with an inductive process. In the following we will denote by $\langle \cdot \rangle$ the saturation of a subsheaf of a vector bundle.

First define $U_1$ to be the line subbundle of $E$ generated by the section $s$, then clearly $\deg U_1 \geq 0$ so the bound is satisfied and also if $\phi(U_1) \subseteq U_1\otimes K$ then we can stop the construction and the filtration $0 \subset U_1 \subset E$ satisfies the properties.

If $\phi(U_1) \nsubseteq U_1 \otimes K$, $\phi$ induces a nonzero map $U_1 \rightarrow E/U_1 \otimes K$. Take:
\begin{equation*}
S_2=\left \langle \im( \phi: U_1 \rightarrow E/U_1 \otimes K) \right\rangle \otimes K^{-1}.
\end{equation*}

Then $S_2$ is a subbundle of $E/U_1$ and hence there exists a unique subbundle $U_2$ of $E$ containing $U_1$ and such that $S_2$ is the image of $U_2$ with respect to the projection $E \rightarrow E/U_1$. By construction we have a bound $\deg U_2 \geq \deg U_1 +2-2g \geq 2-2g$ and also if $\phi(U_2) \subseteq U_2 \otimes K$ then the filtration $0 \subset U_1 \subset U_2 \subset E$ satisfies the properties.

The construction can be carried through inductively. Assume that we have constructed $0 \subset U_1 \subset U_2 \subset \dots \subset U_p$ satisfying $\rk U_i=i$, $\deg U_i \geq i(i-1)(1-g)$ and that there are nonzero maps $U_i/U_{i-1} \rightarrow U_{i+1}/U_i \otimes K$ for $1 \leq i < p$. Then we claim that $\phi(U_p) \nsubseteq U_p \otimes K$ if and only if $\phi^p(s) \notin U_p \otimes K^p$. To prove this observe that $\phi(U_p) \subseteq U_p \otimes K$ if and only if $\phi$ induces the zero map $U_p \rightarrow E/U_p \otimes K$. Since by construction $\phi(U_{p-1}) \subseteq U_p \otimes K$ then $U_p \rightarrow E/U_p \otimes K$ is zero if and only if $U_p/U_{p-1} \rightarrow E/U_p \otimes K$. Consider the iterated applications of $\phi$:
\begin{equation*}
U_1 \rightarrow U_2/U_1 \otimes K \rightarrow \dots \rightarrow U_p/U_{p-1} \otimes K^{p-1} \rightarrow E/U_p \otimes K^p.
\end{equation*}

From this we see that since the composition of the first $p-1$ applications is nonzero and all the quotients are rank one except for the last, then $U_p/U_{p-1} \otimes K^{p-1} \rightarrow E/U_p \otimes K^p$ is zero if and only if $U_1 \rightarrow E/U_p \otimes K^p$ is zero and this is equivalent to $\phi^p(s) \in U_p \otimes K^p$.

Now that we proved the claim, assume $\phi(U_p) \nsubseteq U_p \otimes K$. Then this is equivalent to $\phi^p(s) \notin U_p \otimes K^p$ and hence consider:
\begin{equation*}
S_{p+1}=\left \langle \im( \phi^p: U_1 \rightarrow E/U_p \otimes K) \right\rangle \otimes K^{-p}
\end{equation*}
which is a subbundle of $E/U_p$. Then there exists a unique subbundle $U_{p+1}$ of $E$ containing $U_p$ and such that $S_{p+1}$ is the image of $U_{p+1}$ with respect to the projection $E \rightarrow E/U_p$. Then $U_{p+1}$ also satisfies the recursive hypotheses and the construction can proceed.

Note that this sequence of steps will terminate either when we get $U_r$ that has rank $r$ and hence $U_r=E$ or when we reach a $\phi$-invariant subbundle of $E$.

The statement about the minimality of $U_l$ is clear because it is generically generated by $s ,\phi(s),\dots, \phi^{p-1}(s)$.
\end{proof}
\end{theorem}
\begin{rmk}
Note that, by construction, the filtration above is unique for each triple $(E,\phi,s)$.
\end{rmk}
\begin{defn}[$U$-filtration]
Let $(E,\phi,s)$ be a triple. We will call \emph{$U$-filtration} the filtration whose existence was proved in the above theorem.\end{defn}
\begin{theorem}
\label{extremal}
\begin{itemize}
%\item[(a)] $\mathcal{M}_\sigma^{r,d}$ is non-empty if and only if $d \geq r(r-1)(g-1)$.\\
\item[(i)] Assume that $\sigma$ is very close to $0$ (i.e. smaller than the first critical value), then $\sigma$-stability for a triple $(E,\phi,s)$ implies the semistability of $(E,\phi)$ and so we have an Abel-Jacobi map:
\begin{equation*}
AJ: \mathcal{M}_\sigma^{r,d} \rightarrow \mathcal{M}^{r,d}.
\end{equation*}
For $d$ large enough (e.g. $d > r(2g-1)+(r-1)^2(2g-2) $) for any semistable Higgs bundle $(E,\phi)$ we have $H^1(E)=0$ and therefore $AJ$ is a projective bundle over the stable part of $\mathcal{M}^{r,d}$.\\
\item[(ii)] For $\sigma>(r-1)d+r(r-1)(r-2)(g-1)$ and a $\sigma$-stable triple $(E,\phi,s)$ corresponding to a pair $(\mathcal{F},s)$ the following three equivalent conditions are realized:
\begin{itemize}
\item there are no $\phi$-invariant subbundles of $E$ which contain the section
\item $s, \phi(s), \dots, \phi^{r-1}(s)$ generically generate $E$
\item $s$ as a map $\mathcal{O}_{\P X} \rightarrow \mathcal{F}$ has zero dimensional cokernel.
\end{itemize}
\end{itemize}
\begin{proof}
%The proof of (a) will be in the last chapter \marginpar{reference}.\\
For (i) it is proved in \cite[corollary 3.4]{nitsure1991moduli} that if $d > r(2g-1)+(r-1)^2(2g-2) $ then any semistable Higgs bundle of degree $d$ and rank $r$ has vanishing $H^1$. Therefore on the stable locus of $\mathcal{M}^{r,d}$ the fibers of $AJ$ are projective spaces of constant dimension and so $AJ$ is a projective bundle.

For (ii) let us first show that the three conditions are equivalent. The first and second are equivalent because the minimal $\phi$-invariant subbundle of $E$ that contains $s$ is $U_l$ constructed above and it is generated by $s, \phi(s), \dots, \phi^{l-1}(s)$. The second and the third are equivalent since the subsheaf of $E$ generated by $s, \phi(s), \dots, \phi^{r-1}(s)$ is the pushforward of the image of $\mathcal{O}_X \rightarrow \mathcal{F}$.

To prove that the three conditions hold for $\sigma>(r-1)d+r(r-1)(r-2)(g-1)$ we use the $U$-filtration of the triple. Saying that the three conditions are satisfied is the same as saying that the $U$-filtration is of full length $r=\rk E$. Therefore assume that the filtration stops at the step $l < r$. We have the bound $\deg U_l \geq l(l-1)(1-g)$. Therefore:
\begin{equation*}
(l-1)(1-g)+\frac{\sigma}{l} \leq \frac{\deg U_l+\sigma}{l}\leq \frac{d+\sigma}{r}
\end{equation*}
where the first inequality follows from the bound on $\deg U_l$, while the second follows from assuming that the triple is $\sigma$-(semi)stable. Rearranging we get:
\begin{equation*}
\sigma \leq d \frac{l}{r-l}+r(g-1)\frac{l(l-1)}{r-l}\leq d(r-1)+r(r-1)(r-2)(g-1),
\end{equation*}
where the second inequality follows from $l < r$. This contradicts the bound on $\sigma$.
\end{proof}
\end{theorem}
\begin{defn}
We will refer to $\sigma$-stability as $\varepsilon$-stability, when $\sigma$ is very close to zero, and as $\infty$-stability, when $\sigma$ is very large (i.e. larger than the bound in the last proposition).
\end{defn}

From the previous proposition we see that $\infty$-stability implies that the section $s$ is a cyclic vector for $\phi$. Observe that, as opposed to Bradlow pairs, crossing the last critical value for triples will not cause the moduli space to be empty. Rather, the moduli space is constant (and non-empty) after $\sigma=(r-1)d+r(r-1)(r-2)(g-1)$. The cyclic property of $s$ and $\phi$ will have important consequences later, when we study the fixed points of a $\C^*$-action on the moduli space of triples.

For the moment we conclude by observing that the condition of $s: \mathcal{O}_X \rightarrow \mathcal{F}$ having zero dimensional cokernel is the one mentioned by R. Pandharipande and R. Thomas in \cite{pandharipande2009curve}. The only difference is that our ambient variety is a surface, namely $X$ and in \cite{pandharipande2009curve} the ambient variety is a Calabi-Yau threefold and the pairs are used for curve counting purposes.\\
In the case of surfaces however, the moduli space of $\infty$-stable triples satisfies a nice property, following from \cite[Proposition B.8]{pandharipande2010stable} using the BNR correspondence.
\begin{theorem}
\label{relhilb}
The Hitchin map:
\begin{equation*}
\chi^{r,d}_\infty: \mathcal{M}_\infty^{r,d} \rightarrow \mathcal{A}^r
\end{equation*}
is the relative Hilbert scheme of $d+r(r-1)(g-1)$ points over the family of spectral curves $\mathcal{A}^r$.
\begin{proof}
Follows from \cite[Proposition B.8]{pandharipande2010stable}, the BNR correspondence and the relation between $d'$ and $d$ we found in observation \ref{degs}.
\end{proof}
\end{theorem}
\section{Deformation theory of Bradlow-Higgs triples}
\label{deform}
Using \v{C}ech cohomology we can compute the tangent space to $\mathcal{M}^{r,d}_\sigma$ at $(E,\phi,s)$. Note that a similar analysis has been carried out in \cite[proposition 4.12]{le1995faisceaux}, under the point of view of coherent systems of $\P X$.
\begin{theorem}
\label{deformations}
Let $(E,\phi,s)$ be a $\sigma$-stable Bradlow-Higgs triple, then the tangent space at $(E,\phi,s)$ is given by the first cohomology $\H^1(E,\phi,s)$ of the complex:
\begin{equation*}
C^0(\End E)\rightarrow C^1(\End E) \oplus C^0(K\End E) \oplus C^0(E) \rightarrow C^1(E) \oplus C^1(K \End E)
\end{equation*}
where the first map is $$p(k)=(dk, [k,\phi],k \cdot s)$$ and the second one is $$q(\tau, \nu, \gamma)=(\tau \cdot s+d \gamma, [\tau,\phi]+d \nu).$$
Furthermore, $\H^0(E,\phi,s)=0$.

The same result can be obtained from the hypercohomology of the complex:
\begin{align*}
\End E &\rightarrow K \End E \oplus E\\
f & \mapsto ([f,\phi],f(s))
\end{align*}
from which we can also deduce the long exact sequence
\begin{align}
0 \rightarrow H^0(\End E) \rightarrow H^0(K \End E\oplus E) \rightarrow T_{(E,\phi,s)} \rightarrow \\ \nonumber
\rightarrow H^1(\End E) \rightarrow H^1(K \End E\oplus E) \rightarrow \H^2 \rightarrow 0 
\end{align}
%and that, for high degree and rank 2, $\dim \mathcal{M}_\sigma^{n,d}=d+6(g-1)+1$.
\begin{proof}
If we want to understand maps $\Spec \C[\varepsilon]/(\varepsilon^2) \rightarrow \mathcal{M}^{r,d}_\sigma$ such that the image of the closed point is $(E,\phi,s)$ then we can choose an open cover $\{U_\alpha\}$ over which $E$ trivializes. So the open sets $\Spec \C[\varepsilon]/(\varepsilon^2) \times U_\alpha$ constitute a trivializing open cover for any deformation $(\tilde{E}, \tilde{\phi}, \tilde{s})$ of our triple.

We can write the transition functions for $E$ as $1+\varepsilon \tau_{\alpha \beta}$ for some $\tau \in C^1(\End E)$, the section can be expressed as $s+\varepsilon \gamma_\alpha$ for $\gamma \in C^0(E)$ and finally the Higgs field as $\phi+\varepsilon \nu_\alpha$ for $\nu \in C^0(K \otimes \End E)$.

If we write the compatibility conditions necessary to obtain a deformation we get:
\begin{align*}
&(1+\varepsilon \tau_{\alpha \beta})(s+\varepsilon \gamma_\beta)=s+\varepsilon \gamma_\alpha\\
&(1+\varepsilon \tau_{\alpha \beta})(\phi+\varepsilon \nu_\beta)(1+\varepsilon \tau_{\alpha \beta})^{-1}=\phi+\varepsilon \nu_\alpha.
\end{align*}
The equations then yield $\tau \cdot s+d \gamma=0$ and $[\tau,\phi]+d\nu=0$ which is precisely the condition $q(\tau,\nu, \gamma)=0$.

Since $\sigma$-stable triples have no automorphisms we can assume that isomorphic deformations only come from the change of trivializations on the open cover $U_\alpha$. In other words we obtain isomorphic deformations (whose cocycles will be distinguished by an apex) if and only if there exists a $k \in C^0(\End E)$ such that
\begin{align*}
&(1+\varepsilon k_\alpha)(1+\varepsilon \tau'_{\alpha \beta})(1-\varepsilon k_\alpha)=1+\varepsilon \tau_{\alpha \beta}\\
&(1+\varepsilon k_\alpha)(s+\varepsilon \gamma'_\alpha)=s+\varepsilon \gamma_\alpha\\
&(1+\varepsilon k_\alpha)(\phi+\varepsilon \nu'_\alpha)(1-\varepsilon k_\alpha)=\phi+\varepsilon \nu_\alpha.
\end{align*}

These three equations are equivalent to requiring that $(\tau'-\tau, \nu'-\nu, \gamma'-\gamma)$ lies in the image of $p$.

Our complex fits into an exact sequence of complexes:
\footnotesize
\begin{center}
\begin{tikzcd}
& 0 \arrow{r} \arrow{d} & C^0(\End E) \arrow{r}{=} \arrow{d}{p} & C^0(\End E)  \arrow{d} \\
& C^0(E) \arrow{r}{incl} \arrow{d}{d} & C^1(\End E) \oplus C^0(K\End E) \oplus C^0(E) \arrow{r}{proj} \arrow{d}{q} & C^1(\End E) \oplus C^0(K\End E)  \arrow{d} \\
& C^1(E) \arrow{r}{incl} & C^1(E) \oplus C^1(K\End E) \arrow{r}{proj} & C^1(K\End E) 
\end{tikzcd}
\end{center}
\normalsize
The right column is the complex associated with the deformations of the Higgs bundle $(E, \phi)$ and its cohomology will be denoted by $\H^i([\cdot, \phi])$.

The cohomology of the deformation complex for our moduli problem will be denoted by $\H^i(E,\phi,s)$ or simply $\H^i$ when no confusion can arise.

This yields an exact sequence of cohomology groups:
\begin{align}
\label{tan}
0 \rightarrow \H^0(E,\phi,s) \rightarrow \H^0([\cdot, \phi]) \rightarrow H^0(E) \rightarrow T_{(E,\phi,s)} \rightarrow \\ \nonumber
\rightarrow \H^1([\cdot,\phi]) \rightarrow H^1(E) \rightarrow \H^2(E,\phi,s) \rightarrow \H^2([\cdot, \phi]) \rightarrow 0.
\end{align}

Note that, of course, we obtain the same result by considering the long exact sequence of hypercohomologies associated to the exact sequence of complexes:
\begin{center}
\begin{tikzcd}
& 0\arrow{r} & 0 \arrow{r} \arrow{d} & \End E \arrow{r} \arrow{d} & \End E \arrow{r} \arrow{d} & 0\\
& 0\arrow{r} & E \arrow{r} & K\End E \oplus E \arrow{r} & K\End E \arrow{r} & 0\\
\end{tikzcd}
\end{center}
The right vertical arrow is the commutator with $\phi$, while the central vertical arrow is the map $f \mapsto ([f,\phi],f(s))$.

The map $\H^0([\cdot, \phi]) \rightarrow H^0(E)$ is evaluation of an endomorphism of $E$ on the section. $\H^0(E,\phi,s)$ is therefore the kernel of such map, i.e. $\H^0(E,\phi,s)$ consists of the endomorphisms of $E$ that commute with $\phi$ and that annihilate the section $s$ and we already know that, by $\sigma$-stability, these are zero. From this we deduce that $\H^0(E,\phi,s)=0$.
\end{proof}
\end{theorem}
Recall that the deformation complex for a Higgs bundle is self dual and in particular
$$\dim \H^2([\cdot, \phi])=\dim \H^0([\cdot, \phi]) \geq 1.$$

Since $\H^2$ surjects onto $\H^2([\cdot, \phi])$, we also have $\dim \H^2 \geq 1$.\\
This allows to get
$$\dim T_{(E,\phi,s)} \mathcal{M}^{r,d}_\sigma=d+(2r^2-r)(g-1)+\dim \H^0([\cdot,\phi]).$$

Note that stability of the underlying Higgs bundle is an open condition and therefore we know that $\dim \H^0([\cdot,\phi])$ is 1 on an open subset of $\mathcal{M}^{r,d}_\sigma$, as we already noted in proposition \ref{opendense}.

We can make a further step in the understanding of the deformations.
\begin{prop}
\label{h2}
Let $(E,\phi,s)$ be a $\sigma$-stable triple and let us denote by $\H^2$ the hypercohomology group in \ref{deformations}. Then $(\H^2)^*$ is the kernel of the following map:
\begin{align*}
H^0(\End E) \oplus H^0(K E^*) & \rightarrow H^0(K \End E)\\
(\alpha, \beta) &\mapsto [\alpha,\phi]+\beta \otimes s.
\end{align*}
In particular, $\dim \H^2 \geq \dim \H^2([\cdot, \phi])$.
\begin{proof}
It is enough to note that, by Serre's duality, $(\H^2)^*$ is the zeroth hypercohomology group of the dual of the deformation complex for $(E,\phi,s)$, i.e.:
\begin{align*}
\End E \oplus K E^* &\rightarrow K \End E\\
(\alpha, \beta) &\mapsto [\alpha,\phi]+\beta\otimes s
\end{align*}
and therefore $(\H^2)^*$ is the kernel of the previous map applied to the global sections.\\
It is worth noting that by $\beta \otimes s$ we mean the map $E \rightarrow E \otimes K$ obtained by tensoring $\beta: E \rightarrow K$ with $s: \mathcal{O} \rightarrow E$.\\
The last statement follows, for example, from the exact sequence \ref{tan}. It can also be deduced directly from the characterization of $(\H^2)^*$ which clearly contains all pairs $(\alpha, \beta)$ with $\alpha$ commuting with $\phi$ and $\beta=0$.
\end{proof}
\end{prop}
\begin{rmk}
\label{adhm}
From \cite[definition 2.1]{diaconescu2012moduli} we see that the data $(E,\phi,s,\alpha,\beta)$ with $(E,\phi,s)$ a Bradlow-Higgs triple, $\alpha \in H^0(\End E)$, $\beta \in H^0(K E^*)$ such that:
\begin{equation*}
[\alpha,\phi]+\beta \otimes s=0
\end{equation*}
define an \emph{ADHM sheaf} on the curve $C$. If $(E,\phi,s)\in \mathcal{M}_\sigma^{r,d}$ is a singular point, then there exist $\alpha$ and $\beta$ not both zero for which $(E,\phi,s,\alpha,\beta)$ is an ADHM sheaf.

For a triple $(E,\phi,s)$ we can deduce from \cite[lemma 2.5]{diaconescu2012moduli} that if $\dim \H^2 > \dim \H^2([\cdot,\phi])$ then there exists a subbundle of $E$ which is proper, $\phi$-invariant and contains $s$. Therefore the following holds:
\begin{itemize}
\item the triples $(E,\phi,s)$ for which the spectral curve is integral are always smooth points
\item if there are no $\phi$-invariant subbundles that contain the section then $\H^2(E,\phi,s)^*$ has the same dimension as the space of endomorphisms of $(E,\phi)$.
\end{itemize}
\end{rmk}
We can prove more precise results once we concentrate on specific combinations of $r,d$ and $\sigma$. Let us first state a technical lemma.
\begin{lemma}
\label{unstableext}
Let $E$ be a rank two vector bundle that is given by an extension:
\begin{equation*}
0 \rightarrow E_2 \rightarrow E \rightarrow E_1 \rightarrow 0
\end{equation*}
with $\deg E_2 > \deg E_1$. Denote by $i$ the inclusion $E_2 \rightarrow E$ and by $p$ the projection $E \rightarrow E_1$.\\
If the extension is nonsplit, we have an isomorphism:
\begin{align*}
\C \oplus H^0( E_1^*E_2) &\rightarrow H^0(\End E)\\
(\lambda, g) & \mapsto \lambda \cdot \id_E+igp.
\end{align*}

If, instead, the extension is split, then there is an isomorphism:
\begin{align*}
\C \oplus \C \oplus H^0( E_1^*E_2) &\rightarrow H^0(\End E)\\
(\lambda_1,\lambda_2, g) & \mapsto \begin{pmatrix} \lambda_1 & 0\\0& \lambda_2 \end{pmatrix}+igp.
\end{align*}
\begin{proof}
Let us first assume the extension is nonsplit. Then we have an exact sequence:
\begin{equation*}
0 \rightarrow \mathcal{O} \rightarrow E_2^* E \rightarrow E_2^*E_1 \rightarrow 0
\end{equation*}
and since $H^0(E_2^* E_1)=0$ we get $H^0(E_2^* E)=\C$. Consider then the exact sequence:
\begin{equation*}
0 \rightarrow E_1^*E_2 \rightarrow E_1^* E \rightarrow \mathcal{O} \rightarrow 0
\end{equation*}
where the second map is defined by taking $f: E_1 \rightarrow E$ and sending it to $pf : E_1 \rightarrow E_1$.

Note that the induced map on the global sections $H^0(E_1^* E) \rightarrow H^0(\mathcal{O})$ has to be zero because if not, we would have a map $f: E_1 \rightarrow E$ such that $p f = \id_{E_1}$ and this implies that the sequence is split. Therefore $H^0(E_1^* E_2) \cong H^0(E_1^* E)$.

Lastly, consider the exact sequence:
\begin{equation*}
0 \rightarrow E_1^*E \rightarrow \End E \rightarrow E_2^* E \rightarrow 0
\end{equation*}
where the second map is defined by sending $f : E \rightarrow E$ to $fi: E_2 \rightarrow E$. Note that the map $ H^0(\End E) \rightarrow H^0(E_2^* E) \cong \C$ is surjective because clearly it does not send $f=\id_E$ to zero. If we trace back the identifications we made we get the statement.\\
The proof for the split case is easier because $H^0(\End E)=H^0(\mathcal{O} \oplus \mathcal{O} \oplus E_1^*E_2 \oplus E_2^* E_1)=\C \oplus \C \oplus H^0(E_1^* E_2)$ and the statement follows immediately.
\end{proof}
\end{lemma}
\begin{cor}
\label{smoothr2}
Assume that $r=2$ and that $(E,\phi,s)$ is a triple whose underlying Higgs bundle is stable and of degree $d$. Then:
\begin{itemize}
\item[(i)] if $d<0$ then $(E,\phi,s)$ is a smooth point in $\mathcal{M}_\varepsilon^{2,d}$,
\item[(ii)] if $d>4g-4$ is odd then $(E,\phi,s)$ is a smooth point in $\mathcal{M}_\varepsilon^{2,d}$.
\end{itemize}
In particular $\mathcal{M}_\varepsilon^{2,d}$ is smooth for $d<0$ and for $d>4g-4$ odd.
\begin{proof}
Here we use the characterization in proposition \ref{h2}. Assume that $E$ is itself stable. Then $H^0(\End E)$ consists only of scalar multiples of the identity. Therefore $[\alpha,\phi]=0=\beta \otimes s$ and so $\beta =0$. So we deduce that $\dim \H^2=1$ which implies $(E,\phi,s)$ is a smooth point.

Assume that $E$ is unstable of rank 2. We have an extension:
\begin{equation*}
0 \rightarrow E_2 \rightarrow E \rightarrow E_1 \rightarrow 0
\end{equation*}
with $\deg E_2 > \deg E_1$ and denote by $i$ the inclusion and $p$ the projection. Here we are in the hypotheses of lemma \ref{unstableext}. Since $E_2$ is destabilizing but $(E,\phi)$ is stable, then $p \phi i$ cannot be zero. Suppose that the above extension is nonsplit. Pick $(\alpha, \beta) \in (\H^2)^*$. If $\alpha \in H^0(\End E)$ we know that $\alpha = \lambda \cdot \id_E+igp$ for $g \in H^0(E_1^*E_2^*)$. Therefore $$[\alpha, \phi]=igp \phi - \phi igp.$$

Since $$p [\alpha,\phi] i=pigp \phi i - p\phi igpi=0,$$ we know that
$$0=p(\beta \otimes s)=\beta i \otimes p(s)$$
and so either $\beta i=0$ or $p(s)=0$. In the first case:
$$igp\phi i=[\alpha, \phi] i = (\beta \otimes s) i =\beta i \otimes s=0$$
and since $p \phi i \neq 0$ we get $ig=0$ and hence $g=0$. In particular $[\alpha,\phi]=0=\beta \otimes s$ and hence $\alpha$ is scalar multiple of the identity while $\beta=0$. This implies $\dim \H^2=1$ and the triple is a smooth point.

If instead $p(s)=0$ then from
$$-p \phi igp=p[\alpha,\phi]=\beta \otimes p(s)=0$$
and since $p \phi i \neq 0$ we have $gp =0$ and hence $g=0$. So once again $\alpha$ is a multiple of the identity and $\beta=0$ so $\dim \H^2=1$ and we get again a smooth point.

We are therefore left with the case of the split extension. Write $s=s_1+s_2$ with $s_i \in H^0(E_i)$ and
$$\phi=\begin{pmatrix}
\phi_{11} & \phi_{21}\\
\phi_{12} & \phi_{22}
\end{pmatrix}.$$

In this case we must have $\phi_{21} \neq 0$ for the stability. Pick $(\alpha,\beta) \in (\H^2)^*$ and write
$$\alpha=\begin{pmatrix}
a & 0\\
c & d
\end{pmatrix},$$
$\beta= \beta_1+\beta_2$. By expanding $[\alpha,\phi]+\beta \otimes s$ we get the following relations:
\begin{align*}
&\phi_{21} c = \beta_1 \otimes s_1\\
&(d-a) \phi_{21} = \beta_2 \otimes s_1\\
&c \phi_{11}-\phi_{22} c+(d-a) \phi_{12} =- \beta_1 \otimes s_2\\
&-c \phi_{21} = \beta_2 \otimes s_2
\end{align*}
Note that if $d<0$ then $\deg E_1 <0$ and so $s_1=0$ so we deduce that $c=0$ and $a=d$ and hence $\alpha$ is a scalar multiple of the identity and $\beta=0$ so $\dim \H^2=1$.

If instead $d>4g-4$ then $d_2>2g-2$ and so $\beta_2: E_2 \rightarrow K$ has to be zero. From the above relations we deduce once again that $\alpha$ is a scalar multiple of the identity and $\beta=0$ so that $\dim \H^2=1$.

Assertions $(i)$ and $(ii)$ as well as the last statement follow from the previous remarks.
\end{proof}
\end{cor}
Let us conclude the section with some remarks. From theorem \ref{extremal} part (i) it follows that the Abel-Jacobi map $\mathcal{M}_\varepsilon^{r,d} \rightarrow \mathcal{M}^{r,d}$ is a projective bundle when $d$ and $r$ are coprime and $d >r(2g-1)+(r-1)^2(2g-2)$. In particular, under these hypotheses $\mathcal{M}_\varepsilon^{r,d}$ is smooth.

For rank 2 it easy to find a singular point in $\mathcal{M}_\varepsilon^{2,d}$ even when $d$ is odd and $0 \leq d \leq 4g-4$. We can build it as follows. Pick $d_1$ and $d_2$ integers satisfying $d_1 \geq 0$ and $d_2 \leq 2g-2$. Pick divisors $S$ and $B$ on $C$ of degrees $d_1$ and $2g-2-d_2$ respectively. Define $D=S+B$.

With these data we define a Bradlow-Higgs triple. Let $E_1$ be the line bundle associated to $S$ and $s_1$ a section of $E_1$ (unique up to scaling by a constant) whose divisor is $S$. Let $E_2$ be defined by $K(-B)$. Then $\deg E_1=d_1$ and $\deg E_2=d_2$. Define $E=E_1\oplus E_2$ and $\phi: E_2 \rightarrow E_1 \otimes K$ be the map associated to the divisor $D$. Lastly let $\beta_2: E_2 \rightarrow K$ be the map associated to the divisor $B$. Clearly we will have $\beta_2 \otimes s_1=\lambda\phi$ as maps $E_2 \rightarrow E_1 \otimes K$ for some constant $\lambda\neq 0$. Then we define $\alpha \in H^0(\End E)$ and $\beta: E \rightarrow K$ as:
\begin{align*}
\alpha &= \begin{pmatrix}
\lambda & 0\\
0 & 0
\end{pmatrix}\\
\beta &= \begin{pmatrix}
0 & \beta_2
\end{pmatrix}.
\end{align*}

It is immediate to check that $(E,\phi)$ is stable and $[\alpha,\phi]=\beta \otimes s$ and since $\alpha$ and $\beta$ are both nonzero we see that $\dim \H^2 > \dim \H^2([\cdot,\phi])$. It is also easy to verify that for this particular triple $\dim \H^2=2$.
We can also produce singular points for the other moduli spaces $\mathcal{M}_\sigma^{r,d}$.

The idea is to produce a split Higgs bundle such that the two summands are not too distant in degree and have both a nonzero section so that the resulting triple will be stable.

Let us build an example for $r=2$, $d$ odd and $\sigma>1$. Consider $E=L \oplus M$ with $\deg M=\deg L+1$, $\phi=diag(\alpha, \beta)$ with $\alpha, \beta \in H^0(K)$ distinct and nonzero. Finally take $s=u \oplus v$ with $u$ and $v$ nonzero sections of $L$ and $M$ respectively. Assume that $v \neq c(u)$ for all maps $c: L \rightarrow M$. Then $L$ and $M$ are clearly the only $\phi$-invariant sub bundles of $E$ and none of them contains the section $s$. If we have $\deg E=d=2k+1$ then $\deg L=k$ and $\deg M=k+1$ so we have:
\begin{align*}
\deg L < \frac{\deg E+\sigma}{2} \qquad \text{ and }\qquad \deg M < \frac{\deg E+\sigma}{2}
\end{align*}
provided that $\sigma>1$. This proves that under our assumption this triple is stable. Moreover any diagonal endomorphism of $E$ will commute with $\phi$ so that we have $\dim \H^2([\cdot, \phi])=\dim \H^0([\cdot, \phi]) \geq 2$ and so $\dim \H^2 \geq 2$ and the triple is a singular point for our moduli space.

Later, we will also construct explicitly fixed points for the $\C^*$-action that are singular in $\mathcal{M}_\sigma^{2,d}$ for $\sigma$ after the first wall.
\section{Properness of the Hitchin map}
In this brief section we will prove the following theorem.
\begin{theorem}
The Hitchin maps
\begin{equation*}
\chi^{r,d}_\sigma: \mathcal{M}_\sigma^{r,d} \rightarrow \mathcal{A}^r
\end{equation*}
are all proper.
\end{theorem}

In the subsequent sections we will see some consequences of this.

Recall from section \ref{secbnr} that $\mathcal{M}^{r,d}_\sigma$ is an open subset of $\text{Syst}_{\P X,\sigma}(r[C],d+r(1-g),1)$ where $\P X$ is the smooth projective surface obtained by adding a divisor at infinity to $T^* C$. We also have a map:
\begin{align*}
\mathcal{S}^{r,d}_\sigma: \text{Syst}_{\P X,\sigma}(r[C],d+r(1-g),1) &\rightarrow \mathcal{H}^r\\
(\mathcal{F},s) & \mapsto \supp{F}
\end{align*}
where $\mathcal{H}^r$ is a Hilbert scheme of curves in $\P X$ whose is Hilbert polynomial is determined by $r$. Observe that, even though we don't include it in the notation, both $\text{Syst}_{\P X,\sigma}(r[C],d+r(1-g),1)$ and $\mathcal{S}^{r,d}_\sigma$ depend on the choice of some polarization $L$ on $\P X$, as explained in section \ref{secbnr}.

It is clear that $\mathcal{S}^{r,d}_\sigma$ is a proper map. In fact, according to \cite[th\'{e}or\`eme 4.11]{le1995faisceaux}, $\text{Syst}_{\P X,\sigma}(r[C],d+r(1-g),1)$ is a projective variety and, according to \cite[theorem 1.1.2]{gottsche1994hilbert}, $\mathcal{H}^r$ is also projective. Using \cite[theorem 4.9]{hartshorne1977algebraic} we see that projective morphisms are proper therefore the composition
\begin{align*}
\text{Syst}_{\P X,\sigma}(r[C],d+r(1-g),1) \rightarrow \mathcal{H}^r \rightarrow \{pt\}
\end{align*}
is proper and $\mathcal{H}^r \rightarrow \{pt\}$ is proper and hence separated. According to \cite[corollary 4.8]{hartshorne1977algebraic}, $\mathcal{S}^{r,d}_\sigma$ is then proper.

Recall that $\mathcal{M}^{r,d}_\sigma$ is the locus of $\text{Syst}_{\P X,\sigma}(r[C],d+r(1-g),1)$ where the support of the sheaves does not intersect the divisor at infinity of $\P X$. As we observed in \ref{secbnr}, if we restrict to such a locus then the stability condition on $\text{Syst}_{\P X,\sigma}(r[C],d+r(1-g),1)$ does not actually depend on the choice of the polarization $L$ on $\P X$. Since the condition defining $\mathcal{M}^{r,d}_\sigma$ relies exclusively on the support of the underlying sheaves we see that the Hitchin map $\chi^{r,d}_\sigma$ is given by the following base change diagram:
\begin{center}
\begin{tikzcd}
& \mathcal{M}^{r,d}_\sigma \arrow{r}{\chi_\sigma^{r,d}} \arrow{d}& \mathcal{A}^r\arrow{d} \\
& \text{Syst} \arrow{r}{\mathcal{S}^{r,d}_\sigma} & \mathcal{H}^r
\end{tikzcd}
\end{center}
and, according to \cite[corollary 4.8]{hartshorne1977algebraic}, $\chi^{r,d}_\sigma$ is then a proper map.
\section{The $\C^*$-action}
As in the case of the moduli space of stable Higgs bundles, the moduli spaces $\mathcal{M}^{r,d}_\sigma$ admit a $\C^*$-action that scales the Higgs field.
\begin{defn}
We define an action:
\begin{align*}
\C^* \times \mathcal{M}_\sigma^{r,d} &\rightarrow \mathcal{M}_\sigma^{r,d}\\
(\lambda, (E,\phi,s)) &\mapsto \lambda \cdot (E,\phi,s)=(E,\lambda\phi,s).
\end{align*}
\end{defn}
As with the moduli space of Higgs bundles we can use the $\C^*$-action and the properness of the Hitchin map to decompose the $\mathcal{M}_\sigma^{r,d}$ into attracting sets. We will also examine the fixed point loci as $\sigma$ varies. Special attention is dedicated to the case of rank 2.

We also have the following theorem.
\begin{theorem}
For all $(E,\phi,s) \in \mathcal{M}_\sigma^{r,d}$ the limit
\begin{equation*}
\lim_{\lambda \rightarrow 0} \lambda \cdot (E,\phi,s)
\end{equation*}
exists in $\mathcal{M}_\sigma^{r,d}$.
\begin{proof}
We have proper maps
\begin{equation*}
\chi^{r,d}_\sigma: \mathcal{M}_\sigma^{r,d} \rightarrow \mathcal{A}^r
\end{equation*}
that are also equivariant with respect to the action of $\C^*$ on $\mathcal{A}^r$ we described in \ref{higgs}. The claimed existence of the limits follows from the valuative criterion for properness and the fact that $\C^*$ acts with positive weights on $\mathcal{A}^r$.
\end{proof}
\end{theorem}
\subsection{The fixed point locus}
Let us first understand which points in $\mathcal{M}_\sigma^{r,d}$ are fixed by the action. The proof of the following is completely analogous to the characterizations of the fixed points in $\mathcal{M}^{r,d}$
\begin{prop}
\label{fixedbht}
Suppose that $(E,\phi,s)$ is a $\sigma$-(semi)stable pair in $\mathcal{M}_\sigma^{r,d}$ such that $\lambda \cdot (E,\phi,s) \cong (E,\phi,s)$ for some $\lambda \in \C^*$ that is not a root of unity, then $E=E_1 \oplus \dots \oplus E_m$, $\phi(E_i) \subseteq E_{i-1} \otimes K$, $\phi(E_1)=0$ and $s \in E_i$ for some $1 \leq i \leq m$.
\begin{proof}
If we don't have the section the proof is the same as in \cite[Lemma 4.1]{simpson1992higgs}.\\
$\lambda \cdot (E,\phi,s) \cong (E,\phi,s)$ implies that there is an automorphism $f$ of $E$ such that $f(s)$ is a multiple of $s$ and $f \phi = \lambda \phi f$. Since the coefficients of the characteristic polynomial are sections of $\mathcal{O}_C$, we know that they are constants. Therefore also the eigenvalues of $f$ are constant and we have the decomposition of $E$ according to the generalized eigenspaces $E=\oplus_{\mu} E_\mu$ where $E_\mu=\ker (f-\mu)^n$. Now, $(f-\lambda \mu)^n \phi=\lambda^n \phi (f-\mu)^n$ so, in particular, $\phi(E_\mu) \subset E_{\lambda \mu}$. Ordering the eigenvalues of $f$ and assuming that $\lambda$ is not a root of unity, we get the statement about $\phi$. Furthermore, by definition of the $E_i$, we know that $f=\sum_i f_i$ with $f_i \in \Aut(E_i)$. This implies that, if $f(s)=k \cdot s$ for some $k \in \C^*$ then whenever we can write $s=\sum_i s_i$ with $s_i \in H^0(E_i)$ then each of the $f_i$ scales $s_i$ by a different constant. This clearly implies that only one of the $s_i$ can be nonzero.
\end{proof}
\end{prop}
Therefore we have
\begin{cor}
The fixed points of the $\C^*$-action on $\mathcal{M}_\sigma^{r,d}$ are exactly the points of the form above.
\begin{proof}
We simply need to check that every such point is fixed. But this is clear because, if $(E,\phi,s)$ is in the form above, then $f$ defined by $f_{|E_i}=\lambda^{i-1} \id_{E_i}$ is an isomorphism between $(E,\lambda\phi,s)$ and $(E,\phi,s)$.
\end{proof}
\end{cor}
\subsection{Structure of the attracting sets}
The results by Bia\l{}ynicki-Birula can only be applied for the few $\mathcal{M}_\sigma^{r,d}$ that are smooth. For example from proposition \ref{smoothr2} we know that when $r=2$ and $\sigma$ is small, we should have $d<0$ or $d$ odd and $>4g-4$ in order to have smooth moduli spaces. In the other cases in fact the moduli spaces are not smooth. Nevertheless, since we still have a proper Hitchin map, we have the decomposition of $\mathcal{M}_\sigma^{r,d}$ into attracting loci.
\begin{defn}
Given a partition $\underline{r}=(r_1, \dots ,r_m)$ of $r$, a partition $\underline{d}=(d_1, \dots, d_m)$ and $1 \leq k\leq m$ we denote by
\begin{equation*}
F_{\underline{r},\sigma}^{\underline{d},k}=\left \{
\begin{array}{c}
(E,\phi,s)\in \mathcal{M}_\sigma^{r,d} \text{ such that } E=E_1 \oplus \dots \oplus E_m\\
\text{with } \deg E_i=d_i, \rk{E_i}=r_i\\
\phi(E_i) \subseteq E_{i-1} \otimes K, \phi(E_1)=0\\
s \in H^0(E_k)
\end{array} \right \}
\end{equation*}
\end{defn}
\begin{defn}
Let us fix $F_{\underline{r},\sigma}^{\underline{d},k}$. We define two locally closed subsets by:
\begin{align*}
F_{\underline{r},\sigma}^{\underline{d},k+}&=\{(E,\phi,s) \in \mathcal{M}_\sigma^{r,d}\text{ such that } \lim_{\lambda \rightarrow 0}(E,\phi,s) \in F_{\underline{r},\sigma}^{\underline{d},k}\}\\
F_{\underline{r},\sigma}^{\underline{d},k-}&=\{(E,\phi,s) \in \mathcal{M}_\sigma^{r,d}\text{ such that } \lim_{\lambda \rightarrow \infty}(E,\phi,s) \in F_{\underline{r},\sigma}^{\underline{d},k}\}
\end{align*}
\end{defn}
Then, similar to the situation for $\mathcal{M}^{r,d}$ described in section \ref{semipr}, we have:
\begin{prop}
Let $(E,\phi,s)$ be a point in $F_{\underline{r},\sigma}^{\underline{d},k}$, for $\sigma$ not a critical value. With the previous notations, the weight $0$ part of the $\C^*$-action on $T_{(E,\phi,s)} \mathcal{M}_\sigma^{r,d}$ is given by the first hypercohomology of the complex:
\begin{equation*}
\bigoplus_{i=1}^{m} \Hom(E_i,E_i) \rightarrow \bigoplus_{i=2}^{m} \Hom(E_i,E_{i-1} \otimes K) \oplus E_k,
\end{equation*}
the positive weight part is the first hypercohomology of:
\begin{equation*}
\bigoplus_{i<j} \Hom(E_i,E_j) \rightarrow \bigoplus_{i\leq j} \Hom(E_i,E_{j} \otimes K) \oplus \bigoplus_{i>k} E_i
\end{equation*}
while the negative part is the first hypercohomology of:
\begin{equation*}
\bigoplus_{i>j}^{m} \Hom(E_i,E_j) \rightarrow \bigoplus_{i>j+1}^{m} \Hom(E_i,E_{j} \otimes K) \oplus \bigoplus_{i<k} E_i
\end{equation*}
\begin{proof}
Since $(E,\phi,s)$ is $\sigma$-stable, it only has scalar automorphisms. In particular the group of isomorphisms between $(E,\lambda \phi, s)$ and $(E,\phi,s)$ is also $\C^*$. Among these automorphisms, there is a unique one that preserves the section, which is exactly the one that acts with weight 0 on $E_k$.

In other words, suppose we have a fixed point $E=E_1 \oplus \dots \oplus E_m$, $\phi$ as usual and $s \in H^0(E_k)$. Then the canonical automorphism that we want is $f:E \rightarrow E$ such that $f_{|E_i}=\lambda^{i-k} \id_{E_i}$. Recall that the Zariski tangent space to $(E,\phi,s)$ is generated by the kernel of the map
\begin{align*}
C^1(\End E) \oplus C^0(K\End E) \oplus C^0(E) &\rightarrow C^1(E) \oplus C^1(K \End E)\\
(\tau, \nu, \gamma) & \mapsto (\tau \cdot s+d \gamma, [\tau,\phi]+d \nu).
\end{align*}

The $\C^*$-action induced on the Zariski tangent space is given by:
\begin{align*}
\lambda(\tau,0,0)=\lambda^{j-i}(\tau,0,0) &\text{ if } \tau \in C^1(E_i^* E_j)\\
\lambda(0,\nu,0)=\lambda^{j-i+1}(0,\nu,0) &\text{ if } \nu \in C^0(E_i^* E_j K)\\
\lambda(0,0,\gamma)=\lambda^{i-k}(0,0,\gamma) &\text{ if } \gamma \in C^0(E_i).\\
\end{align*}

In order to decompose the tangent complex for $(E,\phi,s)$ we can therefore consider $E_i^*E_j$ acted on with weight $i-j$ (observe that this does not depend on $k$), $E_i$ with weight $i-k$ and $K$ with weight 1 and we obtain the above complexes.
\end{proof}
\end{prop}
Later we will use the complexes to compute the dimension of the positive weight part of the $\C^*$-action on the tangent space at the fixed points.
\begin{rmk}
\label{smoothres}
An important remark is that, even though \cite[Theorem 4.1]{bialynicki1973some} does not apply to the singular $\mathcal{M}_\sigma^{r,d}$, it is still possible to apply them to the subvarieties that contain the smooth fixed points and all the points whose limit as $\lambda \rightarrow 0$ is one of the smooth fixed points. Note that these subvarieties are clearly smooth, since the dimension of the Zariski tangent space can only go up in dimension with specialization and it is already minimal at the smooth fixed points. As we will see later, some of the attracting sets lie entirely in these smooth subvarieties.
\end{rmk}
If the fixed point is smooth then we get the dimension of the affine fiber of $F_{\underline{r},\sigma}^{\underline{d},k+}\rightarrow F_{\underline{r},\sigma}^{\underline{d},k}$ at that point, by the previous remark.
\subsection{The rank 2 case}
\label{rank2disc}
We now specialize on the case $r=2$. First of all, we remark that the set of critical values for $\sigma$, when $r=2$ is, as for Bradlow pairs, the set of positive integers with the same parity as $d$.

First let us comment on the case of $d < 0$.
\begin{prop}
\label{negdegrk2}
Let $d < 0$. Then $\mathcal{M}_\infty^{2,d}=\mathcal{M}_\varepsilon^{2,d}$ and both are non-empty iff $d \geq 2-2g$.
\begin{proof}
Let $(E,\phi,s) \in \mathcal{M}_\varepsilon^{2,d}$ and consider the subbundle $S=\left \langle s \right \rangle$ generated by the section $s$. Then certainly $\deg S \geq 0 > d$ so that $S$ is stricly destabilizing for $E$ and cannot be fixed by $\phi$. This implies $s$ is cyclic for $\phi$ as we are in rank 2.\\
Viceversa, assume $(E,\phi,s) \in \mathcal{M}_\infty^{2,d}$. Suppose there exists a rank one Higgs subbundle $(L,\psi)$ of $(E, \phi)$ such that $\deg L \geq d/2$. Since $s$ is cyclic for $\phi$, we know that $\s \notin H^0(L)$. If we denote again by $S$ the subbundle generated by $s$, we would get a nonzero map $S \rightarrow E / L$. This would imply that $0 \leq \deg S \leq d - \deg L \leq d/2$ but this is impossible.

We now prove the statement about non-emptiness. If $d \geq 0$ then we know there exist $\sigma$-stable pairs and so both moduli spaces are non-empty (but they are not equal).\\
For $d<0$ the two moduli spaces coincide and if we have at least a stable triple $(E,\phi,s)$ there is a nonzero map $S \rightarrow E/S \otimes K$, since $s$ is cyclic for $\phi$. But then $0 \leq \deg S \leq d- \deg S +2g-2$ which implies $0 \leq 2 \deg S \leq d+2g-2$ or $d \geq 2-2g$.\\
Viceversa, note that for any point $P \in C$ there is an embedding
\begin{align*}
\mathcal{M}_\infty^{2,d} &\rightarrow \mathcal{M}_\infty^{2,d+2}\\
(E,\phi,s) &\mapsto (E(P),\phi,s(P))
\end{align*}
so that it is enough to prove that $\mathcal{M}_\infty^{2,d}$ is non-empty for $d=2-2g$ and $d=3-2g$.\\
For $d=2-2g$ take $E=\mathcal{O} \oplus K^*$, $s$ a constant section of $\mathcal{O}$ and $\phi$ annihilating $K^*$ and restricting to the constant map $\mathcal{O} \rightarrow K^* \otimes K=\mathcal{O}$. This is a triple in $\mathcal{M}_\infty^{2,2-2g}$.\\
For $d=3-2g$ pick $P \in C$ and define $E=\mathcal{O} \oplus K^*(P)$, $s$ a constant section of $\mathcal{O}$ and $\phi$ annihilating $K^*(P)$ and restricting to the only nonzero map $\mathcal{O} \rightarrow K^* \otimes K(P)=\mathcal{O}(P)$. This is instead a triple in $\mathcal{M}_\infty^{2,3-2g}$.
\end{proof}
\end{prop}
Note that this proposition also implies that if $(E,\phi,s) \in \mathcal{M}_\varepsilon^{2,d}$ for $d <0$ then $(E,\phi)$ is strictly stable. %Furthermore, $\mathcal{M}_\varepsilon^{2,d}$ is smooth for all $2-2g \leq d < 0$ (not only for $d$ odd).

We also prove a more precise version of \cite[corollary 3.4]{nitsure1991moduli} for rank $2$.
\begin{prop}
\label{sharpnits}
Let $(E,\phi)$ be a semistable Higgs bundle of rank 2 and degree $d$. If $d \geq 6g-5$ then $H^1(E)=0$.
\begin{proof}
If $E$ is semistable itself then there are no non-zero maps $E \rightarrow K$ since $\mu(E) > 3g-2 >2g-2 = \mu(K)$, therefore $H^1(E)=H^0(KE^*)^*=0$. Otherwise assume that $L \subset E$ is the maximal destabilizing line subbundle of $E$. Then $\deg L >d/2$ but since $(E,\phi)$ is stable there is also a non-zero map $L \rightarrow E/L \otimes K$ and so $\deg L \leq d/2 +g-1$. In particular we also get $1-g+d/2 \leq \deg E/L < d/2$. If $d \geq 6g-4$ then we have $\deg L > d/2 > 3g-2 >2g-2$ and $\deg E/L > 3g-3+1-g=2g-2$ and therefore $H^1(L)=H^1(E/L)=0$ which implies $H^1(E)=0$.
\end{proof}
\end{prop}
\begin{rmk}
\label{nits}
Note that the estimate in \cite[corollary 3.4]{nitsure1991moduli} would say that $H^1(E)=0$ for $d > 4g-2 + 2g-2=6g-4$ i.e. $d \geq 6g-3$ which leaves out a value of $d$ that we will need in the following chapters. It is also reasonable to think that the gap between the sharp estimate and the estimate in \cite[corollary 3.4]{nitsure1991moduli} will widen as the rank increases.
\end{rmk}

We can prove the following.
\begin{theorem}
\label{fpr2}
Let $d\geq 2-2g$ be an integer and $\sigma>0$ different from a critical value. Then we can classify the components of the fixed point locus of $\mathcal{M}_\sigma^{2,d}$ as follows:
\begin{itemize}
\item[(i)] if $d \geq 0$ and $\sigma < d$ then one of the components of the fixed points for the $\C^*$-action is $F_{(2),\sigma}^{(d),1}=M_\sigma^{2,d}$, i.e. the moduli space of $\sigma$-stable Bradlow pairs embedded as triples with zero Higgs field. If $d <0$ then there are no $\sigma$-stable Bradlow pairs and so this component is empty.
\item[(ii)] if there exists and integer $m$ such that $\max\{0,d/2 +1-g\} \leq m < \frac{d-\sigma}{2}$, then there exist components:
\begin{align*}
F_{(1,1),\sigma}^{(d_1,d_2),1}&=\left \{
\begin{array}{c}
(E,\phi,s)\in \mathcal{M}_\sigma^{2,d} \text{ such that } E=E_1 \oplus E_2\\
\text{with } \deg E_i=d_i, \rk{E_i}=1\\
\phi(E_2) \subseteq E_{1} \otimes K, \phi(E_1)=0\\
s \in H^0(E_1)
\end{array} \right \} \cong \\
&\cong S^{d_1}( C) \times S^{d_1-d_2+2g-2}( C).
\end{align*}
Here $d_1$ and $d_2$ are integers satisfying $d_1+d_2=d$ and one of the following equivalent inequalities:
\begin{align*}
&\max\{0,d/2 +1-g\} \leq d_1 < \frac{d-\sigma}{2} \\
&\frac{d+\sigma}{2} < d_2 \leq \min\{d,d/2+g-1\}\\
&\max\{0,2g-2-d\}\leq d_1-d_2+2g-2 < 2g-2-\sigma \text{ same parity as $d$}.
\end{align*}
\item[(iii)] There exist components:
\begin{align*}
F_{(1,1),\sigma}^{(d_1,d_2),2}&=\left \{
\begin{array}{c}
(E,\phi,s)\in \mathcal{M}_\sigma^{2,d} \text{ such that } E=E_1 \oplus E_2\\
\text{with } \deg E_i=d_i, \rk{E_i}=1\\
\phi(E_2) \subseteq E_{1} \otimes K, \phi(E_1)=0\\
s \in H^0(E_2)
\end{array} \right \} \cong\\
&\cong S^{d_2}( C) \times S^{d_1-d_2+2g-2}( C).
\end{align*}
Here $d_1$ and $d_2$ are integers satisfying $d_1+d_2=d$ and one of the following equivalent inequalities:
\begin{align*}
&d/2 +1-g \leq d_1 < \min \left \{\frac{d+\sigma}{2},d+1 \right \} \\
&\max \left \{-1,\frac{d-\sigma}{2} \right \}< d_2 \leq d/2+g-1\\
& 0\leq d_1-d_2+2g-2 < 2g-2+\min\{\sigma,d+1\} \text{ same parity as $d$}.
\end{align*}
\end{itemize}
\begin{proof}
Point (i) simply follows by noting that $M^{2,d}$ is nonempty if and only if $d \geq 0$ and $\sigma < d$ and $(E,s) \in M^{2,d}$ if and only if $(E,0,s) \in \mathcal{M}_\sigma^{2,d}$. Clearly, every $\sigma$-stable triple whose Higgs field vanishes is fixed by the $\C^*$-action.

For (ii), recall from proposition \ref{fixedbht} that for the partition $(1,1)$ we must have fixed points of the form above with $s \in H^0(E_1)$ if $k=1$. Furthermore, note that $E_1$ is the only $\phi$-invariant subbundle in this case. $s \in H^0(E_1)$ implies $d_1 \geq 0$, $\sigma$-stability implies $d_1<(d-\sigma)/2$ and the existence of a nonzero map $E_2 \rightarrow E_1 \otimes K$ implies $d_1-d_2+2g-2 \geq 0$. Putting all together and using $d_1+d_2=d$ we get the claimed set of inequalities.

For (iii) we know that $s \in E_2$ when $k=2$. Here $s \in H^0(E_2)$ forces $d_2 \geq 0$, $\sigma$-stability implies $d_1 < (d+\sigma)/2$ and as before $d_1-d_2+2g-2 \geq 0$. Using $d_1+d_2=d$ we get the claimed set of inequalities.
\end{proof}
\end{theorem}
From now on we will call \emph{non-split points} the points in $F_{(2),\sigma}^{(d),1}$, \emph{split type 1} the points in $F_{(1,1),\sigma}^{(d_1,d_2),1}$ and \emph{split type 2} the points in $F_{(1,1),\sigma}^{(d_1,d_2),2}$.

Note that split type 2 points always exist for all combinations of $d \geq 2-2g$ and $\sigma >0$ and nonsplit fixed points exist only for $0<\sigma<d$. Split type 1 points are more complicated, in fact they exist only if $\sigma < d$ and $\sigma < 2g-3$ for $d$ odd or $\sigma < 2g-2$ for $d$ even. In particular, for $d >2g-2$ this implies that there is a critical value for $\sigma$ after which the type 1 split points cease to exist.

We can briefly analyze the dimension of $\H^2$ for the deformations of the fixed points. Recall the characterization of $(\H^2)^*$ we gave in proposition \ref{h2}.

For a nonsplit point $(E,0,s)$, clearly $[\alpha,\phi]=0$ and so $\beta=0$. Therefore $\dim \H^2=\dim H^0( \End E)$ and so $(E,0,s)$ is smooth if and only if $E$ is simple, which happens for example if $E$ is stable.

For the split fixed points instead, we interpret $\phi$ as a matrix
\begin{equation*}
\begin{pmatrix}
0 & \phi\\
0 & 0
\end{pmatrix}
\end{equation*}
and an endomorphism of $E_1 \oplus E_2$ has the form
\begin{equation*}
\begin{pmatrix}
a & b\\
c & d
\end{pmatrix}
\end{equation*}
with $a \in \C$ and $c \in \Hom(E_1,E_2)$, $b \in \Hom(E_2,E_1)$. We have $\beta=\beta_1+\beta_2$ with $\beta_i : E_i \rightarrow K$ and $s=s_1+s_2$ with $s_i \in H^0(E_i)$.
Then:
\begin{equation*}
[\alpha,\phi]=\begin{pmatrix}
-\phi c & (a-d)\phi\\
0 & c \phi
\end{pmatrix}
\end{equation*}
and
\begin{equation*}
\beta \otimes s=\begin{pmatrix}
\beta_1 s_1 & \beta_2 s_1\\
\beta_1 s_2 & \beta_2 s_2
\end{pmatrix}.
\end{equation*}

If $s_2=0$, i.e. for a type 1 fixed point, we get $c=0$, $\beta_1=0$ and $\beta_2 s_1=(a-d) \phi$. Note that since $d_2-d_1>0$ for a type 1 fixed point, we will always have $b=0$. In particular
\begin{equation*}
\dim \H^2=1+
\begin{cases}
1 & \text{if  } \text{div} s_1 + \text{div} \beta_2 = \text{div} \phi\\
0 & \text{otherwise}
\end{cases}
\end{equation*}
note that in the second case $\H^2 \cong \H^2([\cdot,\phi])$ and in the first $\dim \H^2 = \dim \H^2([\cdot,\phi])+1$.

If $s_1=0$, i.e. for a type 2 fixed point, we always get $c=0$ and $a=d$, therefore $[\alpha, \phi]=0$ and $\beta=0$ so $\dim \H^2 = \dim \H^2([\cdot,\phi])$.

We can summarize the discussion in the following proposition.
\begin{prop}
\label{deffixed}
A nonsplit fixed point is smooth if and only if the underlying vector bundle is simple.

A split type 1 fixed point $(E_1 \oplus E_2, \phi: E_2 \rightarrow E_1 \otimes K, s_1)$ is smooth if and only if $ \text{div}\phi-\text{div} s_1 $ is not an effective divisor. In other words the singular points in the component $S^{d_1}( C) \times S^{d_1-d_2+2g-2}( C)$ are those in the image of the map:
\begin{align*}
S^{d_1}( C) \times S^{-d_2+2g-2}( C) & \rightarrow S^{d_1}( C) \times S^{d_1-d_2+2g-2}( C)\\
(D_1, D_2) & \mapsto (D_1, D_1+D_2)
\end{align*}
which of course is empty when $d_2>2g-2$.

A split type 2 fixed point is smooth if and only if the underlying Higgs bundle is simple.
\end{prop}
We can also say something about the structure of the attracting sets, even though it is not always possible to conclude that they are affine fibrations over the fixed point loci since we are not inside smooth ambient varieties.

The fiber of the limit map $F_{(2),\sigma}^{(d),1+} \rightarrow F_{(2),\sigma}^{(d),1}$ over $(E,0,s)$ is clearly $H^0(K \otimes \End E)$. Observe that they are still affine spaces but not of constant dimension as
$$\dim H^0(K \otimes \End E)= \dim H^1(\End E)=4g-4+\dim H^0(\End E)$$
and this can vary if $E$ is not simple. Note that over the locus of $M_\sigma^{2,d} \cong F_{(2),\sigma}^{(d),1}$ where $E$ is stable, we have an affine fibration of dimension $4g-3$.

For fixed points of split type 1 we can compute the dimension of the fibers of $F_{(1,1),\sigma}^{(d_1,d_2),1+} \rightarrow F_{(1,1),\sigma}^{(d_1,d_2),1}$ over the smooth fixed points. We can decompose the deformation complex into weighted parts:
\begin{equation*}
\Hom(E_2,E_2)\oplus\Hom(E_1,E_1) \rightarrow \Hom(E_2,E_1 \otimes K) \oplus E_1
\end{equation*}
for the weight 0 part,
\begin{equation*}
\Hom(E_2,E_1)\rightarrow 0
\end{equation*}
for the negative weight part and
\begin{equation*}
\Hom(E_1,E_2) \rightarrow \Hom(E_1,E_1 \otimes K) \oplus \Hom(E_2,E_2 \otimes K) \oplus \Hom(E_1,E_2 \otimes K) \oplus E_2
\end{equation*}
for the positive weight part. The dimension of the weight 0 part, which clearly corresponds to the tangent space to the component of the fixed points, is $2d_1-d_2+2g-2$ equal to the dimension of $S^{d_1}( C) \times S^{d_1-d_2+2g-2} ( C)$. The dimension of the negative weight part is $d_2-d_1+g-1$ and finally the dimension of the positive weight part is $1+4(g-1)+d_2+1-g$.

For fixed points of split type 2 we can still decompose the deformation complex into weighted parts:
\begin{equation*}
\Hom(E_2,E_2)\oplus\Hom(E_1,E_1) \rightarrow \Hom(E_2,E_1 \otimes K) \oplus E_2
\end{equation*}
for the weight 0 part,
\begin{equation*}
\Hom(E_2,E_1)\rightarrow E_1
\end{equation*}
for the negative weight part and
\begin{equation*}
\Hom(E_1,E_2) \rightarrow \Hom(E_1,E_1 \otimes K) \oplus \Hom(E_2,E_2 \otimes K) \oplus \Hom(E_1,E_2 \otimes K)
\end{equation*}
for the positive weight part. Here type 2 split fixed points can be singular and it is actually easy to see that the positive part of the tangent space contains entirely the part that gives the extra dimension at the singular points. More specifically, the dimension of the weight 0 part is $d_1+2g-2$ equal to the dimension of $S^{d_2}( C) \times S^{d_1-d_2+2g-2} ( C)$. The dimension of the negative weight part is $d_2$ and finally the dimension of the positive weight part is $1+4(g-1)+\dim \Hom(E_2,E_1)$. As we already said, this gives the dimension of the fibers of $F_{(1,1),\sigma}^{(d_1,d_2),2+} \rightarrow F_{(1,1),\sigma}^{(d_1,d_2),2}$ over the smooth fixed points.

As for rank $2$ Higgs bundles we can prove that the decomposition of $\mathcal{M}_\sigma^{2,d}$ into attracting sets is the same as the decomposition according to the Harder-Narasimhan type of the Bradlow pair underlying the triple. We will call \emph{Hitchin stratification} the stratification of $\mathcal{M}_\sigma^{2,d}$ into attracting sets and \emph{Shatz stratification} the stratification of $\mathcal{M}_\sigma^{2,d}$ according to the Harder-Narasimhan type of the underlying pair.
\begin{theorem}
\label{shatz}
Let $\sigma$ be different from a critical value.\\
We have:
\begin{equation*}
F_{(2),\sigma}^{(d),1+}=
\left \{
\begin{array}{c}
(E,\phi,s) \text{ such that } (E,s) \text{ is $\sigma$-stable}
\end{array} \right \}
\end{equation*}
and the limit map $F_{(2),\sigma}^{(d),1+} \rightarrow F_{(2),\sigma}^{(d),1}$ takes $(E,\phi,s)$ to $(E,0,s)$.

Otherwise, the underlying Bradlow pair is $\sigma$-unstable and:
\begin{equation*}
F_{(1,1),\sigma}^{(d_1,d_2),1+}=
\left \{
\begin{array}{c}
(E,\phi,s) \text{ where } E \text{ is defined by}\\
0 \rightarrow E_2 \rightarrow E \rightarrow E_1 \rightarrow 0\\
s \notin H^0(E_2), \deg E_i =d_i, d_2 > (d+\sigma)/2\\
\text{and } \phi \text{ does not preserve } E_2
\end{array} \right \}
\end{equation*}
and the limit map $F_{(1,1),\sigma}^{(d_1,d_2),1+} \rightarrow F_{(1,1),\sigma}^{(d_1,d_2),1}$ sends $(E,\phi,s)$, described as above, to $(E_2\oplus E_1, p_{E_1} \phi i_{E_2}, p_{E_1}(s))$ where $i_{E_2}$ and $p_{E_1}$ are the inclusions and projections in the description of $E$ as an extension.

The other possibility is:
\begin{equation*}
F_{(1,1),\sigma}^{(d_1,d_2),2+}=
\left \{
\begin{array}{c}
(E,\phi,s) \text{ where } E \text{ is defined by}\\
0 \rightarrow E_2 \rightarrow E \rightarrow E_1 \rightarrow 0\\
s \in H^0(E_2), \deg E_i =d_i, d_2 > (d-\sigma)/2\\
\text{and } \phi \text{ does not preserve } E_2
\end{array} \right \}
\end{equation*}
and the limit map $F_{(1,1),\sigma}^{(d_1,d_2),2+} \rightarrow F_{(1,1),\sigma}^{(d_1,d_2),2}$ sends $(E,\phi,s)$, described as above, to $(E_2\oplus E_1, p_{E_1} \phi i_{E_2}, s)$ where $i_{E_2}$ and $p_{E_1}$ are the inclusions and projections in the description of $E$ as an extension.
\begin{proof}
The statement about $F_{(2),\sigma}^{(d),1+}$ is clear because, if $(E,s)$ is a $\sigma$-stable Bradlow pair, then $(E,\phi,s)$ will be a $\sigma$-stable triple for every $\phi \in H^0(K \End E)$. Furthermore, the entire family $\lambda \mapsto (E,\lambda \phi,s)$ consists of $\sigma$-stable triples, even for $\lambda=0$. Therefore the limit of $\lambda \cdot (E,\phi, s)$ as $\lambda \rightarrow 0$ is $(E,0,s)$ in this case.

For the other two types of cells we can use a \v{C}ech cohomology description of the triple to understand what the $\C^*$-action does. We spell all the details for $F_{(1,1),\sigma}^{(d_1,d_2),1+}$ and for the other cells the proof is analogous.

Let $(E,\phi,s)$ be defined by an extension:
\begin{equation*}
0 \rightarrow E_2 \rightarrow E \rightarrow E_1 \rightarrow 0,
\end{equation*}
call $i_{E_2}$ and $p_{E_1}$ the inclusions and projections and assume $d_2 > (d+\sigma)/2$, $s \notin H^0(E_2)$ and $\phi(E_2) \nsubseteq E_2 \otimes K$.

Then $(E,s)$ is not $\sigma$-stable and $E_2$ is the maximal destabilizing of $(E,s)$. If there were an $L \subset E$ that is $\phi$-invariant and contains the section, then clearly the map $L \rightarrow E \rightarrow E_1$ is nonzero and so $\deg L \leq \deg E_1 < (d-\sigma)/2$. This implies that every such triple is actually $\sigma$-stable.

Choose a fine enough affine cover $U_\alpha$ of $C$ such that $E$, $E_2$ and $E_1$ all trivialize. Suppose that $(E,\phi,s)$ is given by the \v{C}ech cohomology data $(U_\alpha, g_{\alpha \beta},\phi_\alpha,s_\alpha)$ relative to the open cover $U_\alpha$, where:
\begin{equation*}
g_{\alpha\beta}=
\begin{pmatrix}
e_{1\alpha\beta} & 0\\
\tau_{\alpha\beta} & e_{2\alpha\beta}
\end{pmatrix}
\end{equation*}
are the transition functions of $E$, where $e_i$ are the transition functions for $E_i$ and $\tau \in C^1( E_1^* E_2)$,
\begin{equation*}
\phi_\alpha=
\begin{pmatrix}
a_{\alpha} & b_{\alpha}\\
c_{\alpha} & d_{\alpha}
\end{pmatrix}
\end{equation*}
for the Higgs field,
\begin{equation*}
s_\alpha=
\begin{pmatrix}
s_{1\alpha}\\
s_{2\alpha}
\end{pmatrix}
\end{equation*}
for $s$ where $s_i$ is the part of the section lying in $E_i$.

Since these data define a triples, we have the compatibility relations:
\begin{align*}
& g_{\alpha \beta}=g_{\beta \alpha}^{-1}\\
& g_{\alpha \beta} g_{\beta \gamma}=g_{\alpha \gamma}\\
& s_\alpha = g_{\alpha \beta} s_{\beta}\\
& g_{\alpha \beta} \phi_\beta = \phi_\alpha g_{\alpha \beta}.
\end{align*}

Now consider the triple $(E,\lambda \phi, s)$ for $\lambda \in \C^*$, which is represented by the data $(U_\alpha, g_{\alpha \beta},\lambda \phi_\alpha,s_\alpha)$. We want to show that $(E,\lambda\phi,s)$ is also represented by the data $(U_\alpha, g^{\lambda}_{\alpha \beta},\phi^{\lambda}_\alpha,s^{\lambda}_\alpha)$ where:
\begin{equation*}
g^{\lambda}_{\alpha\beta}=
\begin{pmatrix}
e_{1\alpha\beta} & 0\\
\lambda \tau_{\alpha\beta} & e_{2\alpha\beta}
\end{pmatrix},
\end{equation*}
\begin{equation*}
\phi^{\lambda}_\alpha=
\begin{pmatrix}
\lambda a_{\alpha} & b_{\alpha}\\
\lambda^2 c_{\alpha} & \lambda d_{\alpha}
\end{pmatrix}
\end{equation*}
and 
\begin{equation*}
s^{\lambda}_\alpha=
\begin{pmatrix}
s_{1\alpha}\\
\lambda s_{2\alpha}
\end{pmatrix}.
\end{equation*}

First of all the compatibility relations for $(U_\alpha, g_{\alpha \beta},\lambda \phi_\alpha,s_\alpha)$ imply the compatibility relations for $(U_\alpha, g^{\lambda}_{\alpha \beta},\phi^{\lambda}_\alpha,s^{\lambda}_\alpha)$ so that the second set of data actually defines a triple.

Using the change of trivialization (independent of $\alpha$):
\begin{equation*}
\gamma_\alpha=
\begin{pmatrix}
\lambda & 0\\
0 & 1
\end{pmatrix} \in C^0(\GL_2)
\end{equation*}
we can see that $(U_\alpha, g_{\alpha \beta},\lambda \phi_\alpha,s_\alpha)$ and $(U_\alpha, g^{\lambda}_{\alpha \beta},\phi^{\lambda}_\alpha,s^{\lambda}_\alpha)$ are equivalent set of data and so they define the same triple $(E,\lambda \phi, s)$. What we need to verify explicitly are the relations:
\begin{align*}
& \gamma_\alpha^{-1}g_{\alpha \beta} \gamma_\beta = g^{\lambda}_{\alpha \beta}\\
& \gamma_\alpha^{-1} \lambda \phi_\alpha \gamma_\alpha =\phi^\lambda_\alpha\\
& \gamma_\alpha^{-1} \lambda s_\alpha=s_\alpha^{\lambda}
\end{align*}
which are an immediate check. Note that in the last relation we are allowed to change $s_\alpha$ by $\lambda s_\alpha$ because two triples which are identical except for the sections which differ by multiplication by a constant are always isomorphic.

The data $(U_\alpha, g^{\lambda}_{\alpha \beta},\phi^{\lambda}_\alpha,s^{\lambda}_\alpha)$ produce a $\C^*$-family of $\sigma$-stable triples whose limit as $\lambda \rightarrow 0$ is given by the data:
\begin{equation*}
\begin{pmatrix}
e_{1\alpha\beta} & 0\\
0 & e_{2\alpha\beta}
\end{pmatrix},
\end{equation*}
\begin{equation*}
\begin{pmatrix}
0 & b_\alpha\\
0 & 0
\end{pmatrix},
\end{equation*}
and
\begin{equation*}
\begin{pmatrix}
s_{1 \alpha}\\
0
\end{pmatrix}
\end{equation*}
that correspond to the triple $(E_1\oplus E_2, p_{E_1} \phi i_{E_2}, p_{E_1}(s))$ which lies in $F_{(1,1),\sigma}^{(d_1,d_2),1}$.

Observe also that:
\begin{equation*}
F_{(1,1),\sigma}^{(d_1,d_2),1}\subset
\left \{
\begin{array}{c}
(E,\phi,s) \text{ where } E \text{ is defined by}\\
0 \rightarrow E_2 \rightarrow E \rightarrow E_1 \rightarrow 0\\
s \notin H^0(E_2), \deg E_i =d_i, d_2>(d+\sigma)/2\\
\text{and } \phi \text{ does not preserve } E_2
\end{array} \right \}
\end{equation*}
as is immediately checked by definition. This proves that
\begin{equation*}
F_{(1,1),\sigma}^{(d_1,d_2),1+}\supset
\left \{
\begin{array}{c}
(E,\phi,s) \text{ where } E \text{ is defined by}\\
0 \rightarrow E_2 \rightarrow E \rightarrow E_1 \rightarrow 0\\
s \notin H^0(E_2), \deg E_i =d_i, d_2>(d+\sigma)/2\\
\text{and } \phi \text{ does not preserve } E_2
\end{array} \right \}.
\end{equation*}

The same kind of argument allows to prove
\begin{equation*}
F_{(1,1),\sigma}^{(d_1,d_2),2+}\supset
\left \{
\begin{array}{c}
(E,\phi,s) \text{ where } E \text{ is defined by}\\
0 \rightarrow E_2 \rightarrow E \rightarrow E_1 \rightarrow 0\\
s \in H^0(E_2), \deg E_i =d_i, d_1 > (d-\sigma)/2\\
\text{and } \phi \text{ does not preserve } E_2
\end{array} \right \}
\end{equation*}
and that the limit map has the claimed form.

To conclude we observe that both the decomposition into attracting set and the suggested decomposition cover $\mathcal{M}^{2,d}_\sigma$. The first one because limits for the $\C^*$-action as $\lambda \rightarrow 0$ always exists, the second one because it exhausts all the possible Harder-Narasimhan types of the underlying Bradlow pair of a triple in $\mathcal{M}^{2,d}_\sigma$. Therefore we can conclude that all the inclusions we proved are actually equalities.
\end{proof}
\end{theorem}

%% file: Motives.tex
\chapter{Motivic invariants and wall-crossing}
%
%\begin{quote}
%{\small
%Summary of chapter contents\\
%\begin{itemize}
%\item initial stability for odd degree
%\item initial stability for even degree
%\item Wall crossing (???Proj construction of the flip loci, stratification of the flip loci)
%\item Motives of the flip loci
%\item Direct computation of $[\mathcal{M}_\infty^{2,d}]$
%\item Few comments on the direct computation of the single cells
%\item Comments of recent methods applied to triples
%\item Corollary about motives of Hilbert schemes and Jacobians of weird curves.
%\end{itemize}
%}
%\end{quote}
%
\section{Wall crossing for rank $2$}
In this section we examine what happens when the stability parameter for rank 2 Bradlow-Higgs triples is modified. First we examine the geometry of the loci that change when crossing a critical value, the so called \emph{flip loci}. Second, we compute their class in the Grothendieck ring of varieties.
\subsection{The flip loci}
First of all note that, in the rank 2 case, the possible critical values for $\sigma$ are the integers greater than or equal to 1, with the same parity as the degree $d$.

As we discussed in section \ref{Bradlow}, the moduli spaces will be modified by adding and subtracting objects that are defined using extensions of Higgs bundles. Fix a critical value $\bar\sigma$. Let us start with a definition.
\begin{defn}[Flip loci]
Let $\bar\sigma$ be a critical value. We denote by $\mathcal{W}_{\bar\sigma}^{d,+}$ the locally closed subvariety of $\mathcal{M}_{\bar\sigma_-}^{2,d}$ consisting of those triples that are $\bar\sigma_-$-stable but not $\bar\sigma_+$-stable.

Analogously we denote by $\mathcal{W}_{\bar\sigma}^{d,-}$ the locally closed subvariety of $\mathcal{M}_{\bar\sigma_+}^{2,d}$ consisting of those triples that are $\bar\sigma_+$-stable but not $\bar\sigma_-$-stable.
\end{defn}
Before going into the details let us give some useful definitions.
\begin{defn}
Given a degree $d \geq 1$ and $\bar\sigma$ a critical value for $d$ we define:
\begin{align*}
p^d_{\bar\sigma}: S^{(d-\bar\sigma)/2}( C) &\rightarrow J^{(d-\bar\sigma)/2}( C)\\
D \mapsto &\mathcal{O}(D),
\end{align*}
and
\begin{equation*}
X^d_{\bar\sigma}:=p_{\bar\sigma}^{d*}T^*J^{(d-\bar\sigma)/2}( C)\times T^*J^{(d+\bar\sigma)/2}( C).
\end{equation*}
Finally, we denote by:
\begin{equation*}
q^d_{\bar\sigma}: X_{\bar\sigma} \times C \rightarrow X_{\bar\sigma}
\end{equation*}
the proper projection on the first factor.
\end{defn}
The variety $X^d_{\bar\sigma}$ has a cumbersome definition but is actually simple to understand. In fact it is a product of the total spaces of two vector bundles but they are both trivial, therefore $$X^d_{\bar\sigma}\cong  S^{(d-\bar\sigma)/2}( C) \times J^{(d+\bar\sigma)/2}( C) \times H^0(K) \times H^0(K).$$
It is also clear from the definition of $X^d_{\bar\sigma}$ that it is the moduli space of pairs of rank one Higgs bundles, plus a section of the first one.

First we treat the case of wall crossing occurring when the parameter is increasing, i.e. for the flip locus $\mathcal{W}_{\bar\sigma}^{d,+}$. Looking at definition \ref{sigmast}, we see that the inequality involving the section becomes stronger while the other one becomes weaker.

We are looking for triples $(E,\phi,s)$ that are $\bar\sigma_-$-stable but not $\bar\sigma_+$-stable. These two conditions together impose the existence of a subobject that is preserved by $\phi$ and is just barely $\bar\sigma_+$-destabilizing.
More precisely, we want triples $(E,\phi,s)$ such that there is a $\phi$-invariant subbundle $L \subset E$, $s \in H^0(L)$ and $\deg L=(d-\bar\sigma)/2$. In other words, we need to have an extension of Higgs bundles:
\begin{equation*}
0 \rightarrow (L,\psi_1) \rightarrow (E,\phi) \rightarrow (M, \psi_2) \rightarrow 0
\end{equation*}
such that $s \in H^0(L)$. This, in particular, implies that $L=\mathcal{O}(D)$ where $D$ is the effective divisor associated to the section $s$.

Let us first establish a result about extensions of Higgs bundles of rank 1.
\begin{lemma}
\label{exts}
Let $(L,\psi_1)$ and $(M, \psi_2)$ be two rank one Higgs bundles. Then there exists a vector space, that we denote by $\H^1((M,\psi_2),(L,\psi_1))$ that parametrizes all extensions of Higgs bundles starting with $(L, \psi_1)$ and ending with $(M,\psi_2)$. Such a vector space is the first hypercohomology of the complex:
\begin{align*}
M^*   L &\rightarrow M^*   L   K\\
f & \mapsto f \psi_2 - \psi_1 f
\end{align*}
and therefore fits into the long exact sequence:
\begin{align}
0 \rightarrow \H^0((M,\psi_2),(L,\psi_1)) \rightarrow H^0(M^*   L) \rightarrow H^0(M^*   L   K) \rightarrow \H^1((M,\psi_2),(L,\psi_1)) \rightarrow\nonumber\\ \rightarrow H^1(M^*   L) \rightarrow H^1(M^*   L   K) \rightarrow \H^2((M,\psi_2),(L,\psi_1)) \rightarrow 0. \label{eq:exext}
\end{align}

The zero class in $\H^1((M,\psi_2),(L,\psi_1))$ corresponds to a split extension of Higgs bundles and if two classes differ by the multiplication by a nonzero scalar, then they define isomorphic Higgs bundles.
\begin{proof}
Choose an open cover $\{U_\alpha\}$ of $C$ on which $L$, $M$ and $K$ are all trivial. Let us denote by $\tau \in C^1(M^*   L)$ the off diagonal portion of the transition matrix for our vector bundle $E$. If $l_{\alpha\beta}$ and $m_{\alpha\beta}$ are the transition functions for $L$ and $M$ respectively then the transition functions for $E$ can be written as
\begin{equation*}
\begin{pmatrix}
l_{\alpha\beta} & \tau_{\alpha\beta}\\
0 & m_{\alpha\beta}
\end{pmatrix}.
\end{equation*}

Let then $\nu \in C^0(M^*   L   K)$ be the off diagonal part of the Higgs field. So that the restriction of the Higgs field to the open subset $U_\alpha$ will be:
\begin{equation*}
\begin{pmatrix}
\psi_{1\alpha} & \nu_{\alpha}\\
0 & \psi_{2\alpha}
\end{pmatrix}.
\end{equation*}

After the notation is fixed, we can write the compatibility relation between $\tau$ and $\nu$:
\begin{equation*}
\begin{pmatrix}
l_{\alpha\beta} & \tau_{\alpha\beta}\\
0 & m_{\alpha\beta}
\end{pmatrix}
\begin{pmatrix}
\psi_{1\beta} & \nu_{\beta}\\
0 & \psi_{2\beta}
\end{pmatrix}=
\begin{pmatrix}
\psi_{1\alpha} & \nu_{\alpha}\\
0 & \psi_{2\alpha}
\end{pmatrix}
\begin{pmatrix}
l_{\alpha\beta} & \tau_{\alpha\beta}\\
0 & m_{\alpha\beta}
\end{pmatrix}
\end{equation*}
which, after taking into account the trivializations of the bundles, becomes $d \nu=\psi_1 \tau-\tau \psi_2$, i.e. the kernel of the following map:
\begin{align*}
C^0(M^*   L   K) \oplus C^1(M^*   L) &\rightarrow C^1(M^*   L   K)\\
(\nu, \tau) &\mapsto d\nu -\tau \psi_2+\psi_1 \tau.
\end{align*}

Note that with this procedure we produce two isomorphic Higgs bundles $(E,\psi)$ and $(E',\psi')$ (that of course still have $(L,\psi_1)$ as a subobject and $(M,\psi_2)$ as a quotient) if and only if there is $\gamma \in C^0(M^*   L)$ that yields an equivalence between the \v{C}ech data representing $E$ and $E'$. The matrices of such an equivalence will look like:
\begin{equation*}
\begin{pmatrix}
1 & \gamma_{\alpha}\\
0 & 1
\end{pmatrix}.
\end{equation*}

The compatibility relation with the respective transition matrices of $E$ and $E'$ therefore becomes
\begin{equation*}
\begin{pmatrix}
1 & \gamma_{\alpha}\\
0 & 1
\end{pmatrix}
\begin{pmatrix}
l_{\alpha\beta} & \tau'_{\alpha\beta}\\
0 & m_{\alpha\beta}
\end{pmatrix}
=
\begin{pmatrix}
l_{\alpha\beta} & \tau_{\alpha\beta}\\
0 & m_{\alpha\beta}
\end{pmatrix}
\begin{pmatrix}
1 & \gamma_{\beta}\\
0 & 1
\end{pmatrix}
\end{equation*}
while for the Higgs fields we get:
\begin{equation*}
\begin{pmatrix}
1 & \gamma_{\alpha}\\
0 & 1
\end{pmatrix}
\begin{pmatrix}
\psi_{1\alpha} & \nu'_{\alpha}\\
0 & \psi_{2\alpha}
\end{pmatrix}=
\begin{pmatrix}
\psi_{1\alpha} & \nu_{\alpha}\\
0 & \psi_{2\alpha}
\end{pmatrix}
\begin{pmatrix}
1 & \gamma_{\alpha}\\
0 & 1
\end{pmatrix}.
\end{equation*}

The relations can be written in a compact form as $\tau'-\tau=d \gamma$ and $\nu'-\nu=\psi_1 \gamma-\gamma \psi_2$. These relations can also be expressed as the image of the map
\begin{align*}
C^0(M^*   L) &\rightarrow C^0(M^*   L   K) \oplus C^1(M^*   L)\\
\gamma &\mapsto (\psi_1 \gamma-\gamma \psi_2, d\gamma).
\end{align*}
This clearly completes the assertion about the hypercohomology, from which also follows the statement about the long exact sequence.

The zero class corresponds to the data
\begin{equation*}
\begin{pmatrix}
l_{\alpha\beta} & 0\\
0 & m_{\alpha\beta}
\end{pmatrix}, \quad
\begin{pmatrix}
\psi_{1\alpha} & 0\\
0 & \psi_{2\alpha}
\end{pmatrix}
\end{equation*}
which clearly represent the Higgs bundle $(L,\psi_1) \oplus (M,\psi_2)$. Changing $\{\tau_{\alpha\beta}\}$ and ${\nu_{\alpha}}$ with $\{\lambda\tau_{\alpha\beta}\}$ and ${\lambda\nu_{\alpha}}$, for some $\lambda \neq 0$, corresponds to acting with the change of trivialization:
\begin{equation*}
\theta_\alpha=
\begin{pmatrix}
\lambda & 0\\
0 & 1
\end{pmatrix}
\end{equation*}
and therefore does not change $(E,\phi)$.
\end{proof}
\end{lemma}
We have the following:
\begin{prop}
Fix a degree $d$ and a critical value $\bar\sigma$. The locus $\mathcal{W}_{\bar\sigma}^{d,+}$ consists of those triples $(E,\phi,s)$ for which $(E,\phi)$ fits into a non-split extension:
\begin{equation*}
0 \rightarrow (L,\psi_1) \rightarrow (E,\phi) \rightarrow (M, \psi_2) \rightarrow 0
\end{equation*}
where $\deg L =(d-\bar\sigma)/2$, $\deg M =(d+\bar\sigma)/2$ and $s \in H^0(L)$.

Furthermore, there is a projective map $\pi_{\bar\sigma}^{d,+}: \mathcal{W}_{\bar\sigma}^{d,+} \rightarrow X^d_{\bar\sigma}$ whose fibers are projective spaces (but not of constant dimension).
\begin{proof}
Consider a triple $(E,\phi,s)$ that fits in the non-split extension:
\begin{equation*}
0 \rightarrow (L,\psi_1) \rightarrow (E,\phi) \rightarrow (M, \psi_2) \rightarrow 0
\end{equation*}
with $s \in H^0(L)$, $\deg L=(d-\bar\sigma)/2$ and $\deg M=(d+\bar\sigma)/2$. Assume there is a $\phi$-invariant subbundle $L'$. Then if $s \in H^0(L')$ we get $L=L'$ and hence $\deg L' =(d-\bar\sigma)/2<(d-\bar\sigma_-)/2$. Otherwise assume $s \notin H^0(L')$. This implies there is a nonzero map $L' \rightarrow M$ and hence $\deg L' \leq (d+\bar\sigma)/2$. The equality can never occur because otherwise the map would be an isomorphism and the sequence would be split. This implies that all such triples are $\bar\sigma_-$-stable. It is clear by construction that they are not $\bar\sigma_+$-stable.

Viceversa, every triple in $\mathcal{W}_{\bar\sigma}^{d,+}$ must be $\bar\sigma_-$-stable and $\bar\sigma_+$-unstable so it must have a $\phi$-invariant subbundle containing the section of degree $(d-\bar\sigma)/2$.

The map
\begin{align*}
\pi_{\bar\sigma}^{d,+}: \mathcal{W}_{\bar\sigma}^{d,+} &\rightarrow X^d_{\bar\sigma}\\
(E,\phi,s) & \mapsto (L,s,M,\psi_1,\psi_2)
\end{align*}
is defined by sending a triple $(E,\phi,s)$ to the canonical subbundle $L$ containing the section (which has degree $(d-\bar\sigma)/2$ and is $\phi$-invariant) and to the canonical quotient $M$, both of which inherit a Higgs field from $(E,\phi)$. The fibers of such a map are, by lemma \ref{exts}, projectivized extension spaces and their dimension can vary, as we will see.
\end{proof}
\end{prop}

In lemma \ref{exts} we gave an interpretation only for the first hypercohomology of the complex as a space of extensions of Higgs bundles. We can give an interpretation for the 0th and the 2nd hypercohomology as well. In fact we can think of $\H^0((M,\psi_2),(L,\psi_1))$ the subspace of $H^0(M^* L)$ of morphisms $f: M \rightarrow L$ that satisfy $f \psi_2=\psi_1 f$, i.e. the morphisms of Higgs bundles $(M, \psi_2) \rightarrow (L, \psi_1)$. Furthermore, by using Serre duality on the starting complex we get the dual complex:
\begin{align*}
L^* M &\rightarrow L^* M   K\\
f & \mapsto f \psi_2 - \psi_1 f
\end{align*}
which is the same complex but with $L$ and $M$ interchanged. So we see that
$$\H^2((M,\psi_2),(L,\psi_1))$$
is dual to $$\H^0((L,\psi_1),(M,\psi_2)).$$
Finally from the long exact sequence we get
$$\chi((M,\psi_2),(L,\psi_1))=\chi(M^*   L)-\chi(M^*   L   K)=2-2g.$$

Since for $ \mathcal{W}_{\bar\sigma}^{d,+}$ we always have $\deg L < \deg M$, we also have $H^0(M^*   L)=0$ and also $\H^0((M,\psi_2),(L,\psi_1))=0$ because of the injection $\H^0((M,\psi_2),(L,\psi_1)) \rightarrow H^0(M^*   L)$. However $\H^2((M,\psi_2),(L,\psi_1))$ might very well jump. This, as we mentioned, means that the fibers of $\pi_{\bar\sigma}^{d,+}$ can vary in dimension.

Since the final goal will be to compute the motives of our moduli spaces, let us now examine the stratification of $X^d_{\bar\sigma}$ that is induced by the dimension of the fibers of $\pi_{\bar\sigma}^{d,+}$.
\begin{defn}
\label{defstr}
Let $\bar\sigma$ be a critical value, then we define locally closed subsets:
\begin{equation*}
S_{{\bar\sigma},i}^{d,+}:=\{(L,s,\psi_1, M, \psi_2) \in X^d_{\bar\sigma} | \dim \H^2((M,\psi_2),(L,\psi_1)) \geq i\}.
\end{equation*}
\end{defn}

Note that $S_{{\bar\sigma},i+1}^{d,+} \subset S_{{\bar\sigma},i}^{d,+}$ and $S_{{\bar\sigma},0}^{d,+}=X^d_{\bar\sigma}$. Let us fix a point $(L,s,\psi_1, M, \psi_2) \in X^d_{\bar\sigma}$. First of all, since $\dim \P\H^1((M,\psi_2),(L,\psi_1))=2g-2+\dim \H^2((M,\psi_2),(L,\psi_1))-1$, we see that $S_{{\bar\sigma},i}^{d,+}$ is exactly the stratification that agrees with the dimension of the fibers of $\pi_{\bar\sigma}^{d,+}$, after shifting the indices. In other words, the fiber of $\pi_{\bar\sigma}^{d,+}$ over each point of $S_{{\bar\sigma},i}^{d,+}\setminus S_{{\bar\sigma},i+1}^{d,+}$ is exactly $\C \P^{2g-3+i}$.
Furthermore, $\H^2((M,\psi_2),(L,\psi_1))$ is dual to $\H^0((L,\psi_1),(M,\psi_2))$ and we can easily compute the dimension of the latter. In fact $L$ and $M$ are both line bundles and so if $\psi_1 \neq \psi_2$ there are certainly no maps $L \rightarrow M$ that commute with the Higgs fields.

On the other hand, if $\psi_1=\psi_2$, then $\H^0((L,\psi_1),(M,\psi_2))=H^0(L^* M)$. From this we can deduce that $S_{\bar\sigma,1}^{d,+}$ is contained in the locus of $X^d_{\bar\sigma}$ where the Higgs fields agree, furthermore, when $\psi_1=\psi_2$, $\dim \H^0((L,\psi_1),(M,\psi_2))= \dim H^0(L^* M)$ and this is positive if and only if $L^* M=\mathcal{O}(D')$; therefore $M=L(D')$. We will soon use these considerations to compute the motive of $\mathcal{W}_{\bar\sigma}^{d,+}$.

We can now examine the second family involved in the wall crossing, which is the one appearing for decreasing parameter. Fix again a critical value $\bar\sigma$. From the definition of stabilty \ref{sigmast} we see that we should look for those triples $(E,\phi,s)$ such that there exists a $\phi$-invariant line subbundle $M \subset E$ with $s \notin H^0(M)$ and $\deg M=(d+\bar\sigma)/2$ so that the triple will be barely $\bar\sigma_-$-unstable. The fact that $M$ does not contain the section is equivalent to the section $s$ projecting on a nonzero section $\bar{s}$ of the quotient $L=E/M$.

Here we invert the notations for the subobject and the quotient, both for consistency with the degrees and for the fact that, for us, $L$ is the line bundle that comes with a section.\\
We will need a lemma about extensions as in the previous case.
\begin{lemma}
\label{extssec}
Let $(L,\psi_1)$ and $(M, \psi_2)$ be two rank one Higgs bundles such that $\bar{s} \in H^0(L)$ is a nonzero section. Then there exists a vector space, that we denote by
$$\widetilde\H^1((L, \bar{s},\psi_1),(M,\psi_2))$$
that parametrizes all extensions
\begin{equation*}
0 \rightarrow (M, \psi_2) \rightarrow (E,\phi) \rightarrow (L,\psi_1) \rightarrow 0
\end{equation*}
together with a section $s \in H^0(E)$ that projects onto $\bar{s}$. This vector space is the first hypercohomology of the complex:
\begin{align*}
L^*   M &\rightarrow L^*   M   K \oplus M\\
f & \mapsto (f \psi_1 - \psi_2 f,f(\bar{s})).
\end{align*}
The zero class corresponds to the split extension of Higgs bundles together with a trivial lift of the section, meaning that the section is the given one in the quotient. Again, if two classes in $\widetilde\H^1((L, \bar{s},\psi_1),(M,\psi_2))$ differ by the multiplication by a nonzero scalar, then they define isomorphic triples.
\begin{proof}
The proof is similar to the \v{C}ech cohomology proof that we outlined in lemma \ref{exts}, however a bit of care is needed because now subobject and quotient are inverted.

Let us denote by $\tau \in C^1(L^*   M)$ the off diagonal portion of the transition matrix for our vector bundle $E$. If $l_{\alpha\beta}$ and $m_{\alpha\beta}$ are the transition functions for $L$ and $M$ respectively then the transition functions for $E$ can be written as
\begin{equation*}
\begin{pmatrix}
m_{\alpha\beta} & \tau_{\alpha\beta}\\
0 & l_{\alpha\beta}
\end{pmatrix}.
\end{equation*}

Let then $\nu \in C^0(L^*   M   K)$ be the off diagonal part of the Higgs field. So that the restriction of the Higgs field to the open subset $U_\alpha$ will be:
\begin{equation*}
\begin{pmatrix}
\psi_{2\alpha} & \nu_{\alpha}\\
0 & \psi_{1\alpha}
\end{pmatrix}.
\end{equation*}

Let also $\theta \in C^0(M)$ be the lift of $\bar{s}$ restricted to $M$ so that on the open set $U_\alpha$ the section is written:
\begin{equation*}
\begin{pmatrix}
\theta_{\alpha}\\
\bar{s}_{\alpha}
\end{pmatrix}.
\end{equation*}

We can write the compatibility relation between $\tau$ and $\nu$ that becomes $d \nu=\psi_2 \tau-\tau \psi_1$. Also, the compatibility between $\theta$ and $\tau$ becomes $d \theta = \tau \bar{s}$. The two relations together are the kernel of the following map:
\begin{align*}
C^0(L^*   M   K) \oplus C^0(M) \oplus C^1(L^*   M) & \rightarrow C^1(L^*   M   K) \oplus C^1(M)\\
(\nu, \theta, \tau) &\mapsto (d\nu -\tau \psi_1+\psi_2 \tau, d\theta -\tau \bar{s}).
\end{align*}

We produce two isomorphic extensions if and only if there is $\gamma \in C^0(L^*   M)$ that produces an equivalence between the \v{C}ech data representing them. The matrices of such an equivalence will look like:
\begin{equation*}
\begin{pmatrix}
1 & \gamma_{\alpha}\\
0 & 1
\end{pmatrix}.
\end{equation*}

The compatibility relation with the respective transition matrices of $E$ and $E'$ therefore becomes  $\tau'-\tau=d \gamma$ while for the Higgs fields is $\nu'-\nu=\psi_2 \gamma-\gamma \psi_1$ and for the sections is $\theta'-\theta=\gamma \bar{s}$. These relations can also be expressed as the image of the map
\begin{align*}
C^0(L^*   M) &\rightarrow C^0(L^*   M   K) \oplus C^0(M) \oplus C^1(L^*   M)\\
\gamma &\mapsto (\psi_2 \gamma-\gamma \psi_1, d\gamma, \gamma \bar{s}).
\end{align*}

This completes the assertion about the hypercohomology. The last statement is obtained in the same way as in the proof of lemma \ref{exts}.
\end{proof}
\end{lemma}
We have the following.
\begin{prop}
Fix a degree $d$ and a critical value $\bar\sigma$. The locus $\mathcal{W}_{\bar\sigma}^{d,-}$ consists of those triples $(E,\phi,s)$ for which $(E,\phi)$ fits into an extension:
\begin{equation*}
0 \rightarrow (M,\psi_2) \rightarrow (E,\phi) \rightarrow (L, \psi_1) \rightarrow 0
\end{equation*}
where $\deg L =(d-\bar\sigma)/2$, $\deg M =(d+\bar\sigma)/2$ and $s \in H^0(E)$ projects to a nonzero $\bar{s} \in H^0(L)$. The condition on $(E,\phi,s)$ is that either the extension is nonsplit or it is split but $s \neq \bar{s}$.

Furthermore, there is a projective map $\pi_{\bar\sigma}^{d,-}: \mathcal{W}_{\bar\sigma}^{d,-} \rightarrow X^d_{\bar\sigma}$ whose fibers are projective spaces.
\begin{proof}
Consider a triple $(E,\phi,s)$ that fits in the extension:
\begin{equation*}
0 \rightarrow (M,\psi_2) \rightarrow (E,\phi) \rightarrow (L, \psi_1) \rightarrow 0
\end{equation*}
with $s \notin H^0(M)$, $\deg L=(d-\bar\sigma)/2$ and $\deg M=(d+\bar\sigma)/2$. Call $\bar{s}$ the projection of $s$ on $L$ and assume that either the extension is nonsplit or that it is split and $\bar{s} \neq s$. Take a $\phi$-invariant line subbundle $L'$ of $E$. Then, if $s \in H^0(L')$ then certainly $L'$ is not contained in $M$ so there is a nonzero map $L'\rightarrow L$ and therefore $\deg L' \leq \deg L = (d-\bar\sigma)/2$.

Observe that if $\deg L' = \deg L$ then $L=L'$ and the extension has to be split. In this case however the section is contained entirely in $L$ which is a contradiction. Therefore it has to be $\deg L' < \deg L$ and hence $\deg L' < (d-\bar\sigma_+)/2$. In any case $\deg L'$ is either contained in $M$, and then $\deg L' \leq \deg M < (d+\bar\sigma)/2$ or it has a nonzero map $L' \rightarrow L$ so again $\deg L' < (d+\bar\sigma_+)/2$. This proves that $(E,\phi,s)$ is $\bar\sigma_+$-stable. Clearly $(E,\phi,s)$ is also $\bar\sigma_-$-unstable because of $M$.

Vice versa, every triple in $\mathcal{W}_{\bar\sigma}^{d,-}$ has to be $\bar\sigma_-$-unstable and $\bar\sigma_+$-stable so it has to be of the above form.

The map
\begin{align*}
\pi_{\bar\sigma}^{d,-}: \mathcal{W}_{\bar\sigma}^{d,-} &\rightarrow X^d_{\bar\sigma}\\
(E,\phi,s) & \mapsto (L,\bar{s},M,\psi_1,\psi_2)
\end{align*}
is defined by sending a triple $(E,\phi,s)$ to the canonical $\bar\sigma_-$-destabilizing subbundle $M$ and to the pair $(L,\bar{s})$ consisting of the quotient $E/M$ and the projection of $s$ onto the quotient. Both $L$ and $M$ are endowed with Higgs fields because $M$ is $\phi$-invariant.

The fibers of $\pi_{\bar\sigma}^{d,-}$ are the projectivized hypercohomology spaces of lemma \ref{extssec}.
\end{proof}
\end{prop}
In order to compute $\dim \widetilde\H^1((L, \bar{s},\psi_1),(M,\psi_2))$ we can use the following exact sequence of complexes:
\begin{center}
\begin{tikzcd}
& 0\arrow{r} & 0 \arrow{r} \arrow{d} & L^*M \arrow{r}\arrow{d} & L^*M \arrow{r}  \arrow{d} & 0\\
& 0\arrow{r} & KL^*M \arrow{r} & KL^*M \oplus M \arrow{r} & M \arrow{r}& 0.
\end{tikzcd}
\end{center}

The 0th hypercohomology of the first complex of course vanishes, furthermore the same is true for the 0th hypercohomology of the last complex because the map $L^*M \rightarrow M$ is injective on the global sections. This in particular implies that the 0th hypercohomology of the middle complex also vanishes because it is between two zeroes in the long exact sequence of hypercohomologies associated to the exact sequence of complexes we wrote.

The rest of the sequence is:
\begin{align*}
0 \rightarrow& H^0(KL^*M) \rightarrow \widetilde{\H}^1(M,\bar{s},\psi_2, L, \psi_1) \rightarrow  H^0(M\mathcal{O}_D) \rightarrow \\
\rightarrow& H^1(KL^*M) \rightarrow \widetilde{\H}^2(M,\bar{s},\psi_2, L, \psi_1) \rightarrow 0.
\end{align*}
where $D$ is the divisor associated to $\bar{s}$. Since $H^1(KL^*M) \cong H^0(M^*L)^*$ and, in our situation, $\deg M > \deg L$, we see that $H^1(KL^*M)=0$. Therefore $\widetilde{\H}^2(M,\bar{s},\psi_2, L, \psi_1)=0$ as well and we can compute
$$\dim \widetilde{\H}^1(M,\bar{s},\psi_2, L, \psi_1)=\dim H^0(KL^*M)+\dim H^0(M\mathcal{O}_D)=g-1+(d+\bar\sigma)/2.$$

For the map $\pi_{\bar\sigma}^{d,-}: \mathcal{W}_{\bar\sigma}^{d,-} \rightarrow X^d_{\bar\sigma}$ we find that the fibers are actually of constant dimension so that $\mathcal{W}_{\bar\sigma}^{d,-}$ is a projective bundle over $\rightarrow X^d_{\bar\sigma}$.
%
%\begin{rmk}
%Similar to remark \ref{fam1} we can still use the $\text{\textbf{Proj}}$ construction to build the second flip locus.
%
%Take once again the universal triple $(\mathscr{L},\bar{\mathcal{S}},\Psi_1)$ and the universal Higgs bundle $(\mathscr{M},\Psi_2)$ on $X^d_{\bar\sigma} \times C$. Form the complex $\mathscr{M}^*   \mathscr{L} \oplus \mathscr{K} \mathscr{M}^* \rightarrow \mathscr{M}^*   \mathscr{L}   \mathscr{K}$ on $X_{\bar\sigma} \times C$ defined as the dual complex to the one in proposition \ref{extssec}. Define the variety:
%\begin{equation*}
%\widetilde{\mathcal{W}}^{d,-}_{\bar\sigma}:=\text{\textbf{Proj}}\left ( \text{Sym}\left (R^1 q^d_{\bar\sigma*}(\mathscr{M}^*   \mathscr{L} \oplus \mathscr{K} \mathscr{M}^* \rightarrow \mathscr{M}^*   \mathscr{L}   \mathscr{K} \right ) \right ).
%\end{equation*}
%We can then conclude again by the properties of $\textbf{Proj}$.
%\end{rmk}
%
To sum up, both families have a map to $X^d_{\bar\sigma}$. However, in the context of the first family, the divisor yields a line bundle that is the subobject $L$ of the extension we produce, while the Jacobian part yields the quotient $M$. For the second family the roles are exchanged, meaning that the quotient comes with a section, while the subobject is parametrized by the Jacobian part.

Another difference is that, while the fibers of the map $\pi_{\bar\sigma}^{d,+} : \mathcal{W}_{\bar\sigma}^{d,+} \rightarrow X^d_{\bar\sigma}$ can vary in dimension, the fibers of $\pi_{\bar\sigma}^- : \mathcal{W}_{\bar\sigma}^{d,-} \rightarrow X^d_{\bar\sigma}$ are of constant dimension.
\subsection{Motives of the flip loci}
We can start by handling the computation of $[\mathcal{W}_{\bar\sigma}^{d,+}]$.
\begin{prop}
\label{wplus}
We have the following motivic equality:
\begin{equation*}
[\mathcal{W}_{\bar\sigma}^{d,+}]=\L^{2g} \cdot [\C\P^{2g-3}]\cdot [S^{(d-{\bar\sigma})/2}( C)] \cdot [J( C)]+\L^{3g-2} \cdot [S^{(d-{\bar\sigma})/2}( C)] \cdot [S^{\bar\sigma}( C)].
\end{equation*}
\begin{proof}
We have a map $\pi_{\bar\sigma}^{d,+}: \mathcal{W}^{d,+}_{\bar\sigma} \rightarrow X^d_{\bar\sigma}$ and we can stratify $X^d_{\bar\sigma}$ according to the dimension of the fibers of this map. The strata $S_{{\bar\sigma},i}^{d,+}$ appear in definition \ref{defstr}.
We get a preliminary formula
\begin{equation*}
[\mathcal{W}_{\bar\sigma}^{d,+}]=\sum_{i=0}^{\bar\sigma}[S_{{\bar\sigma},i}^{d,+}\setminus S_{{\bar\sigma},i+1}^{d,+}]\cdot [\C\P^{2g-3+i}].
\end{equation*}

For $i>0$, we can identify $S_{{\bar\sigma},i}^{d,+}$ with the locus of points $(L,s,\psi_1,M, \psi_2)$ of $X^d_{\bar\sigma}$ such that $\psi_1=\psi_2$ and $M=L(D')$ with $\dim H^0(\mathcal{O}(D')) \geq i$. For $i=0$ the same holds but there is an extra disjoint part coming from the open locus of $X^d_{\bar\sigma}$ where $\psi_1 \neq \psi_2$.

Let us consider the maps:
\begin{align*}
S^{(d-{\bar\sigma})/2}( C) \times S^{\bar\sigma}( C)\times H^0(K) &\rightarrow S^{(d-{\bar\sigma})/2}( C) \times J^{\bar\sigma}( C)\times H^0(K)\\
(L,s,D',\psi) &\mapsto (L,s, L(D'),\psi)\\
S^{(d-{\bar\sigma})/2}( C) \times J^{\bar\sigma}( C)\times H^0(K) & \rightarrow X^d_{\bar\sigma}\\
(L,s,A,\psi) & \mapsto (L,s, ,\psi, L \otimes A, \psi).
\end{align*}

Clearly the composition of the two is the map $(L,s,D',\psi) \mapsto (L,s,\psi, L(D'),\psi)$ and the second one is injective.
Since we have an alternative description of $S_{{\bar\sigma},i}^{d,+}$ for $i >0$ we can compute:
\begin{equation*}
[S_{{\bar\sigma},i}^{d,+}\setminus S_{{\bar\sigma},i+1}^{d,+}]=[S^{(d-{\bar\sigma})/2}( C)]\cdot[V^{\bar\sigma}_{i}]\cdot[H^0(K)].
\end{equation*}
where $V^{k}_{i}=\{A \in J^{k}( C) | \dim H^0(A) = i\}$.

For $i=0$, there is an extra contribution to $S_{{\bar\sigma},0}^{d,+}\setminus S_{{\bar\sigma},1}^{d,+}$ coming from the locus of $X_{\bar\sigma}^d$ where $\psi_1 \neq \psi_2$. Therefore:
\begin{align*}
[S_{{\bar\sigma},0}^{d,+}\setminus S_{{\bar\sigma},1}^{d,+}]=[S^{(d-{\bar\sigma})/2}( C)]\cdot[V^{\bar\sigma}_{0}]\cdot[H^0(K)]+\\
[S^{(d-{\bar\sigma})/2}( C)]\cdot[J^{(d+{\bar\sigma})/2}( C)]\cdot \left([H^0(K)]^2-[H^0(K)]\right).
\end{align*}

Putting all together we get:
\begin{align*}
[\mathcal{W}_{\bar\sigma}^{d,+}]&=[\C \P^{2g-3}]\cdot[S^{(d-{\bar\sigma})/2}( C)]\cdot[J^{(d+{\bar\sigma})/2}( C)]\cdot \left([H^0(K)]^2-[H^0(K)]\right) +\\
&+[S^{(d-{\bar\sigma})/2}( C)]\cdot[H^0(K)] \cdot \sum_{i=0}^{\bar\sigma} [V^{\bar\sigma}_{i}]\cdot [\C \P^{2g-3+i}].
\end{align*}

First, we try to evaluate the second part of the sum. With a simple algebraic trick we get:
\begin{align*}
\sum_{i=0}^{\bar\sigma} &[V^{\bar\sigma}_{i}]\cdot [\C \P^{2g-3+i}]=[V^{\bar\sigma}_{0}]\cdot [\C \P^{2g-3}]+\\
&+\sum_{i=1}^{\bar\sigma} [V^{\bar\sigma}_{i}]\cdot ([\C \P^{2g-3+i}]-[\C\P^{i-1}]) + \sum_{i=1}^{\bar\sigma} [V^{\bar\sigma}_{i}]\cdot [\C\P^{i-1}].
\end{align*}

Now, clearly
$$[\C \P^{2g-3+i}]-[\C\P^{i-1}]=\L^i (1+ \dots+\L^{2g-3})=\L^i\cdot [\C\P^{2g-3}].$$
Recall that the canonical map:
\begin{align*}
S^{\bar\sigma}( C) &\rightarrow J^{\bar\sigma} ( C)\\
D &\mapsto \mathcal{O}(D)
\end{align*}
has $V^{\bar\sigma}_{1}$ as image and the fiber over $L \in J^{\bar\sigma} ( C)$ is $\P H^0(L)$. This implies that:
\begin{equation*}
\sum_{i=1}^{\bar\sigma} [V^{\bar\sigma}_{i}]\cdot [\C\P^{i-1}] = [S^{\bar\sigma} ( C)].
\end{equation*}

Putting everything in the formula we get:
\begin{align*}
\sum_{i=0}^{\bar\sigma} [V^{\bar\sigma}_{i}]\cdot [\C \P^{2g-3+i}]=[S^{\bar\sigma} ( C)]+\sum_{i=0}^{\bar\sigma} [V^{\bar\sigma}_{i}]\cdot [\C \P^{2g-3}] \cdot \L^i.
\end{align*}

With tricks similar to the ones we already used we find:
\begin{align*}
\sum_{i=0}^{\bar\sigma} [V^{\bar\sigma}_{i}]\cdot \L^i&=\sum_{i=1}^{\bar\sigma} [V^{\bar\sigma}_{i}]\cdot (\L^i-1) + \sum_{i=0}^{\bar\sigma} [V^{\bar\sigma}_{i}]=\\
&=\sum_{i=1}^{\bar\sigma} [V^{\bar\sigma}_{i}]\cdot (\L^i-1)+J^{\bar\sigma} ( C)
\end{align*}
and
\begin{align*}
\sum_{i=1}^{\bar\sigma} [V^{\bar\sigma}_{i}]\cdot (\L^i-1)=(\L-1)\cdot \sum_{i=1}^{\bar\sigma} [V^{\bar\sigma}_{i}]\cdot [\C\P^{i-1}] = (\L-1) \cdot [S^{\bar\sigma} ( C)].
\end{align*}

Recalling that $[H^0(K)]=\L^{g}$ and that the Jacobian varieties of any degree of the curve are all isomorphic we can find a final formula for:
\begin{align*}
&[\mathcal{W}_{\bar\sigma}^{d,+}]=[\C \P^{2g-3}]\cdot [S^{(d-{\bar\sigma})/2}( C)] \cdot[J^{(d+{\bar\sigma})/2}( C)]\cdot \left([H^0(K)]^2-[H^0(K)]\right) +\\
&+[S^{(d-{\bar\sigma})/2}( C)]\cdot[H^0(K)] \cdot \sum_{i=0}^{\bar\sigma} [V^{\bar\sigma}_{i}]\cdot [\C \P^{2g-3+i}]=\\
&=\L^g\cdot(\L^g-1) \cdot [\C \P^{2g-3}]\cdot[S^{(d-{\bar\sigma})/2}( C)]\cdot[J^d( C)]+\\
&+\L^g\cdot [S^{(d-{\bar\sigma})/2}( C)] \cdot \left( [S^{\bar\sigma} ( C)]+(\L-1) \cdot [\C\P^{2g-3}] \cdot [S^{\bar\sigma} ( C)] + [\C\P^{2g-3}] \cdot [J^{\bar\sigma} ( C)]\right)=\\
&=\L^{2g} \cdot [\C\P^{2g-3}]\cdot [S^{(d-{\bar\sigma})/2}( C)] \cdot [J( C)]+\L^{3g-2} \cdot [S^{(d-{\bar\sigma})/2}( C)] \cdot [S^{\bar\sigma}( C)].
\end{align*}
\end{proof}
\end{prop}

As for $\mathcal{W}_{\bar\sigma}^{d,-}$, we can easily compute the motive because the fibers of the map $\pi_\sigma^{d,-}$ are projective spaces of constant dimension. Therefore we get:
\begin{prop}
We have the motivic equality:
\begin{equation*}
[\mathcal{W}_{\bar\sigma}^{d,-}]=\L^{2g}\cdot[S^{(d-{\bar\sigma})/2}( C)]\cdot[J^{(d+{\bar\sigma})/2}( C)]\cdot [\C\P^{(d+\bar\sigma)/2+g-2}].
\end{equation*}
\end{prop}
\subsection{Interaction between the flip loci and the attracting sets}
One of the nice features of the moduli spaces of Bradlow-Higgs triples is that they combine the wall-crossing with the presence of the $\C^*$-action. In this section we study how the flip loci intersect the various attracting sets.
\begin{defn}
\label{weirdl}
We will denote by $S\mathcal{W}^{d,+}_{\bar\sigma}$ the subvariety of $\mathcal{W}^{d,+}_{\bar\sigma}$ containing those triples in $\mathcal{M}^{2,d}_{\bar\sigma_-}$ whose underlying Bradlow pair is $\bar\sigma_-$-stable.
$NS\mathcal{W}^{d,+}_{\bar\sigma}$ will instead denote the subvariety of $S\mathcal{W}^{d,+}_{\bar\sigma}$ containing the triples whose underlying Bradlow pair is not $\bar\sigma_-$-stable.
Analogously for $S\mathcal{W}^{d,-}_{\bar\sigma}$ and $NS\mathcal{W}^{d,-}_{\bar\sigma}$ but in relation to $\bar\sigma_+$-stability.

We define $B^{d,+}_{\bar\sigma}$ to be the locally closed subvariety containing those triples in $F_{(2),\bar\sigma_-}^{(d),1+}$ such that the underlying Bradlow pair is contained in $\P W^{d,+}_{\bar\sigma}$. 
We also define $B^{d,-}_{\bar\sigma}$ as the locally closed subvariety containing those triples in $F_{(2),\bar\sigma_+}^{(d),1+}$ such that the underlying Bradlow pair is contained in $\P W^{d,-}_{\bar\sigma}$.
\end{defn}
Let us first remark that when $d <0 $ there are no $\sigma$-stable Bradlow pairs, for any $\sigma$. Therefore many of the subvarieties we examine in the following are empty for $d <0$.
We can prove two propositions summarizing the relations and properties of the loci defined above.
\begin{prop}
\label{weirdexpl1}
Let $\bar\sigma$ be a critical value. Then:
\begin{itemize}
\item[(i)] $\mathcal{W}^{d,+}_{\bar\sigma}=S\mathcal{W}^{d,+}_{\bar\sigma} \sqcup NS\mathcal{W}^{d,+}_{\bar\sigma}$.
\item[(ii)] $S\mathcal{W}^{d,+}_{\bar\sigma}=F_{(2),\bar\sigma_-}^{(d),1+} \cap \mathcal{W}^{d,+}_{\bar\sigma}$ and in particular $S\mathcal{W}^{d,+}_{\bar\sigma} \subset F_{(2),\bar\sigma_-}^{(d),1+}$.
\item[(iii)] The canonical map $\mathcal{W}^{d,+}_{\bar\sigma} \rightarrow X^d_{\bar\sigma}$ restricts to a map $NS\mathcal{W}^{d,+}_{\bar\sigma} \rightarrow X^d_{\bar\sigma}$ whose fiber over $(L,s,M,\psi_1,\psi_2)$ is $\P H^0(M^* L K)$.
\item[(iv)] $\mathcal{W}^{d,-}_{\bar\sigma}=S\mathcal{W}^{d,-}_{\bar\sigma} \sqcup NS\mathcal{W}^{d,-}_{\bar\sigma}$.
\item[(v)] $S\mathcal{W}^{d,-}_{\bar\sigma}=F_{(2),\bar\sigma_+}^{(d),1+} \cap \mathcal{W}^{d,-}_{\bar\sigma}$ and in particular $S\mathcal{W}^{d-}_{\bar\sigma} \subset F_{(2),\bar\sigma_+}^{(d),1+}$.
\item[(vi)] The canonical map $\mathcal{W}^{d,-}_{\bar\sigma} \rightarrow X^d_{\bar\sigma}$ restricts to a map $NS\mathcal{W}^{d,-}_{\bar\sigma} \rightarrow X^d_{\bar\sigma}$ whose fiber over $(L,\bar{s},M,\psi_1,\psi_2)$ is $\P H^0(L^* M K)$.
\item[(vii)] $NS\mathcal{W}^{d,+}_{\bar\sigma} \subset F_{(1,1),\bar\sigma_-}^{((d-\bar\sigma)/2,(d+\bar\sigma)/2),1+}$.
\item[(viii)] $NS\mathcal{W}^{d,-}_{\bar\sigma} \subset F_{(1,1),\bar\sigma_+}^{((d+\bar\sigma)/2,(d-\bar\sigma)/2),2+}$.
\end{itemize}
\begin{proof}
(i) follows from the fact that if $(E,\phi,s) \in \mathcal{W}^{d,+}_{\bar\sigma}$ then it is $\bar\sigma_-$-stable. Obviously, this can happen either because the underlying Bradlow pair is already $\bar\sigma_-$-stable or because the underlying Bradlow pair is unstable but the triple is $\bar\sigma_-$-stable anyway.

For (ii) it is enough to recall that $F_{(2),\bar\sigma_-}^{(d),1+}$ contains all the triples in $\mathcal{M}_{\bar\sigma_-}^{2,d}$ for which the underlying Bradlow pair is itself $\bar\sigma_-$-stable.

To prove (iii) consider a triple in $(E,\phi,s) \in \mathcal{W}^{d,+}_{\bar\sigma}$ whose subobject is $(L,\psi_1)$, with $s \in H^0(L)$ and whose quotient is $(M,\psi_2)$. From lemma \ref{exts} we know that there is an exact sequence:
\begin{align*}
& \H^0((M,\psi_2),(L,\psi_1)) \rightarrow H^0(M^* \otimes L) \rightarrow H^0(M^* \otimes L \otimes K) \rightarrow \H^1((M,\psi_2),(L,\psi_1)) \rightarrow\\
&\rightarrow H^1(M^* \otimes L) \rightarrow H^1(M^* \otimes L \otimes K) \rightarrow \H^2((M,\psi_2),(L,\psi_1)) \rightarrow 0
\end{align*}
and the triple is represented, up to isomorphism, by a class in $\P \H^1((M,\psi_2),(L,\psi_1))$.

Also, the underlying Bradlow pair is $\bar\sigma_-$-stable if and only if the class representing the triple does not lie in $$\ker( \H^1((M,\psi_2),(L,\psi_1)) \rightarrow H^1(M^* \otimes L)).$$

Since $\deg M > \deg L$ we know that $H^0(M^* L)=0$ so that $$H^0(M^*  L  K)=\ker( \H^1((M,\psi_2),(L,\psi_1)) \rightarrow H^1(M^*  L)).$$

From this we deduce that the projectivization of the previous kernel is the fiber of $NS\mathcal{W}^{d,+}_{\bar\sigma} \rightarrow X^d_{\bar\sigma}$.

Statements (iv) and (v) follow as (i) and (ii). As for (vi), let $(E,\phi,s) \in \mathcal{W}^{d,-}_{\bar\sigma}$ and $(M,\psi_2)$ be the subobject, while $(L,\psi_1)$ be the quotient and $s \notin H^0(M)$ projects to $\bar{s} \in H^0(M)$. From lemma \ref{extssec} we know there is a different exact sequence:
\begin{align*}
0 \rightarrow H^0(KL^*M) \rightarrow \widetilde\H^1((L, \bar{s},\psi_1),(M,\psi_2)) \rightarrow  H^0(M\mathcal{O}_D) \rightarrow 0
\end{align*}
where $D$ is the divisor associated to $\bar{s}$. The underlying Bradlow pair, $(E,s)$ is $\bar\sigma_+$-stable if and only if the class representing $(E,\phi,s)$ in $\P  \widetilde\H^1((L, \bar{s},\psi_1),(M,\psi_2))$ is not in $$\ker( \widetilde\H^1((L, \bar{s},\psi_1),(M,\psi_2)) \rightarrow H^0(M\mathcal{O}_D))=H^0(KL^*M).$$

Therefore the canonical map $\mathcal{W}^{d,-}_{\bar\sigma} \rightarrow X^d_{\bar\sigma}$ also restricts to a canonical map
$$NS\mathcal{W}^{d,-}_{\bar\sigma} \rightarrow X^d_{\bar\sigma}$$
whose fiber over $(L,\bar{s},M,\psi_1, \psi_2)$ is $\P H^0(KL^*M)$.\\
For (vii) consider a triple in $NS\mathcal{W}^+_{\bar\sigma}$. The underlying Bradlow pair is split, the maximal $\bar\sigma_-$-destabilizing subbundle does not contain the section and is of degree $(d+\bar\sigma)/2$. So the Harder-Narasimhan type of the underlying pair matches with $F_{(1,1),\bar\sigma_-}^{((d-\bar\sigma)/2,(d+\bar\sigma)/2),1+}$.

Statement (viii) is proved as (vii).
\end{proof}
\end{prop}
\begin{prop}
\label{weirdexpl}
Let $\bar\sigma$ be a critical value. Then:
\begin{itemize}
\item[(i)] $B^{d,+}_{\bar\sigma}$ is the inverse image of $\P W^{d,+}_{\bar\sigma}$ with respect to the limit map $F_{(2),\bar\sigma_-}^{(d),1+} \rightarrow F_{(2),\bar\sigma_-}^{(d),1}\cong M^{2,d}_{\bar\sigma_-}$.
\item[(ii)] $S\mathcal{W}^{d,+}_{\bar\sigma} \subset B^{d,+}_{\bar\sigma}$ and $B^{d,+}_{\bar\sigma} \setminus S\mathcal{W}^{d,+}_{\bar\sigma} \subset F_{(1,1),\bar\sigma_+}^{((d+\bar\sigma)/2,(d-\bar\sigma)/2),2+}$. The complement can be described as:
\begin{equation*}
\left \{
\begin{array}{c}
(E,\phi,s) \text{ such that }\\
E=L \oplus M\\
s \in H^0(L), \deg L =(d-\bar\sigma)/2, \deg M=(d+\bar\sigma)/2\\
\text{and } \phi \text{ does not preserve } L
\end{array} \right \}.
\end{equation*}
and therefore contains $NS\mathcal{W}^{d,-}_{\bar\sigma}$. Furthermore, if $\bar\sigma > 2g-2$ then the locally closed subvariety we just described is equal to $NS\mathcal{W}^{d,-}_{\bar\sigma}$.
\item[(iii)] $B^{d,-}_{\bar\sigma}$ is the inverse image of $\P W^{d,-}_{\bar\sigma}$ with respect to the limit map $F_{(2),\bar\sigma_+}^{(d),1+} \rightarrow F_{(2),\bar\sigma_+}^{(d),1}\cong M^{2,d}_{\bar\sigma_+}$.
\item[(iv)] $S\mathcal{W}^{d,-}_{\bar\sigma} \subset B^{d,-}_{\bar\sigma} $ and $B^{d,-}_{\bar\sigma} \setminus S\mathcal{W}^{d,-}_{\bar\sigma} \subset F_{(1,1),\bar\sigma_-}^{((d-\bar\sigma)/2,(d+\bar\sigma)/2),1+}$. The complement can be described as:
\begin{equation*}
\left \{
\begin{array}{c}
(E,\phi,s) \text{ such that }\\
E=M \oplus L\\
s \in H^0(L), \deg L =(d-\bar\sigma)/2, \deg M=(d+\bar\sigma)/2\\
\text{and } \phi \text{ does not preserve } M
\end{array} \right \}.
\end{equation*}
and therefore contains $NS\mathcal{W}^{d,+}_{\bar\sigma}$. Furthermore, if $\bar\sigma > 2g-2$ then the locally closed subvariety we just described is empty and therefore $NS\mathcal{W}^{d,+}_{\bar\sigma}$ is empty as well.
\end{itemize}
\begin{proof}
(i) follows by the identification $F_{(2),\bar\sigma_-}^{(d),1}\cong M^{2,d}_{\bar\sigma_-}$.

For (ii), the first inclusion follows from the definitions. Note that $B^{d,+}_{\bar\sigma} \setminus S\mathcal{W}^{d,+}_{\bar\sigma}$ consists of triples that are $\bar\sigma_-$-stable and whose underlying pair is also $\bar\sigma_-$-stable. However those triples in $B^{d,+}_{\bar\sigma}$ that are $\bar\sigma_+$-unstable are contained in $S\mathcal{W}^{d,+}_{\bar\sigma}$ and by definition the triples in $B^{d,+}_{\bar\sigma}$ whose underlying Bradlow pair is $\bar\sigma_+$-stable are also contained in $S\mathcal{W}^{d,+}_{\bar\sigma}$. Therefore $B^{d,+}_{\bar\sigma} \setminus S\mathcal{W}^{d,+}_{\bar\sigma}$ consists of triples that are both $\bar\sigma_-$ and $\bar\sigma_+$-stable. The underlying Bradlow pairs are however $\bar\sigma_-$-stable but $\bar\sigma_+$-unstable. Since the maximal destabilizing has degree $(d-\bar\sigma)/2$ and contains the section, it will fit the Harder-Narasimhan type of $F_{(1,1),\bar\sigma_+}^{((d+\bar\sigma)/2,(d-\bar\sigma)/2),2+}$. For the second part we see immediately that the complement has to consist of triples whose underlying pair is both $\bar\sigma_+$-unstable and $\bar\sigma_-$-unstable and so has to be split and all triples in $NS\mathcal{W}^{d,-}_{\bar\sigma}$ are of this form. The statement about the case $\bar\sigma > 2g-2$ follows because, since $\deg M^*LK <0$ then any Higgs field assigned to such a split pair has to preserve the subobject not containing the section (denoted $M$) and therefore any triple in
$$F_{(1,1),\bar\sigma_+}^{((d+\bar\sigma)/2,(d-\bar\sigma)/2),2+} \setminus (B^{d,+}_{\bar\sigma} \setminus S\mathcal{W}^{d,+}_{\bar\sigma})$$
will also lie in $NS\mathcal{W}^{d,-}_{\bar\sigma}$ proving the converse inclusion.

Statement (iii) is proved by using the identification $F_{(2),\bar\sigma_+}^{(d),1}\cong M^{2,d}_{\bar\sigma_+}$ and statement (iv) has an analogous proof to (ii).
\end{proof}
\end{prop}
We conclude the section by giving a more intuitive idea of the previous proposition. Let us examine the flip locus $\mathcal{W}^{d,+}_{\bar\sigma}$ for instance. Then it can be decomposed as $\mathcal{W}^{d,+}_{\bar\sigma}=S\mathcal{W}^{d,+}_{\bar\sigma} \sqcup NS\mathcal{W}^{d,+}_{\bar\sigma}$. $S\mathcal{W}^{d,+}_{\bar\sigma}$ is the intersection of the flip locus with the attracting set $F_{(2),\bar\sigma_-}^{(d),1+}$ while $NS\mathcal{W}^{d,+}_{\bar\sigma}$ is the intersection of the flip locus with the attracting set $F_{(1,1),\bar\sigma_-}^{((d-\bar\sigma)/2,(d+\bar\sigma)/2),1+}$. After $\mathcal{W}^{d,+}_{\bar\sigma}$ has been removed from the moduli space some of the limit points are also removed, but then the limits for the triples in $B^{d,+}_{\bar\sigma} \setminus S\mathcal{W}^{d,+}_{\bar\sigma}$ are replaced by new limit points after the critical value is passed and in fact, for example, $B^{d,+}_{\bar\sigma} \setminus S\mathcal{W}^{d,+}_{\bar\sigma} \subset F_{(1,1),\bar\sigma_+}^{((d+\bar\sigma)/2,(d-\bar\sigma)/2),2+}$ indicating that the limit is in a different attracting set.
\section{The case of low stability parameter}
\subsection{Odd degree}
\label{odddeg}
As we saw in corollary \ref{smoothr2}, $\mathcal{M}^{2,d}_\varepsilon$ is smooth for $d<0$ and for $d$ odd bigger than $4g-4$. In particular in this case $\mathcal{M}^{2,d}_\varepsilon$ is semiprojective and we can use the Bia\l{}ynicki-Birula stratification to compute the motive. We start by introducing a notation:
\begin{defn}
Let $d$ be an odd integer. Denote by $I_1^o(d)$ the set of pairs of integers $(d_1,d_2)$ for which $F_{(1,1),\varepsilon}^{(d_1,d_2),1}$ is a nonempty fixed point component of $\mathcal{M}_\varepsilon^{2,d}$. Analogously, let $I_2^o(d)$ denote the set of pairs of integers $(d_1,d_2)$ for which $F_{(1,1),\varepsilon}^{(d_1,d_2),2}$ is a nonempty fixed point component of $\mathcal{M}_\varepsilon^{2,d}$.
\end{defn}
From theorem \ref{fpr2}, we see that $I_1^o(d)$ is empty if $d <1$ and if $d \geq 1$ consists of the pairs $(d_1,d_2)$ satisfying $d_1+d_2=d$ and the following three equivalent sets of inequalities:
\begin{align*}
&\max \left \{0,\frac{d+1}{2}+1-g\right \} \leq d_1\leq \frac{d-1}{2} \\
&\frac{d+1}{2} \leq d_2 \leq \min \left \{ \frac{d-1}{2}+g-1, d\right \} \\
&\max\{1,2g-2-d\} \leq d_1-d_2+2g-2 \leq 2g-3 \text{ only odd values}.
\end{align*}

$I_2^o(d)$ instead consists of the pairs of integers $(d_1,d_2)$ satisfying $d_1+d_2=d$ and the following three equivalent sets of inequalities:
\begin{align*}
&\frac{d+1}{2}+1-g \leq d_1\leq \min \left \{d,\frac{d-1}{2} \right\} \\
&\max \left \{0,\frac{d+1}{2} \right\} \leq d_2 \leq \frac{d-1}{2}+g-1 \\
&1 \leq d_1-d_2+2g-2 \leq \min\{ d+2g-2,2g-3\} \text{ only odd values}.
\end{align*}

From proposition \ref{deffixed} we know that every nonsplit fixed point is smooth. Since the dimension of the Zariski tangent space can only increase by specialization, we see that $F_{(2),\varepsilon}^{(d),1+}$ always lies in the smooth part of $\mathcal{M}^{2,d}_\varepsilon$. The same is true for split type 2 fixed points and their corresponding attracting sets. In order to compute the motive of those type 1 attracting sets containing singular points we need to work a bit more. Recall from theorem \ref{fpr2} that we have the following description of the type 1 attracting loci:
\begin{equation*}
F_{(1,1),\varepsilon}^{(d_1,d_2),1+}=
\left \{
\begin{array}{c}
(E,\phi,s) \text{ where } E \text{ is defined by}\\
0 \rightarrow E_2 \rightarrow E \rightarrow E_1 \rightarrow 0\\
s \notin H^0(E_2), \deg E_i =d_i\\
\text{and } \phi \text{ does not preserve } E_2
\end{array} \right \}
\end{equation*}
for $(d_1,d_2) \in I_1^o(d)$. We will make use of the following definition.
\begin{defn}
Let $(d_1,d_2)\in I_1^o(d)$. Let us denote by $NSPF_{(1,1),\varepsilon}^{(d_1,d_2),1+}$ the locus of $F_{(1,1),\varepsilon}^{(d_1,d_2),1+}$ where the underlying Bradlow pair lies in $\P W^{d,-}_{d_2-d_1}$ and by $SPF_{(1,1),\varepsilon}^{(d_1,d_2),1+}$ the complement, i.e. the locus where the pair is split.
\end{defn}
With the above definition we clearly have:
\begin{equation*}
NSPF_{(1,1),\varepsilon}^{(d_1,d_2),1+} = B_{d_2-d_1}^{d,-} \cap F_{(1,1),\varepsilon}^{(d_1,d_2),1+}
\end{equation*}
and also:
\begin{equation*}
F_{(1,1),\varepsilon}^{(d_1,d_2),1+} = B_{d_2-d_1}^{d,-} \cap F_{(1,1),\varepsilon}^{(d_1,d_2),1+} \sqcup SPF_{(1,1),\varepsilon}^{(d_1,d_2),1+}.
\end{equation*}

In order to compute the motives we need the following.
\begin{prop}
\begin{equation*}
B_{d_2-d_1}^{d,-} = NSPF_{(1,1),\varepsilon}^{(d_1,d_2),1+} \sqcup S \mathcal{W}^{d,-}_{d_2-d_1} = B_{d_2-d_1}^{d,-} \cap F_{(1,1),\varepsilon}^{(d_1,d_2),1+} \sqcup S \mathcal{W}^{d,-}_{d_2-d_1}.
\end{equation*}
\begin{proof}
The second equality is immediate from the previous observation. From part (iv) of proposition \ref{weirdexpl} we see that $S \mathcal{W}^{d,-}_{d_2-d_1} = B_{d_2-d_1}^{d,-}$, also by definition if $(E,\phi,s) \in NSPF_{(1,1),\varepsilon}^{(d_1,d_2),1+}$ then $(E,s) \in \P W^{d,-}_{d_2-d_1}$ and so by definition \ref{weirdl} $(E,\phi,s) \in B_{d_2-d_1}^{d,-}$. This proves one of the inclusions.

For the reverse inclusion, again from definition \ref{weirdl} we have $(E,\phi,s) \in B_{d_2-d_1}^{d,-}$ if and only if $(E,s) \in \P W^{d,-}_{d_2-d_1}$. If this is the case, then either $\phi$ preserves the subobject of degree $d_2$ that does not contain the section, and then $(E,\phi,s) \in S \mathcal{W}^{d,-}_{d_2-d_1}$ or $\phi$ does not preserve such subobject and then $(E,\phi,s) \in NSPF_{(1,1),\varepsilon}^{(d_1,d_2),1+}$. This proves the second inclusion and also the fact that the union is disjoint.
\end{proof}
\end{prop}
Now we can write the motivic identity:
\begin{equation*}
[F_{(1,1),\varepsilon}^{(d_1,d_2),1+}]=[B_{d_2-d_1}^{d,-}]-[S \mathcal{W}^{d,-}_{d_2-d_1}]+[SPF_{(1,1),\varepsilon}^{(d_1,d_2),1+}].
\end{equation*}
Therefore it suffices to compute the previous three motives. Since we will use $[B_{\bar\sigma}^{d,-}]$ in other sections, we can compute it now.
\begin{prop}
\label{bminus}
Let $\bar\sigma$ be a critical value. Then we have:
\begin{align*}
[B_{\bar\sigma}^{d,-}]&=\L^{1+4(g-1)} [S^{(d-\bar\sigma)/2}( C)] \left ( (\L-1)[S^{(d+\bar\sigma)/2}( C)]+[J^{(d+\bar\sigma)/2}( C)] \right ) + \\&+\L^{1+4(g-1)} \frac{\L^{(d-\bar\sigma)/2}-\L}{\L-1} [S^{(d-\bar\sigma)/2}( C)] \left ( (\L-1)[S^{\bar\sigma}( C)]+[J^{\bar\sigma}( C)] \right ).
\end{align*}
\begin{proof}
Recall that $(E,\phi,s) \in B_{\bar\sigma}^{d,-}$ if and only if $(E,s) \in \P W^{d,-}_{\bar\sigma}$. Therefore consider the composition of maps:
\begin{align*}
B_{\bar\sigma}^{d,-} &\rightarrow \P W^{d,-}_{\bar\sigma} \rightarrow S^{(d-\bar\sigma)/2}( C) \times J^{(d+\bar\sigma)/2}( C)\\
(E,\phi,s) & \mapsto (E,s) \mapsto (L,\bar{s}, M)
\end{align*}
where the second map is defined by using the fact that, if $(E,s) \in \P W^{d,-}_{\bar\sigma}$, then $E$ is an extension:
\begin{equation*}
0 \rightarrow M \rightarrow E \rightarrow L \rightarrow 0
\end{equation*}
with $s \notin H^0(M)$ projecting to $\bar{s} \in H^0(L)$ and $\deg M =(d+\bar\sigma)/2$, $\deg L =(d-\bar\sigma)/2$.

As explained in \cite[proposition 3.3]{thaddeus1994stable}, over $(L,\bar{s}, M)$ the fiber of
$$\P W^{d,-}_{\bar\sigma} \rightarrow S^{(d-\bar\sigma)/2}( C) \times J^{(d+\bar\sigma)/2}( C)$$
is $\P H^0(M \mathcal{O}_D)$ where $D$ is the divisor associated to $\bar{s}$ i.e. $\C\P^{(d-\bar\sigma)/2-1}$. Inside $\P H^0(M \mathcal{O}_D)$ there is the locus corresponding to those pairs whose underlying vector bundle is split. Such a locus is the projectivization of $$\ker \left ( H^0(M \mathcal{O}_D) \rightarrow H^1(L^*M) \right) \cong \quotient{H^0(M)}{H^0(L^*M)}.$$

Over the points of such a kernel the vector bundle will be a split sum and so the fiber of $B_{\bar\sigma}^{d,-} \rightarrow \P W^-_{\bar\sigma}$ is $H^0(K \End E)$ whose dimension, according to lemma \ref{unstableext}, is $2+4(g-1)+\dim H^0(L^*M)$, over the complement instead the fiber is still $H^0(K \End E)$ but the dimension is now $1+4(g-1)+\dim H^0(L^*M)$.

We introduce the stratification:
\begin{equation*}
Z^{d, \bar \sigma}_{i,j}:= \{ (L,\bar{s}, M) \in S^{(d-\bar\sigma)/2}( C) \times J^{(d+\bar\sigma)/2}( C): \dim H^0(L^*M)=i, \dim H^0(M)=j \}.
\end{equation*}

With this definition and the previous observations we can compute the motive of $B_{\bar\sigma}^{d,-}$:
\begin{align*}
[B_{\bar\sigma}^{d,-}]&=\sum_{i,j} [Z^{d, \bar \sigma}_{i,j}] \left ( [\C\P^{j-i-1}] \L^{2+4(g-1)+i} + ([\C\P^{(d-\bar\sigma)/2-1}] - [\C\P^{j-i-1}]) \L^{1+4(g-1)+i} \right )=\\
&=\frac{\L^{1+4(g-1)}}{\L-1} \sum_{i,j} [Z^{d, \bar \sigma}_{i,j}] \left ( \L^{j+1}-\L^{i+1} + \L^{(d-\bar\sigma)/2+i} - \L^{j} \right )
\end{align*}

Now observe that:
\begin{align*}
\sum_{i,j}& [Z^{d, \bar \sigma}_{i,j}] \L^j = [S^{(d-\bar\sigma)/2}( C)] \sum_j [V_j^{(d+\bar\sigma)/2}]\L^j =\\
&=[S^{(d-\bar\sigma)/2}( C)] \left ( (\L-1)[S^{(d+\bar\sigma)/2}( C)]+[J^{(d+\bar\sigma)/2}( C)] \right ).
\end{align*}

In order to compute
\begin{equation*}
\sum_{i,j} [Z^{d, \bar \sigma}_{i,j}] \L^i
\end{equation*}
we can argue as follows. Consider the isomorphism:
\begin{align*}
S^{(d-\bar\sigma)/2}( C) \times J^{(d+\bar\sigma)/2}( C) & \rightarrow S^{(d-\bar\sigma)/2}( C) \times J^{\bar\sigma}( C)\\
(L,\bar{s},M) & \mapsto (L,\bar{s},L^*M).
\end{align*}

The stratum $\cup_j Z^{d, \bar \sigma}_{i,j}$ in $S^{(d-\bar\sigma)/2}( C) \times J^{(d+\bar\sigma)/2}( C) $ is then isomorphic to $S^{(d-\bar\sigma)/2}( C) \times V_i^{\bar\sigma}$ in $S^{(d-\bar\sigma)/2}( C) \times J^{\bar\sigma}( C)$. Therefore:
\begin{align*}
\sum_{i,j}& [Z^{d, \bar \sigma}_{i,j}] \L^i = [S^{(d-\bar\sigma)/2}( C)] \sum_i [ V_i^{\bar\sigma}] \L^i = [S^{(d-\bar\sigma)/2}( C)] \left ( (\L-1)[S^{\bar\sigma}( C)]+[J^{\bar\sigma}( C)] \right ).
\end{align*}

Putting all together, we get
\begin{align*}
[B_{\bar\sigma}^{d,-}]&=\L^{1+4(g-1)} [S^{(d-\bar\sigma)/2}( C)] \left ( (\L-1)[S^{(d+\bar\sigma)/2}( C)]+[J^{(d+\bar\sigma)/2}( C)] \right ) + \\&+\L^{1+4(g-1)} \frac{\L^{(d-\bar\sigma)/2}-\L}{\L-1} [S^{(d-\bar\sigma)/2}( C)] \left ( (\L-1)[S^{\bar\sigma}( C)]+[J^{\bar\sigma}( C)] \right ).
\end{align*}
\end{proof}
\end{prop}
\begin{prop}
We have the identity:
\begin{equation*}
[SPF_{(1,1),\varepsilon}^{(d_1,d_2),1+}]=\L^{3g-1+d_2-d_1} [S^{d_1} ( C)] [S^{d_1-d_2+2g-2}( C)].
\end{equation*}
\begin{proof}
Recall that $SPF_{(1,1),\varepsilon}^{(d_1,d_2),1+}$ consists of those pairs in the form $(E_1\oplus E_2,s,\phi)$ where $\deg E_i=d_i$, $s \in H^0(E_1)$ and $\phi$ does not preserve $E_2$. First of all, note that the moduli space of such split pairs is
\begin{equation*}
S^{d_1} ( C) \times J^{d_2} (C ).
\end{equation*}

Second, we can always write the Higgs field into matrix form:
\begin{equation*}
\phi = 
\begin{pmatrix}
\phi_{11} & \phi_{21}\\
\phi_{12} & \phi_{22}
\end{pmatrix}
\end{equation*}
where
\begin{equation*}
\phi_{21} : E_2 \rightarrow E_1 \otimes K
\end{equation*}
has to be nonzero, and
\begin{equation*}
\phi_{12} : E_1 \rightarrow E_2 \otimes K.
\end{equation*}

Since the endomorphisms of the underlying pair are just diagonal (possibly non scalar) then we can see that $(E_1 \oplus E_2, s , \phi)$ and $(E_1 \oplus E_2, s , \phi')$ are isomorphic if and only if there exists $\lambda \in \C^*$ such that:
\begin{equation*}
\phi' = 
\begin{pmatrix}
\phi'_{11} & \phi'_{21}\\
\phi'_{12} & \phi'_{22}
\end{pmatrix}=
\begin{pmatrix}
\phi_{11} & \lambda\phi_{21}\\
\lambda^{-1}\phi_{12} & \phi_{22}
\end{pmatrix}.
\end{equation*}

Let's consider the map:
\begin{equation*}
SPF_{(1,1),\varepsilon}^{(d_1,d_2),1+} \rightarrow S^{d_1} ( C) \times J^{d_2} (C ) \times H^0(K)^2
\end{equation*}
forgetting the off-diagonal parts of the Higgs field. Then, by what we said, the fiber of such a map over $(E_1,s,E_2)$ will be:
\begin{equation*}
\quotient{\left (H^0(KE_1 E_2^*) \setminus \{0\} \right ) \times H^0(K E_2 E_1^*)}{ \C^*}.
\end{equation*}

This allows us to compute the motive of $SPF_{(1,1),\sigma}^{(d_1,d_2),1+}$ by stratifying $S^{d_1} ( C) \times J^{d_2} (C )$ according to the dimension of $H^0(KE_1 E_2^*)$. Note in fact that the dimension of $H^0(K E_2 E_1^*)$ is constantly equal to $g-1+d_2-d_1$. Call $V^{d_1-d_2+2g-2}_i$ the stratum of $J^{d_1-d_2+2g-2}$ where the dimension of the global sections of the line bundle is equal to $i$. Then we have:
\begin{align*}
[SPF_{(1,1),\varepsilon}^{(d_1,d_2),1+}]=&\sum_i \left (\L^{2g} [S^{d_1} ( C)] [V^{d_1-d_2+2g-2}_i] \L^{d_2-d_1+g-1} [\C\P^{i-1} ]\right )=\\
=&\L^{3g-1+d_2-d_1} [S^{d_1} ( C)] [S^{d_1-d_2+2g-2}( C)].
\end{align*}
\end{proof}
\end{prop}
Lastly, we have the following.
\begin{prop}
Let $\bar\sigma$ be a critical value. The following motivic equalities hold:
\begin{align*}
[NS\mathcal{W}_{\bar\sigma}^{d,-}]&=\L^{2g} [S^{(d-\bar\sigma)/2} ( C)] [J ( C)] [\C\P^{g-2+\bar\sigma}]\\
[S\mathcal{W}_{\bar\sigma}^{d,-}]&=\L^{2g} [S^{(d-\bar\sigma)/2} ( C)] [J ( C)] \cdot ([\C\P^{(d+\bar\sigma)/2+g-2}]-[\C\P^{g-2+\bar\sigma}]).
\end{align*}
\begin{proof}
From part (vi) of proposition \ref{weirdexpl1} we know that there exists a map $NS\mathcal{W}_{\bar\sigma}^{d,-} \rightarrow X^d_{\bar\sigma}$ whose fiber over $(L,\bar{s},M,\psi_1,\psi_2)$ is $\P H^0(KL^*M)$. Note that the dimension of the fibers is constant and is always $g-2+\bar\sigma$. Therefore:
\begin{equation*}
[NS\mathcal{W}_{\bar\sigma}^{d,-}]=\L^{2g} [S^{(d-\bar\sigma)/2} ( C)] [J ( C)] [\C\P^{g-2+\bar\sigma}]
\end{equation*}
and
\begin{equation*}
[S\mathcal{W}_{\bar\sigma}^{d,-}]=[\mathcal{W}_{\bar\sigma}^{d,-}]-[NS\mathcal{W}_{\bar\sigma}^{d,-}]=\L^{2g} [S^{(d-\bar\sigma)/2} ( C)] [J ( C)] \cdot ([\C\P^{(d+\bar\sigma)/2+g-2}]-[\C\P^{g-2+\bar\sigma}]).
\end{equation*}
\end{proof}
\end{prop}
We can now have a formula for the type 1 attracting sets.
\begin{prop}
\begin{align*}
&[F_{(1,1),\varepsilon}^{(d_1,d_2),1+}]=\L^{3g-2+d_2} [S^{d_1}( C)] [S^{d_1-d_2+2g-2}( C)]+\\
&+(\L^{4g-2}-\L^{4g-3}) [S^{d_1}( C)] \left ( [S^{d_2}( C)]-[J( C)] \frac{\L^{d_2+1-g}-1}{\L-1} \right ).
\end{align*}
\begin{proof}
We can compute:
\begin{align*}
&[F_{(1,1),\varepsilon}^{(d_1,d_2),1+}]=[B_{d_2-d_1}^{d,-}]-[S \mathcal{W}^{d,-}_{d_2-d_1}]+[SPF_{(1,1),\varepsilon}^{(d_1,d_2),1+}]=\\
&=\L^{1+4(g-1)} [S^{d_1}( C)] \left ( (\L-1)[S^{d_2}( C)]+[J( C)] \right ) + \\
&+\L^{1+4(g-1)} \frac{\L^{d_1}-\L}{\L-1} [S^{d_1}( C)] \left ( (\L-1)[S^{d_2-d_1}( C)]+[J( C)] \right )+\\
&-\L^{2g} [S^{d_1} ( C)] [J ( C)] \cdot ([\C\P^{d_2+g-2}]-[\C\P^{g-2+d_2-d_1}])+\\
&+\L^{3g-1+d_2-d_1} [S^{d_1} ( C)] [S^{d_1-d_2+2g-2}( C)]=\\
&=(\L^{4g-2}-\L^{4g-3})  [S^{d_1}( C)]  [S^{d_2}( C)]+\L^{4g-3}  [S^{d_1}( C)]  [J( C)]+\\
&+(\L^{4g-3+d_1}-\L^{4g-2})  [S^{d_1}( C)]  [S^{d_2-d_1}( C)]+\L^{4g-3} \frac{\L^{d_1}-\L}{\L-1}[S^{d_1}( C)]  [J( C)]+\\
&-\L^{3g-1+d_2-d_1}\frac{\L^{d_1}-1}{\L-1}[S^{d_1}( C)]  [J( C)]+\\
&+ \L^{3g-1+d_2-d_1} [S^{d_1} ( C)] [S^{d_1-d_2+2g-2}( C)]=\\
&=(\L^{4g-2}-\L^{4g-3})  [S^{d_1}( C)]  [S^{d_2}( C)]+\\
&+ (\L^{4g-3}-\L^{3g-1+d_2-d_1})\frac{\L^{d_1}-1}{\L-1}[S^{d_1}( C)]  [J( C)]+\\
&+ [S^{d_1}( C)]\left ( (\L^{4g-3+d_1}-\L^{4g-2})  [S^{d_2-d_1}( C)] + \L^{3g-1+d_2-d_1} [S^{d_1-d_2+2g-2}( C)] \right).
\end{align*}
From Serre's duality we have the identity:
\begin{equation*}
[S^{d_2-d_1}( C)]=\L^{d_2-d_1-g+1}[S^{d_1-d_2+2g-2}( C)]+\frac{\L^{d_2-d_1-g}-1}{\L-1}[J( C)]
\end{equation*}
for $0 \leq d_2-d_1 \leq 2g-2$. Using this we can proceed in the computation:
\begin{align*}
&[F_{(1,1),\varepsilon}^{(d_1,d_2),1+}]=(\L^{4g-2}-\L^{4g-3}) [S^{d_1}( C)]  [S^{d_2}( C)]+ \\
&+(\L^{4g-3}-\L^{3g-1+d_2-d_1})\frac{\L^{d_1}-1}{\L-1}[S^{d_1}( C)]  [J( C)]+\\
&+ (\L^{4g-3+d_1}-\L^{4g-2})\frac{\L^{d_2-d_1+1-g}-1}{\L-1}[S^{d_1}( C)]  [J( C)]+\\
&+(\L^{4g-3+d_1}-\L^{4g-2})\L^{d_2-d_1+1-g}[S^{d_1}( C)][S^{d_1-d_2+2g-2}( C)]+\\
&+\L^{3g-1+d_2-d_1} [S^{d_1}( C)] [S^{d_1-d_2+2g-2}( C)]=\\
&=(\L^{4g-2}-\L^{4g-3}) [S^{d_1}( C)]  [S^{d_2}( C)]+ \\
&-(\L^{4g-2}-\L^{4g-3})\frac{\L^{d_2+1-g}-1}{\L-1}  [S^{d_1}( C)]  [J( C)]+\\
&+\L^{3g-2+d_2} [S^{d_1}( C)] [S^{d_1-d_2+2g-2}( C)]=\\
&=\L^{3g-2+d_2} [S^{d_1}( C)] [S^{d_1-d_2+2g-2}( C)]+\\
&+(\L^{4g-2}-\L^{4g-3}) [S^{d_1}( C)] \left ([S^{d_2}( C)]-[J( C)] \frac{\L^{d_2+1-g}-1}{\L-1} \right).
\end{align*}
\end{proof}
\end{prop}
\begin{rmk}
Observe that, according to proposition \ref{deffixed}, the attracting set $F_{(1,1),\varepsilon}^{(d_1,d_2),1+}$ lies in the smooth part if and only if $d_2>2g-2$. In that case we can apply \cite[Theorem 4.1]{bialynicki1973some} and deduce that $F_{(1,1),\varepsilon}^{(d_1,d_2),1+} \rightarrow F_{(1,1),\varepsilon}^{(d_1,d_2),1}$ is a smooth affine fibration of rank $3g-2+d_2$ (see section \ref{rank2disc}). In this case we have $[S^{d_2}( C)]=[J( C)][\C\P^{d_2-g}]$ and hence:
\begin{align*}
&[F_{(1,1),\varepsilon}^{(d_1,d_2),1+}]= \L^{3g-2+d_2} [S^{d_1}( C)] [S^{d_1-d_2+2g-2}( C)]
\end{align*}
which is the same as computing $F_{(1,1),\varepsilon}^{(d_1,d_2),1+}$ as an affine fibration.
\end{rmk}
We can summarize the section in the following.
\begin{theorem}
\label{motodd}
For $d \geq 0 $ odd and $0<\varepsilon<1$ we have:
\begin{align*}
[\mathcal{M}_\varepsilon^{2,d}]=&\L^{1+4(g-1)}[M_\varepsilon^{2,d}]+\sum_{(d_1,d_2) \in I_1^o(d)} \L^{1+3(g-1)+d_2}[S^{d_1}( C)][S^{d_1-d_2+2g-2}]+\\
&+\sum_{(d_1,d_2) \in I_2^o(d)}\L^{1+4(g-1)} [S^{d_2}( C)][S^{d_1-d_2+2g-2}( C)]+\\
&+\sum_{(d_1,d_2) \in I_1^o(d)} (\L^{4g-2}-\L^{4g-3}) [S^{d_1}( C)] \left ([S^{d_2}( C)]-[J( C)] \frac{\L^{d_2+1-g}-1}{\L-1}\right).
\end{align*}
\end{theorem}
\subsection{The moduli spaces for even degree}
For even degree $d$, the situation is quite different due to the presence of strictly semistable Higgs bundles. In fact, for $0<\varepsilon<2$ and even $d$, it is clear that for every $\varepsilon$-stable triple $(E,\phi,s)$ the underlying Higgs bundle $(E,\phi)$ is semistable. Therefore, there is still an Abel-Jacobi map $\mathcal{M}_\varepsilon^{2,d} \rightarrow \mathcal{M}^{2,d}$.

We can say something a bit more precise for $\varepsilon$-stability for even degree. It is easy to see that $(E,\phi,s)$ is $\varepsilon$-stable if and only if $(E,\phi)$ is semistable and if $L \subset E$ is a $\phi$-invariant line subbundle of $E$ with $\deg L = \deg E /2$, then $s \notin H^0(L)$.

Furthermore, the fibers of the Abel-Jacobi map are more complicated than just projectivized spaces of global sections because of the fact that the equivalence relation for the moduli spaces is isomorphism for triples and S-equivalence for Higgs bundles.\\
However, it is still true that above the stable locus of $\mathcal{M}^{2,d}$ the Abel-Jacobi map behaves exactly as in the case of odd degree.

Let us outline the strategy for the computation of $[\mathcal{M}_\varepsilon^{2,d}]$.
\begin{defn}
Let $0<\varepsilon<2$ and $d$ be an even integer. We define two locally closed subvarieties of $\mathcal{M}_\varepsilon^{2,d}$ by:
\begin{equation*}
 \mathscr{X}^d_1=\left \{\begin{array}{c}
(E,\phi,s) \in \mathcal{M}_\varepsilon^{2,d} \text{ such that}\\
(E,s) \text{ is $\varepsilon$-stable and}\\
E \text{ is strictly semistable}
\end{array} \right \}
\end{equation*}
and
\begin{equation*}
 \mathscr{X}^d_2=\left \{\begin{array}{c}
(E,\phi,s) \text{ such that $E$ fits into an exact sequence}\\
0 \rightarrow L \rightarrow E \rightarrow M \rightarrow 0\\
\text{with $\deg L=\deg M$, $s \in H^0(L)$ and}\\
\phi \text{ does not preserve } L
\end{array} \right \}.
\end{equation*}
\end{defn}
First we examine the attracting sets of $\mathcal{M}_\varepsilon^{2,d}$ more closely. For the following, refer to theorem \ref{fpr2} and  \ref{rank2disc}.

The attracting set corresponding to nonsplit fixed points is $F_{(2),\varepsilon}^{(d),1+}$ and is nonempty if and only if $d \geq 0$. There is a limit map $F_{(2),\varepsilon}^{(d),1+} \rightarrow M_\varepsilon^{2,d}$ whose fiber over a pair $(E,s)$ is $H^0(K  \End E)$. We can write:
\begin{equation*}
M_\varepsilon^{2,d}=M_{\varepsilon,st}^{2,d} \sqcup M_{\varepsilon,ss}^{2,d}
\end{equation*} 
where the first subvariety corresponds to the locus of pairs having underlying stable vector bundle, while the second one corresponds to the locus of pairs having strictly semistable vector bundle. Then $F_{(2),\varepsilon}^{(d),1+} \rightarrow M_\varepsilon^{2,d}$ is an affine fibration over $M_{\varepsilon,st}^{2,d}$ with fibers of dimension $4g-3$. The inverse image of $F_{(2),\varepsilon}^{(d),1+} \rightarrow M_\varepsilon^{2,d}$ over $M_{\varepsilon,ss}^{2,d}$ is instead precisely $\mathscr{X}^d_1$ defined above.

For the split cells we can introduce a notation as for the odd degree case.
\begin{defn}
Let $d \geq 0$ be an even integer. Denote by $I_1^e(d)$ the set of pairs of integers $(d_1,d_2)$ for which $F_{(1,1),\varepsilon}^{(d_1,d_2),1}$ is a nonempty fixed point component of $\mathcal{M}_\varepsilon^{2,d}$. Analogously, let $I_2^e(d)$ denote the set of pairs of integers $(d_1,d_2)$ for which $F_{(1,1),\varepsilon}^{(d_1,d_2),2}$ is a nonempty fixed point component of $\mathcal{M}_\varepsilon^{2,d}$
\end{defn}
From theorem \ref{fpr2}, we see that $I_1^e(d)$ is nonempty if and only if $d \geq 2$ and consists of pairs of integers $(d_1,d_2)$ satisfying $d_1+d_2=d$ and the following three equivalent sets of inequalities:
\begin{align*}
&\max \left \{0, d/2+1-g\right \} \leq d_1\leq d/2-1 \\
&d/2+1 \leq d_2 \leq \min \left \{ d, d/2+g-1\right \} \\
&\max\{0,2g-2-d\} \leq d_1-d_2+2g-2 \leq 2g-4 \text{ only even values}
\end{align*}
and $I_2^e(d)$ is the set of pairs of integers $(d_1,d_2)$ satisfying $d_1+d_2=d$ and the following three equivalent sets of inequalities:
\begin{align*}
&d/2+1-g \leq d_1\leq \min\{ d/2,d\} \\
&\max\{0,d/2\} \leq d_2 \leq d/2+g-1 \\
&0 \leq d_1-d_2+2g-2 \leq \min\{2g-2,2g-2+d\} \text{ only even values}.
\end{align*}

As we proved in theorem \ref{shatz}, the decomposition into attracting cells coincides with the decomposition according to the Harder-Narasimhan type of the underlying Bradlow pair. From this we see that $F_{(1,1),\varepsilon}^{(d/2,d/2),2+}$ corresponds to the Harder-Narasimhan type $(d/2,d/2)$ i.e. the destabilizing subobject of the pair $(E,s)$ is a line bundle of degree $d/2$ that contains the section $s$. Therefore $F_{(1,1),\varepsilon}^{(d/2,d/2),2+}=\mathscr{X}^d_2$.

If we denote by $\tilde{I}_2^e(d)=I_2^e(d) \setminus \{(d/2,d/2)\}$ then the previous discussion allows us to write the following preliminary formula for $d \geq 0$. Note that the computation of the type 1 attracting set is exactly the same as in section \ref{odddeg}.
\begin{align*}
[\mathcal{M}_\varepsilon^{2,d}]&=\L^{1+4(g-1)}  [M^{2,d}_{\varepsilon,st}]+\sum_{(d_1,d_2) \in I_1^e(d)} \L^{1+3(g-1)+d_2} [S^{d_1}( C)][S^{d_1-d_2+2g-2}( C)] +\\
&+\sum_{(d_1,d_2) \in \tilde{I}_2^e(d)}\L^{1+4(g-1)} [S^{d_2}( C)][S^{d_1-d_2+2g-2}( C)]+\\
&+[\mathscr{X}^d_1]+[\mathscr{X}^d_2]+\\
&+\sum_{(d_1,d_2) \in I_1^e(d)} (\L^{4g-2}-\L^{4g-3}) [S^{d_1}( C)] ([S^{d_2}( C)]-[J( C)][\C\P^{d_2-g}]).
\end{align*}

%The case $2-2g \leq d < 0$ will be dealt in the next chapter\marginpar{reference}.\\
In order to complete the above formula, we are left with computing $[M^{2,d}_{\varepsilon,st}]$, $[M^{2,d}_{\varepsilon,ss}]$, $[\mathscr{X}^d_1]$ and $[\mathscr{X}^d_2]$. Let us first state some technical lemmas that we will need in the computation.
\begin{lemma}
\label{dim_end_even}
Let $E$ be a rank 2 vector bundle over $C$ fitting in the following exact sequence:
\begin{equation*}
0 \rightarrow L \rightarrow E \rightarrow M \rightarrow 0
\end{equation*}
with $L$ and $M$ line bundles of the same degree.

Then:
\begin{equation*}
\dim H^0(\End E)= \begin{cases}
1 & \text{if $L \neq M$ and the extension is nonsplit}\\
2 & \text{if $L \neq M$ and the extension is split}\\
2 & \text{if $L = M$ and the extension is nonsplit}\\
4 & \text{if $L = M$ and the extension is split}
\end{cases}
\end{equation*}
\begin{proof}
The statement is clear in the last case, where $\End E = \mathcal{O}^4$ and in the second case where $\End E = \mathcal{O}^2\oplus L^*M \oplus M^*L$ since $M^*L$ and $L^*M$ have global sections if and only if $L=M$.

Suppose now that the extension is nonsplit and that $L \neq M$. Consider the exact sequence:
\begin{equation*}
0 \rightarrow \mathcal{O} \rightarrow E L^* \rightarrow ML^* \rightarrow 0.
\end{equation*}

Since $H^0(L^*M)=0$, the long exact sequence obtained by applying $H^0$ to the above short exact sequence will yield $H^0(EL^*)=H^0( \mathcal{O})=\C$. Consider now:
\begin{equation*}
0 \rightarrow LM^* \rightarrow E M^* \rightarrow \mathcal{O} \rightarrow 0.
\end{equation*}

Again, $H^0(LM^*)=0$ therefore we get an injection $H^0(EM^*) \rightarrow H^0( \mathcal{O})$. Such a map however sends $\gamma : M \rightarrow E$ to $p \gamma$ where $p: E \rightarrow M$ is the projection in the exact sequence for $E$.

If $H^0(EM^*) \neq 0$ then there would be a map $\gamma : M \rightarrow E$ such that $p \gamma \neq 0$ i.e. a splitting of the sequence which is a contradiction. Therefore $H^0(EM^*) = 0$.\\
We can conclude this case by observing that there is a further exact sequence:
\begin{equation*}
0 \rightarrow EM^* \rightarrow \End E \rightarrow EL^* \rightarrow 0
\end{equation*}
yielding the long exact sequence:
\begin{equation*}
0 \rightarrow H^0(EM^*) \rightarrow H^0(\End E) \rightarrow H^0(EL^*) \rightarrow \dots
\end{equation*}

Note that the map $H^0(\End E) \rightarrow H^0(EL^*)$ sends $f: E \rightarrow E$ to $fi : L \rightarrow E$ where $i : L \rightarrow E$ is the inclusion map in the exact sequence for $E$. Since $\id_E$ is sent to $i : L \rightarrow E$ we conclude that the map is surjective and therefore $H^0(\End E)=\C$.

When the sequence is nonsplit and $L=M$, from the sequence:
\begin{equation*}
0 \rightarrow LM^* \rightarrow E M^* \rightarrow \mathcal{O} \rightarrow 0,
\end{equation*}
arguing as before, we deduce that $H^0(EM^*)=\C$. However this also implies that $H^0(EL^*)= \C$ and therefore from the sequence
\begin{equation*}
0 \rightarrow EM^* \rightarrow \End E \rightarrow EL^* \rightarrow 0
\end{equation*}
we can deduce that $H^0(\End E)= \C^2$ because $H^0( \End E ) \rightarrow H^0(EL^*)$ is still surjective.
\end{proof}
\end{lemma}
\begin{lemma}
Let $n \geq 0$ be an integer and let $\Gr(2,n+2)$ denote the Grassmannian variety of vector spaces of dimension two inside $\C^{n+2}$. Then:
\begin{equation*}
[\Sym^2 ( \C\P^n)]=[\Gr (2,n+2)]=\frac{(\L^{n+2}-1)(\L^{n+1}-1)}{(\L^2-1)(\L-1)}.
\end{equation*}
\begin{proof}
From \cite[Lemma 4.4]{gottsche2000motive} we deduce, by putting $X=\{ pt \}$, that $\Sym^n(\L^m)=\L^{mn}$ in the ring of motives. Therefore:
\begin{align*}
[\Sym^2 ( \C\P^n)]&= \sum_{0 \leq i < j \leq n} \L^{i+j} + \sum_{j=0}^n \L^{2j}=\sum_{j=1}^n \L^j \left ( \frac{\L^j-1}{\L-1}\right ) +  \frac{\L^{2n+2}-1}{\L^2-1}=\\
&= \sum_{j=0}^n \frac{\L^{2j}}{\L-1}- \sum_{j=0}^n \frac{\L^{j}}{\L-1}+\frac{\L^{2n+2}-1}{\L^2-1}=\\
&=\frac{\L^{2n+2}-1}{(\L^2-1)(\L-1)}-\frac{\L^{n+1}-1}{(\L-1)^2}+\frac{\L^{2n+2}-1}{\L^2-1}=\\
&=\frac{\L^{2n+2}-1-\L^{n+2}+\L-\L^{n+1}+1+\L^{2n+3}-\L^{2n+2}-\L+1}{(\L^2-1)(\L-1)}=\\
&=\frac{(\L^{n+2}-1)(\L^{n+1}-1)}{(\L^2-1)(\L-1)}.
\end{align*}

To compute the motive of $\Gr(2,n+2)$ we can view it as a global quotient of:
\begin{equation*}
\{(u,v) \in \C^{n+2} : u \neq v, u \neq 0 \neq v\}
\end{equation*}
by $GL_2$. Therefore:
\begin{align*}
[\Gr(2, n+2)]&=\frac{(\L^{n+2}-1)(\L^{n+2}-\L)}{(\L^2-1)(\L^2-\L)}=\frac{(\L^{n+2}-1)(\L^{n+1}-1)}{(\L^2-1)(\L-1)}
\end{align*}
and the result follows.
\end{proof}
\end{lemma}
\begin{lemma}
Let $E$ be a rank two vector bundle fitting in the nonsplit exact sequence:
\begin{equation*}
0 \rightarrow L \rightarrow E \rightarrow M \rightarrow 0
\end{equation*}
with $L$ and $M$ line bundles and $\deg L = \deg M$.

Then, the locus of $H^0(K \End E)$ consisting of the Higgs fields preserving the subobject $L$ is a vector subspace of $H^0(K \End E)$ of dimension $3g-1$.
\begin{proof}
Let us call $i$ and $p$ the inclusion and projection in the exact sequence defining $E$. First observe that $\phi \in H^0(K \End E)$ preserves $L$ if and only if $p \phi i=0$. This means that the locus we are looking for is the kernel of the map:
\begin{align*}
H^0(K \End E) &\rightarrow H^0(KL^*M)\\
\phi & \mapsto p \phi i.
\end{align*}

Furthermore, such a map factors as:
\begin{align*}
H^0(K \End E) &\rightarrow H^0(KE^*M)\\
\phi & \mapsto p \phi
\end{align*}
followed by
\begin{align*}
H^0(K E^*M) &\rightarrow H^0(KL^*M)\\
\gamma & \mapsto \gamma i.
\end{align*}

Suppose that $L \neq M$. As we computed in \ref{dim_end_even}, $\dim H^0(K \End E)=1+4(g-1)$, $\dim H^0(K E^* L)=2g-1$, $\dim H^0(KE^* M)=2g-2$. From the exact sequence:
\begin{equation*}
0 \rightarrow H^0(K E^*L) \rightarrow H^0( K \End E) \rightarrow H^0(K E^* M) \rightarrow \dots
\end{equation*}
we see that $\dim \im (H^0( K \End E) \rightarrow H^0(K E^* M))=2g-2=\dim H^0(K E^* M)$. Therefore $H^0( K \End E) \rightarrow H^0(K E^* M)$ is surjective and this implies:
\begin{align*}
&\dim \im (H^0(K \End E) \rightarrow H^0(KL^*M))=\dim \im (H^0(K E^*M) \rightarrow H^0(KL^*M)).
\end{align*}

Furthermore, from the sequence:
\begin{equation*}
0 \rightarrow H^0(K) \rightarrow H^0( K E^*M) \rightarrow H^0(K L^* M) \rightarrow \dots
\end{equation*}
we deduce that $\dim \im (H^0(K E^*M) \rightarrow H^0(KL^*M))=g-2$.

To conclude, we have:
\begin{align*}
&\dim \ker (H^0(K \End E) \rightarrow H^0(KL^*M) )=\\
&=1+4(g-1) - \dim \im (H^0(K E^*M) \rightarrow H^0(KL^*M))=3g-1.
\end{align*}

In the case when $L=M$ the proof works in the same way but $\dim H^0(K \End E)=2+4(g-1)$, $\dim H^0(K E^* L)=2g-1$, $\dim H^0(KE^* M)=2g-1$. Here $H^0( K \End E) \rightarrow H^0(K E^* M)$ is still surjective and $\dim \im (H^0(K E^*M) \rightarrow H^0(KL^*M))=g-1$ from which we deduce that the dimension of our kernel is again $3g-1$.
\end{proof}
\end{lemma}
Let us start by defining a decomposition of  $ \mathscr{X}^d_1$ and $M_{\varepsilon,ss}^{2,d}$ that we will use in the computation of the motive.
\begin{defn}
Let $(E,\phi,s)$ be a triple in $ \mathscr{X}^d_1$. Then there are three possibilities for $E$: it is a nonsplit extension of two line bundles of the same degree, it is a split extension of two different line bundles of the same degree or it is a split extension of two copies of the same line bundle. Accordingly we define:
\begin{equation*}
\mathscr{U}^d_1=\left \{\begin{array}{c}
(E,\phi,s) \text{ such that $E$ fits into a nonsplit exact sequence}\\
0 \rightarrow L \rightarrow E \rightarrow M \rightarrow 0\\
\text{with $\deg L=\deg M$, $s \notin H^0(L)$ and}\\
\phi \in H^0(K\End E) \text{ arbitrary}
\end{array}\right \}
\end{equation*}
covering the first case,
\begin{equation*}
\mathscr{V}^d_1=\left \{\begin{array}{c}
(E,\phi,s) \text{ such that } E=L \oplus M \\
\text{with $\deg L=\deg M$, $s \notin H^0(L)$ and $s \notin H^0(M)$, $L \neq M$ and}\\
\phi \in H^0(K\End E) \text{ arbitrary}
\end{array}\right \}
\end{equation*}
covering the second and finally
\begin{equation*}
\mathscr{W}^d_1=\left \{\begin{array}{c}
(E,\phi,s) \text{ such that } E=L \oplus L \\
\text{$s= u\oplus v$ with $u$, $v$ linearly independent and}\\
\phi \in H^0(K\End E) \text{ arbitrary}
\end{array}\right \}
\end{equation*}
covering the third.

Similarly for $M_{\varepsilon,ss}^{2,d}$ we define:
\begin{equation*}
\widetilde{\mathscr{U}}^d_1=
\left \{\begin{array}{c}
(E,\phi,s) \text{ such that $E$ fits into a nonsplit exact sequence}\\
0 \rightarrow L \rightarrow E \rightarrow M \rightarrow 0\\
\text{with $\deg L=\deg M$, $s \notin H^0(L)$}
\end{array}\right \}
\end{equation*}
covering the first case,
\begin{equation*}
\widetilde{\mathscr{V}}^d_1=\left \{\begin{array}{c}
(E,\phi,s) \text{ such that } E=L \oplus M \\
\text{with $\deg L=\deg M$, $L \neq M$, $s \notin H^0(L)$ and $s \notin H^0(M)$ }
\end{array}\right \}
\end{equation*}
covering the second and finally
\begin{equation*}
\widetilde{\mathscr{W}}^d_1=\left \{\begin{array}{c}
(E,\phi,s) \text{ such that } E=L \oplus L \\
\text{$s= u\oplus v$ with $u$, $v$ linearly independent}
\end{array}\right \}
\end{equation*}
covering the third.
\end{defn}
\begin{rmk}
A few comments are in order for the previous definitions. Suppose that $0 < \varepsilon <2$. For a nonsplit exact sequence:
\begin{equation*}
0 \rightarrow L \rightarrow E \rightarrow M \rightarrow 0
\end{equation*}
with $\deg L=\deg M$ it is clear that $(E,s)$ is $\varepsilon$-stable if and only if $s \notin L$. For $E=L \oplus M$ with $L \neq M$ instead, since $\Aut E = (\C^*)^2$, we see that $(E,s)$ is $\varepsilon$-stable if and only if $s$ is not concentrated in $L$ nor $M$, which is the same as saying that $s=u \oplus v$ with $0 \neq u \in H^0(L)$ and $0 \neq v \in H^0(M)$. Finally, when $E=L \oplus L$ then $\Aut E = \GL_4(\C)$ so $(E,s)$ is $\varepsilon$-stable if and only if $s= u \oplus v$ with $u, v \in H^0(L)$ generating a dimension 2 subspace.
\end{rmk}
%Recall also the definition of the Brill-Noether loci \marginpar{reference}, that will be used in the following.
%
\begin{prop}
\label{motiveU}
We have the following motivic equalities:
\begin{align*}
[\mathscr{U}^d_1]&=\L^{2+4(g-1)}\sum_{j \geq 1} [V_j^{d/2}] \left ( \frac{\L^j-1}{\L-1} \right )\left ( \frac{\L^{d/2}-1}{\L-1}-\frac{\L^{j-1}-1}{\L-1}\right )+\\
&+\L^{1+4(g-1)}\sum_{j \geq 0} [V_j^{d/2}] \left ( [S^{d/2}( C)]-\frac{\L^j-1}{\L-1} \right )\left ( \frac{\L^{d/2}-1}{\L-1}-\frac{\L^{j}-1}{\L-1}\right ).
\end{align*}
and
\begin{align*}
[\widetilde{\mathscr{U}}^d_1]&=\sum_{j \geq 1} [V_j^{d/2}] \left ( \frac{\L^j-1}{\L-1} \right )\left ( \frac{\L^{d/2}-1}{\L-1}-\frac{\L^{j-1}-1}{\L-1}\right )+\\
&+\sum_{j \geq 0} [V_j^{d/2}] \left ( [S^{d/2}( C)]-\frac{\L^j-1}{\L-1} \right )\left ( \frac{\L^{d/2}-1}{\L-1}-\frac{\L^{j}-1}{\L-1}\right ).
\end{align*}
\begin{proof}
Let $(E,s) \in \widetilde{\mathscr{U}}^d_1$ and $E$ be defined by a nonsplit exact sequence:
\begin{equation*}
0 \rightarrow L \rightarrow E \rightarrow M \rightarrow 0.
\end{equation*}

Then, by Lemma \ref{dim_end_even}, we know that $\dim H^0( \End E)=1$ if $L \neq M$ and $2$ if $L=M$. Note that in both cases, if $p: E \rightarrow M$ is the projection map in the exact sequence defining $E$, $p(s) \in H^0(M)$ is preserved by the action of invertible global endomorphisms of $E$. In particular this implies that we have a well defined map:
\begin{align*}
\widetilde{\mathscr{U}}^d_1 &\rightarrow S^{d/2} ( C) \times J^{d/2} ( C)\\
(E,s) & \mapsto (M,p(s),L).
\end{align*}

As we can deduce easily from the exact sequence:
\begin{equation*}
0 \rightarrow H^0(M^*L) \rightarrow H^0(L) \rightarrow H^0(L \mathcal{O}_D ) \rightarrow H^1(M^*L) \rightarrow H^1(L) \rightarrow 0
\end{equation*}
where $D$ is the divisor of $p(s)$, the fiber of the previous map above $(M,p(s),L)$ is $$\C\P^{d/2-1} \setminus \C\P^{\dim H^0(L) - \dim H^0(M^*L) -1}.$$

Recall that in this case $\P H^0(L \mathcal{O}_D)$ parametrizes isomorphism classes of pairs $(E,s)$ that have $L$ as a subobject, $M$ as a quotient and the section $s$ projects to the one fixed for $M$.

We stratify $S^{d/2} ( C) \times J^{d/2} ( C)$ according to the two dimensions appearing in the fibers of the above map. Define:
\begin{equation*}
Z_{ij}^d:=\{ (M, p(s), L) \in S^{d/2} ( C) \times J^{d/2} ( C): \dim H^0(M^*L) =i,\ \dim H^0(L)=j \}.
\end{equation*}

Note that $\dim H^0(M^*L)$ can only be 1 if $L=M$ and 0 if $L \neq M$. In particular:
\begin{equation*}
Z_{0j}^d:=\{ (M, p(s), L) \in S^{d/2} ( C) \times J^{d/2} ( C): L \neq M,\ \dim H^0(L)=j \}
\end{equation*}
and
\begin{equation*}
Z_{1j}^d:=\{ (M, p(s), L) \in S^{d/2} ( C) \times J^{d/2} ( C): L=M,\ \dim H^0(L)=j \}.
\end{equation*}

Furthermore it is clear that:
\begin{equation*}
[Z_{0j}^d]+[Z_{1j}^d]=[S^{d/2}( C)][V_j^{d/2}]
\end{equation*}
and
\begin{equation*}
[Z_{1j}^d]=[V_j^{d/2}][\C\P^{j-1}].
\end{equation*}

From these considerations we can deduce that:
\begin{align*}
[\widetilde{\mathscr{U}}^d_1]&=\sum_{j \geq 2} [Z_{1j}^d] \left ( [\C\P^{d/2-1}]-[\C\P^{j-2}] \right ) + [Z_{11}^d][\C\P^{d/2-1}]+\\
&+\sum_{j \geq 1} [Z_{0j}^d] \left ( [\C\P^{d/2-1}]-[\C\P^{j-1}]\right ) + [Z_{00}^d][\C\P^{d/2-1}]=\\
&=\sum_{j \geq 1} [V_j^{d/2}] \left ( \frac{\L^j-1}{\L-1} \right )\left ( \frac{\L^{d/2}-1}{\L-1}-\frac{\L^{j-1}-1}{\L-1}\right )+\\
&+\sum_{j \geq 0} [V_j^{d/2}] \left ( [S^{d/2}( C)]-\frac{\L^j-1}{\L-1} \right )\left ( \frac{\L^{d/2}-1}{\L-1}-\frac{\L^{j}-1}{\L-1}\right ).
\end{align*}

Instead, since $\dim H^0( K\End E)=1+4(g-1)$ if $L \neq M$ and $2+4(g-1)$ if $L=M$ we can look at the following map as well:
\begin{align*}
\mathscr{U}^d_1&\rightarrow \widetilde{\mathscr{U}}^d_1\\
(E,\phi,s) & \mapsto (E,s)
\end{align*}
and deduce
\begin{align*}
[\mathscr{U}^d_1]&=\L^{2+4(g-1)}\sum_{j \geq 2} [Z_{1j}^d] \left ( [\C\P^{d/2-1}]-[\C\P^{j-2}] \right ) +\L^{2+4(g-1)} [Z_{11}^d][\C\P^{d/2-1}]+\\
&+\L^{1+4(g-1)}\sum_{j \geq 1} [Z_{0j}^d] \left ( [\C\P^{d/2-1}]-[\C\P^{j-1}]\right ) +\L^{1+4(g-1)} [Z_{00}^d][\C\P^{d/2-1}]=\\
&=\L^{2+4(g-1)}\sum_{j \geq 1} [V_j^{d/2}] \left ( \frac{\L^j-1}{\L-1} \right )\left ( \frac{\L^{d/2}-1}{\L-1}-\frac{\L^{j-1}-1}{\L-1}\right )+\\
&+\L^{1+4(g-1)}\sum_{j \geq 0} [V_j^{d/2}] \left ( [S^{d/2}( C)]-\frac{\L^j-1}{\L-1} \right )\left ( \frac{\L^{d/2}-1}{\L-1}-\frac{\L^{j}-1}{\L-1}\right ).
\end{align*}
\end{proof}
\end{prop}
\begin{rmk}
Note the similarities in the computation of $[\mathscr{U}^d_1]$ and $[B^{d,-}_{\bar\sigma}]$ in proposition \ref{bminus}.
\end{rmk}
\begin{prop}
\label{motiveV}
We have the following motivic equalities:
\begin{align*}
[\mathscr{V}^d_1]&=\L^{2+4(g-1)}\left ( [\Sym^2( S^{d/2} ( C) )]-\sum_{j \geq 1} [V_j^{d/2}] [\Sym^2 ( \C\P^{j-1})] \right ).
\end{align*}
and
\begin{align*}
[\widetilde{\mathscr{V}}^d_1]&=[\Sym^2( S^{d/2} ( C) )]-\sum_{j \geq 1} [V_j^{d/2}] [\Sym^2 ( \C\P^{j-1})].
\end{align*}
\begin{proof}
Note that for all $(E,s) \in \widetilde{\mathscr{V}}^d_1$ we have $\dim H^0(K \End E)=2+4(g-1)$ therefore $[\mathscr{V}^d_1]=\L^{2+4(g-1)}[\widetilde{\mathscr{V}}^d_1]$. Second, by the description of $\widetilde{\mathscr{V}}^d_1$ it is easy to see that it is isomorphic to $\Sym^2 ( S^{d/2} ( C) ) \setminus \Gamma$ where $\Gamma$ is the locus of divisors $(D_1, D_2) \in \Sym^2 ( S^{d/2} ( C) )$ such that $ \mathcal{ O }(D_1) \cong \mathcal{O} (D_2 )$. To understand $\Gamma$, we look at the map
\begin{align*}
\Sym^2 ( S^{d/2} ( C) )  & \rightarrow \Sym^2 ( J^{d/2} ( C) )\\
(D_1, D_2) &\mapsto (\mathcal{ O }(D_1), \mathcal{O} (D_2 ) ).
\end{align*}

$\Gamma$ is the inverse image of the diagonal in $\Sym^2 ( J^{d/2} ( C) )$ with respect to the previous map. In particular there is a map:
\begin{align*}
\Gamma  & \rightarrow J^{d/2} ( C) \\
(D_1, D_2) &\mapsto \mathcal{ O }(D_1) = \mathcal{O} (D_2 ) .
\end{align*}
whose fiber over $L \in J^{d/2} (C )$ is $\Sym^2 ( \P H^0( L) )$. Therefore we get:
\begin{align*}
[\widetilde{\mathscr{V}}^d_1]&=[\Sym^2( S^{d/2} ( C) )]-\sum_{j \geq 1} [V_j^{d/2}] [\Sym^2 ( \C\P^{j-1})].
\end{align*}
\end{proof}
\end{prop}
\begin{prop}
\label{motiveW}
We have the following motivic equalities:
\begin{align*}
[\mathscr{W}^d_1]&=\L^{4g}\sum_{j\geq 2} [V_j^{d/2}] [\Gr(2,j)]
\end{align*}
and
\begin{align*}
[\widetilde{\mathscr{W}}^d_1]&=\sum_{j\geq 2} [V_j^{d/2}] [\Gr(2,j)].
\end{align*}
\begin{proof}
Since for all $(E,s) \in \widetilde{\mathscr{W}}^d_1$ we have $\dim H^0( K \End E)= 4g$ then we have $[\mathscr{W}^d_1]=\L^{4g}[\widetilde{\mathscr{W}}^d_1]$.

Furthermore, by the definition of $\widetilde{\mathscr{W}}^d_1$ we have a map:
\begin{align*}
\widetilde{\mathscr{W}}^d_1 & \rightarrow J^{d/2}\\
(L \oplus L,s) &\mapsto L
\end{align*}
and, since $s= u \oplus v$ with $u, v \in H^0(L)$ linearly independent, the fiber of such a map above $L$ is $\Gr (2, H^0(L))$, the Grassmannian of two dimensional subspaces of $H^0(L)$. The result follows.
\end{proof}
\end{prop}
\begin{prop}
We have the following motivic equalities:
\begin{align*}
[\mathscr{X}^d_1]&=\L^{2+4(g-1)} [S^{d/2}( C)] (\L^{d/2-1}-1)+\L^{1+4(g-1)} [J^{d/2} ( C)] [S^{d/2}( C)] [\C\P^{d/2-1}] +\\
&-\L^{1+4(g-1)} [S^{d/2}( C)] ^2 + \L^{2+4(g-1)} [\Sym^2( S^{d/2} ( C) )]
\end{align*}
and
\begin{align*}
[M_{\varepsilon,ss}^{2,d}]&=[J^{d/2} (C )][S^{d/2}( C)] [\C\P^{d/2-1}]-[S^{d/2}( C)]^2+[\Sym^2( S^{d/2} ( C) )].
\end{align*}
\begin{proof}
We use propositions \ref{motiveU}, \ref{motiveV} and \ref{motiveW} to compute directly:
\begin{align*}
&[\mathscr{X}^d_1]=[\mathscr{U}^d_1]+[\mathscr{V}^d_1]+[\mathscr{W}^d_1]=\\
&=\L^{2+4(g-1)}\sum_{j \geq 1} [V_j^{d/2}] \left ( \frac{\L^j-1}{\L-1} \right )\left ( \frac{\L^{d/2}-1}{\L-1}-\frac{\L^{j-1}-1}{\L-1}\right )+\\
&+\L^{1+4(g-1)}\sum_{j \geq 0} [V_j^{d/2}] \left ( [S^{d/2}( C)]-\frac{\L^j-1}{\L-1} \right )\left ( \frac{\L^{d/2}-1}{\L-1}-\frac{\L^{j}-1}{\L-1}\right )+\\
&+\L^{2+4(g-1)}\left ( [\Sym^2( S^{d/2} ( C) )]-\sum_{j \geq 1} [V_j^{d/2}] [\Sym^2 ( \C\P^{j-1}] \right )+\\
&+\L^{4g}\sum_{j\geq 2} [V_j^{d/2}] [\Gr(2,j)]=\\
&=\L^{2+4(g-1)}[\C\P^{d/2-1}]\sum_{j \geq 1} [V_j^{d/2}] [\C\P^{j-1}]-\L^{2+4(g-1)}\sum_{j \geq 1} [V_j^{d/2}] \left ( \frac{\L^j-1}{\L-1} \right )\left ( \frac{\L^{j-1}-1}{\L-1}\right )+\\
&+\L^{1+4(g-1)}\sum_{j \geq 0} [V_j^{d/2}] [S^{d/2}( C)] [\C\P^{d/2-1}]-\L^{1+4(g-1)}\sum_{j \geq 0} [V_j^{d/2}] [S^{d/2}( C)] [\C\P^{j-1}]+\\
&-\L^{1+4(g-1)}\sum_{j \geq 0} [V_j^{d/2}][\C\P^{d/2-1}][\C\P^{j-1}]+\L^{1+4(g-1)}\sum_{j \geq 0} [V_j^{d/2}] \left ( \frac{\L^{j}-1}{\L-1}\right )^2+\\
&+\L^{2+4(g-1)}[\Sym^2( S^{d/2} ( C) )]-\L^{2+4(g-1)}\sum_{j \geq 1} [V_j^{d/2}] [\Gr(2,j+1)]+\\
&+\L^{4g}\sum_{j\geq 2} [V_j^{d/2}] [\Gr(2,j)]=\\
&=\L^{2+4(g-1)}[\C\P^{d/2-1}][S^{d/2}( C)]+\L^{1+4(g-1)}[J^{d/2} ( C)] [S^{d/2}( C)] [\C\P^{d/2-1}]+\\
&-\L^{1+4(g-1)}[S^{d/2}( C)]^2-\L^{1+4(g-1)}[S^{d/2}( C)][\C\P^{d/2-1}]+\L^{2+4(g-1)}[\Sym^2( S^{d/2} ( C) )]+\\
&+\L^{1+4(g-1)} \sum_{j \geq 1} [V_j^{d/2}] \bigg ( -\L  \left ( \frac{\L^j-1}{\L-1} \right )\left ( \frac{\L^{j-1}-1}{\L-1}\right ) + 
\left ( \frac{\L^{j}-1}{\L-1}\right )^2 +\\
&- \L \frac{(\L^{j+1}-1)(\L^{j}-1)}{(\L^2-1)(\L-1)} + \L^3 \frac{(\L^{j}-1)(\L^{j-1}-1)}{(\L^2-1)(\L-1)} \bigg )=\\
&=\L^{1+4(g-1)}[\C\P^{d/2-1}][S^{d/2}( C)](\L-1)+\L^{1+4(g-1)} [J^{d/2} ( C)] [S^{d/2}( C)] [\C\P^{d/2-1}] +\\
&-\L^{1+4(g-1)} [S^{d/2}( C)] ^2 + \L^{2+4(g-1)} [\Sym^2( S^{d/2} ( C) )]+\\
&+\L^{1+4(g-1)} \sum_{j \geq 1} [V_j^{d/2}] [\C\P^{j-1}](1-\L)=\\
&=\L^{2+4(g-1)} [S^{d/2}( C)] (\L^{d/2-1}-1)+\L^{1+4(g-1)} [J^{d/2} ( C)] [S^{d/2}( C)] [\C\P^{d/2-1}] +\\
&-\L^{1+4(g-1)} [S^{d/2}( C)] ^2 + \L^{2+4(g-1)} [\Sym^2( S^{d/2} ( C) )].
\end{align*}
Furthermore:
\begin{align*}
[\widetilde{\mathscr{X}}^d_1]&=[\widetilde{\mathscr{U}}^d_1]+[\widetilde{\mathscr{V}}^d_1]+[\widetilde{\mathscr{W}}^d_1]=\\
&=[\C\P^{d/2-1}]\sum_{j \geq 1} [V_j^{d/2}] [\C\P^{j-1}]-\sum_{j \geq 1} [V_j^{d/2}] \left ( \frac{\L^j-1}{\L-1} \right )\left ( \frac{\L^{j-1}-1}{\L-1}\right )+\\
&+\sum_{j \geq 0} [V_j^{d/2}] [S^{d/2}( C)] [\C\P^{d/2-1}]-\sum_{j \geq 0} [V_j^{d/2}] [S^{d/2}( C)] [\C\P^{j-1}]+\\
&-\sum_{j \geq 0} [V_j^{d/2}][\C\P^{d/2-1}][\C\P^{j-1}]+\sum_{j \geq 0} [V_j^{d/2}] \left ( \frac{\L^{j}-1}{\L-1}\right )^2+\\
&+[\Sym^2( S^{d/2} ( C) )]-\sum_{j \geq 1} [V_j^{d/2}] [\Gr(2,j+1)]+\sum_{j\geq 2} [V_j^{d/2}] [\Gr(2,j)]=\\
&=[S^{d/2}( C)] [\C\P^{d/2-1}]+[J^{d/2} (C )][S^{d/2}( C)] [\C\P^{d/2-1}]+\\
&-[S^{d/2}( C)] [\C\P^{d/2-1}]-[S^{d/2}( C)]^2+[\Sym^2( S^{d/2} ( C) )]+\\
&+\sum_{j \geq 1} \frac{[V_j^{d/2}] (L^j-1)}{(\L^2-1)(\L-1)} \left ( (\L^j-\L^{j-1})(\L+1)-(\L^{j+1}-1)+(\L^{j-1}-1)\right )=\\
&=[J^{d/2} (C )][S^{d/2}( C)] [\C\P^{d/2-1}]-[S^{d/2}( C)]^2+[\Sym^2( S^{d/2} ( C) )].
\end{align*}
\end{proof}
\end{prop}
\begin{prop}
We have the following motivic equality:
\begin{align*}
[\mathscr{X}^d_2]&=[S^{d/2}( C)][J^{d/2}( C)][\C\P^{g-2}]\L^{1+4(g-1)}+[S^{d/2}( C)]\L^{3g-2}(\L^{2g-2}+\L^g-1).
\end{align*}
\begin{proof}
Recall the description:
\begin{equation*}
 \mathscr{X}^d_2=\left \{\begin{array}{c}
(E,\phi,s) \text{ such that $E$ fits into an exact sequence}\\
0 \rightarrow L \rightarrow E \rightarrow M \rightarrow 0\\
\text{with $\deg L=\deg M$, $s \in H^0(L)$ and}\\
\phi \text{ does not preserve } L
\end{array} \right \}.
\end{equation*}

We can divide the computation of the motive according to whether or not the vector bundle underlying the triple is split. In both cases it is relevant to distinguish when the quotient and the subobject are the same or different.

More precisely:
\begin{align*}
\mathscr{X}^d_2&=\left \{\begin{array}{c}
(E,\phi,s) \text{ such that $E=L\oplus L$}\\
s = u \oplus 0 \text{ with $u \in H^0(L)$}\\
\phi \text{ does not preserve}\\
\text{ the copy of $L$ containing $s$}
\end{array} \right \} \sqcup
\left \{\begin{array}{c}
(E,\phi,s) \text{ such that $E=L\oplus M$}\\
L \neq M, s \in H^0(L), \deg L = \deg M\\
\phi \text{ does not preserve $L$}
\end{array} \right \} \sqcup\\
& \sqcup
\left \{\begin{array}{c}
(E,\phi,s) \text{ with $E$ fitting}\\
\text{ in a nonsplit exact sequence}\\
0 \rightarrow L \rightarrow E \rightarrow L \rightarrow 0\\
\text{$s \in H^0(L)$ and}\\
\phi \text{ does not preserve the subobject}
\end{array} \right \} \sqcup
\left \{\begin{array}{c}
(E,\phi,s) \text{ with $E$ fitting}\\
\text{ in a nonsplit exact sequence}\\
0 \rightarrow L \rightarrow E \rightarrow M \rightarrow 0\\
\text{$L \neq M$, $s \in H^0(L)$ and}\\
\phi \text{ does not preserve $L$}
\end{array} \right \}.
\end{align*}

Let us name the four strata $\mathscr{Y}^d_i$ for $i=1, \dots, 4$, in the order they appear. We have a map:
\begin{align*}
\mathscr{Y}^d_1&\rightarrow S^{d/2} ( C)\\
(L \oplus L, \phi, u \oplus 0) &\mapsto (L, u \oplus 0).
\end{align*}

The datum of $L$ and $u$ essentially recovers the pair but then $(E,s)$ will have automorphisms. More precisely $\Aut(E,s)=(\C^*)^2 \times \C$ given by the invertible upper triangular matrices. These automorphisms act by conjugation on the set of possible Higgs fields:
\begin{equation*}
\phi=\begin{pmatrix}
\psi_1 & \vartheta \\
\gamma & \psi_2
\end{pmatrix}
\end{equation*}
where the entries of the matrix are in $H^0(K)$ and $\gamma \neq 0$. A quick computation shows that the stabilizers of the action are only the scalar multiples of the identity. Therefore the motive of the fibers of:
\begin{align*}
\mathscr{Y}^d_1 &\rightarrow S^{d/2} ( C)
\end{align*}
is:
\begin{equation*}
\left [\frac{H^0(K)^3 \times (H^0(K) \setminus \{0\})}{\C^* \times \C} \right ]=\frac{\L^{3g}(\L^g-1)}{\L(\L-1)}
\end{equation*}
and
\begin{equation*}
[\mathscr{Y}^d_1]=[S^{d/2}( C)]\frac{\L^{3g}(\L^g-1)}{\L(\L-1)}.
\end{equation*}

For the computation of $[\mathscr{Y}_2^d]$ let us first note that the locus of:
\begin{equation*}
S^{d/2}( C) \times J^{d/2}( C)
\end{equation*}
defined by the triples $(L,s,M)$ where $L = M$ is the graph $\Gamma^d$ of the map:
\begin{align*}
S^{d/2}( C) &\rightarrow J^{d/2}( C)\\
(L,s) & \mapsto L
\end{align*}
and in particular it is isomorphic to $S^{d/2}( C)$.

We have a map
\begin{align*}
\mathscr{Y}^d_2&\rightarrow S^{d/2} ( C)\times J^{d/2}( C) \setminus \Gamma^d\\
(L \oplus M, \phi, s) &\mapsto (L, s, M).
\end{align*}

As before, the datum of $(L,s,M)$ determines the pair $(E,s)$ but $\Aut(E,s)=(\C^*)^2$. The set of Higgs fields not preserving $L$ is $H^0(K)^2 \times H^0(KM^*L) \times ( H^0(KL^*M) \setminus \{0\})$ and the action of $\Aut(E,s)$ again has only scalar multiples of the identity as stabilizers. Therefore the motive of the fibers of
\begin{align*}
\mathscr{Y}^d_2&\rightarrow S^{d/2} ( C)\times J^{d/2}( C) \setminus \Gamma^d
\end{align*}
is
\begin{equation*}
\left [\frac{H^0(K)^2 \times H^0(KM^*L) \times ( H^0(KL^*M) \setminus \{0\})}{\C^*} \right ]=\frac{\L^{3g-1}(\L^{g-1}-1)}{\L-1}
\end{equation*}
and so
\begin{equation*}
\left [ \mathscr{Y}_2^d \right ]=\left ([S^{d/2}( C)][J^{d/2}( C)]-[S^{d/2}( C)] \right )\frac{\L^{3g-1}(\L^{g-1}-1)}{\L-1}.
\end{equation*}
For $\mathscr{Y}_3^d$ we have a map:
\begin{align*}
\mathscr{Y}^d_3 &\rightarrow S^{d/2} ( C) \times \P\Ext^1( \mathcal{O},\mathcal{O})\\
(E,\phi,s) & \mapsto (L,s,  [E] )
\end{align*}
remembering $(L,s)$ and the class of the extension. Note that the target is a trivial projective bundle of rank $g-1$ because it is the pullback to $S^{d/2} ( C)$ of the constant bundle $\P\Ext^1( \mathcal{O},\mathcal{O}) \rightarrow \{pt\}$. As we already observed earlier, $\Aut(E,s)=\C^* \times \C$ generated by nonzero multiples of the identity and multiples of $i \circ p$ where $i$ and $p$ are the inclusion and the projection in the exact sequence defining $E$. The locus of possible Higgs fields is isomorphic to $\C^{2+4(g-1)} \setminus \C^{3g-1}$ and the conjugation action of $\Aut(E,s)$ has scalar multiples of the identity as stabilizers. Therefore:
\begin{equation*}
[\mathscr{Y}^d_3]=[S^{d/2}( C)][\C\P^{g-1}]\left ( \L^{1+4(g-1)} - \L^{3g-2} \right ).
\end{equation*}

Finally, for $\mathscr{Y}^d_4$, there is a map:
\begin{align*}
\mathscr{Y}^d_4 &\rightarrow \mathcal{B}^d\\
(E,\phi,s) & \mapsto (L,s, M, [E] )
\end{align*}
where $\mathcal{B}^d$ is the projective bundle on $S^{d/2} ( C)\times J^{d/2}( C) \setminus \Gamma^d$ defined by the pulling back the projective bundle on $J^0( C)$ whose fiber over $A$ is $\P H^1(A)$, with respect to the map:
\begin{align*}
S^{d/2} ( C)\times J^{d/2}( C) \setminus \Gamma^d &\rightarrow J^0( C)\\
(L,s,M) & \mapsto M^*L. 
\end{align*}

Once again, $\mathcal{B}^d$ will remember the pair $(E,s)$. In this case the automorphisms of the pairs are just scalar multiples of the identity and therefore we get:
\begin{align*}
[\mathscr{Y}^d_4]=\left ([S^{d/2}( C)][J^{d/2}( C)]-[S^{d/2}( C)] \right )[\C\P^{g-2}]\left ( \L^{1+4(g-1)} - \L^{3g-1} \right ).
\end{align*}

Putting everything together we get:
\begin{align*}
[\mathscr{X}^d_2]&=[S^{d/2}( C)]\frac{\L^{3g}(\L^g-1)}{\L(\L-1)}+\left ([S^{d/2}( C)][J^{d/2}( C)]-[S^{d/2}( C)] \right )\frac{\L^{3g-1}(\L^{g-1}-1)}{\L-1}\\
&+[S^{d/2}( C)][\C\P^{g-1}]\left ( \L^{1+4(g-1)} - \L^{3g-2} \right )+\\
&+\left ([S^{d/2}( C)][J^{d/2}( C)]-[S^{d/2}( C)] \right )[\C\P^{g-2}]\left ( \L^{1+4(g-1)} - \L^{3g-1} \right )=\\
&=[S^{d/2}( C)][\C\P^{g-1}]\left ( \L^{1+4(g-1)} + \L^{3g-1}- \L^{3g-2} \right )+\\
&+\left ([S^{d/2}( C)][J^{d/2}( C)]-[S^{d/2}( C)] \right )[\C\P^{g-2}]\L^{1+4(g-1)}=\\
&=[S^{d/2}( C)][J^{d/2}( C)][\C\P^{g-2}]\L^{1+4(g-1)}+[S^{d/2}( C)]\L^{3g-2}(\L^{2g-2}+\L^g-1).
\end{align*}
\end{proof}
\end{prop}
We can summarize the results of this section as follows.
\begin{theorem}
\label{moteven}
Let $d \geq 0$ be an even integer and $0 < \varepsilon < 2$. Then:
\begin{align*}
[\mathcal{M}_\varepsilon^{2,d}]&=\L^{1+4(g-1)}  [M^{2,d}_\varepsilon]+\sum_{(d_1,d_2) \in I_1^e(d)} \L^{1+3(g-1)+d_2} [S^{d_1}( C)][S^{d_1-d_2+2g-2}( C)] +\\
&+\sum_{(d_1,d_2) \in I_2^e(d)}\L^{1+4(g-1)} [S^{d_2}( C)][S^{d_1-d_2+2g-2}( C)]+\\
&+(\L-1) \L^{1+4(g-1)}[\Sym^2( S^{d/2} ( C) )] +\L^{1+4(g-1)} [S^{d/2}( C)][J^{d/2}( C)][\C\P^{g-2}]+\\
&+[S^{d/2}( C)] \L^{3g-2} \left ( \L^{d/2+g-1}+\L^{2g-2}-1\right )+\\
&+\sum_{(d_1,d_2) \in I_1^e(d)} (\L^{4g-2}-\L^{4g-3}) [S^{d_1}( C)] \left ([S^{d_2}( C)]-[J( C)] \frac{\L^{d_2+1-g}-1}{\L-1} \right).
\end{align*}
\end{theorem}
\section{Direct computation of the motive of $\mathcal{M}^{2,d}_{\infty}$}
Using the content of the previous sections, we can compute the motive of $\mathcal{M}^{2,d}_{\infty}$ by simply starting from the motive of $\mathcal{M}^{2,d}_{\varepsilon}$ and then add and subtract the motive of the flip loci. In other words we have all the ingredients to compute $[\mathcal{M}^{2,d}_{\infty}]$.

It is interesting to note that it is possible to proceed in the reverse direction, i.e. first compute $[\mathcal{M}^{2,d}_{\infty}]$ directly, then add the motives of the flip loci and ultimately deduce $[\mathcal{M}^{2,d}_{\varepsilon}]$. Here we outline the strategy for the computation of $[\mathcal{M}^{2,d}_{\infty}]$.

Recall that $[\mathcal{M}^{2,d}_{\infty}]$ contains only split type 2 attracting sets. These have the form:
\begin{equation*}
F_{(1,1),\infty}^{(d_1,d_2),2+}=
\left \{
\begin{array}{c}
(E,\phi,s) \text{ $\sigma$-stable where } E \text{ is defined by}\\
0 \rightarrow L \rightarrow E \rightarrow M \rightarrow 0\\
s \in H^0(L), \deg L =d_2, \deg M = d_1\\
\text{and } \phi \text{ does not preserve } L
\end{array} \right \}
\end{equation*}
with $d_1$ and $d_2$ satisfying:
\begin{align*}
&d/2 +1-g\leq d_1 \leq d \\
& 0 \leq d_2 \leq d/2+g-1\\
& 0 \leq d_1-d_2+2g-2 \leq 2g-2+d \text{ same parity as $d$}.
\end{align*}

Observe that for $d_2>d_1$, which is equivalent to $d_1<d/2$, $F_{(1,1),\infty}^{(d_1,d_2),2+}=F_{(1,1),\varepsilon}^{(d_1,d_2),2+}$ lies entirely in the smooth part of the moduli space (see proposition \ref{deffixed}) and therefore the motive can be computed using \cite[Theorem 4.1]{bialynicki1973some}. If $d$ is even and $d_1=d_2=d/2$ we know that $F_{(1,1),\infty}^{(d/2,d/2),2+}= \mathscr{X}^d_2$.

To compute the remaining motives we assume that $d_1>d_2$. Similar to the computation of the motive of the type 1 attracting sets, we can further decompose an attracting set according to the properties of the underlying pair. In fact, we know that for pairs of the form:
\begin{equation*}
0 \rightarrow L \rightarrow E \rightarrow M \rightarrow 0
\end{equation*}
with $\deg M-\deg L= k >0$, being in $\P W^{d,+}_k$ is equivalent to the extension being non-split. Therefore we can distinguish between the cases where the underlying pair is non-split and when $E=L\oplus M$.
\begin{defn}
Let $(d_1,d_2)$ satisfy $d_1+d_2=d$, the three inequalities above and $d_1> d_2$. Let us denote by $NSPF_{(1,1),\infty}^{(d_1,d_2),2+}$ the locus of $F_{(1,1),\infty}^{(d_1,d_2),2+}$ where the underlying pair lies in $\P W^{d,+}_{d_1-d_2}$ and by $SPF_{(1,1),\infty}^{(d_1,d_2),2+}$ the locus where the pair is split.
\end{defn}
It is clear from the definitions that we have:
\begin{equation*}
NSPF_{(1,1),\infty}^{(d_1,d_2),2+}=B_{d_1-d_2}^{d,+} \cap F_{(1,1),\infty}^{(d_1,d_2),2+}
\end{equation*}
and:
\begin{equation*}
F_{(1,1),\infty}^{(d_1,d_2),2+} =B_{d_1-d_2}^{d,+} \cap F_{(1,1),\infty}^{(d_1,d_2),2+} \sqcup SPF_{(1,1),\infty}^{(d_1,d_2),2+}.
\end{equation*}

We can also prove the following:
\begin{prop}
\begin{equation*}
B_{d_1-d_2}^{d,+} = NSPF_{(1,1),\infty}^{(d_1,d_2),2+} \sqcup S \mathcal{W}^{d,+}_{d_1-d_2}=B_{d_1-d_2}^{d,+} \cap F_{(1,1),\infty}^{(d_1,d_2),2+} \sqcup S \mathcal{W}^{d,+}_{d_1-d_2}.
\end{equation*}
\begin{proof}
The second equality follows from the previous remarks. Recall that $(E,\phi,s) \in B_{d_1-d_2}^{d,+}$ if and only if $(E,s) \in \P W^{d,+}_{d_1-d_2}$. Since by definition $(E,\phi,s) \in NSPF_{(1,1),\infty}^{(d_1,d_2),2+}$ if and only if $(E,s) \in F_{(1,1),\infty}^{(d_1,d_2),2+}$ and the extension for $(E,s)$ is nonsplit we immediately get that if $(E,\phi,s) \in NSPF_{(1,1),\infty}^{(d_1,d_2),2+}$ then $(E,s) \in \P W^{d,+}_{d_1-d_2}$ and so $(E,\phi,s) \in B_{d_1-d_2}^{d,+}$.

This proves $NSPF_{(1,1),\infty}^{(d_1,d_2),2+} \subset B_{d_1-d_2}^{d,+}$. We already observed in proposition \ref{weirdexpl1} that $S \mathcal{W}^{d,+}_{d_1-d_2} \subset B_{d_1-d_2}^{d,+}$. To conclude it is enough to observe that if $(E,\phi,s) \in NSPF_{(1,1),\infty}^{(d_1,d_2),2+}$ or $(E,\phi,s) \in S \mathcal{W}^{d,+}_{d_1-d_2}$ then $(E,s) \in \P W^{d,+}_{d_1-d_2}$ but in the first case the subobject containing the section is not preserved by $\phi$, since $NSPF_{(1,1),\infty}^{(d_1,d_2),2+} \subset F_{(1,1),\infty}^{(d_1,d_2),2+}$, while in the second case the subobject containing the section is preserved by $\phi$, since $S \mathcal{W}^{d,+}_{d_1-d_2} \subset \mathcal{W}^{d,+}_{d_1-d_2}$. This proves both the reverse inclusion and the fact that the union is disjoint.
\end{proof}
\end{prop}
Therefore we can write the relation:
\begin{equation*}
[F_{(1,1),\sigma}^{(d_1,d_2),2+}]=[B_{d_1-d_2}^{d,+}]-[S \mathcal{W}^{d,+}_{d_1-d_2}]+[SPF_{(1,1),\sigma}^{(d_1,d_2),2+}].
\end{equation*}
\begin{rmk}
Note that it is particularly hard to attempt the computation of $[B_{d_1-d_2}^{d,+}]$ because for $(E,\phi,s) \in B_{d_1-d_2}^{d,+}$ we get a presentation of the pair $(E,s)$ as an extension of line bundles for which the subobject is not destabilizing. This makes the problem particularly hard because we need to know the dimension of $H^0( K \End E)$ for these triples and for that we need to know whether or not $E$ is stable and if not who is the maximal destabilizing subbundle, all of which is not clear from the given presentation.

There is however a way around this. In fact the $B_{d_1-d_2}^{d,+}$ and $B_{d_1-d_2}^{d,-}$ are closely related to the wall-crossing of Bradlow pairs. In fact we have the following relations:
\begin{equation*}
0=[F_{(2),\varepsilon}^{(d),1+}]+\sum_{k=0}^{(d-1)/2} \left ( [B_{2k+1}^{d,-}] -[B_{2k+1}^{d,+}] \right )
\end{equation*}
for odd $d$, and 
\begin{equation*}
0=[F_{(2),\varepsilon}^{(d),1+}]+\sum_{k=0}^{d/2} \left ( [B_{2k}^{d,-}] -[B_{2k}^{d,+}] \right )
\end{equation*}
for even $d$. Both these relations imply that if we only wish to compute the sum of all of the $[B_{d_1-d_2}^{d,+}]$ then it suffices to $[F_{(2),\varepsilon}^{(d),1+}]$, which we computed in the previous sections, and the sum of the $[B_{d_1-d_2}^{d,-}]$ which can be computed using proposition \ref{bminus}.
\end{rmk}
%
%The full formula for $\mathcal{M}^{2,d}_{\infty}$, taking into account the full range for the type 2 cells, will be:
%\begin{align}
%\label{m_inf}
%[\mathcal{M}^{2,d}_{\infty}]&=\sum_{k=1-g}^{-1} [F_{(1,1)}^{(k+(d+1)/2,(d-1)/2-k),2+}]+\sum_{k=}^{(d-1)/2} [F_{(1,1)}^{(k+(d+1)/2,(d-1)/2-k),2+}]=\\ 
%&=\sum_{k=1-g}^{-1} [F_{(1,1)}^{(k+(d+1)/2,(d-1)/2-k),2+}]+\sum_{k=}^{(d-1)/2}  \left ([B_{2k+1}^{d,+}]-[S \mathcal{W}^{d,+}_{2k+1}]+[SPF_{(1,1),\sigma}^{(k+(d+1)/2,(d-1)/2-k),2+}] \right ) \nonumber
%\end{align}
In the following propositions we compute the remaining motives.
\begin{prop}
Let $\bar\sigma$ be a critical value. The following motivic equalities hold:
\begin{align*}
[NS\mathcal{W}_{\bar\sigma}^{d,+}]&=\L^{2g} [S^{(d-\bar\sigma)/2} ( C)] [S^{2g-2-\bar\sigma} ( C)] \text{ if } \bar\sigma < 2g-2 \text{ or } 0 \text{ otherwise}\\
[S\mathcal{W}_{\bar\sigma}^{d,+}]&=\L^{2g} \cdot [\C\P^{2g-3}]\cdot [S^{(d-{\bar\sigma})/2}( C)] \cdot [J( C)]+\L^{3g-2} \cdot [S^{(d-{\bar\sigma})/2}( C)] \cdot [S^{\bar\sigma}( C)]+\\
&-\L^{2g} [S^{(d-\bar\sigma)/2} ( C)] [S^{2g-2-\bar\sigma} ( C)].
\end{align*}
\begin{proof}
We already proved that there is a map
$$NS\mathcal{W}_{\bar\sigma}^{d,+} \rightarrow X_{\bar\sigma}$$
whose fiber over $(L,s,M,\psi_1,\psi_2)$ is $\P H^0(KLM^*)$. Let us define the strata
\begin{equation*}
U_{\bar\sigma,i}^+:=\{(L,s,M) \in S^{(d-\bar\sigma)/2} ( C) \times J^{(d+\bar\sigma)/2}( C) | \dim H^0(KLM^*) = i\}.
\end{equation*}

Consider also the map:
\begin{align*}
S^{(d-\bar\sigma)/2} ( C) \times J^{(d+\bar\sigma)/2}( C) &\rightarrow S^{(d-\bar\sigma)/2}( C) \times J^{2g-2-\bar\sigma}( C)\\
(L,s,M) & \mapsto (L,s, KLM^*)
\end{align*}
which is clearly an isomorphism. Under this isomorphism $U_{\bar\sigma,i}^+$ corresponds to
$$V^{2g-2-\bar\sigma}_{i} \times S^{(d-\bar\sigma)/2} ( C).$$
Therefore we deduce that:
\begin{equation*}
[U_{\bar\sigma,i}^+]=[V^{2g-2-\bar\sigma}_{i} ]\cdot [S^{(d-\bar\sigma)/2} ( C)].
\end{equation*}

So we can compute:
\begin{align*}
[NS\mathcal{W}_{\bar\sigma}^{d,+}]&=\L^{2g} [S^{(d-\bar\sigma)/2} ( C)] \sum_{i=0}^{2g-2-\bar\sigma} [ V^{2g-2-\bar\sigma}_{i} ] \cdot [\C \P^{i-1}]=\\
&=\L^{2g} [S^{(d-\bar\sigma)/2} ( C)] [S^{2g-2-\bar\sigma} ( C)].
\end{align*}

Since $\mathcal{W}_{\bar\sigma}^{d,+}=S\mathcal{W}_{\bar\sigma}^{d,+} \sqcup NS\mathcal{W}_{\bar\sigma}^{d,+}$ we also deduce:
\begin{align*}
[S\mathcal{W}_{\bar\sigma}^{d,+}]=[\mathcal{W}_{\bar\sigma}^{d,+}]-&[NS\mathcal{W}_{\bar\sigma}^{d,+}]=\L^{2g} \cdot [\C\P^{2g-3}]\cdot [S^{(d-{\bar\sigma})/2}( C)] \cdot [J( C)]+\\
&+\L^{3g-2} \cdot [S^{(d-{\bar\sigma})/2}( C)] \cdot [S^{\bar\sigma}( C)]-\L^{2g} [S^{(d-\bar\sigma)/2} ( C)] [S^{2g-2-\bar\sigma} ( C)].
\end{align*}
\end{proof}
\end{prop}
\begin{rmk}
Note that in the previous proposition we used the motive of $\mathcal{W}_{\bar\sigma}^{d,+}$ to compute $[S\mathcal{W}_{\bar\sigma}^{d,+}]$. Observe that it not strictly necessary and we can compute $[S\mathcal{W}_{\bar\sigma}^{d,+}]$ directly by observing that the restriction of $\mathcal{W}_{\bar\sigma}^{d,+} \rightarrow X_{\bar\sigma}^d$ to $S\mathcal{W}_{\bar\sigma}^{d,+}$ has fibers that can be well understood. Namely the fiber over $(L,s,M,\psi_1,\psi_2)$ is $\P\H^1((M,\psi_2),(L,\psi_1)) \setminus \P H^0(KLM^*)$ the dimension of which can be computed by considerations similar to those in proposition \ref{wplus}.
\end{rmk}
\begin{prop}
We have the identity:
\begin{equation*}
[SPF_{(1,1),\sigma}^{(d_1,d_2),2+}]=\L^{2g} [S^{d_2} ( C)] [\C \P^{d_1-d_2+g-2}] \left ( (\L-1) [S^{d_2-d_1+2g-2}( C)] + [J (C )]\right ).
\end{equation*}
\begin{proof}
Recall that $SPF_{(1,1),\sigma}^{(d_1,d_2),2+}$ consists of those pairs in the form $(L\oplus M,s,\phi)$ where $\deg L=d_2$ and $\deg M=d_1$ have the appropriate degrees, $s \in H^0(L)$ and $\phi$ does not preserve $L$. First of all, note that the moduli space of such split pairs is
\begin{equation*}
S^{d_2} ( C) \times J^{d_1} (C ).
\end{equation*}

Second, we can always write the Higgs field into matrix form:
\begin{equation*}
\phi = 
\begin{pmatrix}
\psi_1 & \xi\\
\theta & \psi_2
\end{pmatrix}
\end{equation*}
where
\begin{equation*}
\theta : L \rightarrow M \otimes K
\end{equation*}
has to be nonzero, and
\begin{equation*}
\xi : M \rightarrow L \otimes K.
\end{equation*}

Since the endomorphisms of the underlying pair are just diagonal (possibly non scalar) then we can see that $(L \oplus M, s , \phi)$ and $(L \oplus M, s , \phi')$ are isomorphic if and only if there exists $\lambda \in \C^*$ such that:
\begin{equation*}
\phi' = 
\begin{pmatrix}
\psi'_1 & \xi'\\
\theta' & \psi'_2
\end{pmatrix}=
\begin{pmatrix}
\psi_1 & \lambda \xi\\
\lambda^{-1} \theta & \psi_2
\end{pmatrix}.
\end{equation*}

Let's consider the map:
\begin{equation*}
SPF_{(1,1),\sigma}^{(d_1,d_2),2+} \rightarrow S^{d_2} ( C) \times J^{d_1} (C ) \times H^0(K)^2
\end{equation*}
forgetting the off-diagonal parts of the Higgs field. Then, by what we said, the fiber of such a map over $(L,s,M)$ will be:
\begin{equation*}
\quotient{\left (H^0(KML^*) \setminus \{0\} \right ) \times H^0(KLM^*)}{ \C^*}.
\end{equation*}

This allows us to compute the motive of $SPF_{(1,1),\sigma}^{(d_1,d_2),2+}$ in a familiar way by stratifying $S^{d_2} ( C) \times J^{d_1} (C )$ according to the dimension of $H^0(KLM^*)$. Note in fact that the dimension of $H^0(KML^*)$ is constantly equal to $g-1+d_1-d_2$. Call $V_i$ the stratum of $J^{d_1-d_2+2g-2}$ where the dimension of the global sections of the line bundle is equal to $i$.

Then we have:
\begin{align*}
[SPF_{(1,1),\sigma}^{(d_1,d_2),2+}]=&\sum_i \left (\L^{2g} [S^{d_2} ( C)] [V_i] \L^i \frac{\L^{d_1-d_2+g-1}-1}{\L-1} \right )=\\
=&\L^{2g} [S^{d_2} ( C)] [\C \P^{d_1-d_2+g-2}] \left ( (\L-1) \sum_i [V_i] [\C\P^{i-1}] + \sum_i [V_i]\right )=\\
=&\L^{2g} [S^{d_2} ( C)] [\C \P^{d_1-d_2+g-2}] \left ( (\L-1) [S^{d_2-d_1+2g-2}( C)] + [J (C )]\right ).
\end{align*}
\end{proof}
\end{prop}

%% file: Poles.tex
\chapter{Bradlow-Higgs triples with poles}
In this chapter we will introduce a variation to Bradlow-Higgs triples that is the analogue to the one in section \ref{hbpoles} for Higgs bundles. We will see that the results of the thesis about Bradlow-Higgs triples carry over with almost no modifications.

Let's start with a definition.
\begin{defn}[Bradlow-Higgs triples with poles and $\sigma$-stability]
Let $\gamma > 0$ be an integer and $P$ a point on the smooth projective curve $C$. A triple $(E,\phi,s)$ is said to be a \emph{Bradlow-Higgs $\gamma$-triple} if $(E,\phi)$ is a Higgs $\gamma$-bundle (as in definition \ref{gammahiggs}) and $s \in H^0(C,E)$ is a nonzero section of the underlying vector bundle $E$.

Let $\sigma$ be a positive real number. We say $(E,\phi,s)$ is \emph{$\sigma$-(semi)stable} if, for all proper $\phi$-invariant subbundles $F \subset E$, we have:
\begin{align*}
\frac{\deg(F)}{\rk(F)} \stless \frac{\deg(E)+\sigma}{\rk(E)} \qquad & \text{if } s \notin H^0(C, F)\\
\frac{\deg(F)+\sigma}{\rk(F)} \stless \frac{\deg(E)+\sigma}{\rk(E)} \qquad & \text{if } s \in H^0(C, F)
\end{align*}

We denote by $\mathcal{M}_\sigma^{r,d}(\gamma)$ the corresponding moduli spaces.\
\end{defn}
Many of the results from the previous chapters are still valid with no or small changes. We will discuss the differences in the following.

The BNR correspondence explained in section \ref{secbnr} relates Bradlow-Higgs $\gamma$-triples to coherent systems on the surface $\P X(\gamma)=\P (\mathcal{O}_C \oplus K(\gamma P))$ which is the compactification of the total space of $K (\gamma P)$.

We define a $\gamma$-Hitchin base
$$\mathcal{A}^r(\gamma)=\oplus_{i=0}^{r}H^0(K(\gamma P)^i)$$
which parametrizes spectral curves in $\P X(\gamma)$ and has dimension $r^2(g-1)+ \gamma r(r+1)/2$. Accordingly, there are also modified Hitchin maps:
\begin{align*}
\chi^{r,d}_\sigma(\gamma): \mathcal{M}_\sigma^{r,d}(\gamma) & \rightarrow \mathcal{A}^r(\gamma)\\
(E,\phi,s) & \mapsto \text{char poly}(\phi).
\end{align*}

We can also define a $\C^*$-action on the $\mathcal{M}_\sigma^{r,d}(\gamma)$ by scaling the Higgs field in the usual way. Also, the $\gamma$-Hitchin maps satisfy the same properties as the original one, namely they are equivariant with respect to the $\C^*$-action and they are proper.

The fixed point loci and attracting sets have a similar structure and it is possible to compute the dimension of the positive, negative and zero weight part of the $\C^*$-action restricted to the Zariski tangent space of a fixed point.
\begin{defn}
Given a partition $\underline{r}=(r_1, \dots ,r_m)$ of $r$, a partition $\underline{d}=(d_1, \dots, d_m)$ and $1 \leq k\leq m$ we denote by
\begin{equation*}
F_{\underline{r},\sigma}^{\underline{d},k}(\gamma)=\left \{
\begin{array}{c}
(E,\phi,s)\in \mathcal{M}_\sigma^{r,d}(\gamma) \text{ such that } E=E_1 \oplus \dots \oplus E_m\\
\text{with } \deg E_i=d_i, \rk{E_i}=r_i\\
\phi(E_i) \subseteq E_{i-1} \otimes K(\gamma P), \phi(E_1)=0\\
s \in H^0(E_k)
\end{array} \right \}
\end{equation*}
\end{defn}
\begin{defn}
Let us fix $F_{\underline{r},\sigma}^{\underline{d},k}(\gamma)$. We define two locally closed subsets by:
\begin{align*}
F_{\underline{r},\sigma}^{\underline{d},k+}(\gamma)&=\{(E,\phi,s) \in \mathcal{M}_\sigma^{r,d}(\gamma)\text{ such that } \lim_{\lambda \rightarrow 0}(E,\phi,s) \in F_{\underline{r},\sigma}^{\underline{d},k}(\gamma)\}\\
F_{\underline{r},\sigma}^{\underline{d},k-}(\gamma)&=\{(E,\phi,s) \in \mathcal{M}_\sigma^{r,d}(\gamma)\text{ such that } \lim_{\lambda \rightarrow \infty}(E,\phi,s) \in F_{\underline{r},\sigma}^{\underline{d},k}(\gamma)\}
\end{align*}
\end{defn}
Then, similar to the situation for $\mathcal{M}^{r,d}$ described in section \ref{semipr}, we have:
\begin{prop}
Let $(E,\phi,s)$ be a point in $F_{\underline{r},\sigma}^{\underline{d},k}(\gamma)$, for $\sigma$ not a critical value. With the previous notations, the weight $0$ part of the $\C^*$-action on $T_{(E,\phi,s)} \mathcal{M}_\sigma^{r,d}(\gamma)$ is given by the first hypercohomology of the complex:
\begin{equation*}
\bigoplus_{i=1}^{m} \Hom(E_i,E_i) \rightarrow \bigoplus_{i=2}^{m} \Hom(E_i,E_{i-1} \otimes K(\gamma)) \oplus E_k,
\end{equation*}
the positive weight part is the first hypercohomology of:
\begin{equation*}
\bigoplus_{i<j} \Hom(E_i,E_j) \rightarrow \bigoplus_{i\leq j} \Hom(E_i,E_{j} \otimes K(\gamma)) \oplus \bigoplus_{i>k} E_i
\end{equation*}
while the negative part is the first hypercohomology of:
\begin{equation*}
\bigoplus_{i>j}^{m} \Hom(E_i,E_j) \rightarrow \bigoplus_{i>j+1}^{m} \Hom(E_i,E_{j} \otimes K(\gamma)) \oplus \bigoplus_{i<k} E_i
\end{equation*}
\end{prop}

The remarks about extremal values of $\sigma$ and the $U$-filtration are still valid once we change $K$ with $K(\gamma P)$. There are new estimates on the degrees of the $U$-filtration. Namely
$$\deg U_i \geq i(i-1)(1-g)-\gamma \frac{i(i-1)}{2}.$$

In particular theorem \ref{extremal} is modified as follows:
\begin{theorem}
\begin{itemize}
%\item[(a)] $\mathcal{M}_\sigma^{r,d}$ is non-empty if and only if $d \geq r(r-1)(g-1)$.\\
\item[(i)] Assume that $\sigma$ is very close to $0$, then $\sigma$-stability for a triple $(E,\phi,s)$ implies the semistability of $(E,\phi)$ and so we have an Abel-Jacobi map:
\begin{equation*}
AJ(\gamma): \mathcal{M}_\sigma^{r,d}(\gamma) \rightarrow \mathcal{M}^{r,d}(\gamma).
\end{equation*}
For $d$ large enough (e.g. $d > r(2g-1)+(r-1)^2(2g-2+\gamma) $) for any semistable Higgs $\gamma$-bundle $(E,\phi)$ we have $H^1(E)=0$ and therefore $AJ(\gamma)$ is a projective bundle over the stable part of $\mathcal{M}^{r,d}(\gamma)$.\\
\item[(ii)] For $\sigma>(r-1)d+r(r-1)(r-2)(g-1)+\gamma r(r-1)(r-2)/2$ and a $\sigma$-stable triple $(E,\phi,s)$ corresponding to a pair $(\mathcal{F},s)$ the following three equivalent conditions are realized:
\begin{itemize}
\item there are no $\phi$-invariant subbundles of $E$ which contain the section
\item $s, \phi(s), \dots, \phi^{r-1}(s)$ generically generate $E$
\item $s$ as a map $\mathcal{O}_{\P X(\gamma)} \rightarrow \mathcal{F}$ has zero dimensional cokernel.
\end{itemize}
\end{itemize}
\end{theorem}

Also, theorem \ref{relhilb} modifies as follows.
\begin{theorem}
The Hitchin map:
\begin{equation*}
\chi^{r,d}_\infty(\gamma): \mathcal{M}_\infty^{r,d}(\gamma) \rightarrow \mathcal{A}^r(\gamma)
\end{equation*}
is the relative Hilbert scheme of $d+r(r-1)(g-1)+\gamma r(r-1)/2$ points over the family of spectral curves $\mathcal{A}^r(\gamma)$.
\end{theorem}
Deformation theory deserves a special discussion, as increasing $\gamma$ will make it easier for a Bradlow-Higgs $\gamma$-triple to be a smooth point of the corresponding moduli space.

We still have a complex that computes the dimension of the Zariski tangent space.
\begin{theorem}
Let $(E,\phi,s)$ be a $\sigma$-stable Bradlow-Higgs $\gamma$-triple, then the tangent space at $(E,\phi,s)$ is given by the first cohomology $\H^1(E,\phi,s)$ of the complex:
\begin{equation*}
C^0(\End E)\rightarrow C^1(\End E) \oplus C^0(K(\gamma P)\End E) \oplus C^0(E) \rightarrow C^1(E) \oplus C^1(K(\gamma P) \End E)
\end{equation*}
where the first map is $$p(k)=(dk, [k,\phi],k \cdot s)$$ and the second one is $$q(\tau, \nu, \gamma)=(\tau \cdot s+d \gamma, [\tau,\phi]+d \nu).$$
Furthermore, $\H^0(E,\phi,s)=0$.

The same result can be obtained from the hypercohomology of the complex:
\begin{align*}
\End E &\rightarrow K(\gamma P) \End E \oplus E\\
f & \mapsto ([f,\phi],f(s))
\end{align*}
from which we can also deduce the long exact sequence
\begin{align*}
0 \rightarrow H^0(\End E) \rightarrow H^0(K(\gamma P) \End E\oplus E) \rightarrow T_{(E,\phi,s)} \rightarrow \\ \nonumber
\rightarrow H^1(\End E) \rightarrow H^1(K(\gamma P) \End E\oplus E) \rightarrow \H^2 \rightarrow 0 
\end{align*}
\end{theorem}

The proof of the proposition is completely analogous to proposition \ref{deform}. This also allows us to compute the dimension by noting that
$$\dim T_{(E,\phi,s)} = d+ 2 r^2(g-1)+r(1-g) + r^2 \gamma + \dim \H^2 (E,\phi,s).$$

As in the case of Higgs bundles, we generically have $\H^2(E,\phi,s)=\H^2([\cdot,\phi])=0$. Note that specializing to $\gamma =0$ does not give the formula for the dimension of the moduli space of Bradlow-Higgs triples (compare with \cite[proposition 7.1]{nitsure1991moduli}).

Here $\H^2(E,\phi,s)$ is characterized in a slightly different way.

\begin{prop}
Let $(E,\phi,s)$ be a $\sigma$-stable $\gamma$-triple. Then $(\H^2(E,\phi,s))^*$ is the kernel of the following map:
\begin{align*}
H^0(\End E (-\gamma P) )\oplus H^0(K E^*) & \rightarrow H^0(K \End E)\\
(\alpha, \beta) &\mapsto [\alpha,\phi]+\beta \otimes s.
\end{align*}
In particular, $\dim \H^2(E,\phi,s) \geq \dim \H^2([\cdot, \phi])$.
\end{prop}

This is particularly important to observe because in general $H^0(\End E (-\gamma P) )$ is smaller than $H^0(\End E )$.

Let us discuss the case of rank 2 more in detail. We have the analogue of proposition \ref{negdegrk2}, with the same proof.
\begin{prop}
Let $d < 0$. Then $\mathcal{M}_\infty^{2,d}(\gamma)=\mathcal{M}_\varepsilon^{2,d}(\gamma)$ and both are non-empty iff $d \geq 2-2g-\gamma$.
\end{prop}

We can classify the fixed points as in theorem \ref{fpr2}.

\begin{theorem}
Let $d\geq 2-2g-\gamma$ be an integer and $\sigma>0$ different from a critical value. Then we can classify the components of the fixed point locus of $\mathcal{M}_\sigma^{2,d}(\gamma)$ as follows:
\begin{itemize}
\item[(i)] if $d \geq 0$ and $\sigma < d$ then one of the components of the fixed points for the $\C^*$-action is $F_{(2),\sigma}^{(d),1}(\gamma)=M_\sigma^{2,d}$, i.e. the moduli space of $\sigma$-stable Bradlow pairs embedded as triples with zero Higgs field. If $d <0$ then there are no $\sigma$-stable Bradlow pairs and so this component is empty.
\item[(ii)] if there exists and integer $m$ such that $\max\{0,(d-\gamma)/2 +1-g\} \leq m < \frac{d-\sigma}{2}$, then there exist components:
\begin{align*}
F_{(1,1),\sigma}^{(d_1,d_2),1}(\gamma)&=\left \{
\begin{array}{c}
(E,\phi,s)\in \mathcal{M}_\sigma^{2,d}(\gamma) \text{ such that } E=E_1 \oplus E_2\\
\text{with } \deg E_i=d_i, \rk{E_i}=1\\
\phi(E_2) \subseteq E_{1} \otimes K, \phi(E_1)=0\\
s \in H^0(E_1)
\end{array} \right \} \cong\\
&\cong S^{d_1}( C) \times S^{d_1-d_2+2g-2+\gamma}( C).
\end{align*}
Here $d_1$ and $d_2$ are integers satisfying $d_1+d_2=d$ and one of the following equivalent inequalities:
\begin{align*}
&\max\{0,(d-\gamma)/2+1-g\} \leq d_1 < \frac{d-\sigma}{2} \\
&\frac{d+\sigma}{2} < d_2 \leq \min\{d,(d+\gamma)/2+g-1\}\\
&\max\{0,2g-2+\gamma-d\}\leq d_1-d_2+2g-2+\gamma < 2g-2+\gamma-\sigma \text{ same parity as $d$}.
\end{align*}
\item[(iii)] There exist components:
\begin{align*}
F_{(1,1),\sigma}^{(d_1,d_2),2}(\gamma)&=\left \{
\begin{array}{c}
(E,\phi,s)\in \mathcal{M}_\sigma^{2,d}(\gamma) \text{ such that } E=E_1 \oplus E_2\\
\text{with } \deg E_i=d_i, \rk{E_i}=1\\
\phi(E_2) \subseteq E_{1} \otimes K, \phi(E_1)=0\\
s \in H^0(E_2)
\end{array} \right \} \cong\\
&\cong S^{d_2}( C) \times S^{d_1-d_2+2g-2+\gamma}( C).
\end{align*}
Here $d_1$ and $d_2$ are integers satisfying $d_1+d_2=d$ and one of the following equivalent inequalities:
\begin{align*}
&(d-\gamma)/2 +1-g \leq d_1 < \min \left \{\frac{d+\sigma}{2},d+1 \right \} \\
&\max\left \{-1,\frac{d-\sigma}{2} \right \}< d_2 \leq (d+\gamma)/2+g-1\\
& 0\leq d_1-d_2+2g-2+\gamma < 2g-2+\gamma+\min\{\sigma,d+1\} \text{ same parity as $d$}.
\end{align*}
\end{itemize}
\end{theorem}

We can now discuss the dimension of the Zariski tangent space at the fixed point. Recall that, since limits as $\lambda \rightarrow 0$ for $\C^*$ acting on any $\gamma$-triple always exist, if the dimension of the Zariski tangent space at a fixed point $(E,\phi,s)$ is already minimal, i.e. the fixed point is smooth, then all the triples whose limit as $\lambda \rightarrow 0$ is $(E,\phi,s)$ will also be smooth points. In the following always assume that $\sigma$ is different from a critical value.

Consider a non-split fixed point $(E,0,s) \in \mathcal{M}_\sigma^{2,d}(\gamma)$. Pick $\alpha \in H^0( \End E(-\gamma P))$ and $\beta \in H^0( KE^*)$. Certainly if $(\alpha, \beta) \in \H^2(E,\phi,s)^*$ then $\beta =0$. In particular $\H^2(E,\phi,s)^*=H^0( \End E(-\gamma P))$. Since $(E,s) \in M_\sigma^{2,d}$ is a $\sigma$-stable Bradlow pair, either $E$ is already semistable and then $H^0( \End E(-\gamma P))=0$ or $E$ fits into an extension:
\begin{equation*}
0 \rightarrow M \rightarrow E \rightarrow L \rightarrow 0
\end{equation*}
where $\deg M > \deg L$ and $\deg M < (d+\sigma)/2$. In this case
$$H^0( \End E(-\gamma P))=H^0(L^*M(-\gamma P))$$
regardless of whether or not the extension is split (see proof of lemma \ref{unstableext}). Since $\deg L^*M(-\gamma P) = \sigma - \gamma$, if $\gamma > \sigma$ then all these points are smooth.

A type 1 split fixed point is of the form $E=E_1 \oplus E_2$ with $\phi: E_2 \rightarrow E_1 K(\gamma P)$ and $s \in H^0(E_1)$. For these points we always have $\deg E_2 > \deg E_1$ and so
$$H^0( \End E(-\gamma P))=H^0(E_1^*E_2(-\gamma P)).$$
Pick $\alpha \in H^0( \End E(-\gamma P))$ and $\beta \in H^0( KE^*)$ then:
\begin{equation*}
[\alpha, \phi]=
\begin{pmatrix}
- \phi c & 0 \\
0 & c \phi
\end{pmatrix}
\end{equation*}
for $c \in H^0(E_1^*E_2(-\gamma P))$ and
\begin{equation*}
\beta \otimes s=
\begin{pmatrix}
\beta_1 s & \beta_2 s \\
0 & 0
\end{pmatrix}
\end{equation*}
for $\beta_i \in H^0(K E_i^*)$. If $(\alpha, \beta) \in  \H^2(E,\phi,s)^*$ then $\beta=0$ and $\alpha=0$ and so every type 1 split fixed point is smooth, regardless of $\sigma$.

Finally, a type 2 split fixed point is of the form $E=E_1 \oplus E_2$ with $\phi: E_2 \rightarrow E_1 K(\gamma P)$ and $s \in H^0(E_2)$. Here $H^0( \End E(-\gamma P))=H^0(E_1^*E_2(-\gamma P)) \oplus H^0(E_2^*E_1(-\gamma P)) $. Pick $\alpha \in H^0( \End E(-\gamma P))$ and $\beta \in H^0( KE^*)$ then:
\begin{equation*}
[\alpha, \phi]=
\begin{pmatrix}
- \phi c & 0 \\
0 & c \phi
\end{pmatrix}
\end{equation*}
for $c \in H^0(E_1^*E_2(-\gamma P))$, $b \in H^0(E_2^*E_1(-\gamma P))$ and
\begin{equation*}
\beta \otimes s=
\begin{pmatrix}
0 & 0\\
\beta_1 s & \beta_2 s
\end{pmatrix}
\end{equation*}
for $\beta_i \in H^0(K E_i^*)$. If $(\alpha, \beta) \in  \H^2(E,\phi,s)^*$ then $\beta=0$ and
$$\H^2(E,\phi,s)^*=H^0(E_2^*E_1(-\gamma P)).$$

Since $\deg E_2^*E_1(-\gamma P)=\deg E_1-\deg E_2 -\gamma < -\gamma + \min \{\sigma, d+1\}$, if $\gamma \geq \min\{\sigma ,d+1\}$ then the fixed point is certainly smooth. 

We can summarize the main consequences of the above discussion in the following proposition.
\begin{cor}
$\mathcal{M}_\varepsilon^{2,d}(\gamma)$ is smooth for all $\gamma \geq 1$. In particular $\mathcal{M}_\varepsilon^{2,d}(\gamma)$ is always semiprojective.

If $\gamma > d$ then $\mathcal{M}_\sigma^{2,d}(\gamma)$ is smooth regardless of $\sigma$ as long as it is different from a critical value. In this case then $\mathcal{M}_\sigma^{2,d}(\gamma)$ is semiprojective.
\end{cor}

A further consequence of this is the following theorem.

\begin{theorem}
Let $\gamma \geq 1$ be an integer. Then if $d <0$ we have:
\begin{align*}
[\mathcal{M}_{\varepsilon}^{2,d}(\gamma)]=\sum_{(d_1,d_2) \in I_{2,\varepsilon}(d,\gamma)} \L^{4g-4+3\gamma} [S^{d_2}( C)][S^{d_1-d_2+2g-2+\gamma}( C)]
\end{align*}
and if $d \geq 0$:
\begin{align*}
&[\mathcal{M}_{\varepsilon}^{2,d}(\gamma)]=\L^{4g-4+4 \gamma} [M_{\varepsilon}^{2,d}] + \sum_{(d_1,d_2) \in I_{1,\varepsilon}(d,\gamma)} \L^{3g-3+3\gamma +d_2} [S^{d_1}( C)][S^{d_1-d_2+2g-2+\gamma}( C)]+\\
&+\sum_{(d_1,d_2) \in I_{2,\varepsilon}(d,\gamma)} \L^{4g-4+3\gamma} [S^{d_2}( C)][S^{d_1-d_2+2g-2+\gamma}( C)]
\end{align*}
where $I_{1,\varepsilon}(d,\gamma)$ is the set of pairs of integers $(d_1,d_2)$ satisfying $d_1+d_2=d$ and:
\begin{align*}
&\max\{0,(d-\gamma)/2+1-g\} \leq d_1 < \frac{d}{2} \\
&\frac{d}{2} < d_2 \leq \min\{d,(d+\gamma)/2+g-1\}\\
&\max\{0,2g-2+\gamma-d\}\leq d_1-d_2+2g-2+\gamma < 2g-2+\gamma \text{ same parity as $d$}
\end{align*}
while $I_{2,\varepsilon}(d,\gamma)$ is the set of pairs of integers $(d_1,d_2)$ satisfying $d_1+d_2=d$ and:
\begin{align*}
&(d-\gamma)/2 +1-g \leq d_1 < \min \left \{\frac{d+1}{2},d+1 \right\} \\
&\max\left \{-1,\frac{d-1}{2} \right\}< d_2 \leq (d+\gamma)/2+g-1\\
& 0\leq d_1-d_2+2g-2+\gamma < 2g-2+\gamma+\min\{1,d+1\} \text{ same parity as $d$}.
\end{align*}

Furthermore, if $\gamma > d$:
\begin{align*}
[\mathcal{M}_{\infty}^{2,d}(\gamma)]=\sum_{(d_1,d_2) \in I_{2,\infty}(d,\gamma)} \L^{4g-4+3\gamma} [S^{d_2}( C)][S^{d_1-d_2+2g-2+\gamma}( C)]
\end{align*}
where $I_{2,\infty}(d,\gamma)$ is the set of pairs of integers $(d_1,d_2)$ satisfying $d_1+d_2=d$ and:
\begin{align*}
&(d-\gamma)/2 +1-g \leq d_1 \leq d \\
&0 \leq d_2 \leq (d+\gamma)/2+g-1\\
& 0\leq d_1-d_2+2g-2+\gamma \leq 2g-2+\gamma+d \text{ same parity as $d$}.
\end{align*}
\end{theorem}

In particular, as done in \cite[section 7.2]{hausel2001geometry} we can consider the embeddings:
\begin{equation*}
\mathcal{M}_{\sigma}^{2,d}(\gamma) \rightarrow \mathcal{M}_{\sigma}^{2,d}(\gamma+1)
\end{equation*}
and compute, for $\sigma=\infty$ and $\sigma=\varepsilon$
\begin{equation*}
\lim_{\gamma \rightarrow \infty} P(\mathcal{M}_{\sigma}^{2,d}(\gamma),t).
\end{equation*}
\begin{theorem}
We have:
\begin{align*}
&\lim_{\gamma \rightarrow \infty} P(\mathcal{M}_{\varepsilon}^{2,d}(\gamma),t) = \lim_{\gamma \rightarrow \infty} P(\mathcal{M}_{\infty}^{2,d}(\gamma),t)=\\
&=\frac{(1+t^3)^{2g}(1+t)^{2g}}{(1-t^2)^2(1-t^4)}=P(\C \P^{\infty},t) P(B \overline{\mathcal{G}},t)
\end{align*}
where $B \overline{\mathcal{G}}$ is the classifying space of the group $\overline{\mathcal{G}}$ mentioned in \cite[section 7.2]{hausel2001geometry}.
\begin{proof}
Recall that $\dim \mathcal{M}_{\varepsilon}^{2,d}(\gamma) = \dim \mathcal{M}_{\infty}^{2,d}(\gamma) = d+ 6g-6+4 \gamma$. We can also directly compute the Poincar\'e polynomial of $\mathcal{M}_{\varepsilon}^{2,d}(\gamma)$ since it is semiprojective. Assuming $d \geq 0$, we have:
\begin{align*}
&P(\mathcal{M}_{\varepsilon}^{2,d}(\gamma),t)=P(M_{\varepsilon}^{2,d},t) +\\
&+ \sum_{(d_1,d_2) \in I_{1,\varepsilon}(d,\gamma)} t^{2g-2+2d-4 d_1} P(S^{d_1}( C),t)P(S^{d_1-d_2+2g-2+\gamma}( C),t)+\\
&+\sum_{(d_1,d_2) \in I_{2,\varepsilon}(d,\gamma)} t^{2 d_2} P(S^{d_2}( C),t)P(S^{d_1-d_2+2g-2+\gamma}( C),t).
\end{align*}

Therefore:
\begin{align*}
&\lim_{\gamma \rightarrow \infty} P(\mathcal{M}_{\varepsilon}^{2,d}(\gamma),t)=P(M_{\varepsilon}^{2,d},t) + t^{2g-2} P(S^{\infty}( C),t) \sum_{i=0}^{\lfloor{d/2}\rfloor} t^{2d-4 i} P(S^{i}( C),t)+\\
&+ P(S^{\infty}( C),t) \sum_{i>d/2} t^{2 i} P(S^{i}( C),t)=\\
&=P(M_{\varepsilon}^{2,d},t)+P(S^{\infty}( C),t) \sum_{i\geq 0} t^{2i} P(S^{i}( C),t)+\\
&+P(S^{\infty}( C),t) \left ( \sum_{i=0}^{\lfloor{d/2}\rfloor} (t^{2g-2+2d-4i}-t^{2i})P(S^{i}( C),t) \right ).
\end{align*}

From \cite[remark 6.3]{mozgovoy2013moduli} it is immediate to deduce that:
\begin{equation*}
P(M_{\varepsilon}^{2,d},t)=P(S^{\infty}( C),t) \left ( \sum_{i=0}^{\lfloor{d/2}\rfloor} (t^{2i}-t^{2g-2+2d-4i})P(S^{i}( C),t) \right )
\end{equation*}
and hence:
\begin{align*}
&\lim_{\gamma \rightarrow \infty} P(\mathcal{M}_{\varepsilon}^{2,d}(\gamma),t)=P(S^{\infty}( C),t) \sum_{i\geq 0} t^{2i} P(S^{i}( C),t)=\\
&=\frac{(1+t^3)^{2g}(1+t)^{2g}}{(1-t^2)^2(1-t^4)}
\end{align*}
as we wanted.

For $d<0$, $\mathcal{M}_{\varepsilon}^{2,d}(\gamma)=\mathcal{M}_{\infty}^{2,d}(\gamma)$. Otherwise:
\begin{align*}
&P(\mathcal{M}_{\infty}^{2,d}(\gamma),t)=\sum_{(d_1,d_2) \in I_{2,\infty}(d,\gamma)} t^{2 d_2} P(S^{d_2}( C),t)P(S^{d_1-d_2+2g-2+\gamma}( C),t)
\end{align*}
and so:
\begin{align*}
&\lim_{\gamma \rightarrow \infty} P(\mathcal{M}_{\infty}^{2,d}(\gamma),t)= P(S^{\infty}( C),t) \sum_{i\geq 0} t^{2i} P(S^{i}( C),t)=\\
&=\frac{(1+t^3)^{2g}(1+t)^{2g}}{(1-t^2)^2(1-t^4)}
\end{align*}
as claimed.
\end{proof}
\end{theorem}

We conclude the chapter by examining what happens during the wall-crossing, i.e. what is the structure of the flip loci.

As in the case of Bradlow-Higgs triples we have two flip loci coming with canonical maps:
\begin{equation*}
\pi_{\bar\sigma}^{d,+}(\gamma): \mathcal{W}_{\bar\sigma}^{d,+}(\gamma) \rightarrow X_{\bar\sigma}^{d}(\gamma),
\end{equation*}
\begin{equation*}
\pi_{\bar\sigma}^{d,-}(\gamma): \mathcal{W}_{\bar\sigma}^{d,-}(\gamma) \rightarrow X_{\bar\sigma}^{d}(\gamma)
\end{equation*}
where
\begin{equation*}
X_{\bar\sigma}^{d}(\gamma)=S^{(d-\bar\sigma)/2}( C) \times J^{(d-\bar\sigma)/2}( C) \times H^0(K(\gamma P))^2.
\end{equation*}

Triples $(E,\phi,s) \in \mathcal{W}_{\bar\sigma}^{d,+}(\gamma)$ are characterized as non-split extensions
\begin{equation*}
0 \rightarrow (L,\psi_1) \rightarrow (E,\phi) \rightarrow (M, \psi_2) \rightarrow 0
\end{equation*}
where $\deg L =(d-\bar\sigma)/2$, $\deg M =(d+\bar\sigma)/2$ and $s \in H^0(L)$. Here $\psi_i \in H^0(K(\gamma P))$. For such a triple $\pi_{\bar\sigma}^{d,+}(\gamma)(E,\phi,s)=(L,s,M,\psi_1,\psi_2)$.

If we fix a point $(L,s,M,\psi_1,\psi_2) \in X_{\bar\sigma}^{d}(\gamma)$ then the fiber of $\pi_{\bar\sigma}^{d,+}(\gamma)$ over such a point can be once again characterized using a hypercohomology group $\H^1((M,\psi_2),(L,\psi_1))$ whose dimension is:
\begin{equation*}
\dim \H^1((M,\psi_2),(L,\psi_1))=
\begin{cases}
2g-2+\gamma & \text{if } \psi_1 \neq \psi_2\\
2g-2+\gamma + \dim H^0(L^*M(-\gamma P)) & \text{if } \psi_1=\psi_2.
\end{cases}
\end{equation*}

In particular, if $\gamma > d$, $\mathcal{W}_{\bar\sigma}^{d,+}(\gamma)$ is a $\C \P^{2g-3+\gamma}$-bundle over $X_{\bar\sigma}^{d}(\gamma)$ regardless of $\bar\sigma$.

Triples $(E,\phi,s) \in \mathcal{W}_{\bar\sigma}^{d,-}(\gamma)$ are characterized as extensions
\begin{equation*}
0 \rightarrow (M,\psi_2) \rightarrow (E,\phi) \rightarrow (L, \psi_1) \rightarrow 0
\end{equation*}
where $\deg L =(d-\bar\sigma)/2$, $\deg M =(d+\bar\sigma)/2$ and $s \in H^0(E)$ projects to a nonzero $\bar{s} \in H^0(L)$. We must have that either the extension is non-split or it is split but $s \neq \bar{s}$.

If we fix a point $(L,s,M,\psi_1,\psi_2) \in X_{\bar\sigma}^{d}(\gamma)$ then the fiber of $\pi_{\bar\sigma}^{d,-}(\gamma)$ over such a point can be once again characterized using a hypercohomology group $ \widetilde\H^1((L, \bar{s},\psi_1),(M,\psi_2))$ whose dimension is:
\begin{equation*}
\dim  \widetilde\H^1((L, \bar{s},\psi_1),(M,\psi_2))=(d+\bar\sigma)/2+g-1+\gamma
\end{equation*}

In particular, $\mathcal{W}_{\bar\sigma}^{d,-}(\gamma)$ is always a $\C \P^{(d+\bar\sigma)/2+g-2+\gamma}$-bundle over $X_{\bar\sigma}^{d}(\gamma)$.

%% file: MSMY.tex
\chapter{Hilbert schemes and compactified Jacobians}
\label{MSMY}
%
%\begin{quote}
%{\small
%Summary of chapter contents\\
%\begin{itemize}
%\item Classical Macdonald formula
%\item Macdonald formula for integral curves (perverse filtration, compactified Jacobian)
%\item Macdonald formula for reduced curves (perverse filtration, properties of the families, fine compactified Jacobians
%\item How to get $F=\sum P(\mathcal{M}_\varepsilon^{2,d})$ using the motives
%\item Properties of $F$ (denominator, evaluation at $q=1$, Serre duality)
%\item How to get $PH(\mathcal{M}^{2,1})$ and $G$ (character variety, weight filtration, P=W)
%\item Properties of $G$
%\item Comments and observations about $F-G$ (proof of the known ranges)
%\item Generalizations to higher rank (interesting degree ranges)
%\item ???Twisted Higgs bundles
%\end{itemize}
%}
%\end{quote}
%
\section{Maulik-Yun and Migliorini-Shende-Viviani formulas}
Recently there has been a lot of steady progress in relating the cohomology of the Hilbert scheme of points of a locally planar curve and the cohomology of compactified Jacobians of the same curve.

In \cite{migliorini2011support}, \cite{maulik2014macdonald} and \cite{rennemo2013homology} we find three different approaches to generalize the well known formula by Macdonald relating the cohomology of the symmetric powers of a smooth projective curve and the cohomology of the Jacobian. The last progress made in this direction is \cite{migliorini2015support}.

Let $X$ be a complex projective and locally planar curve. If $X$ is smooth of genus $g$, denote by $S^n(X)$ its $n$-th symmetric power and by $J(X)$ its Jacobian. The classical Macdonald formula states that:
\begin{equation*}
\sum_{n \geq 0} P(S^n(X) ,t) q^n = \frac{(1+qt)^{2g}}{(1-q)(1-qt^2)}=\frac{P(J(X),qt)}{(1-q)(1-qt^2)}.
\end{equation*}

If $X$ is integral then the Macdonald formula has been generalized to this case independently in \cite{migliorini2011support} and \cite{maulik2014macdonald}, later reproved in \cite{rennemo2013homology}. Denote by $X^{[n]}$ the Hilbert scheme of $n$ points on $X$ and by $\overline{J}({X})$ the compactified Jacobian parametrizing degree zero rank one torsion free sheaves on ${X}$. A version of the generalized Macdonald formula can be stated as follows (see for example \cite[formula 1.4]{maulik2014macdonald}):
\begin{equation*}
\bigoplus_{n \geq 0} H^*({X}^{[n]}) q^n = \frac{\oplus_{i=0}^{2 g_a} \Gr_i^P H^*(\overline{J}({X})) q^i}{(1-\Q q)(1-\Q[-2](-1) q)}.
\end{equation*}

The equality is meant to be of (cohomologically) graded vector spaces. Also, $\Gr_i^P$ refers to the $i$-th graded piece with respect to the perverse filtration. This particular filtration on the cohomology of $\overline{J}({X})$ can be defined by deforming ${X}$ in a flat family of curves $\pi:\mathcal{X} \rightarrow \mathcal{B}$ with the following properties (see \cite[section 2.1]{maulik2014macdonald}):
\begin{itemize}
\item $\mathcal{B}$ is irreducible.
\item The fibers of $\pi$ are integral and locally planar.
\item Let $\pi^{[n]}: \mathcal{X}^{[n]} \rightarrow \mathcal{B}$ denote the Hilbert scheme of $n$ points relative to the family $\pi$. Then we assume the total space $\mathcal{X}^{[n]}$ is smooth for all $n \geq 0$.
\item If $b \in \mathcal{B}$ we denote by $\delta(b)$ the delta invariant of $\pi^{-1}(b)$. Then we assume that $codim_\mathcal{B}(\overline{\{b\}}) \geq \delta(b)$ where $\overline{\{b\}}$ is the Zariski closure of the point $b \in \mathcal{B}$.
\end{itemize}

After we have any of these families $\pi:\mathcal{X} \rightarrow \mathcal{B}$ we also get the relative compactified Jacobian $\pi^J:\mathcal{J} \rightarrow \mathcal{B}$ whose total space is smooth as a consequence of the previous axioms. The map $\pi^J$ is proper and the object $\R \pi^J_* \Q_\mathcal{J}$ is filtered by its perverse truncations $\leftidx{^p}{\tau}{^{\leq i}}\R \pi^J_* \Q_\mathcal{J}$. The perverse filtration is thus defined by restricting the images of the natural maps:
\begin{equation*}
H^*\left (\leftidx{^p}{\tau}{^{\leq i}}\R \pi^J_* \Q_\mathcal{J} \right) \rightarrow H^*\left (\R \pi^J_* \Q_\mathcal{J} \right).
\end{equation*}
to the fiber of $\pi$ corresponding to the starting curve ${X}$. In the papers is also explained why the perverse filtration does not depend on the choice of the particular family. Furthermore, a similar version of the formula holds for the cohomology of the total spaces $\mathcal{X}^{[n]}$ and $\mathcal{J} $ with the perverse filtration. In fact in both \cite{migliorini2011support} and \cite{maulik2014macdonald} the formula is obtained from a stronger sheaf theoretic statement, i.e. an equality in $D_c^b(\mathcal{B})[[q]]$ (see \cite[formula 4]{migliorini2011support}):
\begin{equation*}
\bigoplus_{n \geq 0} q^n \R\pi_*^{[n]} \Q= \frac{\bigoplus_i q^i \cdot \R^i \pi_*^J \Q[-i]}{(1-q)(1-q\Q[-2](-1))}.
\end{equation*}

To deduce the statements for a single curve or for the entire family of curves it is then enough to take the stalks or the global cohomology.

For us, the most relevant example of such a family comes from the Hitchin fibration. Regard the Hitchin base $\mathcal{A}^r$ as the base parametrizing the family of spectral curves $\mathcal{C} \rightarrow \mathcal{A}^r$. Then in \cite[proposition 3.3]{maulik2014macdonald} is proved that the restriction of this family to the locus where the spectral curves are integral $\mathcal{C}_{int} \rightarrow \mathcal{A}^r_{int}$ satisfies the properties that are necessary for the generalized Macdonald formula to hold.

Note that the fibers of the Hitchin map $h: \mathcal{M}^{n,d} \rightarrow \mathcal{A}^r$ are indeed (torsors over) compactified Jacobians over $\mathcal{A}^r_{int}$. Also, as we already noted in proposition \ref{relhilb}, the fibers of any of the $\chi_\sigma^{n,d}: \mathcal{M}_\sigma^{n,d} \rightarrow \mathcal{A}^r$ over integral spectral curves are Hilbert schemes and we also proved that the smooth locus of the $\mathcal{M}_\sigma^{n,d}$ is contained in the integral locus. Therefore in the case of $\mathcal{A}^r_{int}$ we know what the analogues of the relative compactified Jacobian and relative Hilbert schemes are. In the family of spectral curves, however, also appear curves that are reducible and even non-reduced.

There has been a recent and even broader generalization of the Macdonald formula to reduced locally planar curves in \cite{migliorini2015support}, even though the assumptions on the family of curves that is allowed are more restrictive. Let us give a brief overview for reduced curves, details can be found in \cite{migliorini2015support}.

Consider a flat family $\mathcal{C} \rightarrow \mathcal{B}$ of reduced locally planar curves containing at least a reducible curve. To generalize the Macdonald formula, several technical assumptions are needed and we highlight the most important ones here:
\begin{itemize}
\item The family has to be \emph{independently broken} \cite[definition 1.10]{migliorini2015support} meaning that we should have a finite set of flat families $\mathcal{C}_S \rightarrow \mathcal{B}_S$ indexed by subsets $S$ of a fixed finite set $V$ satisfying some compatibility properties. The idea behind this definition is that we should be able to distinguish the irreducible components of our curves in families.
%More precisely, if $S' \cup S''=S$ is a disjoint union, then we have an injective map $\mathcal{B}_{S'}\times \mathcal{B}_{S''} \rightarrow \mathcal{B}_S$ and a fibrewise partial normalization $\mathcal{C}_{S'} \times \mathcal{B}_{S''} \bigsqcup \mathcal{B}_{S'} \times \mathcal{C}_{S''} \rightarrow \mathcal{C}_{S| \mathcal{B}_{S'} \times \mathcal{B}_{S''}}$. Furthermore, every decomposition of the curves in the family should appear in some of these partial normalizations and the set $\{\mathcal{C}_S \rightarrow \mathcal{B}_S \}$ should satisfy compatibility relations along iterated disjoint unions.
\item The family has to be \emph{H-smooth} \cite[definition 1.12]{migliorini2015support}, i.e. the relative Hilbert schemes $\mathcal{C}^{[n]}$ should have smooth total space for all $n \geq 0$ and also the families of irreducible components $\mathcal{C}_S \rightarrow \mathcal{B}_S$ should have smooth relative Hilbert schemes.
\end{itemize}

In the case of reduced curves it is also necessary to modify the notion of Jacobian (see \cite[section 2.3]{migliorini2015support}). In fact if a curve has more than one irreducible component, the moduli problem of torsion free sheaves will not be of finite type and then one has to impose a further stability condition. We only mention the main existence result \cite[theorem 2.8]{migliorini2015support}, originally \cite[theorems A and B]{esteves2001compactifying}.
%Roughly speaking, a \emph{general polarization} $\underline{m}$ on a reduced curve $X$ (having $I$ as the index set of its irreducible components) is the datum of one rational number $m_i$ for every irreducible component $X_i$ of $X$ such that for every proper subcurve $Y \subset X$, $Y=\bigcup_{i \in J} X_i$ there is a connected component $Y'$ either of $Y$ or of the complementary subcurve $\bigcup_{i \in I \setminus J} X_i$ for which the sum of the $m_i$ appearing for the irreducible components of $Y'$ is not in $\Z$. Given a general polarization on $X$ we can define a stability condition for rank one pure one dimensional sheaves $\mathcal{I}$ on $X$ by declaring semistable the sheaves for which:
%\begin{equation*}
%\chi(\mathcal{I}_D) \geq \sum_D m_i
%\end{equation*}
%for all proper subcurves $D$ of $X$. Here $\mathcal{I}_D$ denotes the restriction of $\mathcal{I}$ to $D$ modulo its biggest torsion subsheaf, $\chi$ its Euler-Poincar\'e characteristic and $\sum_D m_i$ the sum of the $m_i$ for $i$ in the set of irreducible components of $X$ appearing in $D$. The sheaves for which all previous inequalities are strict are called stable.
%
%For general polarizations on locally planar curves semistability and stability coincide and it is also possible to construct fine moduli spaces of $\underline{m}$-stable rank one pure one dimensional sheaves on $X$, which are called \emph{fine compactified Jacobians}.

\begin{theorem}
Let $X$ be a geometrically connected projective locally planar reduced curve and $\underline{m}$ a general polarization on $X$. Then there exists a projective scheme $\overline{J}_X(\underline{m})$ which is a fine moduli space for rank one pure one dimensional sheaves that are semistable with respect to $\underline{m}$.
\end{theorem}
It is also crucial to note that given a flat family $\pi: \mathcal{C} \rightarrow \mathcal{B}$ of reduced geometrically connected and locally planar curves, we can always construct, up to an \'etale cover of $\mathcal{B}$, a relative fine compactified Jacobian $\pi^J: \overline{J}_\mathcal{C}\rightarrow \mathcal{B}$. A more detailed discussion on the properties of fine compactified Jacobians can be found in \cite{melo2014fine} and \cite{esteves2001compactifying}.

%Before stating the main result we need one last ingredient. Consider the set of bases for an independently broken family $\{\mathcal{B}_S\}$, then for every disjoint union $S=S' \bigsqcup S''$ we identify $\mathcal{B}_{S'}\times \mathcal{B}_{S''}$ as subschemes of $\mathcal{B}_{S}$ using the injective maps in the definition. Define an exponential operator acting on the category of sheaves on $\bigsqcup_S \mathcal{B}_S$ by:
%\begin{equation*}
%\Exp(\mathcal{F})_{|\mathcal{B}_S}=\bigoplus_{S=\bigsqcup S_i} \boxtimes \mathcal{F}_{|B_{S_i}}
%\end{equation*}
%see \cite[formula 1]{migliorini2015support}.

The formula proved in \cite[theorem 1.16]{migliorini2015support} compares once again $\bigoplus_{n=0}^\infty q^n \R \pi_*^{[n]} \Q$ and $\bigoplus_i q^i \cdot \leftidx{^p}{\R}{^i} \pi_*^J \Q [-i]$.

%With these hypotheses the sheaf theoretic result can be stated as follows (see \cite[theorem 1.16]{migliorini2015support}):
%\begin{theorem}
%Let $\pi: \mathcal{C} \rightarrow \mathcal{B}$ be an H-smooth independently broken family with index set $V$, admitting relative fine compactified Jacobians $\overline{\mathcal{J}}_S \rightarrow \mathcal{B}_S$. Denote by $g$ the locally constant function computing the arithmetic genus of the curves in the family. Denote by $\L$ the constant shifted sheaf $\Q[-2](-1)$. The following formula holds in $D_c^b(\bigsqcup_S \mathcal{B}_S)[[q]]$:
%\begin{equation*}
%(q \L)^{1-g} \bigoplus_{n=0}^\infty q^n \R \pi_*^{[n]} \Q = \Exp \left ( (q \L)^{1-g} \frac{\bigoplus_i q^i \cdot \leftidx{^p}{\R}{^i} \pi_*^J \Q [-i]}{(1-q)(1-q \L)}\right ).
%\end{equation*}
%\end{theorem}
%A few remarks are in order. First of all, $\leftidx{^p}{\R}{^i} \pi_*^J \Q [-i] = IC(\bigwedge^i \R^1)[-i]$ where $\R^1$ is the local system obtained by pushing forward $\Q$ from the smooth part of the family of relative compactified Jacobians. Second, it is not mentioned what is the correct choice of general polarization $\underline{m}$ because in fact the formula is true regardless of the choice.

However, contrary to what happened for integral curves, $\bigoplus_{n=0}^\infty q^n \R \pi_*^{[n]} \Q$ has extra summands that are not supported on the full base of the family, but rather on proper closed subvarieties. These extra summands come from partial normalizations of the curves in the family that can be disconnected and therefore clearly cannot appear for integral curves.

Another fundamental difference is that, in the case of integral curves, the Hilbert scheme of points can be characterized as a moduli space of torsion free sheaves with the extra datum of a section. There is then a forgetful map that will send such a pair to the underlying sheaf which is then a torsion free sheaf and so an element of the compactified Jacobian. When the curve is smooth this is the classical Abel-Jacobi map, sending a divisor on the curve to the associated line bundle. For a sufficiently high number of points such a map is even a projective bundle. In particular this implies that for an integral locally planar curve (or a family of such curves) it is equivalent to require that:
\begin{itemize}
\item[(i)] the relative Hilbert schemes have smooth total space for all $n$
\item[(ii)] the relative Hilbert schemes have smooth total space up to $n=2g-1$
\item[(iii)] the relative compactified Jacobian has smooth total space.
\end{itemize}

When the curve instead is reduced, due to stability issues, we do not have a map from the Hilbert scheme of points to the fine compactified Jacobian. It is still true that if the family of curves is H-smooth then the relative fine compactified Jacobian has smooth total space but the converse implication is not true. For more details about the previous remarks see \cite[section 1]{migliorini2015support}.

Let us discuss briefly what happens for our current example: the family of spectral curves. The first difference we see is that in the family there are non reduced curves. This is possibly the biggest issue because, at the moment, it is not clear what the analogue of the Jacobian should be for such curves. The Hilbert scheme of points for a non reduced curve is defined but it is harder to understand. %Later we will see what the wall-crossing approach can say about the motivic invariants of Hilbert schemes of non-reduced curves.\\

%Consider $r,d$ integers such that $(r,d)=1$. It is a well know fact \marginpar{\footnotesize I was told about this by Luca, I can't find a reference though} that the fibers of the Hitchin map $h: \mathcal{M}^{r,d} \rightarrow \mathcal{A}^r$ over reduced curves are fine compactified Jacobians with specific general polarizations. More precisely, if a spectral curve $\widetilde{C}$ is not irreducible then it is the union of integral spectral curves for lower rank problems, say $r_1, \dots, r_k$ such that $r_1+\dots+r_k=r$. The polarization will then assign to the irreducible components that are spectral curves for a rank $r_i$ problem the weight $r_i/r$.

For the relative Hilbert scheme of points instead the situation is more complicated. We already proved that the moduli spaces $\mathcal{M}_\sigma^{r,d}$, together with the Hitchin maps $\chi_\sigma^{r,d}: \mathcal{M}_\sigma^{r,d} \rightarrow \mathcal{A}^r$ provide different extensions of the relative Hilbert schemes of points for the family $\mathcal{C}_{int} \rightarrow \mathcal{A}^r_{int}$ of integral spectral curves to $\mathcal{A}^r \setminus \mathcal{A}^r_{int}$. Only one of these extensions, namely $\mathcal{M}_\infty^{r,d}$ will be a relative Hilbert scheme, but we proved that it is not smooth and the singular locus will intersect the locus where the spectral curves are reduced. If we want some kind of H-smoothness property to be satisfied, then the $\mathcal{M}_\varepsilon^{r,d}$ provide smooth extensions but unfortunately not for all combinations of $r$ and $d$.

Another interesting property that the $\mathcal{M}_\varepsilon^{r,d}$ satisfy is the existence of an Abel-Jacobi map $\mathcal{M}_\varepsilon^{r,d} \rightarrow \mathcal{M}^{r,d}$ that forgets the section of the triple. The fibers of this map are projective spaces for all $d$, not necessarily of the same dimension. For large $d$ however, the map is a projective bundle.
\section{A formula for rank 2 Bradlow-Higgs triples}
At this point the natural question is whether or not we can compare
\begin{equation*}
\bigoplus_{n \geq 0} \R (\chi_\varepsilon^{r,n+\theta(r)})_* IC_{\mathcal{M}_\varepsilon^{r,n+\theta(r)}} \cdot q^n
\end{equation*}
and
\begin{equation*}
\frac{\bigoplus_i \leftidx{^p}{\R}{^i} h^{r,d}_* \Q [-i] \cdot q^i}{(1-q)(1-q\L)}
\end{equation*}
or some variant of the two expressions.

Here $\theta$ is a shift function that only depends on the rank $r$. Recall from proposition \ref{degs} that for a rank on pure one dimensional sheaf $\mathcal{F}$ on $T^*C$ not intersecting the divisor at infinity and $E= \pi_*\mathcal{F}$ we have the relation $\deg E= \deg \mathcal{F} - r(r-1)(g-1)$. Therefore $\theta(r)=-r(r-1)(g-1)$ in the formula above will allow us to relate $\mathcal{M}_\varepsilon^{r,d}$ with the relative Hilbert scheme of $d+ r(r-1)(g-1)$ points.
%For the rank 2 case it is actually possible to compute the cohomology of both sides directly and compare them. We will need a few propositions first.
%

From theorems \ref{motodd} and \ref{moteven} we have the motive of $\mathcal{M}^{2,d}_\varepsilon$. However, the motive will allow to compute the Poincar\'e polynomial only when $\mathcal{M}^{2,d}_\varepsilon$ is smooth, which happens only when either $d<0$ or for $d\geq 4g-5$ odd. Since in general computing intersection cohomology is a very hard task, we will concentrate mainly on the previous degrees.

For the following, denote by $P$ the Poincar\'e polynomial and by $E$ the E-polynomial. We will also use $P^{vir}$ to denote the following specialization of $E$. For a variety $X$:
\begin{equation*}
P^{vir}(X,t)=t^{2 \dim X} E(X, -1/t,-1/t).
\end{equation*}
Note that $P^{vir}(X,t)=P(X,t)$ when the cohomology of $X$ is pure.

Let us also define the following generating functions:
\begin{defn}
\begin{equation*}
F^{sh}(q)= \sum_{n \geq 1-g} \R (\chi_\varepsilon^{2,2n+1})_* (IC_{\mathcal{M}_\varepsilon^{2,2n+1}}) q^{2n+2g-1},
\end{equation*}
\begin{equation*}
F^{mot}(q)= \sum_{n \geq 1-g} [\mathcal{M}_\varepsilon^{2,2n+1}] q^{2n+2g-1},
\end{equation*}
\begin{equation*}
F^{vir}(q,t)= \sum_{n \geq 1-g} P^{vir}(\mathcal{M}_\varepsilon^{2,2n+1},t) q^{2n+2g-1}
\end{equation*}
\end{defn}
We note immediately that $F^{vir}$ is equal to the generating function of the Poincar\'e polynomials for $\deg q \leq 2g-3$ and $\deg q \geq 6g-5$. Also for $2n+1<0$ and $2n+1 \geq 4g-3$ the smoothness of $\mathcal{M}_\varepsilon^{2,2n+1}$ implies that
$$IC_{\mathcal{M}_\varepsilon^{2,2n+1}}=\Q_{\mathcal{M}_\varepsilon^{2,2n+1}}.$$

Using the motives it is possible to compute $F^{mot}$ and $F^{vir}$ explicitly.
\begin{theorem}
We have:
\begin{align*}
&F^{mot}(q)=\L^{4g-3}q^{2g-1}\left (\frac{\L^{g}[J ( C)] Z(C,q^2)}{(\L-1)(1-\L^2 q^2)}-\frac{[J ( C)] Z(C,\L q^2)}{(\L-1)(1- q^2)}\right )+\\
&+\L^{5g-4}q^{4g-4} Z(C,\L q^2) \Gamma(\L^{-1}q^{-1})+\L^{4g-3} Z(C,q^2) \Gamma(q)+\\
&+\Theta(q)
\end{align*}
where
\begin{equation*}
\Gamma(q)=\sum_{j=0}^{g-2} [S^{2j+1}( C)]q^{2j+1}
\end{equation*}
and 
\begin{align*}
\Theta(q)=&\sum_{n=0}^{2g-3} \sum_{i=0}^n (\L^{4g-2}-\L^{4g-3}) [S^{i}( C)] \cdot\\
& \cdot\left ([S^{2n+1-i}( C)]-[J( C)]\frac{\L^{2n+2-i-g}-1}{\L-1} \right)q^{2n+2g-1}.
\end{align*}
Also:
\begin{align*}
F^{vir}(q,t)&= \frac{q^{2g-1}(1+t)^{2g} (1+q^2t^3)^{2g}}{(1-t^2)(1-q^2)(1-q^2t^2)(1-q^2t^4)} - \frac{q^{2g-1}t^{2g}(1+t)^{2g} (1+q^2t)^{2g}}{(1-t^2)(1-q^2)(1-q^2t^2)(1-q^2t^4)}+\\
&+\frac{1}{2} q^{2g-2}t^{4g-4} \left ( \frac{(1+q^2 t)^{2g}}{(1-q^2)(1-q^2t^2)}\right ) \left ( {\frac {( 1+qt )^{2g}}{( 1-qt^2 )( 1-q ) }}- \frac{( 1-qt )^{2g}}{( 1+qt^2)( 1+q ) } \right )+\\
&+ \frac{1}{2} q^{2g-2}t^{4g-4} \left ( \frac{(1+q^2 t)^{2g}}{(1-q^2)(1-q^2t^2)}\right ) \left ( \frac{(1+t)^{2g}}{1-t^2} \left ( \frac{qt^{4-2g}}{1-q^2t^4} - \frac{q}{1-q^2}\right ) \right )+\\
&+\frac{1}{2}\left ( \frac{(1+q^2 t^3)^{2g}}{(1-q^2t^2)(1-q^2t^4)}\right )\left ( {\frac {( 1+qt )^{2g}}{( 1-qt^2 )( 1-q ) }}- \frac{( 1-qt )^{2g}}{( 1+qt^2)( 1+q ) } \right )+\\
&- \frac{q^{2g-1}(1+t)^{2g}}{1-t^2}\left ( \frac{(1+q^2 t^3)^{2g}}{(1-q^2t^2)(1-q^2t^4)}\right ) \left ( \frac{1}{1-q^2} - \frac{t^{2g}}{1-q^2t^4} \right )+\\
&+t^{8g-6} E(\Theta(u),-1/t,-1/t)_{|u=qt^2}.
\end{align*}
\begin{proof}
Recall from theorem \ref{motodd} that if $d >0$ odd then:
\begin{align*}
[\mathcal{M}_\varepsilon^{2,d}]=&\L^{1+4(g-1)}[M_\varepsilon^{2,d}]+\sum_{(d_1,d_2) \in I_1^o(d)} \L^{1+3(g-1)+d_2}[S^{d_1}( C)][S^{d_1-d_2+2g-2}]+\\
&+\sum_{(d_1,d_2) \in I_2^o(d)}\L^{1+4(g-1)} [S^{d_2}( C)][S^{d_1-d_2+2g-2}( C)]+\\
&+\sum_{(d_1,d_2) \in I_1^o(d)} (\L^{4g-2}-\L^{4g-3}) [S^{d_1}( C)] ([S^{d_2}( C)]-[J( C)][\C\P^{d_2-g}])
\end{align*}
and if $d <0 $ odd then:
\begin{align*}
[\mathcal{M}_\varepsilon^{2,d}]=&\sum_{(d_1,d_2) \in I_2^o(d)}\L^{1+4(g-1)} [S^{d_2}( C)][S^{d_1-d_2+2g-2}( C)].
\end{align*}
We can therefore split the computation into 4 main terms:
\begin{align*}
F^{mot}&(q)=\sum_{n \geq 0} \L^{4g-3}[M_\varepsilon^{2,2n+1}] q^{2n+2g-1} +\\
&+\sum_{n \geq 0}q^{2n+2g-1} \left (\sum_{(d_1,d_2) \in I_1^o(2n+1)} \L^{1+3(g-1)+d_2}[S^{d_1}( C)][S^{d_1-d_2+2g-2}] \right ) +\\
&+\sum_{n \geq 1-g}q^{2n+2g-1} \left (\sum_{(d_1,d_2) \in I_2^o(d)}\L^{1+4(g-1)} [S^{d_2}( C)][S^{d_1-d_2+2g-2}( C)]\right )+\\
&+\sum_{n \geq 0} \sum_{(d_1,d_2) \in I_1^o(2n+1)} (\L^{4g-2}-\L^{4g-3}) [S^{d_1}( C)] \cdot \\
&\cdot \left ([S^{d_2}( C)]-[J( C)] \frac{\L^{d_2+1-g}-1}{\L-1} \right ) q^{2n+2g-1}.
\end{align*}

From \cite[Remark 6.3]{mozgovoy2013moduli}, we deduce:
\begin{align*}
\sum_{n \geq 0} [M_\varepsilon^{2,2n+1}] u^{n}=\frac{\L^{g}[J ( C)] Z(C,u)}{(\L-1)(1-\L^2 u)}-\frac{[J ( C)] Z(C,\L u)}{(\L-1)(1- u)}
\end{align*}
where
$$Z(C,u)=\sum_{n \geq 0} [S^n( C)]u^n$$
is the motivic zeta function of $C$. Therefore:
\begin{align*}
&\sum_{n \geq 0} \L^{4g-3}[M_\varepsilon^{2,2n+1}] q^{2n+2g-1}= \L^{4g-3}q^{2g-1}\sum_{n \geq 0} [M_\varepsilon^{2,2n+1}] q^{2n}=\\
&=\L^{4g-3}q^{2g-1}\left (\frac{\L^{g}[J ( C)] Z(C,q^2)}{(\L-1)(1-\L^2 q^2)}-\frac{[J ( C)] Z(C,\L q^2)}{(\L-1)(1- q^2)}\right ).
\end{align*}

We also have:
\begin{align*}
&\sum_{n \geq 0} \left (\sum_{(d_1,d_2) \in I_1^o(2n+1)} \L^{1+3(g-1)+d_2}[S^{d_1}( C)][S^{d_1-d_2+2g-2}] \right )q^{2n+2g-1}=\\
&=\L^{3g-2}q^{2g-1} \sum_{n=0}^{g-2} \sum_{j=g-2-n}^{g-2} [S^{j+n+2-g}( C)][S^{2j+1}] \L^{n-j+g-1} q^{2n}+\\
&+\L^{3g-2}q^{2g-1} \sum_{n\geq g-1} \sum_{j=0}^{g-2} [S^{j+n+2-g}( C)][S^{2j+1}] \L^{n-j+g-1} q^{2n}=\\
&=\L^{3g-2}q^{2g-1} \sum_{n=0}^{g-2} \sum_{j=g-2-n}^{g-2} [S^{j+n+2-g}( C)][S^{2j+1}] \L^{n-j+g-1} q^{2n}+\\
&+\L^{3g-2}q^{2g-1} \sum_{j=0}^{g-2} \sum_{i\geq j+1}  [S^{i}( C)][S^{2j+1}] \L^{i-2j+2g-3} q^{2i-2j+2g-4}=\\
&=\L^{3g-2}q^{2g-1} \sum_{n=0}^{g-2} \sum_{j=g-2-n}^{g-2} [S^{j+n+2-g}( C)][S^{2j+1}] \L^{n-j+g-1} q^{2n}+\\
&+\L^{3g-2}q^{2g-1} \sum_{j=0}^{g-2} [S^{2j+1}] \L^{-2j+2g-3} q^{-2j+2g-4} \sum_{i\geq j+1}  [S^{i}( C)]\L^{i} q^{2i}=\\
&=\L^{3g-2}q^{2g-1} \sum_{n=0}^{g-2} \sum_{j=g-2-n}^{g-2} [S^{j+n+2-g}( C)][S^{2j+1}] \L^{n-j+g-1} q^{2n}+\\
&+\L^{3g-2}q^{2g-1} \sum_{j=0}^{g-2} [S^{2j+1}] \L^{-2j+2g-3} q^{-2j+2g-4} \left(Z(C,\L q^2)- \sum_{i=0}^j  [S^{i}( C)]\L^{i} q^{2i}\right)=\\
&=\L^{5g-4}q^{4g-4} Z(C,\L q^2) \Gamma(\L^{-1}q^{-1}).
\end{align*}

For the other type of split fixed points we get:
\begin{align*}
&\sum_{n \geq 1-g} \sum_{(d_1,d_2) \in I_2^o(d)}\L^{1+4(g-1)} [S^{d_2}( C)][S^{d_1-d_2+2g-2}( C)]q^{2n+2g-1}=\\
&=\L^{4g-3}q^{2g-1} \sum_{n=1-g}^{-1} \sum_{j=0}^{n-1+g}[S^{2j+1}( C)][S^{n-1+g-j}( C)] q^{2n}+\\
&+\L^{4g-3}q^{2g-1} \sum_{n\geq 0} \sum_{j=0}^{g-2}[S^{2j+1}( C)][S^{n-1+g-j}( C)] q^{2n}=\\
&=\L^{4g-3}q^{2g-1} \sum_{n=1-g}^{-1} \sum_{j=0}^{n-1+g}[S^{2j+1}( C)][S^{n-1+g-j}( C)] q^{2n}+\\
&+\L^{4g-3}q^{2g-1} \sum_{j=0}^{g-2} \sum_{i\geq g-1-j} [S^{2j+1}( C)][S^{i}( C)] q^{2i+2j+2-2g}=\\
&=\L^{4g-3}q^{2g-1} \sum_{n=1-g}^{-1} \sum_{j=0}^{n-1+g}[S^{2j+1}( C)][S^{n-1+g-j}( C)] q^{2n}+\\
&+\L^{4g-3}q^{2g-1} \sum_{j=0}^{g-2} [S^{2j+1}( C)] q^{2j+2-2g}\sum_{i\geq g-1-j} [S^{i}( C)] q^{2i}=\\
&=\L^{4g-3}q^{2g-1} \sum_{n=1-g}^{-1} \sum_{j=0}^{n-1+g}[S^{2j+1}( C)][S^{n-1+g-j}( C)] q^{2n}+\\
&+\L^{4g-3}q^{2g-1} \sum_{j=0}^{g-2} [S^{2j+1}( C)] q^{2j+2-2g}\left( Z(C,q^2) - \sum_{i=0}^{g-2-j} [S^{i}( C)] q^{2i}\right )=\\
&=\L^{4g-3} Z(C,q^2) \Gamma(q).
\end{align*}

For the last term we have:
\begin{align*}
&\sum_{n \geq 0}\sum_{(d_1,d_2) \in I_1^o(d)} (\L^{4g-2}-\L^{4g-3}) [S^{d_1}( C)] \left([S^{d_2}( C)]-[J( C)]\frac{\L^{d_2+1-g}-1}{\L-1}\right)q^{2n+2g-1} =\\
&=\sum_{n=0}^{g-2}  \sum_{i=0}^n (\L^{4g-2}-\L^{4g-3}) [S^{i}( C)] \left([S^{2n+1-i}( C)]-[J( C)]\frac{\L^{2n+2-i-g}-1}{\L-1}\right)q^{2n+2g-1}+\\
&+\sum_{n=g-1}^{2g-3} \sum_{i=2n-2g+3}^n (\L^{4g-2}-\L^{4g-3}) [S^{i}( C)]\cdot\\
&\cdot \left([S^{2n+1-i}( C)]-[J( C)]\frac{\L^{2n+2-i-g}-1}{\L-1}\right)q^{2n+2g-1} 
\end{align*}

since for $d_2 \geq 2g-1$,
$$[S^{d_2}( C)]-[J( C)]\frac{\L^{d_2+1-g}-1}{\L-1}=0.$$

For the same reason,
\begin{align*}
&\sum_{n=g-1}^{2g-3} \sum_{i=2n-2g+3}^n (\L^{4g-2}-\L^{4g-3}) [S^{i}( C)] \cdot \\
&\cdot\left([S^{2n+1-i}( C)]-[J( C)]\frac{\L^{2n+2-i-g}-1}{\L-1} \right)q^{2n+2g-1}=\\
=&\sum_{n=g-1}^{2g-3} \sum_{i=0}^n (\L^{4g-2}-\L^{4g-3}) [S^{i}( C)] \left([S^{2n+1-i}( C)]-[J( C)]\frac{\L^{2n+2-i-g}-1}{\L-1}\right)q^{2n+2g-1}
\end{align*}
and therefore
\begin{align*}
&\sum_{n \geq 0}\sum_{(d_1,d_2) \in I_1^o(d)} (\L^{4g-2}-\L^{4g-3}) [S^{d_1}( C)] \left ([S^{d_2}( C)]-[J( C)]\frac{\L^{d_2+1-g}-1}{\L-1}\right)q^{2n+2g-1} =\\
=&\sum_{n=0}^{2g-3} \sum_{i=0}^n (\L^{4g-2}-\L^{4g-3}) [S^{i}( C)] \left([S^{2n+1-i}( C)]-[J( C)]\frac{\L^{2n+2-i-g}-1}{\L-1}\right)q^{2n+2g-1}.
\end{align*}

In order to get the expression for $F^{vir}(q,t)$ we can argue as follows. Since
$$P^{vir}(\mathcal{M}_\varepsilon^{2,2n+1},t)= t^{2 \dim \mathcal{M}_\varepsilon^{2,2n+1}} E(\mathcal{M}_\varepsilon^{2,2n+1},-1/t,-1/t)$$
and $\dim \mathcal{M}_\varepsilon^{2,2n+1}=2n+1+1+6g-6=2n+2+6(g-1)$ we have:
\begin{equation*}
F^{vir}(q,t)=t^{8g-6} E(F^{mot}(u) , -1/t,-1/t)_{|u=qt^2}.
\end{equation*}

Recall that:
\begin{align*}
&E(\L)=xy\\
&E(Z(C,u))=\frac{(1-xu)^g(1-yu)^g}{(1-u)(1-uxy)}\\
&E(J( C))=(1-x)^g(1-y)^g.
\end{align*}

Also, since:
\begin{equation*}
\sum_{n=0}^{2g-2} [S^n( C)] u^n=Z(C,u)- \frac{\L^g[J( C)] u^{2g-1}}{(\L-1)(1-\L u)}+\frac{[J( C)] u^{2g-1}}{(\L-1)(1-u)},
\end{equation*}
we get:
\begin{align*}
&\Gamma(u)=\frac{1}{2}\left ( Z(C,u)-Z(C,-u)\right)-\frac{1}{2}\left (\frac{\L^g[J( C)] u^{2g-1}}{(\L-1)(1-\L u)}+\frac{\L^g[J( C)] u^{2g-1}}{(\L-1)(1+\L u)} \right)+\\
&+\frac{1}{2}\left ( \frac{[J( C)] u^{2g-1}}{(\L-1)(1-u)} + \frac{[J( C)] u^{2g-1}}{(\L-1)(1+u)} \right).
\end{align*}

The last ingredient we need to compute $F^{vir}$ is:
\begin{align*}
&E(\Gamma(u),-1/t,-1/t)=\frac{1}{2} \left ( \frac{(1+u/t)^{2g}}{(1-u)(1-u/t^2)}-\frac{(1-u/t)^{2g}}{(1+u)(1+u/t^2)}\right )+\\
&-\frac{1}{2} \frac{(1+t)^{2g}}{t^{4g-2}(1-t^2)}\left (\frac{u^{2g-1}}{(1-u/t^2)}+\frac{u^{2g-1}}{(1+u/t^2)} \right )+\\
&+\frac{1}{2} \frac{(1+t)^{2g}}{t^{2g-2}(1-t^2)}\left (\frac{u^{2g-1}}{(1-u)}+\frac{u^{2g-1}}{(1+u)} \right )
\end{align*}

Therefore:
\begin{align*}
F^{vir}(q,t)&= \frac{q^{2g-1}(1+t)^{2g} (1+q^2t^3)^{2g}}{(1-t^2)(1-q^2)(1-q^2t^2)(1-q^2t^4)} - \frac{q^{2g-1}t^{2g}(1+t)^{2g} (1+q^2t)^{2g}}{(1-t^2)(1-q^2)(1-q^2t^2)(1-q^2t^4)}+\\
&+\frac{1}{2} q^{2g-2}t^{4g-4} \left ( \frac{(1+q^2 t)^{2g}}{(1-q^2)(1-q^2t^2)}\right ) \left ( {\frac {( 1+qt )^{2g}}{( 1-qt^2 )( 1-q ) }}- \frac{( 1-qt )^{2g}}{( 1+qt^2)( 1+q ) } \right )+\\
&+ \frac{1}{2} q^{2g-2}t^{4g-4} \left ( \frac{(1+q^2 t)^{2g}}{(1-q^2)(1-q^2t^2)}\right ) \left ( \frac{(1+t)^{2g}}{1-t^2} \left ( \frac{qt^{4-2g}}{1-q^2t^4} - \frac{q}{1-q^2}\right ) \right )+\\
&+\frac{1}{2}\left ( \frac{(1+q^2 t^3)^{2g}}{(1-q^2t^2)(1-q^2t^4)}\right )\left ( {\frac {( 1+qt )^{2g}}{( 1-qt^2 )( 1-q ) }}- \frac{( 1-qt )^{2g}}{( 1+qt^2)( 1+q ) } \right )+\\
&- \frac{q^{2g-1}(1+t)^{2g}}{1-t^2}\left ( \frac{(1+q^2 t^3)^{2g}}{(1-q^2t^2)(1-q^2t^4)}\right ) \left ( \frac{1}{1-q^2} - \frac{t^{2g}}{1-q^2t^4} \right )+\\
&+t^{8g-6} E(\Theta(u),-1/t,-1/t)_{|u=qt^2}.
\end{align*}
\end{proof}
\end{theorem}
We can also prove some properties of $F^{mot}$ and $F^{vir}$. Note that, since we have an explicit formula, the proof of these properties is straightforward. However, they can be deduced \emph{a priori} from the geometry of the $\mathcal{M}_\varepsilon^{2,d}$.%and the arguments we will use in the following proposition can actually be generalized to the case of higher rank, for which we don't have a formula.
\begin{prop}
The following holds for $F^{mot}$ and $F^{vir}$:
\begin{itemize}
\item[(i)] There exists a polynomial $Q^{mot}(q)$ such that
\begin{equation*}
F^{mot}(q)=\frac{Q^{mot}(q)}{(1-q^2)(1-q^2\L^2)}
\end{equation*}
and a polynomial $Q^{vir}(q,t)$ such that
\begin{equation*}
F^{vir}(q,t)=\frac{Q^{vir}(q,t)}{(1-q^2)(1-q^2t^4)}
\end{equation*}
\item[(ii)] $Q^{mot}(1)=(1+\L)[\mathcal{M}^{2,1}]$ and $Q^{vir}(1,t)=(1+t^2)P(\mathcal{M}^{2,1},t)$
\item[(iii)] $Q^{vir}(q,t)$ satisfies $$Q^{vir}(q,t)=(qt)^{8g-4}Q^{vir}(q^{-1}t^{-2},t)$$ and $F^{vir}$ satisfies $$F^{vir}(q,t)=(qt)^{8g-8}F^{vir}(q^{-1}t^{-2},t).$$
\end{itemize}
\begin{proof}
\begin{itemize}
\item[(i)] First of all observe that for $n \geq 3g-3$ the Abel-Jacobi map $\mathcal{M}_\varepsilon^{2,2n+1} \rightarrow \mathcal{M}^{2,2n+1}$ is a projective bundle of rank $2n+2-2g$. Also, $\mathcal{M}^{2,2n+1}$ is always isomorphic to $\mathcal{M}^{2,1}$. In particular, for $n \geq 3g-1$ we have:
\begin{equation*}
[\mathcal{M}_\varepsilon^{2,2n+1}]-(1+\L^2)[\mathcal{M}_\varepsilon^{2,2n-1}] + \L^2[\mathcal{M}_\varepsilon^{2,2n-3})]=0
\end{equation*}
and
\begin{equation*}
P(\mathcal{M}_\varepsilon^{2,2n+1},t)-(1+t^4)P(\mathcal{M}_\varepsilon^{2,2n-1},t) + t^4P(\mathcal{M}_\varepsilon^{2,2n-3})=0.
\end{equation*}
Now the property for $F^{mot}$ can be deduced from the fact that:
\begin{align*}
&(1-q^2)(1-q^2\L^2)F^{mot}(q)= Q^{mot}(q)+\\
&+\sum_{n \geq 3g-1} q^{2n+1+2g-2} ( [\mathcal{M}_\varepsilon^{2,2n+1}]-(1+\L^2)[\mathcal{M}_\varepsilon^{2,2n-1}] + \L^2[\mathcal{M}_\varepsilon^{2,2n-3}] ) =\\
&= Q^{mot}(q)
\end{align*}
for some polynomial $Q^{mot}$. And for $F^{vir}$:
\begin{align*}
&(1-q^2)(1-q^2t^4)F^{vir}(q,t)= Q^{vir}(q,t)+\\
&+\sum_{n \geq 3g-1} q^{2n+1+2g-2} ( P(\mathcal{M}_\varepsilon^{2,2n+1},t)-(1+t^4)P(\mathcal{M}_\varepsilon^{2,2n-1},t)+ t^4P(\mathcal{M}_\varepsilon^{2,2n-3},t) ) =\\
&= Q^{vir}(q,t)
\end{align*}
for some polynomial $Q^{vir}$. We can retrieve an explicit expression for $Q^{vir}$ by observing that:
\begin{align*}
&(1-q^2)(1-q^2t^4)F(q,t) = \\
&=q P(\mathcal{M}_\varepsilon^{2,3-2g},t) + q^3 P(\mathcal{M}_\varepsilon^{2,5-2g},t)-q^3 (1+t^4) P(\mathcal{M}_\varepsilon^{2,3-2g},t) + \\
&+\sum_{n=3-g}^{3g-2}q^{2n+1+2g-2} (P(\mathcal{M}_\varepsilon^{2,2n+1},t)-(1+t^4)P(\mathcal{M}_\varepsilon^{2,2n-1},t) + t^4P(\mathcal{M}_\varepsilon^{2,2n-3},t)).
\end{align*}
\item[(ii)] For $Q^{mot}$ we have:
\begin{align*}
&F^{mot}(q)=\sum_{n=1-g}^{3g-4}q^{2n+1+2g-2}  [\mathcal{M}_\varepsilon^{2,2n+1}] + \sum_{n\geq 3g-3}q^{2n+1+2g-2}  [\mathcal{M}_\varepsilon^{2,2n+1}]=\\
&=\sum_{n=1-g}^{3g-4}q^{2n+1+2g-2}   [\mathcal{M}_\varepsilon^{2,2n+1}] +  \frac{[\mathcal{M}^{2,1}]}{\L-1}\left ( \frac{\L^{4g-3}q^{8g-7}}{1-\L^2q^2} - \frac{q^{8g-7} }{1-q^2}\right ).
\end{align*}
Therefore:
\begin{align*}
Q^{mot}(1)&=\left ((1-q^2)(1-\L^2 q^2)F^{mot}(q) \right )_{|q=1}=(1+\L)[\mathcal{M}^{2,1}].
\end{align*}
For $Q^{vir}$, similarly:
\begin{align*}
&F^{vir}(q,t)=\sum_{n=1-g}^{3g-4}q^{2n+1+2g-2}  P^{vir}(\mathcal{M}_\varepsilon^{2,2n+1},t) + \sum_{n\geq 3g-3}q^{2n+1+2g-2}  P(\mathcal{M}_\varepsilon^{2,2n+1},t)=\\
&=\sum_{n=1-g}^{3g-4}q^{2n+1+2g-2}  P^{vir}(\mathcal{M}_\varepsilon^{2,2n+1},t) + \sum_{n\geq 3g-3}q^{2n+1+2g-2}  P(\mathcal{M}^{2,1},t)P(\C\P^{2n+2-2g})=\\
&=\sum_{n=1-g}^{3g-4}q^{2n+1+2g-2}  P^{vir}(\mathcal{M}_\varepsilon^{2,2n+1},t) +  \frac{q^{2g-1}P(\mathcal{M}^{2,1},t)}{1-t^2}\sum_{n\geq 3g-3}q^{2n} (1-t^{4n+4-4g})=\\
&=\sum_{n=1-g}^{3g-4}q^{2n+1+2g-2}  P^{vir}(\mathcal{M}_\varepsilon^{2,2n+1},t) +  \frac{q^{2g-1}P(\mathcal{M}^{2,1},t)}{1-t^2}\left ( \frac{q^{6g-6}}{1-q^2} - \frac{q^{6g-6} t^{8g-8}}{1-q^2t^4}\right ).
\end{align*}
Therefore:
\begin{align*}
Q^{vir}(1,t)&=\left ((1-q^2)(1-q^2t^4)F^{vir}(q,t) \right )_{|q=1} = (1+t^2) P(\mathcal{M}^{2,1},t).
\end{align*}
\item[(iii)] The functional equation for $F^{vir}$ follows from the one for $Q^{vir}$. First of all we prove the following motivic identity, valid for all $n \geq g-1$:
\begin{equation*}
[\mathcal{M}_\varepsilon^{2,2n+1}]-\L^{2n-2g+3}[\mathcal{M}_\varepsilon^{2,4g-4-2n-1}] =[\mathcal{M}^{2,1}][\C\P^{2n-2g+2}].
\end{equation*}
We start by defining:
\begin{equation*}
Y_k^{2n+1}=\{(E,\phi) \in \mathcal{M}^{2,2n+1} : \dim H^0(E)=k\}.
\end{equation*}
Recall that we have a Serre duality isomorphism:
\begin{align*}
\mathcal{M}^{2,2n+1} &\rightarrow \mathcal{M}^{2,4g-4-2n-1}\\
(E,\phi) & \mapsto (K E^*,\phi).
\end{align*}
Since $\dim H^0(E)- \dim H^0(K E^*) = 2n+1+2-2g$ we see that:
\begin{equation*}
[Y_k^{2n+1}]=[Y_{k+2g-3-2n}^{4g-4-2n-1}].
\end{equation*}
This relation in turn implies that:
\begin{align*}
&[\mathcal{M}_\varepsilon^{2,2n+1}]-\L^{2n-2g+3}[\mathcal{M}_\varepsilon^{2,4g-4-2n-1}]=\\
&= \sum_i [Y_i^{2n+1}][\C\P^{i-1}] - \L^{2n-2g+3}\sum_j [Y_j^{4g-4-2n-1}][\C\P^{j-1}]=\\
&=\sum_i [Y_i^{2n+1}][\C\P^{i-1}-  \L^{2n-2g+3} ] \sum_j [Y_{j+2n+3-2g}^{2n+1}][\C\P^{j-1}]=\\
&=\sum_i [Y_i^{2n+1}][\C\P^{i-1}]- \L^{2n-2g+3}  \sum_i [Y_{i}^{2n+1}][\C\P^{i+2g-4-2n}]  =\\
&=\sum_i [Y_i^{2n+1}] \frac{\L^{i}-1-\L^{i}+\L^{2n-2g+3}}{\L-1}=[\mathcal{M}^{2,1}][\C\P^{2n-2g+2}].
\end{align*}

The motivic relation, in turn, implies the following relation between the virtual Poincar\'e polynomials:
\begin{equation*}
P^{vir}(\mathcal{M}_\varepsilon^{2,2n+1},t)-t^{4n+6-4g}P^{vir}(\mathcal{M}_\varepsilon^{2,4g-5-2n},t)=P^{vir}(\mathcal{M}^{2,1},t)P(\C\P^{2n-2g+2},t).
\end{equation*}
To prove the identity for $Q^{vir}$, let us first write:
\begin{equation*}
Q^{vir}(q,t)=\sum_{l=0}^{4g-3} a_l(t) q^{2l+1}
\end{equation*}
for:
\begin{equation*}
a_l(t)=P^{vir}(\mathcal{M}_\varepsilon^{2,2l+3-2g},t)-(1+t^4)P^{vir}(\mathcal{M}_\varepsilon^{2,2l+1-2g},t)+t^4P^{vir}(\mathcal{M}_\varepsilon^{2,2l-1-2g},t).
\end{equation*}

Note that some of the coefficients of $Q^{vir}$ do not include three terms, but the previous formula holds for those coefficients as well since for the appropriate choice of $l$ some of the three moduli spaces appearing in the formula can be empty.\\
The identity $Q^{vir}(q,t)=(qt)^{8g-4}Q^{vir}(q^{-1}t^{-2},t)$ is equivalent to:
\begin{equation*}
a_l(t)=t^{4l+6-8g}a_{4g-3-l}(t).
\end{equation*}

Also note that it is enough to check the previous equation for $l \geq 2g-1$. In fact assume it holds for all $l \geq 2g-1$ and choose $k \leq 2g-2$. Then $4g-3-k \geq 2g-1$ and so:
\begin{equation*}
a_{4g-3-k}(t)=t^{8g-6-4k}a_{k}(t)
\end{equation*}
which is just the equation for $k$.\\
We can compute:
\begin{align*}
a_l(t)-&t^{4l+6-8g}a_{4g-3-l}(t)=\left(P^{vir}(\mathcal{M}_\varepsilon^{2,2l+3-2g},t)-t^{4l+10-8g}P^{vir}(\mathcal{M}_\varepsilon^{2,6g-7-2l},t)\right)+\\
&-(1+t^4)\left (P^{vir}(\mathcal{M}_\varepsilon^{2,2l+1-2g},t) -t^{4l+6-8g} P^{vir}(\mathcal{M}_\varepsilon^{2,6g-5-2l},t)\right )+\\&+t^4 \left ( P^{vir}(\mathcal{M}_\varepsilon^{2,2l-1-2g},t) -t^{4l+2-8g} P^{vir}(\mathcal{M}_\varepsilon^{2,6g-3-2l},t)\right).
\end{align*}
For $2l-1-2g \geq 2g-1$, i.e. for $l \geq 2g$ we can use the identity we proved above, obtaining:
\begin{align*}
a_l(t)-&t^{4l+6-8g}a_{4g-3-l}(t)= P^{vir}(\mathcal{M}^{2,1},t)P(\C\P^{2l-4g+4},t)+\\
&-(1+t^4)P^{vir}(\mathcal{M}^{2,1},t)P(\C\P^{2l-4g+2},t)+t^4 P^{vir}(\mathcal{M}^{2,1},t)P(\C\P^{2l-4g},t)=\\
&=\frac{P^{vir}(\mathcal{M}^{2,1},t)}{1-t^2} (1-t^{4l-8g+10}-(1+t^4)(1-t^{4l-8g+6})+t^4(1-t^{4l-8g+2}))=0.
\end{align*}
Using the same techniques, with a bit of care, we can check the last case $l=2g-1$:
\begin{align*}
a_{2g-1}(t)-&t^{2}a_{2g-2}(t)=\left(P^{vir}(\mathcal{M}_\varepsilon^{2,2g+1},t)-t^{6}P^{vir}(\mathcal{M}_\varepsilon^{2,2g-5},t)\right)+\\
&-(1+t^4)\left (P^{vir}(\mathcal{M}_\varepsilon^{2,2g-1},t) - t^2 P(\mathcal{M}_\varepsilon^{2,2g-3},t)\right )+\\
&+t^4 P^{vir}(\mathcal{M}_\varepsilon^{2,2g-3},t) -t^{2} P^{vir}(\mathcal{M}_\varepsilon^{2,2g-1},t)=\\
&=P^{vir}(\mathcal{M}^{2,1},t) \left( P(\C\P^2,t)-(1+t^4)-t^2\right)=0.
\end{align*}
This completes the proof.
\end{itemize}
\end{proof}
\end{prop}

We now define the second kind of generating functions involved in the formula.
\begin{defn}
Denote by $\text{odd}_q$ the operator that acts on $f \in R[[q]]$ (for some ring $R$) by deleting the even powers of $q$ in $f$. Define:
\begin{equation*}
G(q,t)=\text{odd}_q \left (\frac{PH(\mathcal{M}^{2,1},q,t)}{(1-q)(1-qt^2)} \right )
\end{equation*}
and
\begin{equation*}
G^{sh}(q)=\text{odd}_q \left (\frac{\bigoplus_{i=0}^{8g-6} IC\left (\bigwedge^i\R^1 \right) }{(1-q\Q)(1-q\Q[-2](-1))} \right )
\end{equation*}
where $\R^1$ is the local system $\R^1(h^{2,1}_{sm})_* \Q$ obtained by first restricting $h^{2,1}: \mathcal{M}^{2,1}\rightarrow \mathcal{A}^2$ to the locus of smooth spectral curves and then pushing forward the constant sheaf.
\end{defn}
\begin{rmk}
\label{decthm}
Let us make a couple comments about the previous definition. First of all, the reason we take the odd powers in $G$ and $G^{sh}$ is that we want to compare them with $F^{vir}$ and $F^{sh}$ respectively and they only contains odd powers.

Second, it follows from \cite[theorem 1.1.2]{de2010topology} and \cite[lemma 1.3.5]{de2010topology} that:
\begin{equation*}
\R h^{2,1}_* \Q = \bigoplus_{i=0}^{8g-6} \leftidx{^p}{\R}{^i} h^{2,1}_* \Q[-i] \oplus \mathscr{L} = \bigoplus_{i=0}^{8g-6} IC\left (\bigwedge^i\R^1 \right) \oplus \mathscr{L} 
\end{equation*}
where $\mathscr{L}$ is a complex supported on the locus of reduced curves in $\mathcal{A}^2$. In other words, the summands of $\R h^{2,1}_* \Q$ that are supported on the whole $\mathcal{A}^2$ are exactly the middle extensions of the external powers of $\R^1$, while $\mathscr{L}$ is an extra summand with proper support. Another important remark is that the perversity is actually determined by the index of the external power of $\R^1$. Furthermore, $\mathscr{L}$ will have no global cohomology, as I learned from private communication (\cite{demighein}). These facts together imply that taking the global cohomology of the object:
\begin{equation*}
\bigoplus_{i=0}^{8g-6} \leftidx{^p}{\R}{^i}  h^{2,1}_* \Q[-i] q^i \in D_c^b(\mathcal{A}^2)[[q]]
\end{equation*}
will compute $PH(\mathcal{M}^{2,1},q,t)$.

Last, note that from the results summarized in section \ref{charvar}, we can actually explicitly compute $G$.
\end{rmk}

Observe that $G$ satisfies properties that are similar to the ones of $F^{vir}$.
\begin{prop}
The following holds for $G$:
\begin{itemize}
\item[(i)] There exists a polynomial $V(q,t)$ such that
\begin{equation*}
G(q,t)=\frac{V(q,t)}{(1-q^2)(1-q^2t^4)}
\end{equation*}
\item[(ii)] $V(1,t)=(1+t^2)P(\mathcal{M}^{2,1},t)$
\item[(iii)] $V(q,t)$ satisfies $$V(q,t)=(qt)^{8g-4}V(q^{-1}t^{-2},t)$$ and $G$ satisfies $$G(q,t)=(qt)^{8g-8}G(q^{-1}t^{-2},t)$$
\end{itemize}
\begin{proof}
\begin{itemize}
\item[(i)] Note that, if $U(q,t) \in \Z[t][[q]]$, then
\begin{equation*}
\text{odd}_q U(q,t)= \frac{1}{2} \left ( U(q,t)- U(-q,t)\right ).
\end{equation*}
Therefore
\begin{align*}
G(q,t)&=\frac{1}{2} \left (\frac{PH(\mathcal{M}^{2,1},q,t)}{(1-q)(1-qt^2)} - \frac{PH(\mathcal{M}^{2,1},-q,t)}{(1+q)(1+qt^2)}\right )=\\
&=\frac{1}{2} \left (\frac{(1+q)(1+qt^2) PH(\mathcal{M}^{2,1},q,t)- (1-q)(1-qt^2)PH(\mathcal{M}^{2,1},-q,t)}{(1-q^2)(1-q^2t^4)}\right )
\end{align*}
so
\begin{equation*}
V(q,t)=\frac{1}{2} \left ((1+q)(1+qt^2) PH(\mathcal{M}^{2,1},q,t)- (1-q)(1-qt^2)PH(\mathcal{M}^{2,1},-q,t) \right ).
\end{equation*}
\item[(ii)] We have:
\begin{equation*}
V(1,t)=\frac{1}{2} \left (2(1+t^2) PH(\mathcal{M}^{2,1},1,t) \right )=(1+t^2)P(\mathcal{M}^{2,1},t).
\end{equation*}
\item[(iii)] From \cite[Corollary 1.1.4]{hausel2008mixed} we see that
\begin{align*}
(qt)^{8g-6}PH(\mathcal{M}^{2,1},q^{-1}t^{-2},t)=PH(\mathcal{M}^{2,1},q,t).
\end{align*}
Therefore:
\begin{align*}
(qt)^{8g-4}& V(q^{-1}t^{-2},t)=\frac{1}{2}  ((qt)^{8g-6}(1+q)(1+qt^2) PH(\mathcal{M}^{2,1},q^{-1}t^{-2},t)+\\
&- (qt)^{8g-6}(1-q)(1-qt^2)PH(\mathcal{M}^{2,1},-q^{-1}t^{-2},t) )=\\
&=\frac{1}{2} \left ((1+q)(1+qt^2) PH(\mathcal{M}^{2,1},q,t)- (1-q)(1-qt^2)PH(\mathcal{M}^{2,1},-q,t) \right )=\\
&= V(q,t).
\end{align*}
Clearly the formula for $G$ follows from the formula for $V$.
\end{itemize}
\end{proof}
\end{prop}
%%
%\begin{prop}
%We have:
%\begin{align*}
%F-G&=\frac{1}{4} q^{2g-2} t^{4g-4} \left ( {\frac {( 1+qt )^{2g}}{( 1-qt^2 )( 1-q ) }}- \frac{( 1-qt )^{2g}}{( 1+qt^2)( 1+q ) }  \right ) \left ( {\frac {( 1+qt )^{2g}}{( 1-qt^2 )( 1-q ) }}+ \frac{( 1-qt )^{2g}}{( 1+qt^2)( 1+q ) }  \right )+\\
%&+ \frac{q^{2g-1} t^{2g} (1+t)^{2g}}{(1-t^2)(1-q^2 t^2)} \left ( \frac{(1+q^2 t^3)^{2g}}{(1-q^2 t^4)^2} - \frac{t^{2g-4} (1+q^2 t)^{2g}}{(1-q^2)^2}\right ).
%\end{align*}
%\begin{proof}
%Computation.
%\end{proof}
%\end{prop}
%
\begin{rmk}
Even though for rank 2 it is possible to compute the difference explicitly, it is quite interesting to understand an a priori reason why the two functions should have such a relation. The idea is to compare $F^{sh}$ with $G^{sh}$ and $P(F^{sh}(q),t)$ with $G(q,t)$. $P(F^{sh}(q),t)$ will be the generating function of the Poincar\'e polynomials of the sheaves contained in $F^{sh}$, which are the intersection cohomology sheaves of the $\mathcal{M}_\varepsilon^{2,2n+1}$. Note that, since for $2n+1<0$ and for $2n+1>4g-4$ $\mathcal{M}_\varepsilon^{2,2n+1}$ is smooth, then $F^{vir}$ and $P(F^{sh}(q),t)$ coincide for $\deg q \leq 2g-3$ and for $\deg q \geq 6g-5$.
\end{rmk}

In the next few pages we want to prove that $G$ and $P(F^{sh}(q),t)$ coincide for $\deg q \leq 2g-3$ and for $\deg q \geq 6g-5$. Let us first recall the following result. See \cite[Theorem 3]{chaureport}.
\begin{prop}
\label{dimsupp}
Let $f: X \rightarrow Y$ be a proper map of smooth algebraic varieties whose fibers have constant dimension. Denote by $m$ the dimension of $Y$ and by $n$ the dimension of $X$. Let $\mathscr{S}$ be summand of $f_* \Q_X$ whose support has dimension $s$. Then $s \geq 2m-n$.
\end{prop}
\begin{theorem}
\label{BHMSMY}
$F^{sh}$ and $G^{sh}$ coincide for $\deg q \leq 2g-3$. $P(F^{sh}(q),t)$ and $G$ coincide for $\deg q \leq 2g-3$ and for $\deg q \geq 6g-5$. Furthermore $P(F^{sh}(q),t)-G(q,t)$ is a polynomial with non-negative coefficients.
%Let $F, G \in \Z[t][[q]]$ be defined as above. Then:
%\begin{itemize}
%\item[(i)] $F-G$ is a polynomial
%\item[(ii)] $F-G$ is divisible by $q^{2g-1}$ and has degree at most $6g-7$
%\item[(iii)] $F-G$ has positive coefficients.
%\end{itemize}
\begin{proof}
Let us start by considering the Hitchin map $h: \mathcal{M}^{2,1} \rightarrow \mathcal{A}^2$ and the object $\R h^{2,1}_* \Q \in D_c^b(\mathcal{A}^2)$. Consider also $\mathcal{A}^2_{sm} \subset \mathcal{A}^2$, the locus of smooth spectral curves and call $\R^1$ the local system $\R^1 (h_{sm})_* \Q$ where $h_{sm}$ is the restriction of $h$ to the inverse image of $\mathcal{A}^2_{sm}$.

As we noted in remark \ref{decthm}, we have:
\begin{equation*}
\R h^{2,1}_* \Q = \bigoplus_{i=0}^{8g-6} \leftidx{^p}{\R}{^i} h^{2,1}_* \Q[-i] \oplus \mathscr{L} = \bigoplus_{i=0}^{8g-6} IC\left (\bigwedge^i\R^1 \right) \oplus \mathscr{L} 
\end{equation*}

First of all we prove that:
\begin{equation*}
\bigoplus_{2-2g \leq d <0} q^{d+2g-2} \R (\chi_\varepsilon^{2,d})_* \Q = \left (\frac{\bigoplus_{i=0}^{8g-6} \leftidx{^p}{\R}{^i}  h^{2,1}_* \Q[-i] q^i}{(1-q \Q)(1-q \L)} \right )_{|\deg(q) \leq 2g-3}.
\end{equation*}
which implies that $F^{sh}$ and $G^{sh}$ coincide for $\deg q \leq 2g-3$ and is a bit stronger because it also says something about the even degree case. Recall from proposition \ref{negdegrk2} that $\mathcal{M}_{\varepsilon}^{2,d}=\mathcal{M}_{\infty}^{2,d}$ for $d < 0$.

When we restrict to the locus $\mathcal{A}^2_{int}$ of integral spectral curves, the formula is true since the restriction of $\chi_\sigma^{2,d}$ to $\mathcal{A}^2_{int}$ is always the Hilbert scheme of $d+2g-2$ points relative to the family of integral spectral curves, independently of $\sigma$. The formula then descends from the results in \cite{migliorini2011support}. This means that if the left hand side and the right hand side differ then, since the $\leftidx{^p}{\R}{^i}  h^{2,1}_* \Q[-i] $ are supported on the whole $\mathcal{A}^2$, there must be summands of the left hand side that are supported inside $\mathcal{A}^2\setminus \mathcal{A}^2_{int}$.

For the second step, let us restrict the formula to the locus $\mathcal{A}^2_{red} \subset \mathcal{A}^2$, i.e. the locus of reduced (but possibly reducible) spectral curves. Note that we cannot immediately conclude that the formula holds using the results from \cite{migliorini2015support} because the family of reduced spectral curves does not satisfy the conditions in the paper. However we can argue as follows. Assume that there is a summand with support contained in the closure of $\mathcal{A}^2_{red}$ and denote by $s$ its dimension. Using proposition \ref{dimsupp}, we get:
\begin{equation*}
s \geq 8g-6 - 4g +3 - d -2g +2=2g-1-d.
\end{equation*}
and for $d \leq -1$ we get $s \geq 2g$. Since $\mathcal{A}^2_{red}$ is irreducible, this implies that the support is actually the closure of $\mathcal{A}^2_{red}$.

Note that the generic reduced and reducible curve is nodal, since $\mathcal{A}^2_{red}$ is the image of the following map:
\begin{align*}
\Sym^2 \left ( H^0(K) \right )\setminus \Delta &\rightarrow \mathcal{A}^2\\
(\psi_1, \psi_2) & \mapsto (-\psi_1-\psi_2, \psi_1 \psi_2)
\end{align*}
and two generic sections of $H^0(K)$ intersect in $2g-2$ points.

Using \cite[Lemma 2.13]{migliorini2015support}, we can deduce that if the stalks of
\begin{equation*}
\bigoplus_{2-2g \leq d <0} q^{d+2g-2} \R (\chi_\varepsilon^{2,d})_* \Q
\end{equation*}
and
\begin{equation*}
 \left (\frac{\bigoplus_{i=0}^{8g-6} \leftidx{^p}{\R}{^i}  h^{2,1}_* \Q[-i] q^i}{(1-q \Q)(1-q \L)} \right )_{|\deg(q) \leq 2g-3}
\end{equation*}
were to be different at a generic nodal curve of $\mathcal{A}^2_{red}$ then their weight polynomial would be different. Since we know that the potential support could only be the closure of $\mathcal{A}^2_{red}$, if we verify that the weight polynomials of the stalks at the generic nodal curves are the same, we can conclude that there is no such support. More details about the strategy outlined here can be found in \cite{migliorini2015support}.

Let $\mathcal{C}$ be a generic nodal curve in $\mathcal{A}^2_{red}$ and let $\mathcal{D}$ be an \emph{integral} nodal curve in $\mathcal{A}^2_{int}$. Define:
\begin{equation*}
Z(\mathcal{C})=\sum_{n\geq 0} q^n [\mathcal{C}^{[n]}]
\end{equation*}
the zeta function of the weight polynomials of the Hilbert schemes of points of $\mathcal{C}$, and $Z(\mathcal{D})$ analogously. Recall that if we have a disjoint union of varieties $X \sqcup Y$, then $Z(X \sqcup Y)=Z(X)\cdot Z(Y)$. Define:
\begin{equation*}
Z(N)=\sum_{n\geq 0} q^n [\text{Hilb}_{pt}^n (\text{node})]
\end{equation*}
to be the zeta function of the weight polynomials of the punctual Hilbert schemes of points of the curve $\{xy=0\} \subset \C^2$ supported at the origin. Since both $\mathcal{C}$ and $\mathcal{D}$ can be decomposed as the union of their smooth part (here denoted with a subscript \emph{sm}) and their nodes, and they both have $2g-2$ nodes, we get:
\begin{equation*}
\frac{Z(\mathcal{C})}{Z(\mathcal{D})}=\frac{Z(\mathcal{C}_{sm})\cdot Z(N)^{2g-2}}{Z(\mathcal{D}_{sm}) \cdot Z(N)^{2g-2}}=\frac{Z(\mathcal{C}_{sm})}{Z(\mathcal{D}_{sm})}.
\end{equation*}

If we denote by $\mathcal{C}^\nu$ and $\mathcal{D}^\nu$ the normalization of the two curves, then we get:
\begin{equation*}
Z(\mathcal{C}^\nu)=Z(\mathcal{C}_{sm})\cdot Z(pt)^{4g-4}
\end{equation*}
and analogously for $\mathcal{D}^\nu$. Putting all together and using the fact that $\mathcal{C}^\nu$ is the disjoint union of two copies of $C$ we get:
\begin{equation*}
\frac{Z(\mathcal{C})}{Z(\mathcal{D})}=\frac{Z(\mathcal{C}_{sm})}{Z(\mathcal{D}_{sm})}=\frac{Z(C)^2}{Z(\mathcal{D}^\nu)}.
\end{equation*}

In order to compute the (weight polynomial of the) stalk of $IC\left (\bigwedge^i\R^1 \right)$ at the generic nodal curve in $\mathcal{A}^2_{red}$ we can use the content of section \ref{cks}. We can also avoid the direct computation of the full weight polynomial by using the following idea.

Whenever we write the CKS complexes for curves whose irreducible components have positive genus, we will have $H^1$ of the normalization of the curve appearing only in degree $0$ and all the elements of the external powers of the degree $0$ term in which appear some classes of the $H^1$ of the normalization lie in the kernel of the first map. In our case this implies:
\begin{equation*}
\left( \sum_{i \geq 0} \left [IC(\bigwedge^i\R^1 )_{|\mathcal{C}}[-i] \right] q^i\right )= U(\Sigma) \left( \sum_{i \geq 0} \left [ \bigwedge^i H^1 (C)[-i] \right] q^i\right )^2\qquad \text{mod }q^{2g-2}
\end{equation*}
and
\begin{equation*}
\left( \sum_{i \geq 0} \left [IC(\bigwedge^i\R^1 )_{|\mathcal{D}}[-i] \right] q^i\right )= U(\overline\Sigma) \left( \sum_{i \geq 0} \left [ \bigwedge^i H^1 (\mathcal{D}^\nu)[-i] \right] q^i\right )\qquad \text{mod }q^{2g-2}
\end{equation*}
where $U(\Sigma)$ and $U(\overline\Sigma)$ are defined and examined in \ref{ratcurve}.

Putting all together we have, mod $q^{2g-2}$:
\begin{align*}
Z(\mathcal{D}^\nu)&\left( \sum_{i \geq 0} \left [IC(\bigwedge^i\R^1 )_{|\mathcal{C}}[-i] \right] q^i\right )=U(\Sigma)Z(C)^2Z(\mathcal{D}^\nu)(1-q\Q)^2(1-q \L)^2=\\
&=U(\overline\Sigma)Z(C)^2Z(\mathcal{D}^\nu)(1-q\Q)(1-q \L)=Z(C)^2 \left( \sum_{i \geq 0} \left [IC(\bigwedge^i\R^1 )_{|\mathcal{D}}[-i] \right] q^i\right )=\\
&=Z(C)^2 Z(\mathcal{D})(1-q\Q)(1-q \L)=Z(\mathcal{C})Z(\mathcal{D}^\nu)(1-q\Q)(1-q \L)
\end{align*}
and therefore:
\begin{equation*}
\left( \sum_{i \geq 0} \left [IC(\bigwedge^i\R^1 )_{|\mathcal{C}}[-i] \right] q^i\right )= Z(\mathcal{C})(1-q\Q)(1-q \L) \qquad \text{mod } q^{2g-2}.
\end{equation*}

%Last, we need to see whether or not there can be summands of the left hand side whose support is contained in $\mathcal{A} \setminus \mathcal{A}_{red}$. Let $s$ be the dimension of such a support. Since we assume it is contained in $\mathcal{A} \setminus \mathcal{A}_{red}$ we clearly have $s \leq g$ which is the dimension of the locus of non-reduced curves. On the other hand using proposition \ref{dimsupp}, we get:
%\begin{equation*}
%s \geq 8g-6 - 4g +3 - d -2g +2=2g-1-d.
%\end{equation*}
%Therefore we would have $g \geq s \geq 2g-1-d$, i.e. $d \geq g-1$. This inequality is impossible in the case when $2-2g \leq d < 0$ and so no such summands can exist here. Note that we actually get no such supports for a slightly bigger range than just $2-2g \leq d <0$.
To summarize, we wanted to prove that weight polynomials of the stalks at $\mathcal{C}$ of
\begin{equation*}
\bigoplus_{2-2g \leq d <0} q^{d+2g-2} \R (\chi_\varepsilon^{2,d})_* \Q
\end{equation*}
and
\begin{equation*}
 \left (\frac{\bigoplus_{i=0}^{8g-6} \leftidx{^p}{\R}{^i}  h^{2,1}_* \Q[-i] q^i}{(1-q \Q)(1-q \L)} \right )_{|\deg(q) \leq 2g-3}
\end{equation*}
are the same. In order to do that we related both polynomials to the weight polynomials of the stalks at $\mathcal{D}$ that we know coincide because $\mathcal{D}$ lies in $\mathcal{A}^2_{int}$ where we know the equality holds.

Let us now prove that $P(F^{sh}(q),t)$ and $G$ coincide for high enough degree. Recall that for odd $d \geq 6g-5$ the map $\mathcal{M}_\varepsilon^{2,d} \rightarrow \mathcal{M}^{2,d}$ is a $\C \P^{d+1-2g}$-bundle. This implies:
\begin{equation*}
\R (\chi_\varepsilon^{2,d})_* \Q = \R h^{2,1}_* \Q \otimes H^*(\C \P^{d+1-2g})
\end{equation*}
here $H^*(\C \P^{d+1-2g})$ denotes the constant local system $\Q \oplus \Q[-2] \oplus \dots \oplus \Q[-2(d+1-2g)]$.

Let us compare the coefficients of the left hand side
\begin{equation*}
\bigoplus_{n \geq 1-g} q^{2n+1+2g-2} \R (\chi_\varepsilon^{2,2n+1})_* IC_{\mathcal{M}_\varepsilon^{2,2n+1}}
\end{equation*}
and of the right hand side
\begin{equation*}
\text{odd}_q \left (\frac{\bigoplus_{i=0}^{8g-6} \leftidx{^p}{\R}{^i}  h^{2,1}_* \Q[-i] q^i}{(1-q \Q)(1-q \L)} \right ).
\end{equation*}
for the powers $q^{d+2g-2}$ with odd $d \geq 6g-5$.

The coefficient of $q^{d+2g-2}$ in the first expression is, as we said:
\begin{align*}
\R h^{2,1}_* \Q \otimes H^*(\C \P^{d+1-2g}) &= \left ( \bigoplus_{i=0}^{8g-6} IC\left (\bigwedge^i\R^1 \right)[-i] \oplus \mathscr{L} \right )\otimes \\
&\otimes \left ( \Q \oplus \Q[-2] \oplus \dots \oplus \Q[-2(d+1-2g)]\right ).
\end{align*}

For the coefficient of $q^{d+2g-2}$ in the second expression we have instead:
\begin{align*}
\sum_{i=0}^{8g-6}& IC\left (\bigwedge^i\R^1 \right)[-i] \otimes \left ( \Q \oplus \dots \Q[-2(d+2g-2-i)] \right)=\\
&=\sum_{i=0}^{4g-3} IC\left (\bigwedge^i\R^1 \right)[-i] \otimes \left ( \Q \oplus \dots \Q[-2(d+1-2g)] \right)+\\
&+\sum_{i=0}^{4g-4} IC\left (\bigwedge^i\R^1 \right)[-i] \otimes \left ( \Q[-2(d+2-2g)] \oplus \dots \Q[-2(d+2g-2-i)] \right)+\\
&+\sum_{j=4g-2}^{8g-6} IC\left (\bigwedge^j\R^1 \right)[-j] \otimes \left ( \Q \oplus \dots \Q[-2(d+2g-2-j)] \right).
\end{align*}

Using the fact that (see \cite[equation 1.4.7]{de2010topology}):
\begin{equation*}
IC\left (\bigwedge^j\R^1 \right) \cong IC\left (\bigwedge^{8g-6-j}\R^1 \right)
\end{equation*}
and the substitution $j=8g-6-i$, we can compute:
\begin{align*}
&\sum_{i=0}^{4g-4} IC\left (\bigwedge^i\R^1 \right)[-i] \otimes \left ( \Q[-2(d+2-2g)] \oplus \dots \Q[-2(d+2g-2-i)] \right)=\\
&=\sum_{j=4g-2}^{8g-6} IC\left (\bigwedge^j\R^1 \right)[6+j-8g] \otimes \left ( \Q[-2(d+2-2g)] \oplus \dots \Q[-2(d-6g+4+j)] \right)=\\
&=\sum_{j=4g-2}^{8g-6} IC\left (\bigwedge^j\R^1 \right)[-j] \otimes \left ( \Q[-2(d+2g-1-j)] \oplus \dots \Q[-2(d+1-2g)] \right).
\end{align*}

Putting all together we get:
\begin{align*}
\sum_{i=0}^{8g-6}& IC\left (\bigwedge^i\R^1 \right)[-i] \otimes \left ( \Q \oplus \dots \Q[-2(d+2g-2-i)] \right)=\\
&=\sum_{i=0}^{4g-3} IC\left (\bigwedge^i\R^1 \right)[-i] \otimes \left ( \Q \oplus \dots \Q[-2(d+1-2g)] \right)+\\
&+\sum_{j=4g-2}^{8g-6} IC\left (\bigwedge^j\R^1 \right)[-j] \otimes \left ( \Q[-2(d+2g-1-j)] \oplus \dots \Q[-2(d+1-2g)] \right)+\\
&+\sum_{j=4g-2}^{8g-6} IC\left (\bigwedge^j\R^1 \right)[-j] \otimes \left ( \Q \oplus \dots \Q[-2(d+2g-2-j)] \right)=\\
&=\sum_{i=0}^{4g-3} IC\left (\bigwedge^i\R^1 \right)[-i] \otimes \left ( \Q \oplus \dots \Q[-2(d+1-2g)] \right)+\\
&+\sum_{j=4g-2}^{8g-6} IC\left (\bigwedge^j\R^1 \right)[-j] \otimes \left ( \Q \oplus \dots \Q[-2(d+1-2g)] \right)=\\
&=\sum_{i=0}^{8g-6} IC\left (\bigwedge^i\R^1 \right)[-i] \otimes \left ( \Q \oplus \dots \Q[-2(d+1-2g)] \right).
\end{align*}

Note that this expression only differs from the first by the presence of the factor:
\begin{equation*}
\mathscr{L} \otimes H^*(\C \P^{d+1-2g})
\end{equation*}
which, however, does not have any global cohomology and therefore this proves the equality of the coefficients of $P(F^{sh}(q),t)$ and $G$ for the powers $q^{d+2g-2}$ with odd $d\geq 6g-5$.

Recall also that $F^{vir}-G$ satisfies:
\begin{equation*}
F^{vir}(q,t)-G(q,t)=(qt)^{8g-8}\left (F^{vir}(q^{-1}t^{-2},t)-G(q^{-1}t^{-2},t) \right).
\end{equation*}

In particular, if $F^{vir}$ and $G$ agree for $\deg q$ odd in the range $ 1, \dots, 2g-3$, then they will agree for $\deg q$ odd in the range $6g-5, \dots, 8g-9$. Since $F^{vir}$ and $P(F^{sh}(q),t)$ agree in the previous ranges we can conclude that $P(F^{sh}(q),t)$ and $G$ agree for $\deg q \leq 2g-3$ and $\deg q \geq 6g-5$ as we wanted.

For the last statement about non-negativity recall that $P(F^{sh}(q),t)$ is the global cohomology of:
\begin{equation*}
\bigoplus_{n \geq 1-g} q^{2n+1+2g-2} \R (\chi_\varepsilon^{2,2n+1})_* IC_{\mathcal{M}_\varepsilon^{2,2n+1}}
\end{equation*}
and $G$ is the global cohomology of $G^{sh}$.

Since these two expressions coincide when restricted to $\mathcal{A}^2_{int}$ and the second one only has summands that are supported on the full $\mathcal{A}^2$, we can conclude that $P(F^{sh}(q),t)-G$ is the global cohomology of the summands of $F^{sh}(q)$ that have proper support in $\mathcal{A}^2$ and do not come from $\mathscr{L}$ (if there are any).
\end{proof}
\end{theorem}

%% file: Higher.tex
\chapter{Remarks on higher rank}
%
%
%
%\begin{itemize}
%\item Serre duality
%\item range for projective bundle (is Nitsure correct?)
%\item range for initial agreement with relative Hilbert scheme
%\item generalization of the formulas
%\end{itemize}
Some of the results of the previous section can be generalized for all ranks.
In this chapter we will comment about the behavior of Bradlow-Higgs pairs for higher rank. Namely we will point out several difficulties that arise when trying to compute the motive of the $\mathcal{M}^{r,d}_{\sigma}$, both with similar methods to those we used in this thesis or with methods available for other moduli problems with wall-crossing.

We will also comment on how it could be possible to generalize the results of chapter \ref{MSMY}.

\section{Higher rank motives}
In the case of Bradlow pairs it is still possible, for lower rank, to understand the flip loci explicitly. For rank 3 this has been done in \cite{munoz2008hodge}. It becomes clear soon, though, that a direct approach would not work for Bradlow pairs and so a different approach to the wall-crossing problem was used in \cite{mozgovoy2013moduli}. The case of Bradlow-Higgs triples is more complicated and arguably a direct approach will fail for this problem as well.

It is worth commenting about the approach in \cite{mozgovoy2013moduli} and explain why it will not apply to the case of Bradlow-Higgs triples. We will try to highlight the main points of the strategy, but a fully detailed survey would require a separate thesis, therefore we refer directly to \cite{mozgovoy2013moduli} for all the details.

First of all denote by $\mathcal{A}_0$ the category of coherent sheaves on the smooth projective curve $C$ and by $\mathcal{A}$ the category of triples $(E_0,E_1,s)$ where $E_1 \in \mathcal{A}_0$, $E_0$ is a direct sum of a certain number of copies of $\mathcal{O}_C$ and $s: E_0 \rightarrow E_1$. Both categories admit \emph{Chern characters}, namely group morphisms from the Grothendieck group $K(\mathcal{A}_0)$ to $\Z^2$ and from $K(\mathcal{A})$ to $\Z^3$ sending a coherent sheaf to its rank and degree and a triple $(E_0,E_1,s)$ to $(\rk E_1, \deg E_1, \rk E_0)$ respectively. There are also skew symmetric bilinear forms $\langle \cdot, \cdot \rangle$ on $\Z^2$ and $\Z^3$ defined in terms of simple formulas on the Chern characters (see \cite[section 4]{mozgovoy2013moduli}).

For each of $\mathcal{A}_0$ and $\mathcal{A}$ we define the respective \emph{motivic Hall algebras} which are $\C$-algebras generated by indicator functions of locally closed substacks of the stacks of objects of $\mathcal{A}_0$ and $\mathcal{A}$ respectively. We denote these two algebras as $H(\mathcal{A}_0)$ and $H(\mathcal{A})$ (see \cite[section 5]{mozgovoy2013moduli}). Note that we have an embedding $\mathcal{A}_0 \subset \mathcal{A}$ and accordingly an inclusion of $\C$-algebras. The product of two indicators functions $\id_X$ and $\id_Y$ is the indicator function $\id_Z$ of objects that have a subobject in $Y$ and the quotient by such a subobject lies in $X$.

There is an integration map $I: H(\mathcal{A}) \rightarrow \mathbb{A}$ where $\mathbb{A}$ is a certain completion of the ring $R[x_1,x_2^{\pm 1},x_3]$ and $R$ is the Grothendieck ring of stacks over $\C$. The product on $\mathbb{A}$ is defined as
\begin{equation*}
x^{\alpha} x^{\beta} = (-\L^{1/2})^{\langle \alpha, \beta \rangle} x^{\alpha+\beta}.
\end{equation*}

The map $I$ is defined on the set of indicator functions $\id_X$ where $X$ is a locally closed substack of the stack of objects of $\mathcal{A}$ for which the Chern character is constantly equal to $\alpha \in \Z^3$. We have
\begin{equation*}
I(\id_X)=(-\L^{1/2})^{\langle \alpha, \alpha \rangle} x^{\alpha} [X]
\end{equation*}
and then it is extended by linearity on the whole Hall algebra. It is clearly a linear map and under the hypothesis that $\Ext^2$ groups vanish between the objects in the indicator functions $\id_X$ and $\id_Y$ then $I(\id_X \id_Y)=I(\id_X)I(\id_Y)$ (see \cite[remark 5.1]{mozgovoy2013moduli}).

Define some distinguished elements in $H^0(\mathcal{A})$. Let $u^h(\alpha)$ be the indicator function of semistable vector bundles with character $(\alpha,0)$ and $f^h_\tau(\alpha)$ be the indicator function of $\tau$-semistable triples $(E_0,E_1,s)$ with character $(\alpha,1)$. Let also $u(\alpha)$ and $f_\tau(\alpha)$ be the integration of $u^h(\alpha)$ and $f^h_\tau(\alpha)$ respectively (both up to a factor). Observe that they are essentially the motives of the moduli spaces that we would like to compute. Define
\begin{equation*}
u^h_\tau=1+\sum_{\mu(\alpha)=\tau} u^h(\alpha) \quad f^h_\tau=\sum_\alpha f^h_\tau(\alpha)
\end{equation*}
in the Hall algebra and
\begin{equation*}
u_\tau=1+\sum_{\mu(\alpha)=\tau} u(\alpha) x^\alpha \quad f_\tau=\sum_\alpha f_\tau(\alpha) x^{(\alpha,1)}
\end{equation*}
to be the corresponding generating series in $\mathbb{A}$.

Using the geometry of Bradlow pairs it is proven (see \cite[lemmas 4.11 and 4.13]{mozgovoy2013moduli}) that every $\tau$-semistable triple $(E_0,E_1,s)$ has a canonical $\tau+$ filtration whose quotient is a semistable vector bundle of slope $\tau$ and the subobject is a $\tau+$-stable triple. Also, every $\tau$-semistable triple $(E_0,E_1,s)$ has a canonical $\tau-$ filtration whose subobject is a semistable vector bundle of slope $\tau$ and the subobject is a $\tau-$-stable triple. This implies the following relations in the motivic hall algebra:
\begin{equation*}
f^h_{\tau}=f^h_{\tau+} u_\tau^h \text{ and } f^h_{\tau| \mu<\tau}= u_\tau^h f^h_{\tau-}
\end{equation*}
where $|\mu < \tau$ denotes the truncation of the series to the term for which the slope is less than $\tau$ (see \cite[theorem 5.6]{mozgovoy2013moduli}).

A very important observation at this point is that for $\tau = \infty$ the only $\tau$-semistable triples are those for which $\rk E_1=1$. These correspond to rank 1 Bradlow pairs that in turn can be identified with divisors on the curve. Therefore
\begin{equation*}
f_\infty =x_1 x_3 \sum_{d \geq 0} [S^d( C) ] x_2^d.
\end{equation*}

Also, we have $f^h_\tau(\alpha)=0$ if $\mu(\alpha)>\tau$ which geometrically amounts to say that the moduli spaces of Bradlow pairs are empty if the stability parameter exceeds a prescribed threshold.

Combining all these results it is possible to find a formula (see \cite[theorem 5.6]{mozgovoy2013moduli}):
\begin{equation*}
f_\tau = (u^{-1}_{>\tau} f_\infty u_{\geq \tau})_{|\mu \leq \tau}
\end{equation*}
where
\begin{equation*}
u_{\geq \tau}= \prod_{\tau' \geq \tau} u_{\tau'}
\end{equation*}
and the product is taken in decreasing order of $\tau'$. The formula can then be inverted using the fact that we know the generating functions $u_\tau$. At the end the motive of the moduli spaces of Bradlow pairs is computed in \cite[theorem 6.2]{mozgovoy2013moduli}.

The same approach for Bradlow-Higgs triples will not work as some essential hypotheses are not satisfied. The main issue is that the use of the integration map is conditional to the fact that for the objects we need in the wall-crossing of Bradlow pairs, the $\Ext^2$ will all vanish. This is far from being true in the case of Bradlow-Higgs triples and in lemma \ref{exts} we saw an example of the presence of nonzero $\Ext^2$. In terms of the flip loci this will imply that they are not bundles in general. A positive fact is that for Higgs bundles, $\Ext^2$ is dual to $\Hom$ by Serre's duality.

The second issue is that the moduli spaces of $\sigma$-stable Bradlow-Higgs triples are not empty after a certain value of $\sigma$. This will imply that the analogue of the generating function $f_\infty$ will not be as simple as in the case of Bradlow pairs. Note also that we gave an approach to compute $\mathcal{M}_\infty^{2,d}$ directly but it is not at all trivial already in rank 2.

In conclusion, the analogues of the formulas
\begin{equation*}
f^h_{\tau}=f^h_{\tau+} u_\tau^h \text{ and } f^h_{\tau| \mu<\tau}= u_\tau^h f^h_{\tau-}
\end{equation*}
in the motivic Hall algebra are probably still valid but, even if they can be proved, it will then not be possible to apply an easy integration map to get information about the motives.
\section{Partial Hilbert scheme formula for higher rank}
In this section we will discuss how to possibly generalize the result in theorem \ref{BHMSMY} to higher rank. Consider the following higher rank generating functions.
\begin{defn}
Let $r \geq 2$ be an integer and $1 \leq m \leq r-1$ be an integer coprime with $r$.
\begin{equation*}
F_{r,m}^{sh}(q)= \sum_{\substack{d \geq r(r-1)(1-g)\\ d \text{ mod } r= m}} \R (\chi_\varepsilon^{r,d})_* (IC_{\mathcal{M}_\varepsilon^{r,d}}) q^{d+r(r-1)(g-1)},
\end{equation*}
\begin{equation*}
F_{r,m}^{mot}(q)= \sum_{\substack{d \geq r(r-1)(1-g)\\ d \text{ mod } r= m}} [\mathcal{M}_\varepsilon^{r,d}] q^{d+r(r-1)(g-1)},
\end{equation*}
\begin{equation*}
F_{r,m}^{vir}(q,t)=\sum_{\substack{d \geq r(r-1)(1-g)\\ d \text{ mod } r= m}} P^{vir}(\mathcal{M}_\varepsilon^{2,2n+1},t) q^{d+r(r-1)(g-1)}.
\end{equation*}
\end{defn}

These are clearly the generating functions in theorem \ref{BHMSMY} extended for higher rank. Let us introduce the following notation.
\begin{defn}
Let $r \geq 2$ be an integer and $0 \leq m \leq r-1$ be an integer. We define $\text{del}_{q,r,m}$ be the operator acting on $R[[q]]$ by deleting from power series all terms whose degree $d$ does not have remainder $m$ modulo $r$.
\end{defn}

Let $h^{r,d}:\mathcal{M}^{r,d} \rightarrow \mathcal{A}^r$ denote the usual Hitchin map and $PH(\mathcal{M}^{r,d},q,t)$ denote the perverse hodge polynomial of $\mathcal{M}^{r,d}$. Then we define the second class of generating functions we need.
\begin{defn}
Let $r \geq 2$ be an integer and $1 \leq m \leq r-1$ be an integer coprime with $r$. Define:
\begin{equation*}
G_{r,m}(q,t)=\text{del}_{q,r,m} \left (\frac{PH(\mathcal{M}^{r,m},q,t)}{(1-q)(1-qt^2)} \right )
\end{equation*}
and
\begin{equation*}
G_{r,m}^{sh}(q)=\text{del}_{q,r,m} \left (\frac{ \R h^{r,m}_* \Q }{(1-q\Q)(1-q\Q[-2](-1))} \right )
\end{equation*}
\end{defn}

In this section we will prove a variant of theorem \ref{BHMSMY} for higher rank.
\begin{theorem}
Let $r \geq 2$ be an integer and $1 \leq m \leq r-1$ be an integer coprime with $r$. Then $F_{r,m}^{sh}$ and $G_{r,m}^{sh}$ coincide for $\deg q \leq (r-1)(2g-2)-1$. $P(F_{r,m}^{sh}(q),t)$ and $G_{r,m}$ coincide for $\deg q \leq (r-1)(2g-2)-1$ and for $\deg q > r(2g-1)+(3r^2-5r+2)(g-1)$.

If the following conditions are satisfied:
\begin{itemize}
\item for $d > r(r+1)(g-1)$ then a stable Higgs bundle $(E,\phi)$ of rank $r$ and degree $d$ satisfies $H^1(E)=0$
\item for $d>r^2(g-1)$ coprime with $r$, $\mathcal{M}_\varepsilon^{r,d}$ is smooth
\end{itemize}

then $P(F_{r,m}^{sh}(q),t)$ and $G_{r,m}$ coincide for $\deg q \leq (r-1)(2g-2)-1$ and for $\deg q > (2 r^2-r)(2g-2)$.
\end{theorem}

Recall that the strategy of the proof for rank 2 was to compare the cohomology of $\mathcal{M}_\varepsilon^{2,d}$ for $d < 0$ and for $d > 4g-4$ odd to the perverse filtration on the cohomology of $\mathcal{M}^{2,1}$. Let us summarize what was needed in order to prove the result.

First we saw in proposition \ref{negdegrk2} that $\mathcal{M}_\infty^{2,d}=\mathcal{M}_\varepsilon^{2,d}$ if $2-2g \leq d <0$ and that in this case they are both smooth. To prove theorem \ref{BHMSMY} in this range we then used proposition \ref{dimsupp} to get information about the possible supports in the decomposition theorem and concluded by studying the case of a generic nodal curve in $\mathcal{A}^2_{red}$.

For the range $d > 4g-4$ odd, we used again smoothness of $\mathcal{M}_\varepsilon^{2,d}$ and the fact that for $d \geq 6g-5$ odd the Abel-Jacobi map $\mathcal{M}_\varepsilon^{2,d}\rightarrow \mathcal{M}^{2,d}$ is a projective bundle (see proposition \ref{sharpnits}). In the range $d \geq 6g-5$ odd theorem \ref{BHMSMY} was then easy to prove. For the range $4g-3 \leq d \leq 6g-7$ odd, we then used smoothness of $\mathcal{M}_\varepsilon^{2,d}$ plus Serre's duality to conclude that if the cohomological formula holds for odd values in the range $3-2g, \dots , -1$ then it must also hold for odd values in the range $4g-3, \dots , 6g-7$.

We also used the fact that the extra summand $\mathscr{L}$ in the decomposition theorem for the Hitchin map $\mathcal{M}^{2,1} \rightarrow \mathcal{A}^2$ has no cohomology. A closer look to the argument though will show that we don't need this to show the formula in the ranges $d \geq 6g-5$ and $d <0$.

It is worth discussing which of the previous properties can be generalized to the case of higher rank. Let us first give a brief remark.
\begin{rmk}
\label{cflip}
Note that since $\sigma$-stability is an open condition, the flip loci are closed in each of the $\mathcal{M}_\sigma^{r,d}$. In particular, due to the fact that the Hitchin maps are proper, the restriction of the Hitchin maps to the flip loci is proper as well. It follows that if $(E,\phi,s)$ belongs to a flip locus, then
$$\lim_{\lambda \rightarrow 0} \lambda \cdot (E,\phi,s)$$
also belongs to the same flip locus.
\end{rmk}
We can generalize proposition \ref{negdegrk2} as follows.
\begin{prop}
Let $d < r(r-2)(1-g)$, then $\mathcal{M}_\varepsilon^{r,d} = \mathcal{M}_\infty^{r,d}$. Both are non-empty if and only if $d \geq r(r-1)(1-g)$.
\begin{proof}
The statement about non-emptiness follows from the $U$-filtration for $\mathcal{M}_\infty^{r,d}$ and for $\mathcal{M}_\varepsilon^{r,d}$ follows from the previous statement.

Let us first prove that if $d < r(r-2)(1-g)$ then $\mathcal{M}_\varepsilon^{r,d} \subseteq \mathcal{M}_\infty^{r,d}$. Suppose $(E,\phi,s)$ has stable underlying Higgs bundle and assume by contradiction that the $U$-filtration of the triple has length $l < r$. But then $\mu(U_l) \geq (l-1)(1-g )>(r-2)(1-g)> \mu(E)$ and therefore $U_l$ cannot be preserved by $\phi$ which, by construction, is impossible.

To conclude we use remark \ref{cflip}. From it we can deduce that it is enough to prove that if $(E,\phi,s) \in \mathcal{M}_\infty^{r,d}$ is a fixed point for the $\C^*$-action, then $(E,\phi,s) \in  \mathcal{M}_\varepsilon^{r,d}$. Recall from proposition \ref{fixedbht} and from the fact that $s$ is cyclic for $\phi$ that $(E,\phi,s)$ has to be of the form $E= E_1 \oplus \dots \oplus E_r$ for line bundles $E_i$, $s \in E_r$, $\phi(E_i) \subseteq E_{i-1} \otimes K $ and $\phi(E_1)=0$. Note also that all the maps $E_i \rightarrow E_{i-1} \otimes K$ induced by $\phi$ have to be non-zero. This implies that $d_{r-i} = \deg E_{r-i} \geq 2i (1-g)$.

Note that the only proper $\phi$-invariant subbundles of $E$ are of the form $\oplus_{i=1}^k E_i$ for some $1 \leq k <r$. Assume by contradiction that any of those are destabilizing for $(E,\phi)$.

This means that for some $1 \leq k < r$ we have:
\begin{equation*}
\frac{1}{k} \sum_{i=1}^k d_i \geq \frac{1}{r} \sum_{i=1}^r d_i
\end{equation*}
which is equivalent to:
\begin{equation*}
(r-k)  \sum_{i=1}^k d_i \geq k \sum_{i=k+1}^r d_i.
\end{equation*}

This implies
\begin{equation*}
(r-k)  \sum_{i=1}^r d_i \geq r \sum_{i=k+1}^r d_i = r \sum_{i=0}^{r-k-1} d_{r-i} \geq r(r-k)(r-k-1)(1-g)
\end{equation*}
which contradicts $d < r(r-2)(1-g)$.
\end{proof}
\end{prop}

From remark \ref{adhm} we also deduce the following.
\begin{cor}
If $d < r(r-2)(1-g)$ then $\mathcal{M}_\varepsilon^{r,d} = \mathcal{M}_\infty^{r,d}$ is smooth when non-empty.
\begin{proof}
It follows from the fact that in this case there are no proper $\phi$-invariant subbundles that contain the section and so the Zariski tangent space has minimal dimension, since the only endomorphisms that commute with $\phi$ are scalars.
\end{proof}
\end{cor}

Assume $d < r(r-2)(1-g)$. In the case $r \geq 2$ proposition \ref{dimsupp} applied to the Hitchin map $\mathcal{M}_\infty^{r,d} \rightarrow \mathcal{A}^r$ allows to conclude that any properly supported summand of the pushforward of the constant sheaf $\Q$ should have support with dimension at least:
\begin{equation*}
2+2r^2(g-1)-1-r^2(g-1)-d-r(r-1)(g-1) \geq 2+r(r-1)(g-1).
\end{equation*}

Since
\begin{equation*}
\dim \mathcal{A}^r_{red} = \max_{1 \leq k \leq r-1} \{ \dim \mathcal{A}^k+\dim \mathcal{A}^{r-k}\}=2+(g-1)(r^2-2r+2)
\end{equation*}
and $2+r(r-1)(g-1) \geq 2+(g-1)(r^2-2r+2)$ we can already conclude that proper supported summands, if there are any, should have support equal to the closure of $\mathcal{A}^r_{red}$.

It is also immediate to check that the argument in theorem \ref{BHMSMY} about the comparison of weight polynomials for nodal curves can be applied again. Just note that in this case we should let $\mathcal{C}$ and $\mathcal{D}$ be generic nodal curves in $\mathcal{A}^r_{red}$ and $\mathcal{A}^r_{int}$. Here lemma \ref{ratcurve} should be modified as follows and the rest of the argument is the same.

\begin{lemma}
Denote by $\Sigma$ a curve with two rational components meeting in $(r-1)(2g-2)$ simple nodes and by $\overline\Sigma$ an integral curve with $(r-1)(2g-2)$ simple nodes whose normalization is isomorphic to $\C\P^1$. Define:
\begin{equation*}
U(\Sigma)=\sum_{n = 0}^{(r-1)(2g-2)-1} q^n [CKS^n [-n]]
\end{equation*}
and
\begin{equation*}
U(\overline\Sigma)=\sum_{n= 0}^{(r-1)(2g-2)-1} q^n [\overline{CKS}^n [-n]]
\end{equation*}
where $[CKS^n]$ is the weight polynomial of the $n$-th CKS complex associated to $\Sigma$ and $[\overline{CKS}^n]$ is the analogous object for $\overline\Sigma$.
Then:
\begin{equation*}
U(\overline\Sigma)=(1-q\Q)(1-q\L)U(\Sigma) \qquad \text{mod } q^{(r-1)(2g-2)}.
\end{equation*}
\end{lemma}

This concludes the generalization of theorem \ref{BHMSMY} from the range $d<0$ for rank 2 to the range $d < r(r-2)(1-g)$ for rank $r$.

Now, using \cite[corollary 3.4]{nitsure1991moduli} we can immediately deduce the same result for $d> r(2g-1)+(r-1)^2(2g-2)$ coprime with $r$ since we know that in this case the Abel-Jacobi map is a projective bundle and the argument in theorem \ref{BHMSMY} carries through. As we already pointed out, the presence of extra summands (even with non-trivial global cohomology) in the decomposition theorem for the Hitchin map $\mathcal{M}^{r,d} \rightarrow \mathcal{A}^r$ would not cause problems at this point.

Serre's duality would allow to deduce the cohomological statement from theorem \ref{BHMSMY} for $d$ coprime with $r$ in the range $r^2(g-1) < d \leq r(r+1)(g-1)$ from the statement we already proved for $r(r-1)(1-g) \leq d < r(r-2)(1-g)$, provided $\mathcal{M}_\varepsilon^{r,d}$ is smooth for $r$, $d$ coprime and $d>r^2(g-1)$.

Ultimately, we would be left with $d$ coprime with $r$ in the range $r(r-2)(1-g) < d < r(2g-1)+(r-1)^2(2g-2)$. It would be reasonable to guess that the sharp estimate for $d$ for which $H^1$ of a semistable Higgs bundle of degree $d$ and rank $r$ vanishes is $d>r(r-2)(1-g)$ rather than the one in \cite[corollary 3.4]{nitsure1991moduli}, but we do not attempt to prove it here.